\documentclass{article}
\usepackage{amsmath}
\usepackage{amsfonts}
\usepackage{amscd}

\def\BbR{\mathbb{R}}
\def\BbN{\mathbb{N}}

\def\BbZ{\mathbb{Z}}

\def\e{\epsilon}
\def\a{\alpha}
\def\b{\beta}
\def\c{\gamma}
\def\p{\psi}
\def\vp{\varphi}

\def\s{\sigma}

\def\num#1{1,\ldots,#1}
\def\word#1#2{{#1}_1\ldots{#1}_{#2}}

\def\A{\mathcal{A}}
\def\M{\mathcal{M}}
\def\P{\mathcal{P}}

\def\C{\mathcal{C}}
\def\H{\mathcal{H}}

\def\E{\mathcal{E}}
\def\U{\mathcal{U}}
\def\D{\mathcal{D}}
\def\I{\mathcal{I}}

\def\B{\mathcal{B}}
\def\R{\mathcal{R}}
\def\F{\mathcal{F}}
\def\G{\mathcal{G}}
\def\Q{\mathcal{Q}}

\def\T{\mathcal{T}}
\def\O{\mathcal{O}}

\def\SS{\Sigma}
\def\LL{\Lambda}
\def\GG{\Gamma}
\def\gen{\underset{{\rm GE}}{\sim}}

\def\bb{\mathbf b}

\def\bp{\mathbf p}

\def\CH{\C\H}

\def\diam#1{{\rm diam}(#1)}

\def\inte#1{{\rm int}(#1)}

\def\bp{\mathbf p}

\def\sd#1#2{#1\backslash#2}

\usepackage{amsthm}

\theoremstyle{plain}
\newtheorem{Th}{Theorem}[section]
\newtheorem{Lem}[Th]{Lemma}
\newtheorem{Cor}[Th]{Corollary}
\newtheorem{Prop}[Th]{Proposition}

\theoremstyle{definition}
\newtheorem{Def}[Th]{Definition}
\newtheorem*{Not}{Notation}
\newtheorem{Ex}[Th]{Example}
\newtheorem{Assum}[Th]{Assumption}

\theoremstyle{remark}
\newtheorem*{Rem}{Remark}

\def\thm{\begin{Th}}
\def\endthm{\end{Th}}
\def\lemma{\begin{Lem}}
\def\endlemma{\end{Lem}}
\def\cor{\begin{Cor}}
\def\endcor{\end{Cor}}
\def\prop{\begin{Prop}}
\def\endprop{\end{Prop}}
\def\definition{\begin{Def}}
\def\enddefinition{\end{Def}}
\def\remark{\begin{Rem}}
\def\endremark{\end{Rem}}
\def\example{\begin{Ex}}
\def\endexample{\end{Ex}}
\def\demo{\begin{proof}}
\def\enddemo{\end{proof}}

\def\notation{\begin{Not}}
\def\endnotation{\end{Not}}
\def\assumption{\begin{Assum}}
\def\endassumption{\end{Assum}}
\usepackage{graphicx}
\usepackage{curves}
\usepackage{epic}
\usepackage{eepic}

\def\O{\mathcal{O}}
\def\d{\delta}

\def\tT{\widetilde{T}}

\def\qs{\underset{{\rm QS}}{\sim}}

\def\od{\overline{d}}
\def\orho{\overline{\rho}}
\def\bl{\underset{\text{BL}}{\sim}}
\def\ac{\underset{\text{AC}}{\sim}}
\def\X{\mathcal{X}}
\def\tpi{\widetilde{\pi}}
\def\wT{\widetilde{T}}
\def\Mod{{\rm Mod}}

\def\tvp{\widetilde{\vp}}
\def\wtg{\widetilde{g}}
\def\wth{\widetilde{h}}
\def\oN{\overline{N}}
\def\uN{\underline{N}}
\def\oI{\overline{I}}
\def\uI{\underline{I}}
\def\oR{\overline{R}}
\def\uR{\underline{R}}
\def\oE{\overline{\mathcal{E}}}
\def\uE{\underline{\mathcal{E}}}
\def\M{\mathcal{M}}
\def\oM{\overline{\M}}
\def\uM{\underline{\M}}
\def\oOmega{\overline{\Omega}}
\def\uR{\underline{R}}
\def\oR{\overline{R}}
\def\ud{\underline{d}}
\def\od{\overline{d}}
\def\inn#1{\langle{#1}\rangle}

\usepackage{euscript}
\renewcommand{\theequation}{\arabic{section}.\arabic{equation}}

\begin{document}
\begin{center}
{\Large \bf Weighted partition of a compact metrizable space, 
its hyperbolicity\\
and Ahlfors regular conformal dimension}\\
by\\
{\large Jun Kigami}\\
{Graduate School of Informatics}\\
{Kyoto University}\\
{e-mail: kigami@i.kyoto-u.ac.jp}
\end{center}
\begin{abstract}
Successive divisions of compact metric spaces appear in many different areas of mathematics such as the construction of self-similar sets, Markov partitions associated with hyperbolic dynamical systems, dyadic cubes associated with a doubling metric space. The common feature in these is to divide a space into a finite number of subsets, then divide each subset into finite pieces and repeat this process again and again. In this paper we generalize such successive divisions and call them partitions. Given a partition, we consider the notion of a weight function  assigning  a ``size'' to each piece of the partition. Intuitively we believe that a partition and a weight function should provide a ``geometry'' and an ``analysis'' on the space of our interest. We are going to pursue this idea in three parts. In the first part, the metrizability of a weight function, i.e. the existence of a metric ``adapted to'' a given weight function, is shown to be equivalent to the Gromov hyperbolicity of the graph associated with the weight function. In the second part, the notions like bi-Lipschitz equivalence, Ahlfors regularity, the volume doubling property and quasisymmetry will be shown to be equivalent to certain properties of weight functions. In particular, we find that quasisymmetry and the volume doubling property are the same notion in the world of weight functions. In the third part, a characterization of the Ahlfors regular conformal dimension of a compact metric space is given as the critical index $p$ of $p$-energies associated with the partition and the weight function corresponding to the metric.
\end{abstract}

\tableofcontents

\section{Introduction}\label{INT}
\begin{figure}
\centering
\setlength{\unitlength}{20mm}

\begin{picture}(4.3,2.5)(0.3,-0.5)
\linethickness{1pt}
\thinlines

\put(4,1.8){\circle*{0.1}}
\put(3.5,1.2){\circle*{0.1}}
\put(4.5,1.2){\circle*{0.1}}
\put(3.25,0.6){\circle*{0.1}}
\put(3.75,0.6){\circle*{0.1}}
\put(4.25,0.6){\circle*{0.1}}
\put(4.75,0.6){\circle*{0.1}}
\put(3.125, 0){\circle*{0.1}}
\put(3.375, 0){\circle*{0.1}}
\put(3.625, 0){\circle*{0.1}}
\put(3.875, 0){\circle*{0.1}}
\put(4.125, 0){\circle*{0.1}}
\put(4.375, 0){\circle*{0.1}}
\put(4.625, 0){\circle*{0.1}}
\put(4.875, 0){\circle*{0.1}}
\drawline(4,1.8)(3.5, 1.2)(3.25, 0.6)(3.125,0)
\drawline(3.5,1.2)(3.75,0.6)(3.875,0)
\drawline(3.75,0.6)(3.625, 0)
\drawline(3.25,0.6)(3.375, 0)
\drawline(4,1.8)(4.5,1.2)(4.75,0.6)(4.875, 0)
\drawline(4.5,1.2)(4.25,0.6)(4.125, 0)
\drawline(4.25,0.6)(4.375, 0)
\drawline(4.75,0.6)(4.625, 0)
\put(4,1.87){\makebox(0,0)[b]{$\phi$}}
\put(3.35,1.22){\makebox(0,0){$0$}}
\put(4.65,1.22){\makebox(0,0){$1$}}
\put(3.09,0.62){\makebox(0,0){$00$}}
\put(3.59,0.62){\makebox(0,0){$01$}}
\put(4.41,0.62){\makebox(0,0){$10$}}
\put(4.9,0.62){\makebox(0,0){$11$}}

\thicklines
\drawline(0, 1.8)(2.5,1.8)
\put(0, 1.8){\circle*{0.03}}
\put(2.5, 1.8){\circle*{0.03}}
\put(1.25,1.87){\makebox(0,0)[b]{\large$K_{\phi}$}}
\put(0, 1.87){\makebox(0,0)[b]{$0$}}
\put(2.5, 1.87){\makebox(0,0)[b]{$1$}}
\drawline(0, 1.22)(2.5, 1.22)
\put(0, 1.22){\circle*{0.03}}
\put(2.5, 1.22){\circle*{0.03}}
\put(1.25, 1.22){\circle*{0.03}}
\put(0.625,1.27){\makebox(0,0)[b]{\large$K_0$}}
\put(1.875,1.27){\makebox(0,0)[b]{\large$K_1$}}
\put(2.5, 1.8){\circle*{0.03}}
\put(1.25, 1.27){\makebox(0,0)[b]{$\frac 12$}}
\drawline(0, 0.62)(2.5, 0.62)
\put(0, 0.62){\circle*{0.03}}
\put(2.5, 0.62){\circle*{0.03}}
\put(1.25, 0.62){\circle*{0.03}}
\put(0.625, 0.62){\circle*{0.03}}
\put(1.875, 0.62){\circle*{0.03}}
\put(0.3125, 0.67){\makebox(0,0)[b]{\large$K_{00}$}}
\put(0.9375, 0.67){\makebox(0,0)[b]{\large$K_{01}$}}
\put(1.5625, 0.67){\makebox(0,0)[b]{\large$K_{10}$}}
\put(2.1825, 0.67){\makebox(0,0)[b]{\large$K_{11}$}}
\put(0.625, 0.67){\makebox(0,0)[b]{$\frac 14$}}
\put(1.875, 0.67){\makebox(0,0)[b]{$\frac 34$}}
\drawline(0, 0)(2.5, 0)
\put(0, 0){\circle*{0.03}}
\put(2.5, 0){\circle*{0.03}}
\put(1.25, 0){\circle*{0.03}}
\put(0.625, 0){\circle*{0.03}}
\put(1.875, 0){\circle*{0.03}}
\put(0.3125, 0){\circle*{0.03}}

\thinlines
\dottedline{0.03}(0.3125, 0)(0.3125, -0.2)
\put(0.3125, 0.05){\makebox(0,0)[b]{$\frac 18$}}
\put(0.9375, 0){\circle*{0.03}}
\dottedline{0.03}(0.9375, 0)(0.9375, -0.2)
\put(0.9375, 0.05){\makebox(0,0)[b]{$\frac 38$}}
\put(1.5625, 0){\circle*{0.03}}
\dottedline{0.03}(1.5625, 0)(1.5625, -0.2)
\put(1.5625, 0.05){\makebox(0,0)[b]{$\frac 58$}}
\put(2.1825, 0){\circle*{0.03}}
\dottedline{0.03}(2.1825, 0)(2.1825, -0.2)
\put(2.1825, 0.05){\makebox(0,0)[b]{$\frac 78$}}

\thinlines
\dottedline{0.03}(0, 1.8)(0, -0.2)
\dottedline{0.03}(2.5, 1.8)(2.5, -0.2)
\dottedline{0.03}(1.25, 1.22)(1.25, -0.2)
\dottedline{0.03}(0.625, 0.62)(0.625, -0.2)
\dottedline{0.03}(1.875, 0.62)(1.875, -0.2)

\put(4.1,-0.5){\makebox(0,0){The associated tree}}
\put(1.25,-0.5){\makebox(0, 0){A partition of $[0, 1]$}}

\end{picture}
\caption{A partition of the unit interval $[0, 1]$ and the associated tree}\label{Unit}
\end{figure}
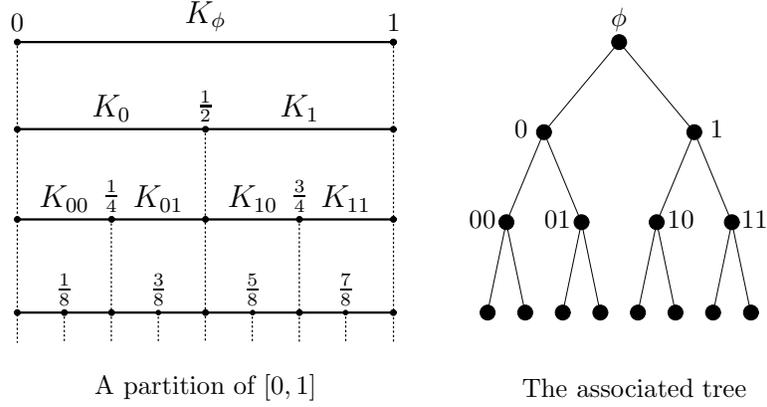

Successive division of a space has played important roles in many areas of mathematics. One of the simplest examples is the binary division of the unit interval $[0, 1]$ shown in Figure~\ref{Unit}. Let $K_{\phi} = [0, 1]$ and divide $K_{\phi}$ in half as $K_0 = [0, \frac 12]$ and $K_1 = [\frac 12, 1]$. Next, $K_0$ and $K_1$ are divided in half again and yield $K_{ij}$ for each $(i, j) \in \{0, 1\}^2$. Repeating this procedure, we obtain $\{K_{\word im}\}_{i_1, \ldots, i_m \in \{0, 1\}}$ satisfying 
\begin{equation}\label{INT.eq10}
K_{\word im} = K_{\word im0} \cup K_{\word im1}
\end{equation}
 for any $m \ge 0$ and $\word im \in \{0, 1\}^m$. In this example, there are two notable properties.\par
The first one is the role of the (infinite) binary tree 
\[
T_b = \{\phi, 0, 1, 00, 01, 10, 11, 000, 001, 010, 011, \ldots\} = \bigcup_{m \ge 0} \{0, 1\}^m,
\]
where $\{0, 1\}^0 = \{\phi\}$. The vertex $\phi$ is called the root or the reference point and $T_b$ is called the tree with the root (or the reference point) $\phi$. Note that the correspondence $\word im \to K_{\word im}$ determines a map from the binary tree to the collection of compact subsets of $[0, 1]$ with the property \eqref{INT.eq10}. \par
Secondly, note that $K_{i_1} \supseteq K_{i_1i_2} \supseteq K_{i_1i_2i_3} \supseteq \ldots$ and
\begin{equation}\label{INT.eq20}
\bigcap_{m \ge 1} K_{\word im}\,\,\text{is a single point}
\end{equation}
for any infinite sequence $i_1i_2\ldots$. (Of course, this is the binary expansion and hence the single point is $\sum_{m \ge 1} \frac{i_m}{2^m}$.) In other words, there is a natural map $\s: \{0, 1\}^{\BbN} \to [0, 1]$ given by 
\[
\s(i_1i_2\ldots) = \bigcap_{m \ge 1} K_{\word im}.
\]
\par

Such a successive division of a compact metric space, which may not be as simple as this one, appears various areas in mathematics. One of the typical examples is a self-similar set in fractal geometry. A self-similar set is a union of finite number of contracted copies of itself.  Then each contracted copy is again a union of contracted copies and so forth. Another example is the Markov partition associated with a hyperbolic dynamical system. See \cite{Adler} for details. Also the division of a metric measure space having the volume doubling property by dyadic cubes  can be thought of as another example of such a division of a space. See Christ\cite{Christ} for example.\par
In general, let $X$ be a compact metrizable topological space with no isolated point. The common properties of the above examples are;\\
(i)\,\,There exists a tree $T$ (i.e. a connected graph without loops) with a root $\phi$. \\
(ii)\,\,
For any vertex $p$ of $T$, there is a corresponding nonempty compact subset of $X$ denoted by $X_p$ and $X = X_{\phi}$.\\
(iii)\,\,
Every vertex $p$ of $T$ except $\phi$ has unique predecessor $\pi(p) \in T$ and 
\begin{equation}\label{INT.eq30}
X_q = \bigcup_{p \in \{p' | \pi(p') = q\}} X_p
\end{equation}
(iv)\,\,
The totality of edges of $T$ is $\{(\pi(q), q)| q \in \sd{T}{\{\phi\}}\}$.\\
(v)\,\,
For any infinite sequence $(p_0, p_1, p_2, \ldots)$ of vertices of $X$ satisfying $p_0 = \phi$ and $\pi(p_{i + 1}) = p_i$ for any $i \ge 1$, 
\begin{equation}\label{INT.eq40}
\bigcap_{i \ge 1} X_{p_i}\,\,\text{is a single point}.
\end{equation}\par
\begin{figure}
\centering
\setlength{\unitlength}{20mm}
\begin{picture}(4.8,2.5)(0.5,-.8)
\linethickness{1pt}
\thinlines

\put(4.1,1.8){\circle*{0.1}}
\put(3.6,0.9){\circle*{0.1}}
\put(4.2,1.1){\circle*{0.1}}
\put(4.8,0.9){\circle*{0.1}}
\put(3.3,0.3){\circle*{0.1}}
\put(3.55,0){\circle*{0.1}}
\put(3.9,-0.3){\circle*{0.1}}
\put(4.45,0.3){\circle*{0.1}}
\put(4.05,0.3){\circle*{0.1}}
\put(5.05,0){\circle*{0.1}}
\put(4.6,-0.3){\circle*{0.1}}

\drawline(4.1,1.8)(3.6,0.9)(3.3,0.3)
\drawline(3.6,0.9)(3.55,0)
\drawline(3.6,0.9)(3.9,-0.3)
\drawline(4.1,1.8)(4.2,1.1)(4.45,0.3)
\drawline(4.2,1.1)(4.05,0.3)
\drawline(4.1,1.8)(4.8,0.9)(4.6,-0.3)
\drawline(4.8,0.9)(5.05,0)
\put(4.1,1.89){\makebox(0,0)[b]{$\phi$}}
\put(3.45,0.9){\makebox(0,0){$1$}}
\put(4.1,1.2){\makebox(0,0){$2$}}
\put(4.95,0.9){\makebox(0,0){$3$}}
\put(3.15,0.3){\makebox(0,0){$11$}}
\put(3.35,0){\makebox(0,0){$12$}}
\put(3.7,-0.3){\makebox(0,0){$13$}}
\put(3.9,0.4){\makebox(0,0){$21$}}
\put(4.6,0.4){\makebox(0,0){$22$}}
\put(4.85,-0.3){\makebox(0,0){$31$}}
\put(5.25,0){\makebox(0,0){$32$}}

\dottedline[\circle*{0.01}]{.05}(3.6,0.9)(4.8,0.9)
\dottedline[\circle*{0.01}]{.05}(3.6,0.9)(4.2,1.1)
\dottedline[\circle*{0.01}]{.05}(4.8,0.9)(4.2,1.1)
\dottedline[\circle*{0.01}]{.05}(3.3,0.3)(3.55,0)
\dottedline[\circle*{0.01}]{.05}(3.9,-0.3)(3.55,0)
\dottedline[\circle*{0.01}]{.05}(3.3,0.3)(4.05,0.3)
\dottedline[\circle*{0.01}]{.05}(4.05,0.3)(3.55,0)
\dottedline[\circle*{0.01}]{.05}(4.6,-0.3)(3.55,0)
\dottedline[\circle*{0.01}]{.05}(3.9,-0.3)(4.6,-0.3)
\dottedline[\circle*{0.01}]{.05}(4.05,0.3)(4.6,-0.3)
\dottedline[\circle*{0.01}]{.05}(4.05,0.3)(4.45,0.3)
\dottedline[\circle*{0.01}]{.05}(4.6,-0.3)(4.45,0.3)
\dottedline[\circle*{0.01}]{.05}(4.6,-0.3)(5.05,0)
\dottedline[\circle*{0.01}]{.05}(5.05,0)(4.45,0.3)

\closecurve(0,1, 0.5,0.2,1.5, 0,2.5,0.2, 3, 1,2.5,1.8, 1.5, 2,0.5,1.8)
\curve(1.5, 0, 1.5, 0.6, 1.4, 0.9, 1.0, 1.6, 0.5, 1.8)
\curve(1.4, 0.9, 1.75,0.9,2.1,1.2,2.5, 1.8)
\curve(2.5, 0.2, 2.2, 0.8,2.1, 1.2)
\curve(1.5, 2, 1.8, 1.7,1.75, 0.9)
\curve(0.5, 0.2, 1, 0.7,1.5, 0.6)
\curve(0, 1, 0.5, 1.1,1, 1.6)
\put(0.5,1.4){\makebox(0,0)[b]{\large $X_{11}$}}
\put(0.7,0.8){\makebox(0,0)[b]{\large $X_{12}$}}
\put(1.1,0.3){\makebox(0,0)[b]{\large $X_{13}$}}
\put(1.45,1.4){\makebox(0,0)[b]{\large $X_{21}$}}
\put(2.05,1.4){\makebox(0,0)[b]{\large $X_{22}$}}
\put(1.9,0.4){\makebox(0,0)[b]{\large $X_{31}$}}
\put(2.6,0.9){\makebox(0,0)[b]{\large $X_{32}$}}
\put(1.5,-0.1){\makebox(0,0)[t]{\large $X = X_{\phi}$}}

\put(4.1,-0.5){\makebox(0,0){The associated graphs $T$}}
\put(4.1,-0.7){\makebox(0,0){Solid lines are vertical edges.}}
\put(4.1,-0.9){\makebox(0,0){Dotted lines are horizontal edges.}}
\put(1.5,-0.5){\makebox(0, 0){A partition of $X$}}
\end{picture}
\caption{A partition and the associated graphs (up to the 2nd stage)}\label{Gen}
\end{figure}
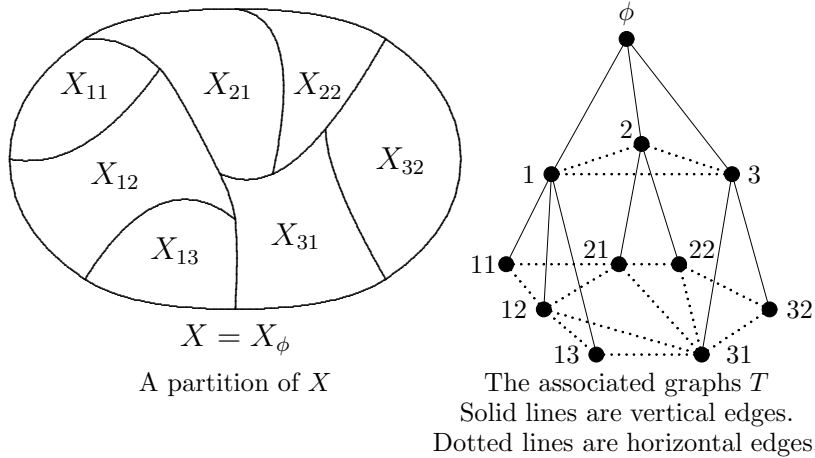

See Figure~\ref{Gen} for illustration of the idea. Note that the properties \eqref{INT.eq30} and \eqref{INT.eq40} correspond to \eqref{INT.eq10} and \eqref{INT.eq20} respectively. In this paper such $\{X_p\}_{p \in T}$ is called a partition of $X$ parametrized by the tree $T$. (We will give the precise definition in Section~\ref{PAS}.) In addition to the ``vertical'' edges , which are the edges of the tree, we provide ``horizontal'' edges to $T$ to describe the combinatorial structure reflecting the topology of $X$ as is seen in Figure~\ref{Gen}. More precisely, a horizontal edge is a pair of $(p, q) \in T \times T$ where $p$ and $q$ have the same distance from the root $\phi$ and $X_p \cap X_q \neq \emptyset$. We call $T$ with horizontal and vertical edges the resolution of $X$ associated with the partition.\par
Another key notion is a weight function on the tree $T$. Note that a metric and a measure give weights of the subsets of $X$. More precisely, let $d$ be a metric on $X$ inducing the original topology of $X$ and let $\mu$ be a Radon measure on $X$ where $\mu(X_p) > 0$ for any $p \in T$. Define $\rho_d:T \to (0, 1]$ and $\rho_{\mu}: T \to (0, 1]$ by
\[
\rho_d(p) = \frac{\diam{X_p, d}}{\diam{X, d}}\quad\text{and}\quad\rho_{\mu}(p) = \frac{\mu(X_p)}{\mu(X)},
\]
where $\diam{A, d}$ is the diameter of $A$ with respect to the metric $d$. Then $\rho_{d}$ (resp. $\rho_{\mu}$) is though of as a natural weight of $X_p$ associated with $d$ (resp. $\mu$).
In both cases where $\# = d$ or $\# = \mu$, the function $\rho_{\#} : T \to (0, \infty)$ satisfy
\begin{equation}\label{INT.eq50}
\rho_{\#}(\pi(p)) \ge \rho_{\#}(p)
\end{equation}
for any $p \in \sd{T}{\{\phi\}}$ and
\begin{equation}\label{INT.eq60}
\lim_{i \to \infty} \rho_{\#}(p_i) = 0
\end{equation}
if $\pi(p_{i + 1}) = p_i$ for any $i  \ge 1$. (To have the second property \eqref{INT.eq60} in case of $\# = \mu$, we must assume that the measure $\mu$ is non-atomic, i.e. $\mu(\{x\}) = 0$ for any $x \in X$.)\par

As we have seen above, given a metric or a measure, we have obtained a weight function $\rho_{\#}$ satisfying \eqref{INT.eq50} and \eqref{INT.eq60}.  In this paper, we are interested in the opposite direction. Namely, given a partition of a compact metrizable topological space parametrized by a tree $T$, we define the notion of weight functions as the collection of functions from $T$ to $(0, 1]$ satisfying the properties \eqref{INT.eq50} and \eqref{INT.eq60}. Then our main object of interest is the space of weight functions including those derived from metrics and measures. Naively we believe that a partition and a weight function essentially determine a ``geometry'' and/or an``analysis'' of the original set no matter where the weight function comes. It may come from a metric, a measure or else.  Keeping this intuition in mind, we are going to develop basic theory of weigh functions in three closely related directions in this paper.\par
The first direction is to study when a weight function is naturally associated with a metric? In brief, our conclusion will be that a power of a weight function is naturally associated with a metric if and only if the rearrangement of the resolution $T$ associated the weight function is Gromov hyperbolic. To be more precise, given a partition $\{X_w\}_{w \in T}$ of a compact metrizable topological space $X$ with no isolated points and a weight function $\rho  : T \to (0, 1]$. In Section~\ref{CNB}, we will define $\d^{\rho}_M(\cdot, \cdot)$, which is called the visual pre-metric associated with $\rho$, in the following way: let $\LL^{\rho}_s$ be the collection of $w$'s  in $T$ where the size $\rho(w)$ is almost $s$. Define a horizontal edge of $\LL^{\rho}_s$ as $(w, v) \in \LL^{\rho}_s$ with $X_w \cap X_v \neq \emptyset$. For $r \in (0, 1)$, the rearranged resolution $\tT^{\rho, r}$ associated with the weight function $\rho$ is defined as the vertices $\cup_{m \ge 0} \LL^{\rho}_{r^m}$ with the vertical edges from the tree structure of $T$ and the horizontal edges of $\LL^{\rho}_{r^m}$. Then the visual pre-metric $\d^{\rho}_M(x, y)$ for $x, y \in X$ is given by the infimum of $s$ where $x$ and $y$ can be connected by an $M$-chain of horizontal edges in $\LL^{\rho}_s$. We think a metric $d$ is naturally associated with the weight function $\rho$ if and only if $d$ and $\LL^{\rho}_M$ are bi-Lipschitz equivalent on $X \times X$. More precisely, we are going to use a phrase``$d$ is adapted to $\rho$'' instead of ``$d$ is naturally associated with $\rho$''. The notion of visual pre-metric is a counterpart of that of visual pre-metric on the boundary of a Gromov hyperbolic metric space, whose detailed account can bee seen in \cite{BuyaloSchr}, \cite{MacTyson} and \cite{KapBen} for example. Now the main conclusion of the first part is Theorem~\ref{HYP.thm20} saying that the hyperbolicity of the rearranged resolution $\tT^{\rho, r}$ is equivalent to the existence of a metric adapted to some power of the weight function. Moreover, if this is the case, the metric adapted to some power of the weight function is shown to be a visual metric in Gromov's sense.\par
The second direction is to establish relationships of various relations between weight functions, metrics and measures. For examples, Ahlfors regularity  and the volume property are relations between measures and metrics.  For $\a > 0$, a measure $\mu$ is $\a$-Ahlfors regular with respect to a metric $d$ if and only if there exist $c_1, c_2 > 0$ such that
\[
c_1r^{\a} \le \mu(B_d(x, r)) \le c_2r^{\a},
\]
where $B_d(x, r) = \{y| y \in X, d(x, y) < r\}$, for any $r \in (0, \diam{X, d}]$ and $x \in X$. See Definition~\ref{VDP.def20} for the precise definition of the volume doubling property. On the other hand, bi-Lipschitz and quasisymmetry are equivalence relations between two metrics. (The precise definitions of bi-Lipschitz equivalence and quasisymmetry are given in Definitions~\ref{BLE.def30} and \ref{INT.def10} respectively.)  Regarding those relations, we are going to claim the following relationships
\begin{equation}\label{INT.eq100}
\text{bi-Lipschitz = Ahlfors regularity = being adapted}
\end{equation}
and
\begin{equation}\label{INT.eq110}
\text{the volume doubling property = quasisymmetry.}
\end{equation}
in the framework of weight functions. To illustrate the first claim more explicitly, let us introduce the notion of bi-Lipschitz equivalence of weight  functions. Two weight functions $\rho_1$ and $\rho_2$ are said to be bi-Lipschitz equivalent if and only if there exist $c_1, c_2 > 0$ such that
\[
c_1\rho_1(p) \le \rho_2(p) \le c_2\rho_1(p)
\]
for any $p \in T$. Now the first claim can be resolved into three parts as follows: let $\rho_1$ and $\rho_2$ be two weight functions.\\
Claim 1: Suppose that $\rho_1 = \rho_{d_1}$ and $\rho_2 = \rho_{d_2}$ for metrics $d_1$ and $d_2$ on $X$. Then $\rho_1$ and $\rho_2$ are bi-Lipschitz equivalent if and only if $d_1$ and $d_2$ are bi-Lipschitz equivalent as metrics.\\
Claim 2: Suppose that $\rho_1 = \rho_d$ and $\rho_2 = \rho_{\mu}$ for a metric $d$ on $X$ and a Radon measure $\mu$ on $X$. Then $(\rho_1)^{\a}$ and $\rho_2$ are bi-Lipschitz equivalent if and only if $\mu$ is $\a$-Ahlfors regularity of $\mu$ with respect to $d$.\\
Claim 3: Suppose that $\rho_1 = \rho_d$ for a metric $d$ on $X$, then $\rho_1$ and $\rho_2$ are bi-Lipschitz equivalent if and only if the metric $d$ is adapted to the weight function $\rho_2$. \\
One can find the precise statement in Theorem~\ref{SMR.thm20} in the case of partitions of $S^2$. The second claim is rationalized in the same manner. See Theorem~\ref{SMR.thm30} for the exact statement in the case of $S^2$ for example.\par
The third direction is a characterization of Ahlfors regular conformal dimension. The Ahlfors regular conformal dimension, AR conformal dimension for short, of a metric space $(X, d)$ is defined as 
\begin{multline*}
\dim_{AR}(X, d) =\\
 \inf\{\a| \text{there exist a metric $\rho$ on $X$ and a Borel regular measure $\mu$ on $X$}\\
\text{ such that $\rho \qs d$ and $\mu$ is $\a$-Ahlfors regular with respect to $\rho$}\},
\end{multline*}
where ``$\rho \qs d$'' means that the two metrics $\rho$ and $d$ are quasisymmetric to each other. In \cite{CarPiag}, Carassco Piaggio has given a characterization of Ahlfors regular conformal dimension in terms of the critical exponent of $p$-combinatorial modulus of discrete path families. In view of the results from the previous part, we have obtained the ways to express the notions of quasisymmetry and Ahlfors regularity in terms of weight functions. So we are going to translate Carassco Piaggio's work into our framework. However, we are going to use the critical exponent of $p$-energy instead of $p$-combinatorial modulus in our work.\footnote{This idea of characterizing AR conformal dimension by $p$-energies was brought to the author by B. Kleiner in a personal communication.} Furthermore, we are going to define the notion of $p$-spectral dimension and present a relation between Ahlfors regular conformal dimension and $p$-spectral dimension. In particular, for $p = 2$, the $2$-spectral dimension has been know to appear in the asymptotic behavior of the Brownian motion and the eigenvalue counting function of the Laplacian on certain fractals like the Sierpinski gasket and the Sierpinski carpet. See \cite{BP}, \cite{BB6} and \cite{KL1} for example. For the Sierpinski carpet, we will show that the $2$-spectral dimension gives an upper bound of Ahlfors regular conformal dimension.\par
One of the ideas behind this study is to approximate a space by a series of graphs. Such an idea has already been explored in association with hyperbolic geometry. For example, in \cite{Elek1} and \cite{BouPaj},  they have constructed an infinite graph whose hyperbolic boundary is homeomorphic to given compact metric space. Their method is first construct a series of coverings of the space, which is a counterpart of our partition, and construct a graph from the series. In \cite{CarPiag}, Carrasco Piaggio has utilized this series of coverings to study Ahlfors regular conformal dimension of the space. His notion of ``relative radius'' essentially corresponds to our weight function. In our framework, the original space is homeomorphic to the analogue of hyperbolic boundary of the resolution $T$ of $X$ even if it is not hyperbolic in the sense of Gromov. See Theorem~\ref{HYP.thm05} for details. In other words, the resolution $T$ of $X$ is a version of hyperbolic filling of the original space $X$. (See \cite{BonkSaks} for the notion of hyperbolic fillings.) In this respect, our study in this paper may be thought of as a theory of weighted hyperbolic fillings.\par
The organization of this paper is as follows. In Section~\ref{SMR}, we give a summary of the main results of this paper in the case of the 2 dimensional sphere as a showcase of the full theory. In Section~\ref{TWR}, we give basic definitions and notations on trees.   Section~\ref{PAS} is devoted to the introduction of partitions and related notions. In Section~\ref{CNB}, we define the notion of weight function and the associated ``visual pre-metric''. We study our first question mentioned above, namely, when a weight function is naturally associated with a (power of) metric in Section~\ref{MAG}. In Section~\ref{HYP}, we are going to relate this question to the hyperbolicity of certain graph associated with a weight function. Section~\ref{BLE} is devoted to justifying the statement \eqref{INT.eq100}. In Sections~\ref{TGF}, \ref{VDP}, \ref{GAE} and \ref{QSY}, we will study the rationalized version of \eqref{INT.eq110} as mathematical statement. In particular, in Section~\ref{VDP}, we introduce the key notion of being ``gentle''.  In Section~\ref{ESS}, we apply our general theory to certain class of subsets of the square and obtain concrete (counter) examples. From Section~\ref{CAM}, we will start arguing a characterization of Ahlfors regular conformal dimension.  From Section~\ref{CAM} to \ref{CAR}, we discuss how to obtain a pair of a metric $d$ and a measure $\mu$ where $\mu$ is $\a$-Ahlfors regular with respect to $d$ for a given order $\a$. The main result of these sections is Theorem~\ref{CAR.thm10}. In Section~\ref{PEN}, we will give a characterization of the Ahlfors regular conformal dimension as a critical index $p$ of $p$-energies. Then we will show the relation of the  Ahlfors regular conformal dimension and $p$-spectral dimension in Section~\ref{APP}. Additionally, we will give another characterization of the Ahlfors regular conformal dimension by $p$-modulus of curve families in Section~\ref{MIT}.  Finally in Section~\ref{LOD}, we present the whereabouts of definitions, notations and conditions appearing in this paper for reader's sake.

\setcounter{equation}{0}
\section{Summary of the main results; the case of $2$-dim. sphere}\label{SMR}

In this section, we summarize our main results in this paper in the case of $2$-dimensional sphere $S^2$ (or the Riemann sphere in other words), which is denoted by $X$ in what follows. We use $d_s$ to denote the standard spherical geodesic metric on $X$. Set
\[
\U = \{A| A \subseteq X, \text{closed}, \inte{A} \neq \emptyset, \text{$\partial{A}$ is homeomorphic to the circle $S^1$.}\}
\]
First we divide $X$ into finite number of subsets $X_1, \ldots, X_{N_0}$ belonging to $\U$, i.e.
\[
X = \bigcup_{i = 1}^{N_0} X_i
\]
We assume that $X_i \cap X_j = \partial{X_i} \cap \partial{X_j}$ if $i \neq j$. Next each $X_i$ is divided into finite number of its subsets $X_{i1}, X_{i2}, \ldots, X_{i N_i} \in \U$ in the same manner as before.  Repeating this process, we obtain $X_{i_1{\ldots}i_k}$ for any $i_1\ldots{i_k}$ satisfying
\begin{equation}\label{SMR.eq10}
X_{i_1{\ldots}i_k} = \bigcup_{j = 1, \ldots, N_{i_1{\ldots}i_k}} X_{i_1{\ldots}i_kj}
\end{equation}
and if $i_1{\ldots}i_k \neq j_1{\ldots}j_k$, then
\begin{equation}\label{SMR.eq20}
X_{i_1{\ldots}i_k} \cap X_{j_1{\ldots}j_k} = {\partial}X_{i_1\ldots{i_k}} \cap {\partial}X_{j_1{\ldots}j_k}.
\end{equation}
Note that \eqref{SMR.eq10} is a counterpart of \eqref{INT.eq30}. Next define
\[
T_k = \{i_1{\ldots}i_k|  i_j \in \{1, \ldots, N_{i_1{\ldots}i_{j - 1}}\}\,\,\text{for any $j = 1, \ldots, k - 1$}\}
\]
for any $k = 0, 1, \ldots$, where $T_0$ is a one point set $\{\phi\}$. Let $T = \cup_{k \ge 0} T_k$. Then $T$ is naturally though of as a (non-directed) tree whose edges are given by the totality of $(i_1{\ldots}i_k, i_1{\ldots}i_k{i_{k + 1}})$. We regard the correspondence $w \in T$ to $X_w \in \U$ as a map from $T$ to $\U$, which is denoted by $\X$. Namely, $\X(w) = X_w$ for any $w \in T$. Note that $\X(\phi) = X$. Define
\[
\SS = \{i_1i_2\ldots | i_1\ldots{i_k} \in T_k\,\,\text{for any $k \ge 0$}\},
\]
which is the ``boundary'' of the infinite tree $T$.\par
Furthermore we assume that for any $i_1i_2\ldots \in \SS$
\[
\bigcap_{k = 1, 2, \ldots} X_{i_1{\ldots}i_k}
\]
is a single point, which is denoted by $\s(i_1i_2\ldots)$. Note that $\s$ is a map from $\SS$ to $X$.  This assumption corresponds to \eqref{INT.eq40} and hence the map $\X$ is a partition of $X$ parametrized by the tree $T$, i.e. it satisfies the conditions (i), (ii), (iii), (iv) and (v) in the introduction. Since $X = \cup_{w \in T_k} X_w$ for any $k \ge 0$, this map $\s$ is surjective. \par
In \cite[Chapter 5]{BonkMeyer}, the authors have constructed ``cell decomposition'' associated with an expanding Thurston map. This ``cell decomposition'' is, in fact, an example of a partition formulated above.\par
Throughout this section, for simplicity, we assume the following conditions (SF) and (TH), where (SF) is called strong finiteness in Definition~\ref{PAS.def200} and (TH) ensures the thickness of every exponential weight function. See Definition~\ref{VDP.def30} for the ``thickness'' of a weight function.\\
(SF)\,\,
\begin{equation}\label{SMR.eq100}
\#(\s^{-1}(x)) < +\infty,
\end{equation}
where $\#(A)$ is the number of elements in a set $A$.\\
(TH)\,\,There exists $m \ge 1$ such that for any $w = i_1{\ldots}i_n \in T$, there exists $v = i_1\ldots{i_n}i_{n + 1}{\ldots}i_{n + m} \in T$ such that $X_v \subseteq \inte{X_w}$.\par
The main purpose of this paper is to describe metrics and measures of $X$ from a given weight assigned to each piece $X_w$ of the partition $\X$.

\definition\label{SMR.def100}
A map $g: T \to (0, 1]$ is called a weight function if and only if it satisfies the following conditions (G1), (G2) and (G3).\\
(G1)\,\,$g(\phi) = 1$\\
(G2)\,\,$g(i_1{\ldots}i_k) \ge g(i_1{\ldots}i_ki_{k + 1})$ for any $i_1{\ldots}i_k \in T$ and $i_1{\ldots}i_ki_{k + 1} \in T$.\\
(G3)
\[
\lim_{m \to 0} \sup_{w \in T_k} g(w) = 0.
\]
Define
\[
\G(T) = \{g| \text{$g: T \to (0, 1]$ is a weight function.}\}
\]
Moreover, we define following conditions (SpE) and (SbE), which represent ``super-exponential'' and ``sub-exponential'' respectively: \\
(SpE)\,\, There exists $\lambda \in (0, 1)$ such that
\[
g(i_1{\ldots}i_ki_{k + 1}) \ge {\lambda}g(i_1{\ldots}i_k)
\]
for any $i_1{\ldots}i_k \in T$ and $i_1{\ldots}i_ki_{k + 1} \in T$. \\
(SbE)\,\,There exist $m \in \BbN$ and $\c \in (0, 1)$ such that 
\[
g(i_1{\ldots}i_ki_{k + 1}{\ldots}i_{k + m}) \le {\c}g(i_1{\ldots}i_k)
\]
for any $i_1{\ldots}i_k \in T$ and $i_1{\ldots}i_ki_{k + 1}{\ldots}i_{k + m} \in T$.\par
Set
\[
\G_e(T) = \{g| \text{$g: T \to (0, 1]$ is a weight function satisfying (SpE) and (SbE).}\}.
\]
\enddefinition

Metrics and measures on $X$ naturally have associated weight functions.

\definition\label{SMR.def120}
Set
\begin{multline*}
\D(X) = \{d| \text{$d$ is a metric on $X$ which produces the original topology of $X$},\\
\text{and}\,\,\diam{X, d} = 1\}
\end{multline*}
and
\begin{multline*}
\M(X) = \{\mu| \text{$\mu$ is a Borel regular probability measure on $X$, $\mu(\{x\}) = 0$}\\
\text{for any $x \in T$ and $\mu(O) > 0$ for any non-empty open set $O \subseteq X$}\}
\end{multline*}
For any $d \in \D(X)$, define $g_d: T \to (0, 1]$ by $g_d(w) = \diam{X_{w}, d}$ and for any $\mu \in \M(X)$, define $g_{\mu}: T \to (0, 1]$ by $g_{\mu}(w) = \mu(X_w)$ for any $w \in T$.
\enddefinition

From Proposition~\ref{PAS.prop00}, we have the following fact.

\prop\label{SMR.prop200}
If $d \in \D(X)$ and $\mu \in \M(X)$, then $g_d$ and $g_{\mu}$ are weight functions.
\endprop

So a metric $d \in \D(X)$ has associated weight function $g_d$. How about the converse direction, i.e. for a given weight function $g$, is there a metric $d$ such that $g = g_d$? To make this question more rigorous and flexible, we define the notion of ``visual pre-metric'' $\d_M^g(\cdot, \cdot)$ associated with a weight function $g$.

\definition\label{SMR.def110}
Let $g \in \G(T)$. Define
\[
\LL_s^g = \{i_1{\ldots}i_k | i_1{\ldots}i_k \in T, g(i_1{\ldots}i_{k - 1}) > s \ge g(i_1{\ldots}i_k)\}
\]
for $s \in (0, 1]$ and 
\begin{multline*}
\d_M^g(x, y) = \inf\{s| \text{there exist $w(1), \ldots, w(M + 1) \in \LL_s^g$ such that}\\
\text{$x \in X_{w(1)}, y \in X_{w(M + 1)}$ and $X_{w(j)} \cap X_{w(j + 1)} \neq \emptyset$ for any $j = 1, \ldots, M$}\}
\end{multline*}
for $x, y \in X$. A weight function is called uniformly finite if and only if 
\[
\sup_{s \in (0, 1], w \in \LL_s^g}\#(\{v | v \in \LL_s^g, X_w \cap X_v \neq \emptyset\}) < +\infty.
\]
\enddefinition

Although $\d_M^g(x, y) \ge 0$, $\d_M^g(x, y) = 0$ if and only if $x = y$ and $\d_M^g(x, y) = \d_M^g(y, x)$, the quantity $\d_M^g$ may not satisfy the triangle inequality in general. The visual pre-metric $\d_M^g(x, y)$ is a counterpart of the visual metric defined in \cite{BonkMeyer}. See Section~\ref{CNB} for details. \par
If the pre-metric $\d_M^g(\cdot, \cdot)$ is bi-Lipschitz equivalent to a metric $d$, we consider $d$ as the metric which is naturally associated with the weight function $g$.

\definition\label{SMR.def130}
Let $M \ge 1$\\
(1)\,\,
A metric $d \in \D(X)$ is said to be $M$-adapted to a weight function $g \in \G(X)$ if and only if there exist $c_1, c_2 > 0$ such that
\[
c_1d(x, y) \le \d_M^g(x, y) \le c_2d(x, y)
\]
for any $x, y \in X$.\\
(2)\,\,
A metric $d$ is said to be $M$-adapted if and only if it is $M$-adapted to $g_d$ and it is said to be adapted if it is $M$-adapted for some $M \ge 1$.\\
(3)\,\,
Define
\begin{align*}
\D_{A, e}(X) &= \{d| d \in \D(X),\,\,\text{$g_d \in \G_e(T)$ and $d$ is adapted.}\}\\
\M_e(X) &= \{\mu| \mu \in \M(X), g_{\mu} \in \G_e(T)\}
\end{align*}
\enddefinition

The value $M$ really matters. See Example~\ref{COM.ex10} for an example.\par
The following definition is used to describe an equivalent condition for the existence of an adapted metric in Theorem~\ref{SMR.thm10}.

\definition\label{SMR.def140}
Let $g \in \G(T)$. For $r \in (0, 1)$, define $\tT^{g, r} = \cup_{m \ge 0} \LL^g_{r^m}$. Define the horizontal edges $E^h_{g, r}$ and the vertical edges $E^v_{g, r}$ of $\tT^{g, r}$ as
\[
E^h_{g, r} = \{(w, v)| w, v \in \LL^g_{r^m}\,\, \text{for some $m \ge 0$}, w \neq v, X_w \cap X_v \neq \emptyset\}
\]
and
\[
E^v_{g, r} = \{(w, v)| w \in \LL^g_{r^m}, v \in  \LL^g_{r^{m + 1}}\,\,\text{for some $m \ge 0$}, X_w \supseteq X_v\}
\]
respectively.
\enddefinition

The following theorem is a special case of Theorem~\ref{HYP.thm20}.

\thm\label{SMR.thm10}
Let $g \in \G(X)$. There exist $M \ge 1$, $\a > 0$ and  a metric $d \in \D(X)$ such that $d$ is $M$-adapted to $g^{\a}$ if and only if the graph $(\tT^{g, r}, E^h_{g, r} \cup E^v_{g, r})$ is Gromov hyperbolic for some $r  \in (0, 1)$. Moreover, if this is the case, then the adapted metric $d$ is a visual metric in the Gromov sense.
\endthm

Next, we define two equivalent relations $\bl$ and $\gen$ on the collection of exponential weight functions. Later, we are going to identify these with known relations according to the types of weight functions.

\definition\label{SMR.def150}
Let $g, h \in \G_e(T)$.\\
(1)\,\,
$g$ and $h$ are said to be bi-Lipschitz equivalent if and only if there exist $c_1, c_2 > 0$ such that
\[
c_1g(w) \le h(w) \le c_2g(w)
\]
for any $w \in T$. We write $g \bl h$ if $g$ and $h$ are bi-Lipschitz equivalent.\\
(2)\,\,
$h$ is said to be gentle to $g$ if and only if there exists $\c > 0$ such that if $w, v \in \LL_s^g$ and $X_w \cap X_v \neq \emptyset$, then $h(w) \le {\c}h(v)$. We write $g \gen h$ if $h$ is gentle to $g$.
\enddefinition

Clearly, $\bl$ is an equivalence relation. On the other hand, the fact that $\gen$ is an equivalence relation is not quite obvious and going to be shown in Theorem~\ref{GAE.thm10}.

\prop\label{SMR.prop10}
The relations $\bl$ and $\gen$ are equivalent relations in $\G_e(T)$. Moreover, if $g \bl h$, then $g \gen h$.
\endprop

Some of the properties of a weight function is invariant under the equivalence relation $\gen$ as follows.

\prop\label{SMR.prop20}
{\rm (1)}\,\,
Being uniformly finite is invariant under the equivalence relation $\gen$, i.e. if $g \in \G_e(T)$ is uniformly finite, $h \in \G_e(T)$ and $g \gen h$, then $h$ is uniformly finite.\\
{\rm (2)}\,\,
The hyperbolicity of $\tT^{g, r}$ is invariant under the equivalence relation $\gen$.
\endprop

The statements (1) and (2) of the above theorem are the special cases of Theorem~\ref{GAE.thm20} and Theorem~\ref{COM.thm30} respectively. \par
The next theorem shows that bi-Lipschitz equivalence of weight functions can be identified with other properties according to types of involved weight functions.

\thm\label{SMR.thm20}
{\rm (1)}\,\,
For $d, \rho \in \D_{A, e}(X)$, $g_d \bl g_{\rho}$ if and only if $d$ and $g$ are bi-Lipschitz equivalent as  metrics.\\
{\rm (2)}\,\,
For $\mu, \nu \in \M(X)$, $g_{\mu} \bl g_{\nu}$ if and only if there exist $c_1, c_2 > 0$ such that 
\[
c_1\mu(A) \le \nu(A) \le c_2\mu(A)
\]
for any Borel set $A \subseteq X$.\\
{\rm (3)}\,\,
For $g \in \G_e(X)$ and $d \in \D_{A, e}(X)$, $g \bl g_d$ if and only if $d$ is $M$-adapted to $g$ for some $M \ge 1$.\\
{\rm (4)}\,\,
For $d \in \D_{A, e}(X)$, $\mu \in \M(X)$ and $\a > 0$, $(g_d)^{\a} \bl g_{\mu}$ and $g_d$ is uniformly finite if and only if $\mu$ is $\a$-Ahlfors regular with respect to $d$, i.e. there exist $c_1, c_2 > 0$ such that
\[
c_1r^{\a} \le \mu(B_d(x, r)) \le c_2r^{\a}
\]
for any $r > 0$ and $x \in X$.
\endthm

The statements (1), (2), (3) and (4) of the above theorem follow from Corollary~\ref{BLE.cor10}, Theorem~\ref{BLE.thm20}, Corollary~\ref{BLE.cor20} and Theorem~\ref{MEM.thm100} respectively.\par
The gentle equivalence relation is identified with ``quasisymmetry'' between metrics and ''volume doubling property'' between a metric and a measure.

\thm\label{SMR.thm30}
{\rm (1)}\,\,
Let $d \in \D_{A, e}(X)$ and $\mu \in \M(X)$. Then $g_{\mu} \in \G_e(T)$, $g_d \gen g_{\mu}$ and $g_d$ is uniformly finite if and only if $\mu$ has the volume doubling property with respect to $d$, i.e. there exists $C > 0$ such that
\[
\mu(B_d(x, 2r)) \le C\mu(B_d(x, r))
\]
for any $r > 0$ and $x \in X$.\\
{\rm (2)}\,\,
For $d \in \D_{A, e}(X)$ and $\rho \in \D(X)$, $d$ is quasisymmetric to $\rho$ if and only if $\rho \in \D_{A, e}(X)$ and $g_d \gen g_{\rho}$.
\endthm

The statement (1) of the above theorem follows from Proposition~\ref{VDP.prop00} and Theorem~\ref{VDP.thm20}-(2). Note that the condition (TH) implies (TH1) appearing in Theorem~\ref{AAA.thm100}. Consequently every exponential weight function is thick by Theorem~\ref{AAA.thm100}. The statement (2) is immediate from Corollary~\ref{QSY.cor10}.\par
In \cite[Section 17]{BonkMeyer}, the authors have shown that the visual metric is quasisymmetric to the chordal metric which is bi-Lipschitz equivalent to the standard geodesic metric $d_S$ on $S^2$ for certain class of expanding Thurston maps. In view of their proof, they have essentially shown a counterpart of the condition given in Theorem~\ref{SMR.thm30}-(2).\par
Next we present a characterization of the Ahlfors regular conformal dimension using the critical index $p$ of $p$-energies.

\definition\label{SMR.def200}
Let $g \in \G_e(T)$ and let $r \in (0, 1)$. For $A \subseteq \LL^g_{r^m}$, $w \in \LL^g_{r^m}$, $M \ge 1$ and $n \ge 0$, define
\[
S^n(A) = \{v | v \in \LL^g_{r^{m + n}}, X_v \subseteq \cup_{w \in A} X_w\}
\]
and
\begin{multline*}
\GG^g_M(w) = \{v| v \in \LL^g_{r^m}, \text{there exists $(v(0), v(1), \ldots, v(M))$ such that}\\
\text{$v(0) = w$ and $(v(i), v(i + 1)) \in E^h_{g, r}$ for any $i = 0, \ldots, M - 1$}\}
\end{multline*}
\enddefinition
The set $S^n(A)$ corresponds to the refinement of $A$ in $\LL^g_{r^{m + n}}$ and
the set $\GG^g_{M}(w)$ is the $M$-neighborhood of $w$ in $\LL^g_{r^m}$. 
\definition\label{SMR.def210}
Let $g \in \G_e(T)$ and let $r \in (0, 1)$. For $p > 0$, $w \in \tT^{g, r}$, $M \ge 1$ and $n \ge 0$, define
\begin{multline*}
\E^g_{M, p, w, n} = \inf\Big\{\sum_{(u, v) \in E^h_{g, r}, u, v \in \LL^g_{r^{m + n}}}|f(u) - f(v)|^p\Big| \\
f: \LL^g_{r^{m + n}} \to \BbR, f|_{S^n(w)} = 1, u|_{\sd{\LL^g_{r^{m + n}}}{S^n(\GG^g_{M}(w)}) }= 0\Big\}
\end{multline*}
and
\[
\E^g_{M, p} = \liminf_{m \to \infty} \sup_{w \in \tT^{g, r}} \E^g_{M, p, w, m}.
\]
\enddefinition

By Theorem~\ref{PEN.thm10}, we have the following characterization of the Ahlfors regular conformal dimension of $(X, d)$ in terms of $\E^g_{M, p}$.

\thm\label{SMR.thm40}
Let $d \in \D_{A, e}(X)$ and set $g = g_d$. Assume that $d$ is uniformly finite and $M$-adapted. If $\E^g_{M, p} = 0$, then there exist $\rho \in \D_{A, e}(X)$ and $\mu \in \M_e(X)$ such that $\mu$ is $p$-Ahlfors regular with respect to $\rho$ and $\rho$ is quasisymmetric to $d$. Moreover, the Ahlfors regular conformal dimension of $(X, d)$ is (finite and) given by $\inf\{p| \E^g_{M, p} = 0\}$.
\endthm

\part{Partitions, weight functions and their hyperbolicity}

\section{Tree with a reference point}\label{TWR}
In this section, we review basic notions and notations on a tree with a reference point.

\definition\label{TWR.def10}
Let $T$ be a countably infinite set and let $\A : T \times T \to \{0, 1\}$ which satisfies $\A(w, v) = \A(v, w)$ and $\A(w, w) = 0$ for any $w, v \in T$. We call the pair $(T, \A)$ a (non-directed) graph with the vertices $T$ and the adjacent matrix $\A$. An element $(u, v) \in T \times T$ is called a edge of $(T, \A)$ if and only if $\A(u, v) = 1$. We will identify the adjacent matrix $\A$ with the collection of edges $\{(u, v)| u, v \in T, \A(u, v) = 1\}$.\newline
(1)\,\,\,
The set $\{v | \A(w, v) = 1\}$ is called the neighborhood of $w$ in $(T, \A)$. $(T, \A)$ is said to be locally finite if the neighborhood of $w$ is a finite set for any $w \in T$.\\
(2)\,\,\,For $w_0, \ldots, w_n \in T$, $(w_0, w_1, \ldots, w_n)$ is called a path between $w_0$ and $w_n$ if $\A(w_i, w_{i + 1}) = 1$ for any $i = 0, 1, \ldots n - 1$. A path $(w_0, w_1, \ldots, w_n)$ is called simple if and only if $w_i \neq w_j$ for any $i, j$ with $0 \le i < j \le n$ and $|i - j| < n$. \\
(3)\,\,\,
$(T, \A)$ is called a (non-directed) tree if and only if there exists a unique simple path between $w$ and $v$ for any $w, v \in T$ with $w \neq v$. For a tree $(T, \A)$, the unique simple path between two vertices $w$ and $v$ is called the geodesic between $w$ and $v$ and denoted by $\overline{wv}$. We write $u \in \overline{wv}$ if $\overline{wv} = (w_0, w_1, \ldots, w_n)$ and $u = w_i$ for some $i$.
\enddefinition

In this paper, we always fix a point in a tree as the root of the tree and call the point the reference point.

\definition\label{TWR.def20}
Let $(T, \A)$ be a tree and let $\phi \in T$. The triple $(T, \A, \phi)$ is called a tree with a reference point $\phi$. \\
(1)\,\,
Define $\pi:T \to T$ by
\[
\pi(w) = \begin{cases}
w_{n - 1} &\text{if $w \neq \phi$ and $\overline{\phi{w}} = (w_0, w_1, \ldots, w_{n - 1}, w_n)$,}\\
\phi &\text{if $w = \phi$}
\end{cases}
\]
 and set $S(w) = \sd{\{v| \A(w, v) = 1\}}{\{\pi(w)\}}$.\\
(2)\,\,
For $w \in T$, we define $|w| = n$ if and only if $\overline{\phi{w}} = (w_0, w_1, \ldots, w_n).$ Moreover, we set $(T)_m = \{w| w \in T, |w| = m\}$.\\
(4)\,\,
An infinite sequence of vertices $(w_0, w_1, \ldots)$ is called an infinite geodesic ray originated from $w_0$ if and only if $(w_0, \ldots, w_n) = \overline{w_0w_n}$ for any $n \ge 0$. Two infinite geodesic rays $(w_0, w_1, \ldots)$ and $(v_0, v_1, \ldots)$ are equivalent if and only if there exists $k \in \BbZ$ such that $w_{n + k} = v_n$ for sufficiently large $n$. An equivalent class of infinite geodesic rays is called an end of $T$. We use $\SS$ to denote the collection of ends of $T$. \\
(5)\,\,
Define $\SS^w$ as the collection of infinite geodesic rays originated from $w \in T$. For any $v \in T$, $\SS^w_v$ is defined as the collection of elements of $\SS^w$ passing through $v$, namely
\[
 \SS^w_v = \{(w, w_1, \ldots)|(w,w_1, \ldots) \in \SS^w, w_n = v\,\,\text{for some $n \ge 1$}\}
\]
\enddefinition

\remark
Strictly, the notations like $\pi$ and $|\cdot|$ should be written as $\pi^{(T, \A, \phi)}$ and $|\cdot|_{(T, \A, \phi)}$ respectively. In fact,  if we will need to specify the tree in question, we are going to use such explicit notations.
\endremark

One of the typical examples of a tree is the infinite binary tree. In the next example, we present a class of trees where $\#(S(w))$ is independent of $w \in T$.

\example\label{TWR.ex10}
Let $N \ge 2$ be an integer. Let $T_m^{(N)} = \{\num N\}^m$ for $m \ge 0$. (We let $T_0^{(N)} = \{\phi\}$, where $\phi$ represents an empty sequence.) We customarily write $(i_1, \ldots, i_m) \in T_m^{(N)}$ as $\word  im$. Define $T^{(N)} = \cup_{m \ge 0} T_m^{(N)}$. Define $\pi: T^{(N)} \to T^{(N)}$ by $\pi(\word i{m}i_{m + 1}) = \word im$ for $m \ge 0$ and $\pi(\phi) = \phi$. Furthermore, define
\[
\A^{(N)}_{w v} = \begin{cases} 1 \quad&\text{if $w \neq v$, and either $\pi(w) = v$ or $\pi(v) = w$,}\\
0 \quad&\text{otherwise.}
\end{cases}
\]
Then $(T^{(N)}, \A^{(N)}, \phi)$ is a locally finite tree with a reference point $\phi$. In particular, $(T^{(2)}, \A^{(2)}, \phi)$ is called the infinite binary tree. 
\endexample

It is easy to see that for any infinite geodesic ray $(w_0, w_1, \ldots)$, there exists a geodesic ray originated from $\phi$ that is equivalent to $(w_0, w_1, \ldots)$. In fact, adding the geodesic $\overline{\phi{w_0}}$ to $(w_0, w_1, \ldots)$ and removing a loop, one can obtain the infinite geodesic ray having required property. This fact shows the following proposition.

\prop\label{TWR.prop10}
There exists a natural bijective map from $\SS$ to $\SS^{\phi}$.
\endprop

Through this map, we always identify the collection of ends $\SS$ and the collection of infinite geodesic rays originated from $\phi$, $\SS^{\phi}$.\par
Hereafter in this paper, we always assume that $(T, \A)$ is a locally finite with a fixed reference point $\phi \in T$. If no confusion can occur, we omit $\phi$ in the notations. For example, we use $\SS$, and $\SS_v$ in place of $\SS^{\phi}$ and $\SS^{\phi}_v$ respectively.

\example\label{TWR.ex20}
Let $N \ge 2$ be an integer. In the case of $(T^{(N)}, \A^{(N)}, \phi)$ defined in Example~\ref{TWR.ex10}, the collection of the ends $\SS$ is $\SS^{(N)} = \{\num N\}^{\BbN} = \{i_1i_2i_3\ldots, | i_j \in \{\num N\}\,\,\text{for any $m \in \BbN$}\}$. With the natural product topology, $\SS^{(N)}$ is a Cantor set, i.e. perfect and totally disconnected.\endexample

\definition\label{TWR.def30}
Let $(T, \A, \phi)$ be a locally finite tree with a reference point $\phi$.\\
(1)\,\,
For $\omega = (w_0, w_1, \ldots) \in \SS$, we define $[\omega]_m$ by $[\omega]_m = w_m$ for any $m \ge 0$. Moreover, let $w \in T$. If $\overline{\phi{w}} = (w_0, w_1, \ldots, w_{|w|})$, then for any $0 \le m \le |w|$, we define $[w]_m = w_m$. For $w \in T$, we define
\[
T_w = \{v| v \in T, w \in \overline{\phi{v}}\}
\]
\\
(2)\,\,
For $w, v \in T$, we define the confluence of $w$ and $v$, $w \wedge v$, by 
\[
w \wedge v = w_{\max\{i| i = 0, \ldots, |w|,  [w]_i = [v]_i\}} 
\]
(3)\,\,
For $\omega, \tau \in \SS$, if $\omega \neq \tau$, we define the confluence of $\omega$ and $\tau$, $\omega \wedge \tau$, by
\[
\omega \wedge \tau = [\omega]_{\max\{m| [\omega]_m = [\tau]_m\}}
\]
(4)\,\,
For $\omega, \tau \in \SS$, we define $\rho_*(\omega, \tau) \ge 0$ by
\[
\rho_*(\omega, \tau) = \begin{cases}
2^{-|\omega \wedge \tau|}\quad&\text{if $\omega \neq \tau$,}\\
\,\,0\quad&\text{if $\omega = \tau$.}
\end{cases}
\]
\enddefinition

It is easy to see that $\rho_*$ is a metric on $\SS$ and $\{\SS_{[\omega]_m}\}_{m \ge 0}$ is a fundamental system of neighborhood of $\omega \in \SS$. Moreover, $\{\SS_v\}_{v \in T}$ is a countable base of open sets. This base of open sets has the following property.

\lemma\label{TWR.lemma10}
Let $(T, \A, \phi)$ be a locally finite tree with a reference point $\phi$. Then for any $w, v \in T$, $\SS_w \cap \SS_v = \emptyset$ if and only if $|w \wedge v| < |w|$ and $|w \wedge v| < |v|$. Furthermore, $\SS_w \cap \SS_v \neq \emptyset$ if and only if $\SS_v \subseteq \SS_w$ or $\SS_w \subseteq \SS_v$.
\endlemma

\demo
If $|w \wedge v| = |w|$, then $w = w \wedge v$ and hence $w \in \overline{\phi{v}}$. Therefore $\SS_v \subseteq \SS_w$. So, $\SS_w \cap \SS_v \neq \emptyset$. Conversely, if $\omega \in \SS_w \cap \SS_v$, then there exist $m, n \ge 0$ such that $w = [\omega]_m$ and $v = [\omega]_m$. It follows that 
\[
w \wedge v = \begin{cases} w \quad&\text{if $m \le n$,}\\
v\quad& \text{if $m \le n$.}
\end{cases}
\]
Hence we see that $|w \wedge v| = |w|$ or $|w \wedge v| = |v|$.
\enddemo

With the help to the above proposition, we may easily verify the following well-known fact. The proof is standard and left to the readers.

\prop\label{TWR.prop20}
If $(T, \A, \phi)$ is a locally finite tree with a reference point $\phi$. Then $\rho_*(\cdot, \cdot)$ is a metric on $\SS$ and $(\SS, \rho)$ is compact and totally disconnected. Moreover, if $\#(S(w)) \ge 2$ for any $w \in T$, then $(\SS, \rho)$ is perfect.
\endprop

By the above proposition, if $\#(S(w)) \ge 2$ for any $w \in T$, then $\SS$ is (homeomorphic to) the Cantor set. 

\section{Partition}\label{PAS}

In this section, we formulate the notion of a partition introduced in Section~\ref{INT} exactly.  A partition is a map from a tree to the collection of nonempty compact subsets of a compact metrizable topological space with no isolated point and it is required to preserve natural hierarchical structure of the tree. Consequently, a partition induces a surjective map from the Cantor set, i.e. the collection of ends of the tree, to the compact metrizable space. \par
Throughout this section, $\T = (T, A, \phi)$ is a locally finite tree with a reference point $\phi$. 

\definition[Partition]\label{PAS.def20}
Let $(X, \O)$ be a compact metrizable topological space having no isolated point, where $\O$ is the collection of open sets,  and let $\C(X, \O)$ be the collection of nonempty compact subsets of $X$. If no confusion can occur, we write $\C(X)$ in place of $\C(X, \O)$.\\
(1)\,\,
A map $K: T \to \C(X, \O)$, where we customarily denote $K(w)$ by $K_w$ for simplicity, is called a partition of $X$ parametrized by $(T, \A, \phi)$ if and only if it satisfies the following conditions (P1) and (P2), which correspond to \eqref{INT.eq30} and \eqref{INT.eq40} respectively. \newline
(P1)\,\,\,
$K_{\phi} = X$ and for any $w \in T$,
\[
K_w = \bigcup_{v \in S(w)} K_v.
\]
(P2)\,\,\,
For any $\omega \in \SS$, $\cap_{m \ge 0}{K_{[\omega]_m}}$ is a single point.\newline
(2)\,\,
Let $K:T \to \C(X, \O)$ be a partition of $X$ parametrized by $(T, \A, \phi)$. Define $O_w$ and $B_w$ for $w \in T$ by
\begin{align*}
O_w &= \sd{K_w}{\Bigg(\bigcup_{v \in \sd{(T)_{|w|}}{\{w\}}} K_v\Bigg)},\\
B_w &= K_w \cap \Bigg(\bigcup_{v \in \sd{(T)_{|w|}}{\{w\}}} K_v\Bigg).
\end{align*}
If $O_w \neq \emptyset$ for any $w \in T$, then the partition $K$ is called minimal.\\
(3)\,\,
Let $K : T \to \C(X, \O)$ be a partition of $X$. Then $(w(1), \ldots, w(m)) \in \cup_{k \ge 0}T^k$ is called a chain of $K$ (or a chain for short if no confusion can occur) if and only if $K_{w(i)} \cap K_{w(i + 1)} \neq \emptyset$ for any $i = 1, \ldots, m - 1$. A chain $(w(1), \ldots, w(m))$ of $K$ is called a chain of $K$ in $\LL \subseteq T$ if $w(i) \in \LL$ for any $i = 1, \ldots, m$. For subsets $A, B \subseteq X$, A chain $(w(1), \ldots, w(m))$ of $K$ is called a chain of $K$ between $A$ and $B$ if and only if $A \cap K_{w(1)} \neq \emptyset$ and $B \cap K_{w(m)} \neq \emptyset$. We use $\CH_K(A, B)$ to denote the collection of chains of $K$ between $A$ and $B$. Moreover, we denote the collections of chains of $K$ in $\LL$ between $A$ and $B$ by $\CH_K^{\LL}(A, B)$.
\enddefinition

As is shown in Theorem~\ref{PAS.thm00}, a partition can be modified so as to be minimal by restricting it to a suitable subtree.\par
The next lemma is an assortment of direct consequences from the definition of the partition.

\lemma\label{PAS.lemma20}
Let $K : T \to \C(X, \O)$ be a partition of $X$ parametrized by $(T, A, \phi)$. \newline
{\rm (1)}\,\,\,
For any $w \in T$, $O_w$ is an open set. $O_v \subseteq O_w$ for any $v \in S(w)$.\newline
{\rm (2)}\,\,\,
$O_w \cap O_v = \emptyset$ if $w, v \in T$ and $\SS_w \cap \SS_v = \emptyset$. 
\newline
{\rm (3)}\,\,\,
If $\SS_w \cap \SS_v = \emptyset$, then $K_w \cap K_v = B_w \cap B_v$.
\endlemma

\demo
(1)\,\,
Note that by (P1), $X = \cup_{w \in (T)_m} K_w$. Hence 
\[
O_w = \sd{K_w}{(\cup_{v \in \sd{(T)_{|w|}}{\{w\}}} K_v)} = \sd{X}{(\cup_{v \in \sd{(T)_{|w|}}{\{w\}}} K_v)}.
\]
The rest of the statement is immediate from the property (P2).\\
(2)\,\,
By Lemma~\ref{TWR.lemma10}, if $u = w \wedge v$, then $|u| < |w|$ and $|u| < |v|$. Let $w' = [w]_{|u| + 1}$ and let $v' = [v]_{|u| + 1}$. Then $w', v' \in S(u)$ and $w' \neq v'$. Since $O_{w'} \subseteq \sd{K_{w'}}{K_{v'}}$, it follows that $O_{w'} \cap O_{v'} = \emptyset$. Using (1), we see $O_w \cap O_v = \emptyset$.\\
(3)\,\,
This follows immediately by (1).
\enddemo

The condition (P2) provides a natural map from the ends of the tree $\SS$ to the space $X$.

\prop\label{PAS.prop10}
Let $K: T \to \C(X, \O)$ be a partition of $X$ parametrized by $(T, \A, \phi)$.\\
{\rm (1)}\,\,\,For $\omega \in \SS$, define $\s(\omega)$ as the single point $\cap_{m \ge 0}{K_{[\omega]_m}}$. Then $\s : \SS \to X$ is continuous and surjective. Moreover. $\s(\SS_w) = K_w$ for any $w \in T$.\newline
{\rm (2)}\,\,\,
The partition $K: T \to \C(X, \O)$ is minimal if and only if $K_w$ is the closure of $O_w$ for any $w \in T$. Moreover, if $K: T \to \C(X, \O)$ is minimal then $O_w$ coincides with the interior of $K_w$.
\endprop

\demo
(1)\,\,
Note that $K_w = \cup_{v \in S(w)} K_v$. Hence if $x \in K_w$, then there exists $v \in S(w)$ such that $x \in K_v$. Using this fact inductively, we see that, for any $x \in X$, there exists $\omega \in \SS$ such that $x \in K_{[\omega]_m}$ for any $m \ge 0$. Since $x \in \cap_{m \ge 0} K_{[\omega]_m}$, (P2) shows that $\s(\omega) = x$. Hence $\omega$ is surjective. At the same time, it follows that $\s(\SS_w) = K_w$. Let $U$ be an open set in $X$. For any $\omega \in \s^{-1}(U)$, $K_{[\omega]_m} \subseteq U$ for sufficiently large $m$. Then $\SS_{[\omega]_m} \subseteq \s^{-1}(U)$. This shows that $\s^{-1}(U)$ is an open set and hence $\s$ is continuous.\newline
(2)\,\,
Let $\overline{O}_w$ be the closure of $O_w$. If $K_w = \overline{O}_w$ for any $w \in T$, then $O_w \neq \emptyset$ for any $w \in T$ and hence $K: T \to \C(X, \O)$ is minimal. Conversely, assume that $K:T \to \C(X, \O)$ is minimal.  By Lemma~\ref{PAS.lemma20}, $\overline{O}_{[\omega]_m} \supseteq \overline{O}_{[\omega]_{m + 1}}$ for any $\omega \in \SS$ and any $m \ge 0$. Hence $\{\s(\omega)\} = \cap_{m \ge 0} K_{[\omega]_m} = \cap_{m \ge 0} \overline{O}_{[\omega]_m} \subseteq \overline{O}_{[\omega]_n}$ for any $n \ge 0$. This yields that $\s(\SS_w) \subseteq \overline{O}_w$. Since $\s(\SS_w) = K_w$, this implies $\overline{O}_w = K_w$.\\
Now if $K$ is minimal, since $O_w$ is open by Lemma~\ref{PAS.lemma20}-(1), it follows that $O_w$ is the interior of $K_w$.
\enddemo

\definition\label{PAS.def200}
 A partition $K: T \to \C(X, \O)$ parametrized by a tree $(T, \A, \phi)$ is called strongly finite if and only if 
\[
\sup_{x \in X} \#(\s^{-1}(x)) < +\infty,
\]
where $\s: \SS \to X$ is the map defined in Proposition~\ref{PAS.prop10}-(1).
\enddefinition

\example\label{PAS.ex10}
Let $(Y, d)$ be a complete metric space and let $\{F_1, \ldots, F_N\}$ be collection of contractions from $(Y, d)$ to itself, i.e. $F_i: Y \to Y$ and
\[
\sup_{x \neq y \in Y} \frac{d(F_i(x), F_i(y))}{d(x, y)} < 1
\]
for any $i = 1, \ldots, N$. Then it is well-known that there exists a unique nonempty compact set $X$ such that
\[
X = \bigcup_{i = 1, \ldots, N} F_i(X).
\]
See \cite[Section~1.1]{AOF} for a proof of this fact for example. $X$ is called the self-similar set associated with $\{F_1, \ldots, F_N\}$. Let $(T^{(N)}, \A^{(N)}, \phi)$ be the tree defined in Example~\ref{TWR.ex10}. For any $\word im \in T$, set $F_{\word im} = F_{i_1}\circ\ldots\circ{F_{i_m}}$ and define $K_w = F_w(X)$. Then $K: T^{(N)} \to \C(X)$ is a partition of $K$ parametrized by $(T^{(N)}, \A^{(N)}, \phi)$. See \cite[Section~1.2]{AOF}. The associated map from $\SS = \{1, \ldots, N\}^{\BbN}$ to $K$ is sometimes called the coding map. To determine if $K$ is minimal or not is known to be rather delicate issue. See \cite[Theorem~1.3.8]{AOF} for example.
\endexample

\begin{figure}
\centering
\setlength{\unitlength}{30mm}

\begin{picture}(4, 2.2)(-1.1, -0.1)
\drawline (-1, 0)(0.8, 0)(0.8, 1.8)(-1, 1.8)(-1, 0)
\drawline(-0.4, 0)(-0.4, 1.8)
\drawline(0.2, 0)(0.2, 1.8)
\drawline(-1, 0.6)(0.8, 0.6)
\drawline(-1, 1.2)(0.8, 1.2)
\multiput(0,0)(0, 0.9){3}{\put(-1,0){\makebox(0,0)[b]{\circle*{0.08}}}}
\put(-1, -0.1){\makebox(0,0)[t]{\large$p_1$}}
\put(-0.1, -0.1){\makebox(0,0)[t]{\large$p_2$}}
\put(0.9, -0.1){\makebox(0,0)[t]{\large$p_3$}}
\put(0.9, 0.9){\makebox(0,0)[l]{\large$p_4$}}
\put(0.9, 1.9){\makebox(0,0)[l]{\large$p_5$}}
\put(-0.1, 1.9){\makebox(0,0)[b]{\large$p_6$}}
\put(-1, 1.9){\makebox(0,0)[r]{\large$p_7$}}
\put(-1.15, 0.9){\makebox(0,0){\large$p_8$}}
\multiput(0,0)(0.9, 0){3}{\put(-1,0){\makebox(0,0)[b]{\circle*{0.08}}}}
\put(0.8,0.9){\makebox(0,0)[b]{\circle*{0.08}}}
\put(0.8,1.8){\makebox(0,0)[b]{\circle*{0.08}}}
\put(-0.1,1.8){\makebox(0,0)[b]{\circle*{0.08}}}
\shade\path(-0.4, 0.6)(0.2,0.6)(0.2,1.2)(-0.4,1.2)(-0.4,0.6)
\put(-0.7,0.3){\makebox(0,0){\large$ 1$}}
\put(-0.1,0.3){\makebox(0,0){\large$ 2$}}
\put(0.5,0.3){\makebox(0,0){\large$ 3$}}
\put(0.5,0.9){\makebox(0,0){\large$ 4$}}
\put(0.5,1.5){\makebox(0,0){\large$ 5$}}
\put(-0.1,1.5){\makebox(0,0){\large$ 6$}}
\put(-0.7,1.5){\makebox(0,0){\large$ 7$}}
\put(-0.7,0.9){\makebox(0,0){\large$ 8$}}

\thinlines

\drawline (1.2, 0)(3, 0)(3, 1.8)(1.2, 1.8)(1.2, 0)
\drawline(1.8, 0)(1.8, 1.8)
\drawline(1.4,0)(1.4, 1.8)
\drawline(1.6,0)(1.6, 1.8)
\drawline(2.0,0)(2.0, 0.6)
\drawline(2.2,0)(2.2, 0.6)
\drawline(2.0,1.2)(2.0, 1.8)
\drawline(2.2,1.2)(2.2, 1.8)
\drawline(2.4, 0)(2.4, 1.8)
\drawline(2.6, 0)(2.6, 1.8)
\drawline(2.8, 0)(2.8, 1.8)
\drawline(1.2, 0.2)(3, 0.2)
\drawline(1.2, 0.4)(3, 0.4)
\drawline(1.2, 0.6)(3, 0.6)
\drawline(1.2,0.8)(1.8,0.8)
\drawline(2.4,0.8)(3,0.8)
\drawline(1.2,1)(1.8,1)
\drawline(2.4,1)(3,1)
\drawline(1.2, 1.2)(3, 1.2)
\drawline(1.2, 1.4)(3, 1.4)
\drawline(1.2, 1.6)(3, 1.6)
\thinlines
\shade\path(1.8,0.6)(2.4,0.6)(2.4,1.2)(1.8,1.2)(1.8,0.6)
\shade\path(1.4,0.2)(1.6,0.2)(1.6,0.4)(1.4, 0.4)(1.4, 0.2)
\shade\path(2,0.2)(2.2,0.2)(2.2,0.4)(2, 0.4)(2, 0.2)
\shade\path(2.6,0.2)(2.8,0.2)(2.8,0.4)(2.6, 0.4)(2.6, 0.2)
\shade\path(1.4,0.8)(1.6,0.8)(1.6,1)(1.4, 1)(1.4, 0.8)
\shade\path(1.4,1.4)(1.6,1.4)(1.6,1.6)(1.4, 1.6)(1.4, 1.4)
\shade\path(2,1.4)(2.2,1.4)(2.2,1.6)(2, 1.6)(2, 1.4)
\shade\path(2.6,1.4)(2.8,1.4)(2.8,1.6)(2.6, 1.6)(2.6, 1.4)
\shade\path(2.6,0.8)(2.8,0.8)(2.8,1)(2.6, 1)(2.6, 0.8)
\put(1.3,0.1){\makebox(0,0){$ {11}$}}
\put(1.5,0.1){\makebox(0,0){$ {12}$}}
\put(1.7,0.1){\makebox(0,0){$ {13}$}}
\put(1.7,0.3){\makebox(0,0){$ {14}$}}
\put(1.7,0.5){\makebox(0,0){$ {15}$}}
\put(1.5,0.5){\makebox(0,0){$ {16}$}}
\put(1.3,0.5){\makebox(0,0){$ {17}$}}
\put(1.3,0.3){\makebox(0,0){$ {18}$}}

\put(0.6, 0){
\put(1.3,0.1){\makebox(0,0){$ {21}$}}
\put(1.5,0.1){\makebox(0,0){$ {22}$}}
\put(1.7,0.1){\makebox(0,0){$ {23}$}}
\put(1.7,0.3){\makebox(0,0){$ {24}$}}
\put(1.7,0.5){\makebox(0,0){$ {25}$}}
\put(1.5,0.5){\makebox(0,0){$ {26}$}}
\put(1.3,0.5){\makebox(0,0){$ {27}$}}
\put(1.3,0.3){\makebox(0,0){$ {28}$}}
}
\put(1.2, 0){
\put(1.3,0.1){\makebox(0,0){$ {31}$}}
\put(1.5,0.1){\makebox(0,0){$ {32}$}}
\put(1.7,0.1){\makebox(0,0){$ {33}$}}
\put(1.7,0.3){\makebox(0,0){$ {34}$}}
\put(1.7,0.5){\makebox(0,0){$ {35}$}}
\put(1.5,0.5){\makebox(0,0){$ {36}$}}
\put(1.3,0.5){\makebox(0,0){$ {37}$}}
\put(1.3,0.3){\makebox(0,0){$ {38}$}}
}
\put(1.2, 0.6){
\put(1.3,0.1){\makebox(0,0){$ {41}$}}
\put(1.5,0.1){\makebox(0,0){$ {42}$}}
\put(1.7,0.1){\makebox(0,0){$ {43}$}}
\put(1.7,0.3){\makebox(0,0){$ {44}$}}
\put(1.7,0.5){\makebox(0,0){$ {45}$}}
\put(1.5,0.5){\makebox(0,0){$ {46}$}}
\put(1.3,0.5){\makebox(0,0){$ {47}$}}
\put(1.3,0.3){\makebox(0,0){$ {48}$}}
}
\put(1.2, 1.2){
\put(1.3,0.1){\makebox(0,0){$ {51}$}}
\put(1.5,0.1){\makebox(0,0){$ {52}$}}
\put(1.7,0.1){\makebox(0,0){$ {53}$}}
\put(1.7,0.3){\makebox(0,0){$ {54}$}}
\put(1.7,0.5){\makebox(0,0){$ {55}$}}
\put(1.5,0.5){\makebox(0,0){$ {56}$}}
\put(1.3,0.5){\makebox(0,0){$ {57}$}}
\put(1.3,0.3){\makebox(0,0){$ {58}$}}
}
\put(0.6, 1.2){
\put(1.3,0.1){\makebox(0,0){$ {61}$}}
\put(1.5,0.1){\makebox(0,0){$ {62}$}}
\put(1.7,0.1){\makebox(0,0){$ {63}$}}
\put(1.7,0.3){\makebox(0,0){$ {64}$}}
\put(1.7,0.5){\makebox(0,0){$ {65}$}}
\put(1.5,0.5){\makebox(0,0){$ {66}$}}
\put(1.3,0.5){\makebox(0,0){$ {67}$}}
\put(1.3,0.3){\makebox(0,0){$ {68}$}}
}
\put(0.0, 1.2){
\put(1.3,0.1){\makebox(0,0){$ {71}$}}
\put(1.5,0.1){\makebox(0,0){$ {72}$}}
\put(1.7,0.1){\makebox(0,0){$ {73}$}}
\put(1.7,0.3){\makebox(0,0){$ {74}$}}
\put(1.7,0.5){\makebox(0,0){$ {75}$}}
\put(1.5,0.5){\makebox(0,0){$ {76}$}}
\put(1.3,0.5){\makebox(0,0){$ {77}$}}
\put(1.3,0.3){\makebox(0,0){$ {78}$}}
}
\put(0.0, 0.6){
\put(1.3,0.1){\makebox(0,0){$ {81}$}}
\put(1.5,0.1){\makebox(0,0){$ {82}$}}
\put(1.7,0.1){\makebox(0,0){$ {83}$}}
\put(1.7,0.3){\makebox(0,0){$ {84}$}}
\put(1.7,0.5){\makebox(0,0){$ {85}$}}
\put(1.5,0.5){\makebox(0,0){$ {86}$}}
\put(1.3,0.5){\makebox(0,0){$ {87}$}}
\put(1.3,0.3){\makebox(0,0){$ {88}$}}
}

\end{picture}
\caption{Partition: the Sierpinski carpet}\label{sc1}
\end{figure}
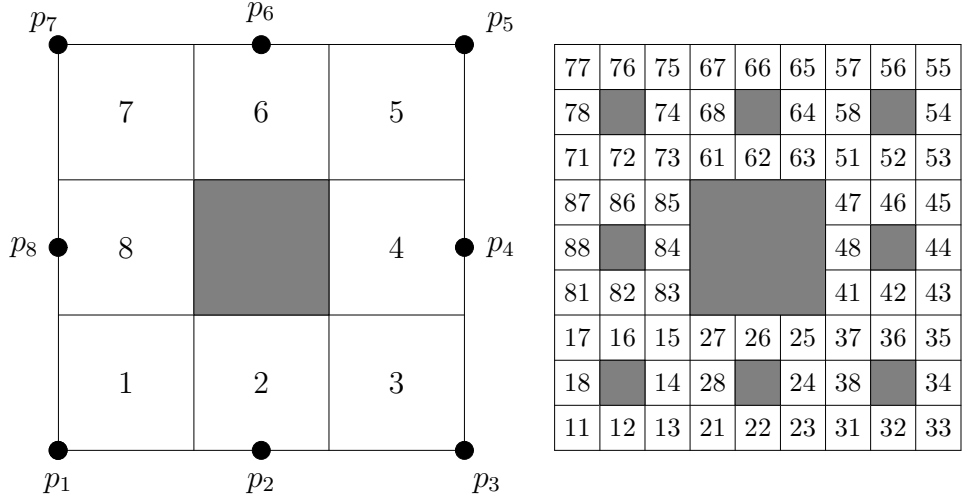

\example[the Sierpinski carpet, Figure~\ref{sc1}]\label{exsc1}
This is the special case of self-similar sets presented in the last example. Let $p_1 = (0, 0)$, $p_2 = (\frac 12, 0)$, $p_3 = (1, 0)$, $p_4 = (1, \frac 12)$, $p_5 = (1, 1)$, $p_6 = (\frac 12, 1)$, $p_7 = (0, 1)$ and $p_8 = (0, \frac 12)$. Set $F_i: [0, 1]^2 \to [0, 1]^2$ for $i = 1, \ldots, 8$ by 
\[
F_i(x) = \frac 13(x - p_i) + p_i
\]
for any $x \in [0, 1]^2$. The unique nonempty compact set $X$ satisfying 
\[
X = \bigcup_{i = 1}^8 F_i(X)
\]
is called the Sierpinski carpet. In this case, the associated tree is $(T^{(8)}, \A^{(8)}, \phi)$. Define $K: T^{(8)} \to \C(X, \O)$ by
\[
K_{i_1\ldots{i_m}} = F_{i_1\ldots{i_m}}(X)
\]
as in Example~\ref{PAS.ex10}. Then $K$ is a partition of $X$ parametrized by the tree $(T^{(8)}, \A^{(8)}, \phi)$. In Figure~\ref{sc1}, $K_{ij}$ is represented by $ij$ for simplicity.
\endexample

Removing unnecessary vertices of the tree, we can always modify the original partition and obtain a minimal one.

\thm\label{PAS.thm00}
Let $K: T \to \C(X, \O)$ be a partition of $X$ parametrized by $(T, \A, \phi)$. There exist $T' \subseteq T$ and $K': T' \to \C(X, \O)$ such that $(T', \A|_{T' \times T'})$ is a tree, $\phi \in T'$, $K'_w \subseteq K_w$ for any $w \in T'$ and $K'$ is a minimal partition of $X$ parametrized by $(T', \A', \phi)$.
\endthm

\demo
We define a sequence $\{T^{(m)}\}_{m \ge 0}$ of subsets of $T$ and $\{K^{(m)}_w\}_{w \in T^{(m)}}$ inductively as follows. First let $T^{(0)} = T$ and $K^{(0)}_w = K_w$ for any $w \in T^{(0)}$.  Suppose we have defined $T^{(m)}$. Define
\[
Q^{(m)} = \Bigg\{w\Bigg| w \in T^{(m)}, K_w \subseteq \bigcup_{v \in (T)_{|w|} \cap T^{(m)}, v \neq w} K_v\Bigg\}.
\]
If $Q^{(m)} = \emptyset$, then set $T^{(m + 1)} = T^{(m)}$ and $K^{(m + 1)}_w = K^{(m)}_w$ for any $w \in T^{(m + 1)}$. Otherwise choose $w^{(m)} \in Q^{(m)}$ so that $|w^{(m)}|$ attains the minimum of $\{|v|: v \in Q^{(m)}\}$. Then define
\[
T^{(m + 1)} = \sd{T^{(m)}}{T_{w^{(m)}}}
\]
and
\[
K^{(m + 1)}_w = \begin{cases}
\cup_{v \in T_w \cap (T)_{|w^{(m)}|} \cap T^{(m + 1)}} K^{(m)}_v&\quad\text{if $w^{(m)} \in T_w$,}\\
K^{(m)}_w&\quad\text{otherwise.}
\end{cases}
\]
In this way, for any $m \ge 0$ and $w \in T^{(m)}$,
\begin{equation}\label{PAS.eq100}
K^{(m)}_w = \bigcup_{v \in S(w) \cap T^{(m)}} K^{(m)}_v.
\end{equation}
Note that $Q^{(m + 1)} \subset \sd{Q^{(m)}}{\{w^{(m)}\}}$. Since $(T)_n$ is a finite set for any $n \ge 0$, it follows that $(T)_n \cap Q^{(m)} = \emptyset$ and $(T^{(m)})_n$ stays the same for sufficiently large $m$. Hence $|w^{(m)}| \to \infty$ as $m \to \infty$ and $(T)_n \cap T^{(m)}$ does not depend on $m$ for sufficiently large $m$. Therefore, letting $T' = \cap_{m \ge 1} T^{(m)}$, we see that $(T', \A|_{T' \times T'})$ is a locally finite tree and $\phi \in T'$. Moreover, note that $K^{(m + 1)}_w \subseteq K^{(m)}_w$ for any $w \in T'$. Hence if we set
\[
K'_w = \bigcap_{m \ge 0} K^{(m)}_w
\]
for any $w \in T'$, then $K'_w$ is nonempty. By \eqref{PAS.eq100}, it follows that
\[
K'_w = \bigcup_{v \in T' \cap S(w)} K'_v
\]
for any $w \in T'$. Thus the map $K': T' \to \C(X, \O)$ given by $K'(w) = K'_w$ is a minimal partition of $X$ parametrized by $(T', \A|_{T' \times T'}, \phi)$.
\enddemo

A partition $K: T \to \C(X, \O)$ induces natural graph structure on $T$. In the rest of this section, we show that  $T$ can be regarded as the hyperbolic filling of $X$ if the induced graph structure is hyperbolic. See \cite{BonkSaks}, for example, about the notion of hyperbolic fillings.

\definition\label{HF.def10}
Let $K: T \to \C(X, \O)$ be a partition. Then define 
\[
E_m^h = \{(w, v)| w, v \in (T)_m, K_w \cap K_v \neq \emptyset\}
\]
and
\[
E^h = \bigcup_{m \ge 0} E_m^h.
\]
An element $(u, v) \in E^h$ is called a horizontal edge associated with $(T, \A, \phi)$ and $K: T \to \C(X, \O)$. The symbol ``h'' in the notation $E_m^h$ and $E^h$ represents the word ``horizontal''. On the contrary, an element $(w, v) \in \A$ is called a vertical edge. Moreover we define
\[
\B(w, v) = \begin{cases} 
\quad1\quad&\text{if $\A(w, v) = 1$ or $(w, v) \in E^h$,}\\
\quad0\quad&\text{otherwise.}
\end{cases}
\]
 The graph $(T, \B)$ is called the resolution of $X$ associated with the partition $K: T \to \C(X, \O)$. We use $d_{(T, \B)}(\cdot, \cdot)$ to denote the shortest path metric, i.e. 
\begin{multline*}
d_{(T, \B)}(w, v) = \min\{n| \text{there exists $(w(1), \ldots, w(n + 1)) \in (\B)^n$}\\\text{$\B(w(i), w(i + 1)) = 1$ for any $i = 1, \ldots, n$}\}
\end{multline*}
\enddefinition

\remark
The horizontal graph $((T)_m, E^h_m)$ is not necessarily connected. More precisely, $((T)_m, E^h_m)$ is connected for any $m \ge 0$ if and only if $X$ is connected. Note that if $X$ is homeomorphic to the Cantor set, then $E^h_m = \emptyset$ for any $m \ge 0$.
\endremark

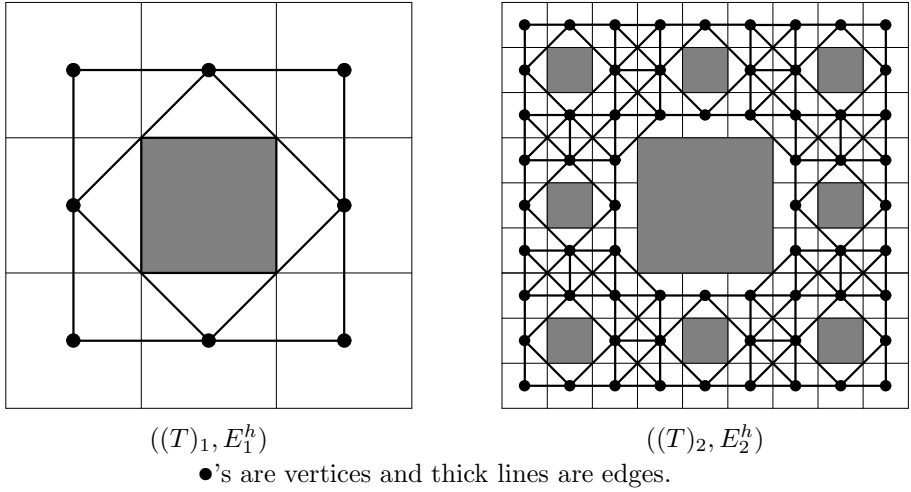
\begin{figure}
\centering
\setlength{\unitlength}{30mm}

\begin{picture}(4, 2)(-1, -0.3)
\drawline (-1, 0)(0.8, 0)(0.8, 1.8)(-1, 1.8)(-1, 0)
\drawline(-0.4, 0)(-0.4, 1.8)
\drawline(0.2, 0)(0.2, 1.8)
\drawline(-1, 0.6)(0.8, 0.6)
\drawline(-1, 1.2)(0.8, 1.2)
\put(-0.7, 0.3){\circle*{0.06}}
\put(-0.7, 0.9){\circle*{0.06}}
\put(-0.7, 1.5){\circle*{0.06}}
\put(-0.1, 0.3){\circle*{0.06}}
\put(-0.1, 1.5){\circle*{0.06}}
\put(0.5, 0.3){\circle*{0.06}}
\put(0.5, 0.9){\circle*{0.06}}
\put(0.5, 1.5){\circle*{0.06}}
\thicklines
\drawline(-0.7, 0.3)(-0.7, 1.5)(0.5, 1.5)(0.5, 0.3)(-0.7, 0.3)
\drawline(-0.7, 0.9)(-0.1, 1.5)(0.5,0.9)(-0.1, 0.3)(-0.7, 0.9)
\shade\path(-0.4, 0.6)(0.2,0.6)(0.2,1.2)(-0.4,1.2)(-0.4,0.6)
\put(-0.1, -0.2){\makebox(0,0)[b]{$((T)_1, E^h_1)$}}

\thinlines

\drawline (1.2, 0)(3, 0)(3, 1.8)(1.2, 1.8)(1.2, 0)
\drawline(1.8, 0)(1.8, 1.8)
\drawline(1.4,0)(1.4, 1.8)
\drawline(1.6,0)(1.6, 1.8)
\drawline(2.0,0)(2.0, 0.6)
\drawline(2.2,0)(2.2, 0.6)
\drawline(2.0,1.2)(2.0, 1.8)
\drawline(2.2,1.2)(2.2, 1.8)
\drawline(2.4, 0)(2.4, 1.8)
\drawline(2.6, 0)(2.6, 1.8)
\drawline(2.8, 0)(2.8, 1.8)
\drawline(1.2, 0.2)(3, 0.2)
\drawline(1.2, 0.4)(3, 0.4)
\drawline(1.2, 0.6)(3, 0.6)
\drawline(1.2, 0.6)(3, 0.6)
\drawline(1.2,0.8)(1.8,0.8)
\drawline(2.4,0.8)(3,0.8)
\drawline(1.2,1)(1.8,1)
\drawline(2.4,1)(3,1)
\drawline(1.2, 1.2)(3, 1.2)
\drawline(1.2, 1.4)(3, 1.4)
\drawline(1.2, 1.6)(3, 1.6)
\thinlines
\shade\path(1.8,0.6)(2.4,0.6)(2.4,1.2)(1.8,1.2)(1.8,0.6)
\shade\path(1.4,0.2)(1.6,0.2)(1.6,0.4)(1.4, 0.4)(1.4, 0.2)
\shade\path(2,0.2)(2.2,0.2)(2.2,0.4)(2, 0.4)(2, 0.2)
\shade\path(2.6,0.2)(2.8,0.2)(2.8,0.4)(2.6, 0.4)(2.6, 0.2)
\shade\path(1.4,0.8)(1.6,0.8)(1.6,1)(1.4, 1)(1.4, 0.8)
\shade\path(1.4,1.4)(1.6,1.4)(1.6,1.6)(1.4, 1.6)(1.4, 1.4)
\shade\path(2,1.4)(2.2,1.4)(2.2,1.6)(2, 1.6)(2, 1.4)
\shade\path(2.6,1.4)(2.8,1.4)(2.8,1.6)(2.6, 1.6)(2.6, 1.4)
\shade\path(2.6,0.8)(2.8,0.8)(2.8,1)(2.6, 1)(2.6, 0.8)

\thicklines
\multiput(1.2, 0)(0.6,0){3}{
\put(0.1,0.1){\circle*{0.04}}
\put(0.1,0.3){\circle*{0.04}}
\put(0.1,0.5){\circle*{0.04}}
\put(0.3,0.5){\circle*{0.04}}
\put(0.5,0.5){\circle*{0.04}}
\put(0.5,0.3){\circle*{0.04}}
\put(0.5,0.1){\circle*{0.04}}
\put(0.3,0.1){\circle*{0.04}}
\drawline(0.1,0.1)(0.1,0.5)(0.5,0.5)(0.5,0.1)(0.1,0.1)
\drawline(0.1,0.3)(0.3,0.5)(0.5,0.3)(0.3,0.1)(0.1,0.3)
}
\put(2.4,0.6){
\put(0.1,0.1){\circle*{0.04}}
\put(0.1,0.3){\circle*{0.04}}
\put(0.1,0.5){\circle*{0.04}}
\put(0.3,0.5){\circle*{0.04}}
\put(0.5,0.5){\circle*{0.04}}
\put(0.5,0.3){\circle*{0.04}}
\put(0.5,0.1){\circle*{0.04}}
\put(0.3,0.1){\circle*{0.04}}
\drawline(0.1,0.1)(0.1,0.5)(0.5,0.5)(0.5,0.1)(0.1,0.1)
\drawline(0.1,0.3)(0.3,0.5)(0.5,0.3)(0.3,0.1)(0.1,0.3)
}
\put(2.4,1.2){
\put(0.1,0.1){\circle*{0.04}}
\put(0.1,0.3){\circle*{0.04}}
\put(0.1,0.5){\circle*{0.04}}
\put(0.3,0.5){\circle*{0.04}}
\put(0.5,0.5){\circle*{0.04}}
\put(0.5,0.3){\circle*{0.04}}
\put(0.5,0.1){\circle*{0.04}}
\put(0.3,0.1){\circle*{0.04}}
\drawline(0.1,0.1)(0.1,0.5)(0.5,0.5)(0.5,0.1)(0.1,0.1)
\drawline(0.1,0.3)(0.3,0.5)(0.5,0.3)(0.3,0.1)(0.1,0.3)
}
\put(1.8,1.2){
\put(0.1,0.1){\circle*{0.04}}
\put(0.1,0.3){\circle*{0.04}}
\put(0.1,0.5){\circle*{0.04}}
\put(0.3,0.5){\circle*{0.04}}
\put(0.5,0.5){\circle*{0.04}}
\put(0.5,0.3){\circle*{0.04}}
\put(0.5,0.1){\circle*{0.04}}
\put(0.3,0.1){\circle*{0.04}}
\drawline(0.1,0.1)(0.1,0.5)(0.5,0.5)(0.5,0.1)(0.1,0.1)
\drawline(0.1,0.3)(0.3,0.5)(0.5,0.3)(0.3,0.1)(0.1,0.3)
}
\multiput(1.2,0.6)(0,0.6){2}{
\put(0.1,0.1){\circle*{0.04}}
\put(0.1,0.3){\circle*{0.04}}
\put(0.1,0.5){\circle*{0.04}}
\put(0.3,0.5){\circle*{0.04}}
\put(0.5,0.5){\circle*{0.04}}
\put(0.5,0.3){\circle*{0.04}}
\put(0.5,0.1){\circle*{0.04}}
\put(0.3,0.1){\circle*{0.04}}
\drawline(0.1,0.1)(0.1,0.5)(0.5,0.5)(0.5,0.1)(0.1,0.1)
\drawline(0.1,0.3)(0.3,0.5)(0.5,0.3)(0.3,0.1)(0.1,0.3)
}
\multiput(0, 0)(1.2,0){2}
{
\multiput(0, 0)(0, 0.6){2}
{
{\multiput(0,0)(0.2, 0){3}{\drawline(1.3, 0.5)(1.3,0.7)}
\multiput(0,0)(0.2,0){2}{\drawline(1.3,0.5)(1.5,0.7)}
\multiput(0,0)(0.2,0){2}{\drawline(1.3,0.7)(1.5,0.5)}}
}
}
\multiput(0,0)(0.6,0){2}
{
\multiput(0,0)(0,1.2){2}
{
\multiput(0,0)(0,0.2){3}{\drawline(1.7,0.1)(1.9,0.1)}
\multiput(0,0)(0,0.2){2}{\drawline(1.7,0.1)(1.9,0.3)}
\multiput(0,0)(0,0.2){2}{\drawline(1.9,0.1)(1.7,0.3)}
}
}
\multiput(0, 0)(0.6,0.6){2}{\drawline(1.9,0.5)(1.7,0.7)}
\multiput(0, 0)(0.6,-0.6){2}{\drawline(1.7,1.1)(1.9,1.3)}
\put(2.1, -0.2){\makebox(0,0)[b]{$((T)_2, E^h_2)$}}
\put(0.9, -0.3){\makebox(0, 0){{\large $\bullet$}'s are vertices and thick lines are edges.}}
\end{picture}
\caption{Horizontal edges: the Sierpinski carpet}\label{sc2}
\end{figure}

In Figure~\ref{sc2}, we present $((T)_1, E^h_1)$ and $((T)_2, E^h_2)$ for the Sierpinski carpet introduced in Example~\ref{exsc1}.

We will show in Lemma~\ref{HYP.lemma20} that if $(\phi, w(1), w(2), \ldots)$ is an infinite geodesic ray in $(T, \B)$ with respect to the metric $d_{(T, \B)}$ starting from $\phi$, then  it coincides with $(\phi, [\omega]_1, [\omega]_2, \ldots)$ for some $\omega \in \SS$. In other words, the collection of geodesic rays of $(T, \B)$ starting from $\phi$ can be identified with $\SS$. The following proposition will be proven in Section~\ref{HYP}.

\prop\label{HF.prop10}
Let $\omega, \tau \in \SS$. If $\sup_{n \ge 1} d_{(T, \B)}([\omega]_n, [\tau]_n) < +\infty$, then $\s(\omega) = \s(\tau)$.
\endprop

By this proposition, whether the resolution $(T, \B)$ is hyperbolic or not, $X$ can be identified with the quotient space of the geodesic rays  under the equivalence relation $\sim$ defined as $\omega \sim \tau$ if and only if $\sup_{n \ge 1} d_{(T, \B)}([\omega]_n, [\tau]_n) < +\infty$. In case of a hyperbolic graph, such a quotient space has been called the hyperbolic boundary of the graph in the framework of Gromov theory of hyperbolic metric spaces. We will give detailed accounts on these points later in Section~\ref{HYP}.\par
In \cite{Elek1}, Elek has constructed a hyperbolic graph whose hyperbolic boundary is homeomorphic to a given compact subset of $\BbR^N$. From our point of view, what he has done is to construct a partition of the compact metric space using dyadic cubes as is seen in the next example. However, the resolution $(T, \B)$ associated with the partition is slightly different from the original graph constructed by Elek. See the details below.
 
\example\label{PAS.ex20}
 Let $X$ be a nonempty compact subset of $\BbR^N$. For simplicity, we assume that $X \subseteq [0, 1]^N$. We are going to construct a partition of $X$ using the dyadic cubes. Let $S_m =\{(m, i_1, \ldots, i_N)| (i_1, \ldots, i_N) \in  \{0, 1, \ldots, 2^m - 1\}^N\}$ and define
\[
C(w) = \prod_{j = 1}^N \Big[\frac {i_j}{2^m}, \frac{i_j + 1}{2^m}\Big]
\]
for $w = (m, i_1, \ldots, i_N) \in S_m$. The collection $\{C(w)| w \in \cup_{m \ge 0} S_m\}$ is called the dyadic cubes. (See \cite{Christ} for example.) Define 
\[
K_{w} = X \cap C(w)
\]
for $w \in \cup_{m \ge 0} S_m$, 
\[
(T)_m = \{w| w \in S_m, K_w \neq \emptyset\}
\]
for $m \ge 0$ and $T = \cup_{m \ge 0} (T)_m$. Moreover, we define $(w, v)  \in \A$ for $(w, v) \in T \times T$ if there exist $m \ge 0$ such that $(w, v) \in (S_m \times S_{m + 1}) \cup (S_{m + 1} \times S_m)$ and $C(w) \supseteq C(v)$ or $C(w) \subseteq C(v)$. Then $(T, \A, \phi)$ is a tree with a reference point $\phi$, where $\phi = (0, 0, \ldots, 0) \in (T)_0$ and the map from $w \in T$ to $K_w \in \C(X, \O)$ is a partition of $X$ parametrized by $(T, \A, \phi)$. The hyperbolic graph constructed by Elek is a slight modification of the resolution $(T, \B)$. In fact, the vertical edges are the same but Elek's graph has more horizontal edges. Precisely set
\[
\widetilde{E}^h_m = \{(w, v)| w, v \in (T)_m, C(w) \cap C(v) \neq \emptyset\}.
\]
and define $\widetilde{\B} = \A \cup (\cup_{m \ge 1} \widetilde{E}^h_m)$. Then Elek's graph coincides with $(T, \widetilde{\B})$. Note that in $(T, \B)$, the horizontal edges are 
\[
E^h_m = \{(w, v)| w, v \in (T)_m, K_w \cap K_v \neq \emptyset\}.
\]
So, $\widetilde{E}^h_m \supseteq E^h_m$ and hence $\widetilde{\B} \supseteq \B$ in general.
 In Example~\ref{HYP.ex10}, we are going to show that $(T, \B)$ is hyperbolic as a corollary of our general framework.
\endexample

\section{Weight function and associated ``visual pre-metric''}\label{CNB}

Throughout this section, $(T, \A, \phi)$ is a locally finite tree with a reference point $\phi$, $(X, \O)$ is a compact metrizable topological space with no isolated point and $K: T \to \C(X, \O)$ is a partition of $X$ parametrized by $(T, \A, \phi)$.\par
In this section, we introduce the notion of a weight function, which assigns each vertex of the tree $T$ a ``size'' or ``weight''. Then, we will introduce a kind of ``balls'' and a ``distance'' of the compact metric space $X$ associated with the weight function.\par

\definition[Weight function]\label{PAS.def10}
A function $g: T \to (0, 1]$ is called a weight function if and only if it satisfies the following conditions (G1), (G2) and (G3):\newline
(G1)\,\,\,
$g(\phi) = 1$\newline
(G2)\,\,\,
For any $w \in T$, $g(\pi(w)) \ge g(w)$\newline
(G3)\,\,\,
$\lim_{m \to \infty} \sup_{w \in (T)_m}g(w) = 0$.\\
We denote the collection of all the weight functions by $\G(T)$.
Let $g$ be a weight function. We define
\[
\LL_s^g = \{w|w \in T, g(\pi(w)) > s \ge g(w)\}
\]
for any $s \in (0, 1]$. $\{\LL_s^g\}_{s \in (0, 1]}$ is called the scale associated with $g$. For $s > 1$, we define $\LL_s^g = \{\phi\}$.
\enddefinition

\remark
To be exact, one should use $\G(T, A, \phi)$ rather than $\G(T)$ as the notation for the collection of all the weight functions because the notion of weight function apparently depends not only on the set $T$ but also the structure of $T$ as a tree. We use, however, $\G(T)$ for simplicity as long as no confusion may occur.
\endremark

\remark
In the case of the partitions associated with a self-similar set appearing in Example~\ref{PAS.ex10}, the counterpart of weight functions was called  gauge functions in \cite{Ki13}. Also $\{\LL_s^g\}_{0 < s \le 1}$ was called the scale associated with the gauge function $g$.
\endremark

Given a weight function $g$, we consider $g(w)$ as a virtual ``size'' or ``diameter'' of $\SS_w$ for each $w \in T$. The set $\LL^g_s$ is the collection of subsets $\SS_w$'s whose sizes are approximately $s$.

\prop\label{CNB.prop10}
Suppose that $g: T \to (0, 1]$ satisfies {\rm (G1)} and {\rm (G2)}.  $g$ is a weight function if and only if 
\begin{equation}\label{CNB.eq10}
\lim_{m \to \infty} g([\omega]_m) = 0
\end{equation}
for any $\omega \in \SS$.
\endprop

\demo
If $g$ is a weight function, i.e.\,(G3) holds, then \eqref{CNB.eq10} is immediate.\par
Suppose that (G3) does not hold, i.e. there exists $\e > 0$ such that
\begin{equation}\label{QSY.eq10}
\sup_{w \in (T)_m} g(w) > \e
\end{equation}
for any $m \ge 0$. Define $Z = \{w | w \in T, g(w) > \e\}$ and $Z_m = (T)_m  \cap Z$. By \eqref{QSY.eq10}, $Z_m \neq \emptyset$ for any $m \ge 0$. Since $\pi(w) \in Z$ for any $w \in Z$, if $Z_{m, n} = \pi^{n - m}(Z_n)$ for any $n \ge m$, where $\pi^k$ is the $k$-th iteration of $\pi$, then $Z_{m, n} \neq \emptyset$ and $Z_{m, n} \supseteq Z_{m, n + 1}$ for any $n \ge m$. Set $Z^*_m = \cap_{n \ge m} Z_{m, n}$. Since $(T)_m$ is a finite set and so is $Z_{m, n}$, we see that $Z^*_m \neq \emptyset$ and $\pi(Z^*_{m + 1}) = Z^*_m$ for any $m \ge 0$. Note that $Z^*_0 = \{\phi\}$. Inductively, we may construct a sequence $(\phi, w(1), w(2), \ldots)$ satisfying $\pi(w(m + 1)) = w(m)$ and $w(m) \in Z^*_m$ for any $m \ge 0$.  Set $\omega = (\phi, w(1), w(2), \ldots)$. Then $\omega \in \SS$ and $g([\omega]_m) \ge \e$ for any $m \ge 0$. This contradicts \eqref{CNB.eq10}.
\enddemo

\prop\label{PAS.prop05}
Let $g: T \to (0, 1]$ be a weight function and let $s \in (0, 1]$. Then
\begin{equation}\label{PAS.eq10}
\bigcup_{w \in \LL_s^g} \SS_w = \SS
\end{equation}
and if $w, v \in \LL_s^g$ and $w \neq v$, then
\[
\SS_w \cap \SS_v = \emptyset.
\]
\endprop

\demo
For any $\omega = (w_0, w_1, \ldots) \in \SS$, $\{g(w_i)\}_{i = 0, 1, \ldots}$ is monotonically non-increasing sequence converging to $0$ as $i \to \infty$. Hence there exists  a unique $m \ge 0$ such that $g(w_{m - 1}) > s \ge g(w_m)$. Therefore, there exists a unique $m \ge 0$ such that $[\omega]_m \in \LL_s^g$. Now \eqref{PAS.eq10} is immediate. Assume $w, v \in \LL_s^g$ and $\SS_v \cap \SS_w \neq \emptyset$. Choose $\omega = (w_0, w_1, \ldots) \in \SS_v \cap \SS_w$. Then there exist $m, n \ge 0$ such that $[\omega]_m = w_m = w$ and $[\omega]_n = w_n = v$. By the above fact, we have $m = n$ and hence $w = v$.
\enddemo

By means of the partition $K: T \to \C(X, \O)$, one can define weight functions naturally associated with metrics and measures on the compact metric space $X$ as follows.

\notation
Let $d$ be a metric on $X$. We define the diameter of a subset $A \subseteq X$ with respect to $d$, $\diam{A, d}$ by $\diam{A, d} = \sup\{d(x, y)| x, y \in A\}$. Moreover, for $x \in X$ and $r > 0$, we set $B_d(x, r) = \{y| y \in X, d(x, y) < r\}$.
\endnotation

\definition\label{PAS.def50}
(1)\,\,
Define
\begin{multline*}
\D(X, \O) = \{d|\text{$d$ is a metric on $X$ inducing the topology $\O$ and}\\
\diam{X, d} = 1\}
\end{multline*}
For $d \in \D(X, \O)$, define $g_d: T \to (0, 1]$ by $g_d(w) = \diam{K_w, d}$ for any $w \in T$.
(2)\,\,
Define
\begin{multline*}
\M_P(X, \O) = \{\mu| \text{$\mu$ is a Radon probability measure on $(X, \O)$}\\
\text{satisfying $\mu(\{x\}) = 0$ for any $x \in X$ and $\mu(K_w) > 0$ for any $w \in T$}\}
\end{multline*}
For $\mu \in \M_P(X, \O)$, define $g_{\mu}: T \to (0, 1]$ by $g_{\mu}(w) = \mu(K_w)$ for any $w \in T$.
\enddefinition

The condition $\diam{X, d} = 1$ in the definition of $\D(X, \O)$ is only for the purpose of normalization. Note that since $(X, \O)$ is compact, if a metric $d$ on $X$ induces the topology $\O$, then $\diam{X, d} < +\infty$.\par

\prop\label{PAS.prop00}
{\rm (1)}\,\,
For any $d \in \D(X, \O)$, $g_d$ is a weight function.\\
{\rm (2)}\,\,
For any $\mu \in \M_P(X, \O)$, $g_{\mu}$ is a weight function.
\endprop

\demo 
(1)\,\,The properties (G1) and (G2) are immediate from the definition of $g_d$. Suppose that there exists $\omega \in \SS$ such that
\begin{equation}\label{QSY.eq15}
\lim_{m \to \infty} g_d([\omega]_m) > 0
\end{equation}
Let $\e$ be the above limit. Since $g_d([\omega]_m) = \diam{K_{[\omega]_m}, d} > \e$, there exist $x_m, y_m \in K_{[\omega]_m}$ such that $d(x_m, y_m) \ge \e$. Note that $K_{[\omega]_m} \supseteq K_{[\omega]_{m + 1}}$ and hence $x_n, y_n \in K_{[\omega]_m}$ if $n \ge m$. Since $X$ is compact, there exist subsequences $\{x_{n_i}\}_{i \ge 1}$, $\{y_{n_i}\}_{i \ge 1}$ converging to $x$ and $y$ as $i \to \infty$ respectively. It follows that $x, y \in \cap_{m \ge 0} K_{[\omega]_m}$ and $d(x, y) \ge \e > 0$. This contradicts (P2). Thus we have shown \eqref{CNB.eq10}. By Proposition~\ref{CNB.prop10}, $g_d$ is a weight function.\\
(2)\,\,
As in the case of metrics, (G1) and (G2) are immediate. Let $\omega \in \SS$. Then $\cap_{m \ge 0} K_{[\omega]_m} = \{\s(\omega)\}$. Therefore, $g_{\mu}([\omega]_m) = \mu(K_{[\omega]_m}) \to 0$ as $m \to \infty$. Hence we verify \eqref{CNB.eq10}. Thus by Proposition~\ref{CNB.prop10}, $g_{\mu}$ is a weight function.
\enddemo

The weight function $g_{d}$ and $g_{\mu}$ are called the weight functions associated with $d$ and $\mu$ respectively. Although the maps $d \to g_d$ and $\mu \to g_{\mu}$ are not injective at all, we sometimes abuse notations and use $d$ and $\mu$ to denote $g_d$ and $g_{\mu}$ respectively. \par
Through a partition we introduce the notion of ``balls'' of a compact metric space associated with a weight function.

\definition\label{PAS.def30}
Let $g : T \to (0, 1]$ be a weight function. \newline
(1)\,\,\,
For $s \in (0, 1], w \in \LL^g_s, M \ge 0$ and $x \in X$, we define 
\begin{multline*}
\LL^{g}_{s, M}(w) = \{v| v \in \LL^g_s, \text{there exists  a chain $(w(1), \ldots, w(k))$ of $K$ in $\LL_s^g$}\\ 
\text{such that $w(1) = w$, $w(k) = v$ and $k \le M + 1$}\}
\end{multline*}
and
\[
\LL^{g}_{s, M}(x) = \bigcup_{w \in \LL_s^g\,\,\text{and}\,\, x \in K_w} \LL^{g}_{s, M}(w).
\]
For $x \in X$, $s \in (0, 1]$ and $M \ge 0$, define
\[
U^g_{M}(x, s) = \bigcup_{w \in \LL^{g}_{s, M}(x)} K_w.
\]
We let $U^g_{M}(x, s) = X$ if $s \ge 1$.
\enddefinition

In Figure~\ref{sc3}, we show examples of $U^g_M(x, s)$ for the Sierpinski carpet introduced in Example~\ref{exsc1}.

The family $\{U^g_M(x, s)\}_{s > 0}$ is a fundamental system of neighborhood of $x \in X$ as is shown in Proposition~\ref{ADD.prop10}. \par
Note that 
\[
\LL^g_{s, 0}(w) = \{w\}\quad\text{and}\quad\LL^g_{s, 1}(w) = \{v | v \in \LL^g_s, K_v \cap K_w \neq \emptyset\}
\]
 for any $w \in \LL^g_s$ and 
 \[
 \LL^g_{s, 0}(x) = \{w| w \in \LL^g_s, x \in K_w\}\quad\text{and}\quad U^g_0(x, s) = \bigcup_{w: w \in \LL^g_s, x \in K_w} K_w
 \]
 for any $x \in X$. Moreover,
 \[
U^g_{M}(x, s) = \{y| y \in X, \text{there exists $(w(1), \ldots, w(M + 1)) \in \CH_K^{\LL_s^g}(x, y)$.}\}
\]

\prop\label{ADD.prop10}
Let $K$ be a partition of $X$ parametrized by $(T, \A, \phi)$ and let $g: T \to (0, 1]$ be a weight function. For any $s \in (0, 1]$ and any $x \in X$, $U_0^g(x, s)$ is a neighborhood of $x$. Furthermore, $\{U^g_M(x, s)\}_{s \in (0, 1]}$ is a fundamental system of neighborhood of $x$ for any $x \in X$.
\endprop

\demo
Let $d$ be a metric on $X$ giving the original topology of $(X, \O)$. Assume that for any $r > 0$, there exists $y \in B_d(x, r)$ such that $y \notin U_0^g(x, s)$. Then there exists a sequence $\{x_n\}_{n \ge 1} \subseteq X$ such that $x_n \to x$ and $x_n \notin U_0^g(x, s)$ for any $n \ge 1$. Since $\LL_s^g$ is a finite set, there exists $w \in \LL_s$ which includes infinite members of $\{x_n\}_{n \ge 1}$. By the closedness of $K_w$, it follows that $x \in K_w$ and $x_n \in K_w \subseteq U_0^g(x, s)$. This contradiction shows that $U_0^g(x, s)$ contains $B_d(x, r)$ for some $r > 0$.\\
Next note $\min_{w \in \LL_s^g} |w| \to \infty$ as $s \downarrow 0$. This along with that fact that $g_d$ is a weight function implies that $\max_{w \in \LL^g_s} \diam{K_w, d} \to 0$ as $s \downarrow 0$. Set $\rho_s = \max_{w \in \LL^g_s} \diam{K_w, d}$. Then $\diam{U^g_M(x, s), d} \le (M + 1)\rho_s \to 0$ as $s \downarrow 0$. This implies that $\cap_{s \in (0, 1]} U_M^g(x, s) = \{x\}$.  Thus $\{U^g_M(x, s)\}_{s \in (0, 1]}$ is a fundamental system of neighborhoods of $x$.
\enddemo

We regard $U^g_M(x, s)$ as a virtual ``ball' of radius $s$ and center $x$. In fact, there exists a kind of ``pre-metric'' $\d_M^g: X \times X \to [0, \infty)$ such that $\d_M^g(x, y) > 0$ if and only if $x \neq y$, $\d_M^g(x, y) = \d_M^g(y, x)$ and
\begin{equation}\label{CNB.eq100}
U_M^g(x, s) = \{y| \d_M^g(x, y) \le s\}.
\end{equation}
As is seen in the next section, however, the pre-metric $\d_M^g$ may not satisfy the triangle inequality in general.

\definition\label{CNB.def10}
Let $M \ge 0$. Define $\d_M^g(x, y)$ for $x, y \in X$ by
\[
\d_M^g(x, y) = \inf\{s| s \in (0, 1], y \in U_M^{g}(x, s)\}.
\]
\enddefinition

\remark
For any $g \in \G(T)$, $M \ge 0$ and $x \in X$, it follows that $\LL^g_{s, M}(x) = \{\phi\}$ and hence $U^g_M(x, 1) = X$. So, $\d^g_M(x, y) \le 1$ for any $x, y \in X$.
\endremark

The pre-metric $\d_M^g$ can be thought of as a counterpart of the ``visual metric'' in the sense of Bonk-Meyer in \cite{BonkMeyer} and the ``visual pre-metric'' in the framework of Gromov hyperbolic metric spaces, whose exposition can be found in \cite{BuyaloSchr} and \cite{MacTyson}. In fact, if certain rearrangement of the resolution $(X, \B)$ is hyperbolic associated with the weight function, then $\d_M^g$ is bi-Lipschitz equivalent to a visual pre-metric in the sense of Gromov. See Theorem~\ref{HYP.thm20} for details.

\prop\label{COM.prop10}
For any $M \ge 0$ and $x, y \in X$,
\begin{equation}\label{COM.eq100}
\d_M^g(x, y) = \min\{s| s \in (0, 1], y \in U_M^g(x, s)\}.
\end{equation}
In particular, \eqref{CNB.eq100} holds for any $M \ge 0$ and $s \in (0, 1]$.
\endprop

\demo
The property (G3) implies that for any $t \in (0, 1]$, there exists $n \ge 0$ such that $\cup_{s \ge t} \LL_s^g \subseteq \cup_{m = 0}^n (T)_m$. Hence $\{(w(1), \ldots, w(M + 1))| w(i) \in \cup_{s \ge t} \LL_s^g\}$ is finite. Let $s_* = \d_M^g(x, y)$. Then there exist a sequence $\{s_m\}_{m \ge 1} \subseteq [s_*, 1]$ and $(w_m(1), \ldots, w_m(M + 1)) \in (\LL^g_{s_m})^{M + 1}$ such that $\lim_{m \to \infty} s_m = s_*$ and $(w_m(1), \ldots, w_m(M + 1))$ is a chain between $x$ and $y$ for any $m \ge 1$. Since $\{(w(1), \ldots, w(M + 1))| w(i) \in \cup_{s \ge s_*} \LL_s^g\}$ is finite, there exists $(w_*(1), \ldots, w_*(M + 1))$ such that $(w_*(1), \ldots, w_*(M + 1)) = (w_m(1), \ldots, w_m(M + 1))$ for infinitely many $m$. For such $m$, we have $g(\pi(w_*(i))) > s_m \ge g(w_*(i))$ for any $i = 1, \ldots, M + 1$. This implies that $w_*(i) \in \LL_{s_*}^g$ for any $i = 1, \ldots, M + 1$ and hence $y \in U^g_M(x, s_*)$. Thus we have shown \eqref{COM.eq100}.
\enddemo

\setcounter{equation}{0}
\section{Metrics adapted to weight function}\label{MAG}

In this section, we consider the first question mentioned in the introduction, which is when a weight function is naturally associated with a metric. Our answer will be given in Theorem~\ref{COM.thm10}.\par
As in the last section, $(T, \A, \phi)$ is a locally finite tree with a reference point $\phi$, $(X, \O)$ is a compact metrizable topological space with no isolated point and $K: T \to \C(X, \O)$ is a partition throughout this section.\par
The purpose of the next definition is to clarify when the virtual balls $U^g_M(x, s)$ induced by a weight function $g$ can be though of as real ``balls'' derived from a metric.

\definition\label{PAS.def35}
Let $M \ge 0$. A metric $d \in \D(X, \O)$ is said to be $M$-adapted to $g$ if and only if there exist $\a_1, \a_2 > 0$ such that
\[
U^{g}_{M}(x, \a_1r) \subseteq B_d(x, r) \subseteq U^g_{M}(x, \a_2r)
\]
for any $x \in X$ and any $r > 0$. $d$ is said to be adapted to $g$ if and only if $d$ is $M$-adapted to $g$ for some $M \ge 0$.
\enddefinition

Now our question is the existence of a metric adapted to a given weight function.
The number $M$ really makes a difference in the above definition. Namely, in Example~\ref{COM.ex10}, we construct an example of a weight function to which no metric is $1$-adapted but some metric is $2$-adapted. 

By \eqref{CNB.eq100}, a metric $d \in \D(X, \O)$ is $M$-adapted to a weight function $g$ if and only if there exist $c_1, c_2 > 0$ such that
\begin{equation}\label{MAG.eq10}
c_1\d_M^g(x, y) \le d(x, y) \le c_2\d_M^g(x, y)
\end{equation}
for any $x, y \in X$. By this equivalence, we may think of a metric adapted to a weight function as a ``visual metric'' associated with the weight function.\par
If a metric $d$ is $M$-adapted to a weight function $g$, then we think of the virtual balls $U^g_M(x, s)$ as the real balls associated with the metric $d$.\par

\begin{figure}
\centering
\setlength{\unitlength}{30mm}

\begin{picture}(4, 2)(-1, -0.3)
\put(-2.2, 0)
{\thicklines
\drawline (1.2, 0)(3, 0)(3, 1.8)(1.2, 1.8)(1.2, 0)
\drawline(1.8, 0)(1.8, 1.8)
\drawline(1.4,0)(1.4, 1.8)
\drawline(1.6,0)(1.6, 1.8)
\drawline(2.0,0)(2.0, 0.6)
\drawline(2.2,0)(2.2, 0.6)
\drawline(2.0,1.2)(2.0, 1.8)
\drawline(2.2,1.2)(2.2, 1.8)
\drawline(2.4, 0)(2.4, 1.8)
\drawline(2.6, 0)(2.6, 1.8)
\drawline(2.8, 0)(2.8, 1.8)
\drawline(1.2, 0.2)(3, 0.2)
\drawline(1.2, 0.4)(3, 0.4)
\drawline(1.2, 0.6)(3, 0.6)
\drawline(1.2, 0.6)(3, 0.6)
\drawline(1.2,0.8)(1.8,0.8)
\drawline(2.4,0.8)(3,0.8)
\drawline(1.2,1)(1.8,1)
\drawline(2.4,1)(3,1)
\drawline(1.2, 1.2)(3, 1.2)
\drawline(1.2, 1.4)(3, 1.4)
\drawline(1.2, 1.6)(3, 1.6)
\thinlines
\shade\path(1.8,0.6)(2.4,0.6)(2.4,1.2)(1.8,1.2)(1.8,0.6)
\shade\path(1.4,0.2)(1.6,0.2)(1.6,0.4)(1.4, 0.4)(1.4, 0.2)
\shade\path(2,0.2)(2.2,0.2)(2.2,0.4)(2, 0.4)(2, 0.2)
\shade\path(2.6,0.2)(2.8,0.2)(2.8,0.4)(2.6, 0.4)(2.6, 0.2)
\shade\path(1.4,0.8)(1.6,0.8)(1.6,1)(1.4, 1)(1.4, 0.8)
\shade\path(1.4,1.4)(1.6,1.4)(1.6,1.6)(1.4, 1.6)(1.4, 1.4)
\shade\path(2,1.4)(2.2,1.4)(2.2,1.6)(2, 1.6)(2, 1.4)
\shade\path(2.6,1.4)(2.8,1.4)(2.8,1.6)(2.6, 1.6)(2.6, 1.4)
\shade\path(2.6,0.8)(2.8,0.8)(2.8,1)(2.6, 1)(2.6, 0.8)
\thinlines
\multiput(1.6,0.2)(0,0.025){8}{\drawline(0,0)(0.4,0)}
\multiput(1.4,0.6)(0,0.025){8}{\drawline(0,0)(0.4,0)}
\multiput(1.4,0.4)(0,0.025){8}{\drawline(0,0)(0.6,0)}
\multiput(1.6,0.2)(0.025,0){8}{\drawline(0,0)(0,0.6)}
\multiput(1.8,0.2)(0.025,0){8}{\drawline(0,0)(0,0.4)}
\multiput(1.4,0.4)(0.025,0){8}{\drawline(0,0)(0,0.4)}
\thicklines
\path (1.6, 0.2)(2.0,0.2)(2.0, 0.6)(1.8, 0.6)(1.8, 0.8)(1.4, 0.8)(1.4, 0.4)(1.6, 0.4)(1.6, 0.2)
\put(1.6667, 0.4667){\circle*{0.03}}
\put(1.6667, 0.35){\makebox(0,0)[b]{{\large $x$}}}
\put(1.95,0.16){\circle*{0.03}}
\put(1.95,0.03){\makebox(0,0)[b]{{\large $y$}}}
\put(2.1, -0.2){\makebox(0,0)[b]{$\Lambda^{\eta}_{\frac 19}$, The meshed region is $U^{\eta}_1(x, \frac 19)$}}
\put(2.1, -0.4){\makebox(0,0)[b]{$\delta^{\eta}_1(x, y) = \frac 13$}}
}

\thicklines

\drawline (1.2, 0)(3, 0)(3, 1.8)(1.2, 1.8)(1.2, 0)
\drawline(1.8, 0)(1.8, 1.8)
\drawline(1.4,0.6)(1.4, 1.2)
\drawline(1.6,0.6)(1.6, 1.2)
\drawline(2.0,0)(2.0, 0.6)
\drawline(2.2,0)(2.2, 0.6)
\drawline(2.0,1.2)(2.0, 1.8)
\drawline(2.2,1.2)(2.2, 1.8)
\drawline(2.4, 0)(2.4, 1.8)
\drawline(2.6, 0.6)(2.6, 1.2)
\drawline(2.8, 0.6)(2.8, 1.2)

\drawline(1.8, 0.2)(2.4, 0.2)
\drawline(1.8, 0.4)(2.4, 0.4)
\drawline(1.2, 0.6)(3, 0.6)
\drawline(1.2,0.8)(1.8,0.8)
\drawline(2.4,0.8)(3,0.8)
\drawline(1.2,1)(1.8,1)
\drawline(2.4,1)(3,1)
\drawline(1.2, 1.2)(3, 1.2)
\drawline(1.8, 1.4)(2.4, 1.4)
\drawline(1.8, 1.6)(2.4, 1.6)
\thinlines
\shade\path(1.8,0.6)(2.4,0.6)(2.4,1.2)(1.8,1.2)(1.8,0.6)
\shade\path(1.4,0.2)(1.6,0.2)(1.6,0.4)(1.4, 0.4)(1.4, 0.2)
\shade\path(2,0.2)(2.2,0.2)(2.2,0.4)(2, 0.4)(2, 0.2)
\shade\path(2.6,0.2)(2.8,0.2)(2.8,0.4)(2.6, 0.4)(2.6, 0.2)
\shade\path(1.4,0.8)(1.6,0.8)(1.6,1)(1.4, 1)(1.4, 0.8)
\shade\path(1.4,1.4)(1.6,1.4)(1.6,1.6)(1.4, 1.6)(1.4, 1.4)
\shade\path(2,1.4)(2.2,1.4)(2.2,1.6)(2, 1.6)(2, 1.4)
\shade\path(2.6,1.4)(2.8,1.4)(2.8,1.6)(2.6, 1.6)(2.6, 1.4)
\shade\path(2.6,0.8)(2.8,0.8)(2.8,1)(2.6, 1)(2.6, 0.8)
\thinlines
\multiput(1.2,0)(0, 0.025){8}{\drawline(0,0)(0.8,0)}
\multiput(1.6,0.2)(0,0.025){8}{\drawline(0,0)(0.4,0)}
\multiput(1.2,0.6)(0,0.025){8}{\drawline(0,0)(0.6,0)}
\multiput(1.2,0.4)(0,0.025){8}{\drawline(0,0)(0.8,0)}
\multiput(1.6,0)(0.025,0){8}{\drawline(0,0)(0,0.8)}
\multiput(1.8,0.0)(0.025,0){8}{\drawline(0,0)(0,0.6)}
\multiput(1.4,0.4)(0.025,0){8}{\drawline(0,0)(0,0.4)}
\multiput(1.2,0)(0.025,0){8}{\drawline(0,0)(0,0.8)}
\multiput(1.4,0)(0.025,0){8}{\drawline(0,0)(0,0.2)}
\multiput(1.2,0.2)(0,0.025){8}{\drawline(0,0)(0.2,0)}
\thicklines
\put(1.6667, 0.4667){\circle*{0.03}}
\put(1.6667, 0.35){\makebox(0,0)[b]{{\large $x$}}}
\put(1.95,0.16){\circle*{0.03}}
\put(1.95,0.03){\makebox(0,0)[b]{{\large $y$}}}
\put(2.1, -0.2){\makebox(0,0)[b]{$\Lambda^g_{\frac 19}$, The meshed region is $U^g_1(x, \frac 19)$}}
\put(2.1, -0.4){\makebox(0,0)[b]{$\delta^g_1(x, y) = \frac 19$}}
\end{picture}
\caption{Visual pre-metrics: the Sierpinski carpet}\label{sc3}
\end{figure}
\example[Figure~\ref{sc3}]\label{PAS.ex100}
Let us consider the case of the Sierpinski carpet introduced in Example~\ref{exsc1}. In this case, the corresponding tree is $(T^{(8)}, \A^{(8)}, \phi)$. Write $T = T^{(8)}$. Define $\eta: T \to (0, 1]$ by $\eta(w) = \frac 1{3^m}$ for any $w \in (T)_m$. Then $\eta$ is a weight function and $\LL^{\eta}_s = (T)_m$ if and only if $\frac 1{3^{m - 1}} > s \ge \frac 1{3^m}$. Let $d_*$ be the (restriction of) Euclidean metric. Then $d_*$ is $1$-adapted to $h$. This can be deduced from the following two observations. First, if $w, v \in (T)_m$ and $K_w \cap K_v \neq \emptyset$, then $\sup_{x \in K_w, y \in K_v} d_*(x, y) \le \frac{2\sqrt{2}}{3^m}$. Second, if $w, v \in (T)_m$ and $K_w \cap K_v = \emptyset$, then $\inf_{x \in K_w, y \in K_v} d_*(x, y) \ge \frac 1{3^m}$. In fact, these two facts implies that
\[
\frac 13\d^{\eta}_1(x, y) \le d_*(x, y) \le 2\sqrt{2}\d^{\eta}_1(x, y)
\] 
for any $x, y \in X$. Next we try another weight function $g$ defined as
\[
g(i_1\ldots{i_m}) = r_{i_1}{\cdots}r_{i_m}
\]
for any $m \ge 0$ and $i_1, \ldots, i_m \in \{1, \ldots, 8\}$, where
\[
r_i = \begin{cases}
\frac 19\quad&\text{if $i$ is odd,}\\
\frac 13\quad&\text{if $i$ is even.}
\end{cases}
\]
Then $\LL^g_{\frac 13} = \{1, \ldots, 8\}$ and 
\[\LL^g_{\frac 13} = \{1, 3, 5, 7\} \cup \{i_1i_2| i_1 \in \{2, 4, 6, 8\}, i_2 \in \{1, \ldots, 8\}\}.
\]
In this case, the existence of an adapted metric is not immediate. However, by \cite[Example~1.7.4]{Ki13}, it follows that $\eta \gen g$. (See Definition~\ref{VDP.def10} for the definition of $\gen$.) By Theorems~\ref{COM.thm30} and \ref{HYP.thm20}, there exists a metric $\rho \in \D(X, \O)$ that is adapted to $g^{\a}$ for some $\a > 0$. Furthermore, Theorem~\ref{QSY.thm10} shows that $\rho$ is quasisymmetric to $d_*$.
\endexample

There is another ``pre-metric'' associated with a weight function.

\definition\label{MAG.def10}
Let $M \ge 0$. Define $D^g_{M}(x, y)$ for $x, y \in X$ by
\[
D^g_{M}(x, y) = \inf\Big\{\sum_{i = 1}^k g(w(i))\Big|1 \le k \le M + 1, (w(1), \ldots, w(k)) \in \CH_K(x, y)\Big\}
\]
\enddefinition

It is easy to see that $0 \le D_M^g(x, y) \le 1$, $D_M^g(x, y) = 0$ if and only if $x = y$ and $D_M^g(x, y) = D_M^g(y, x)$. In fact, the pre-metric $D_M^g$ is equivalent to $\d_M^g$ as follows.

\prop\label{MAG.prop10}
For any $M \ge 0$ and $x, y \in X$,
\[
\d_M^g(x, y) \le D^g_M(x, y) \le (M + 1)\d_M^g(x, y).
\]
\endprop

\demo
Set $s_* = \d_M^g(x, y)$. Using Proposition~\ref{COM.prop10}, we see that there exists a chain $(w(1), \ldots, w(M + 1))$ between $x$ and $y$ such that $w(i) \in \LL_{s_*}^g$ for any $i = 1, \ldots, M + 1$. Then
\[
D_M^g(x, y) \le \sum_{i = 1}^{M + 1} g(w(i)) \le (M + 1)s_*
\]
Next set $d_* = D_M^g(x, y)$. For any $\e > 0$, there exists a chain $(w(1), \ldots, w(M + 1))$ between $x$ and $y$ such that $\sum_{i = 1}^{M + 1} g(w(i)) < d_* + \e$. In particular, $g(w(i)) < d_* + \e$ for any $i = 1, \ldots, M + 1$. Hence for any $i = 1, \ldots, M + 1$, there exists $w_*(i) \in \LL^g_{d_* + \e}$ such that $K_{w(i)} \subseteq K_{w_*(i)}$. Since $(w_*(1), \ldots, w_*(M + 1))$ is a chain between $x$ and $y$, it follows that $\d_M^g(x, y) \le d_* + \e$. Thus we have shown $\d_M^g(x, y) \le D_M^g(x, y)$.
\enddemo

Combining the above proposition with \eqref{MAG.eq10}, we see that $d$ is $M$-adapted to  $g$ if and only if there exist $C_1, C_2 > 0$ such that
\begin{equation}\label{MAG.eq20}
C_1D_M^g(x, y) \le d(x, y) \le C_2D_M^g(x, y)
\end{equation}
for any $x, y \in X$.\par
Next we present another condition which is equivalent to a metric being adapted.

\thm\label{PAS.thm10}
Let $g: T \to (0, 1]$ be a weight function and let $M \ge 0$.  If $d \in \D(X, \O)$, then $d$ is $M$-adapted to $g$ if and only if the following conditions {\rm (ADa)} and $\rm (ADb)_M$ hold:\\
{\rm (ADa)}\,\,\,
There exists $c > 0$ such that $\diam{K_w, d} \le cg(w)$ for any $w \in T$.\newline
{$\rm (ADb)_M$}\,\,\,
For any $x, y \in X$, there exists $(w(1), \ldots, w(k)) \in \CH_K(x, y)$ such that $1 \le k \le M + 1$ and
\[
Cd(x, y) \ge \max_{i = 1, \ldots, k} g(w(i)),
\]
where $C > 0$ is independent of $x$ and $y$.
\endthm

\remark
In \cite[Proposition~8.4]{BonkMeyer}, one find an analogous result in the case of partitions associated with expanding Thurston maps. The condition (ADa) and {$\rm (ADb)_M$} corresponds their conditions (ii) and (i) respectively.
\endremark

\demo
First assuming ${\rm (ADa)}$ and ${\rm (ADb)}_M$, we are going to show \eqref{MAG.eq10}. Let $x, y \in X$. By $\rm (ADb)_M$, there exists a chain $(w(1), \ldots, w(k))$ between $x$ and $y$ such that $1 \le k \le M + 1$ and $Cd(x, y) \ge g(w(i))$ for any $i = 1, \ldots, k$. By (G2), there exists $v(i)$ such that $\SS_{v(i)} \supseteq \SS_{w(i)}$ and $v(i) \in \LL^g_{Cd(x, y)}$. Since $(v(1), \ldots, v(k))$ is a chain in $\LL_{Cd(x, y)}^g$ between $x$ and $y$, it follows that $Cd(x, y) \ge \d_M^g(x, y)$. \par
Next set $t = \d_M^g(x, y)$. Then there exists a chain $(w(1), \ldots, w(M + 1)) \in \CH_K(x, y)$ in $\LL_t^g$. Choose $x_i \in K_{w(i)} \cap K_{w(i + 1)}$ for every $i = 1, \ldots, M$. Then
\begin{multline*}
d(x, y) \le d(x, x_1) + \sum_{i = 1}^{M - 1} d(x_i, x_{i + 1}) + d(x_{M}, y)\\
\le c\sum_{j = 1}^{M + 1} g(w(i)) \le c(M + 1)t = c(M + 1)\d_M^g(x, y).
\end{multline*}
Thus we have \eqref{MAG.eq10}.\par
Conversely, assume that \eqref{MAG.eq10} holds, namely, there exist $c_1, c_2 > 0$ such that $c_1d(x, y) \le \d^g_{M}(x, y) \le c_2d(x, y)$ for any $x, y \in X$. If $x, y \in K_w$, then $w \in \CH_K(x, y)$. Let $m = \min\{k| g(\pi^k(w)) > g(\pi^{k - 1}(w)), k \in \BbN\}$ and set $s = g(w)$. Then $g(\pi^{k - 1}(w)) = s$ and $\pi^{k - 1}(w) \in \LL_s^g$. Since $\pi^{k - 1}(w) \in \CH_K(x, y)$, we have
\[
g(w) = s \ge \d_0^g(x, y) \ge \d_M^g(x, y) \ge c_1d(x, y).
\]
This immediately yields (ADa). \par
Set $s_* = c_2d(x, y)$ for $x, y \in X$. Since $\d^g_M(x, y) \le c_2d(x, y)$, there exists a chain $(w(1), \ldots, w(M + 1))$ in $\LL_{s_*}^g$ between $x$ and $y$. As $g(w(i)) \le s_*$ for any $i = 1, \ldots, M + 1$, we have $\rm (ADb)_M$.
\enddemo

Since ${\rm (ADb)}_M$ implies ${\rm (ADb)}_N$ for any $N \ge M$, we have the following corollary.

\cor\label{PAS.cor10}
Let $g: T \to (0, 1]$ be a weight function. If $d \in \D(X, \O)$ is $M$-adapted to $g$ for some $M \ge 0$, then it is $N$-adapted to $g$ for any $N \ge M$.
\endcor

Recall that a metric $d \in \D(X, \O)$ defines a weight function $g_d$. So one may ask if $d$ is adapted to the weight function $g_d$ or not. Indeed, we are going to give an example of a metric $d \in \D(X, \O)$ which is not adapted to $g_d$ in Example~\ref{ESS.ex20}.

\definition\label{PAS.def60}
Let $d \in \D(X, \O)$. $d$ is said to be adapted if $d$ is adapted to $g_d$.
\enddefinition

\prop\label{PAS.prop40}
Let $d \in \D(X, \O)$. $d$ is adapted if and only if there exists a weight function $g: T \to (0, 1]$ to which $d$ is adapted. Moreover, suppose that $d$ is adapted. If
\[
D^d(x, y) = \inf\{\sum_{i = 1}^{k} g_d(w(k)) | k \ge 1, (w(1), \ldots, w(k)) \in \CH_K(x, y)\}
\]
for any $x, y \in X$, then there exist $c_* > 0$ such that
\[
c_*D^d(x, y) \le d(x, y) \le D^d(x, y)
\]
for any $x, y \in X$.
\endprop

\demo
Necessity direction is immediate. Assume that $d$ is $M$-adapted to a weight function $g$. By (ADa) and $\rm (ADb)_M$, for any $x, y \in X$ there exist $k \in \{1, \ldots, M + 1\}$ and $(w(1), \ldots, w(k)) \in \CH_K(x, y)$ such that
\[
Cd(x, y) \ge \max_{i = 1, \ldots, k} g(w(i)) \ge \frac 1c\max_{i = 1, \ldots, k} g_d(w(i)).
\]
This proves $\rm(ADb)_M$ for the weight function $g_d$. So we verify that $d$ is $M$-adapted to $g_d$. Now, assuming that $d$ is adapted to $g_d$, we see
\[
c_1D^d_M(x, y) \le d(x, y)
\]
by \eqref{MAG.eq20}.
Since $D^d_M(x, y)$ is monotonically decreasing as $M \to \infty$, it follows that
\[
c_1D^d(x, y) \le d(x, y).
\]
On the other hand, if $(w(1), \ldots, w(k)) \in \CH_K(x, y)$, then the triangle inequality yields
\[
d(x, y) \le \sum_{i = 1}^k g_d(w(i)).
\]
Hence $d(x, y) \le D^d(x, y)$.
\enddemo

Let us return to the question on the existence of a metric associated with a given weight function $g$. Strictly speaking, one should try to find a metric adapted to the weight function $g$ itself. In this section, however, we are going to deal with a modified version, i.e. the existence of a metric adapted to $g^{\a}$ for some $\a > 0$. Note that if $g$ is a weight function, then so is $g^{\a}$ and $\d^{g^{\a}}_M = (\d_M^g)^{\a}$.\par
To start with, we present a weak version of ``triangle inequality'' for  the family $\{\d_M^g\}_{M \ge 1}$.

\prop\label{COM.prop20}
\[
\d_{M_1 + M_2 + 1}^g(x, z) \le \max\{\d_{M_1}^g(x, y),\d_{M_2}^g(y, z)\}
\]
\endprop

\demo
Setting $s_* = \max\{\d_{M_1}^g(x, y), \d_{M_2}^g(y, z)\}$, we see that there exist  a chain $(w(1), \ldots, w(M_1 + 1))$ between $x$ and $y$ and a chain $(v(1), \ldots, v(M_2 + 1))$ between $y$ and $z$ such that $w(i), v(j) \in \LL_{s_*}^g$ for any $i$ and $j$. Since $(w(1), \ldots, w(M_1 + 1), v(1), \ldots, v(M_2 + 1))$ is a chain between $x$ and $z$, we obtain the claim of the proposition.
\enddemo

By this proposition, if $\d_M^g(x, y) \le c\d_{2M + 1}(x, y)$ for any $x, y \in X$, then $\d_M^g(x, y)$ is so-called quasimetric, i.e.
\begin{equation}\label{COM.eq200}
\d_M^g(x, y) \le c\big(\d_M^g(x, z) + \d_M^g(z, y)\big)
\end{equation}
for any $x, y, z \in X$. The coming theorem shows that $\d_M^g$ being a quasimetric is equivalent to the existence of a metric adapted to $g^{\a}$ for some $\a$.\par
The following definition and proposition give another characterization of the visual pre-metric $\d_M^g$.

\definition\label{COM.def20}
For $w, v \in T$, the pair $(w, v)$ is said to be $m$-separated with respect to $\LL_s^g$ if and only if whenever $(w, w(1), \ldots, w(k), v)$ is a chain and $w(i) \in \LL_s^g$ for any $i = 1, \ldots, k$, it follows that $k \ge m$.
\enddefinition

\prop\label{MAG.prop20}
For any $x, y \in X$ and $M \ge 1$,
\[
\d_M^g(x, y) = \sup\{s| \text{$(w, v)$ is $M$-separated if $w, v \in \LL_s^g$, $x \in K_w$ and $y \in K_v$}\}.
\]
\endprop

The following theorem gives several equivalent conditions on the existence of a metric adapted to a given weight function $g$. In Theorem~\ref{HYP.thm20}, those condition will be shown to be equivalent to the hyperbolicity of the rearrangement of the resolution $(T, \B)$ associated with the weight function$g$ and the adapted metric is, in fact, a visual metric in the sense of Gromov.  See \cite{BuyaloSchr} and \cite{MacTyson} for details on visual metric in the sense of Gromov.

\thm\label{COM.thm10}
Let $M \ge 1$ and let $g \in \G(T)$. The following four conditions are equivalent:\\
${\rm (EV)}_M$
There exist $\a \in (0, 1]$ and $d \in \D(X, \O)$ such that $d$ is $M$-adapted to $g^{\a}$.\\
${\rm (EV2)}_M$
$\d_M^g$ is a quasimetric, i.e. there exists $c > 0$ such that \eqref{COM.eq200} holds for any $x, y, z \in X$.\\
${\rm (EV3)}_M$
There exists $\c \in (0, 1)$ such that $\c^n\d_M^g(x, y) \le \d_{M + n}^g(x, y)$ for any $x, y \in X$ and $n \ge 1$.\\
${\rm (EV4)}_M$
There exists $\c \in (0, 1)$ such that $\c\d_M^g(x, y) \le \d_{M + 1}^g(x, y)$ for any $x, y \in X$.\\
Moreover, if $K: T \to \C(X, \O)$ is minimal, then all the conditions above are equivalent to the following condition ${\rm (EV5)}_M$.\\
${\rm (EV5)}_M$
There exists $\c \in (0, 1)$ such that if $(w, v) \in \LL_s^g \times \LL_s^g$ is $M$-separated with respect to $\LL_s^g$, then $(w, v)$ is $(M + 1)$-separated with respect to $\LL_{\c{s}}^g$.
\endthm

The symbol ``EV'' in the above conditions ${\rm (EV)}_M$, ${\rm (EV1)}_M$, \ldots, ${\rm (EV5)}_M$  represents ``Existence of a Visual metric''.\par

We use the following lemma to prove this theorem.

\lemma\label{COM.lemma100}
If there exist $\c \in (0, 1)$ and $M \ge 1$ such that $\c\d_M^g(x, y) \le \d_{M + 1}(x, y)$ for any $x, y \in X$, then
\[
\c^n\d_M^g(x, y) \le \d_{M + n}^g(x, y)
\]
for any $x, y \in X$ and $n \ge 1$.
\endlemma

\demo
We use an inductive argument. Assume that
\[
\c^l\d_M^g(x, y) \le \d_{M + l}^g(x, y)
\]
for any $x, y \in X$ and $l = 1, \ldots, n$. Suppose $\d_{M + n + 1}^g(x, y) \le \c^{n + 1}s$. Then there exists a chain $(w(1), \ldots, w(M + n + 2))$ in $\LL_{\c^{n + 1}s}^g$ between $x$ and $y$. Choose any $z \in K_{w(M + n + 1)} \cap K_{w(M + n + 2)}$. Then
\[
\c^n\d_M^g(x, z) \le \d_{M + n}^g(x, z) \le \c^{n + 1}s.
\]
Thus we obtain $\d_M^g(x, z) \le \c{s}$. Note that $\d_0^g(z, y) \le \c^{n + 1}s$. By Proposition~\ref{COM.prop20}, 
\[
\c\d_M^g(x, y) \le \d_{M + 1}^g(x, y) \le \max\{\d_{M + 1}^g(x, z), \d_0^g(z, y)\} \le \c{s}.
\]
This implies $\d_M^g(x, y) \le s$.
\enddemo

\demo[Proof of Theorem~\ref{COM.thm10}]
${\rm (EV)}_M$ $\Rightarrow$ ${\rm (EV4)}_M$:\,\,Since $d$ is $M$-adapted to $g^{\a}$, by Corollary~\ref{PAS.cor10}, $d$ is $M + 1$-adapted to $g^{\a}$ as well. By \eqref{MAG.eq10}, we obtain ${\rm (EV4)}_M$.\\
${\rm (EV3)}_M$ $\Leftrightarrow$ ${\rm (EV4)}_M$:\,\,This is immediate by Lemma~\ref{COM.lemma100}.\\
${\rm (EV3)}_M$ $\Rightarrow$ ${\rm (EV2)}_M$:\,\, Let $n = M + 1$. By Proposition~\ref{COM.prop20}, we have
\[
c_{2M + 1}\d_M^g(x, y) \le \d_{2M + 1}^g(x, y) \le \max\{\d_M^g(x, z), \d_M^g(z, y)\} \le \d_M^g(x, z) + \d_M^g(z, y).
\]
${\rm (EV2)}_M$ $\Rightarrow$ ${\rm (EV)}_M$:\,\,By \cite[Proposition~14.5]{Heinonen}, there exist $c_1, c_2 > 0$, $d \in \D(X, \O)$ and $\a \in (0, 1]$  such that $c_1\d_M^g(x, y)^{\a} \le d(x, y) \le c_2\d_M^g(x, y)^{\a}$ for any $x, y \in X$. Note that $\d_M^{g}(x, y)^{\a} = \d_M^{g^{\a}}(x, y)$.  By \eqref{MAG.eq10}, $d$ is $M$-adapted to $g^{\a}$.\\
${\rm (EV4)}_M$ $\Rightarrow$ ${\rm (EV5)}_M$:\,\,Assume that $w, v \in \LL_s^g$. If $w$ and $v$ are not $(M + 1)$-separated with respect to $\LL_{{\c}s}^g$, then there exist $w(1), \ldots, w(M) \in \LL_{{\c}s}^g$ such that $(w, w(1), \ldots, w(M), v)$ is a chain. Then we can choose $w' \in T_w \cap \LL_{{\c}s}^g$ and $v' \in T_v \cap \LL_{{\c}s}^g$ so that $(w', w(1), \ldots, w(M), v')$ is a chain. Let $x \in O_{w'}$ and let $y \in O_{v'}$. Then $\d_{M + 1}^g(x, y) \le \c{s}$. Hence by ${\rm (EV4)}_M$, $\d^g_M(x, y) \le s$. There exists a chain $(v(1), v(2), \ldots, v(M + 1))$ in $\LL_s^g$ between $x$ and $y$. Since $x \in O_{w'} \subseteq O_w$ and $y \in O_{v'} \subseteq O_v$, we see that $v(0) = w$ and $v(M + 1) = v$. Hence $w$ and $v$ are not $M$-separated with respect to $\LL_s^g$.\\
${\rm (EV5)}_M$ $\Rightarrow$ ${\rm (EV4)}_M$:\,\,Assume that $\d_{M + 1}^g(x, y) \le \c{s}$. Then there exists a chain $(w(1), \ldots, w(M + 2))$ in $\LL_{\c{s}}^g$ between $x$ and $y$. Let $w$ (resp. $v$) be the unique element in $\LL_s^g$ satisfying $w(1) \in T_w$ (resp. $w(M + 2) \in T_v$). Then $(w, v)$ is not $(M + 1)$-separated in $\LL^g_{\c{s}}$. By ${\rm (EV5)}_M$, $(w, v)$ is not $M$-separated in $\LL^g_s$. Hence there exists a chain $(w, v(1), \ldots, v(M - 1), v)$ in $\LL_s^g$. This implies $\d_M^g(x, y) \le s$.
\enddemo

\setcounter{equation}{0}
\section{Hyperbolicity of resolutions and the existence of adapted metrics}\label{HYP}

In this section, we study the hyperbolicity in Gromov's sense of the resolution $(T, \B)$ of a compact metric space $X$. Roughly speaking the hyperbolicity will be shown to be equivalent to the existence of an adapted metric. More precisely, we define the hyperbolicity of a weight function $g$ as that of certain rearranged subgraph of $(T, \B)$ associated with $g$ and show that the weight function $g$ is hyperbolic if and only if there exists a metric adapted to $g^{\a}$ for some $\a > 0$. Furthermore, in such a case, the adapted metric is shown to be a ``visual metric''. See Theorem~\ref{HYP.thm20} for exact statements. Another important point is the ``boundary'' of the resolution $(T, \B)$ is always identified with the original metric space $X$ with or without hyperbolicity of $(T, \B)$ as is shown in Theorem~\ref{HYP.thm05}. Furthermore, we are going to obtain counterparts by Elek~\cite{Elek1} and Lau-Wang~\cite{LauWang} on the constructions of a hyperbolic graph whose hyperbolic boundary is a given compact metric space as by-products of our general framework.\par
Throughout this section, $(T, \A, \phi)$ is a locally finite tree with a reference point $\phi$, $(X, \O)$ is a compact metrizable topological space with no isolated point and $K: T \to \C(X, \O)$ is a partition of $X$ parametrized by $(T, \A, \phi)$. Moreover, $(T, \B)$ is the resolution of $X$ associated with the partition $K: T \to C(X, \O)$.\par
 The first lemma claims that the collection of geodesic rays of $(T, \B)$ starting from $\phi$ equals $\SS$, which is the collection of geodesic rays of the tree $(T, \A)$ starting from $\phi$.

\lemma\label{HYP.lemma20}
If $(w(0), w(1), w(2), \ldots)$ is a geodesic ray from $\phi$ of $(T, \B)$, then $\pi(w(i + 1)) = w(i)$ for any $i = 1, 2, \ldots$.  In other word, all the edges of a geodesic ray from $\phi$ are vertical edges and the collection of geodesic rays of $(T, \B)$ coincides with $\SS$.
\endlemma

\demo
Suppose that $\pi(w(i)) = w(i - 1)$ for any $i = 1, \ldots, n$ and $(w(n), w(n + 1))$ is a horizontal edge. Then $|w(i)| = i$ for any $i = 0, 1, \ldots, n$ and $|w(n + 1)| = n$. Since $d_{(T, \B)}(\phi, w(n + 1)) = n$, the sequence $(\phi, w(1), \ldots, w(n), w(n + 1))$ can not be a geodesic. Hence there exists no horizontal edge in $(w(0), w(1), w(2), \ldots)$.
\enddemo

The following proposition is the restatement of Proposition~\ref{HF.prop10}.

\prop[= Proposition~\ref{HF.prop10}]\label{HYP.prop10}
Let $\omega, \tau \in \SS$. Then 
\[
\sup_{n \ge 1} d_{(T, \B)}([\omega]_n, [\tau]_n) < +\infty
\]
 if and only if $\s(\omega) = \s(\tau)$.
\endprop

To prove the above proposition, we need to study the structure of geodesics of $(T, d_{(T, \B)})$.

\definition\label{HYP.def10}
(1)\,\,Let $w, v \in (T)_m$ for some $m \ge 0$. The pair $(w, v)$ is called horizontally minimal if and only if there exists a geodesic of the resolution $(T, \B)$ between $w$ and $v$ which consists only of horizontal edges. \\
(2)\,\,Let $w \neq v \in T$. Then a geodesic $\bb$ of $(T, \B)$ between $w$ and $v$ is called a bridge if and only if there exist $i, j \ge 0$ and a horizontal geodesic $(v(1), \ldots, v(k))$ such that $\pi^i(w) = v(1), \pi^j(v) = v(k)$ and 
\[
\bb = (w, \pi(w), \ldots, \pi^i(w), v(2), \ldots, v(k  - 1), \pi^j(v), \ldots, \pi(v), v)
\]
 The number $|v(1)|$ is called the height of the bridge. Also $(w, \pi(w), \ldots, \pi^i(w))$, $(v(1), \ldots, v(k))$ and $(\pi^j(v), \ldots, \pi(v), v)$ are called the ascending part, the horizontal part and the descending part respectively.
 \enddefinition

\lemma\label{HYP.lemma10}
For any $w, v \in T$, there exists a bridge between $w$ and $v$.
\endlemma

\demo
Let $(w(1), \ldots, w(m))$ be a geodesic of $(T, \B)$ between $w$ and $v$. Note that there exists no dent, which is a segment $(w(i), \ldots, w(k), w(k + 1))$ satisfying $|w(i + 1)| = |w(i)| + 1$, $|w(i + 1)| = |w(i + 2)| = \ldots = |w(k)|$ and $|w(k + 1)| = |w(k)| - 1$, because applying $\pi$ to the dent, we can reduce the length at least by $2$. Therefore, if $m_* = \min\{|w(i)| : i = 1, \ldots, m\}$, then there exist $i_* < j_*$ such that $\{i: |w(i)| = m_*\} = \{i: i_* \le i \le j_*\}$ and $|w(i)|$ is monotonically nonincreasing on $I_1 = \{i| 1 \le i \le i_*\}$ and monotonically nondecreasing on $I_2 = \{i|  j_*  \le i \le m\}$. On $I_1$, if there exists $(w(i - 1), w(i), w(i + 1))$ such that $|w(i - 1)| = |w(i)|$ and $w(i + 1) = \pi(w(i))$, then we modify this part to $(w(i - 1), \pi(w(i - 1)), w(i + 1))$. After the modification, the resulting sequence is also a geodesic. Similarly, on $I_2$, if there exists $(w(i -1)), w(i), w(i + 1))$ such that $w(i - 1) = \pi(w(i))$ and $|w(i)| = |w(i + 1)|$, then we modify this part to $(w(i - 1), \pi(w(i + 1)), w(i + 1))$. After modification, the resulting sequence is still a geodesic. Iterating those modifications on $I_1$ and $I_2$ repeatedly as many times as possible, we obtain a bridge between $w$ and $v$ in the end.
\enddemo

\demo[Proof of Proposition~\ref{HYP.prop10}]
Write $\omega_m = [\omega]_m$ and $\tau_m = [\tau]_m$ for any $m \ge 0$. Assume that $\sup_{n \ge 1} d_{(T, \B)}(\omega_n, \tau_n) < +\infty$. Let $d \in \D(X, \O)$. Then $g_d$ is a weight function by Proposition~\ref{PAS.prop00}. In particular, $\max_{w \in (T)_m} \diam{K_w, d} \to 0$ as $m \to \infty$. Set $x = \s(\omega)$ and $y = \s(\tau)$. Let $N = \sup_{n \ge 1} d_{(T, \B)}(\omega_n, \tau_n) < +\infty$. Then for any $m$, there exists  a bridge $(\omega_m, \ldots, \omega_{m - n}, \ldots, \tau_{m - n}, \ldots, \tau_m)$ between $[\omega]_m$ and $[\tau]_m$, where $(\omega_{m - n}, \ldots, \tau_{m - n})$ is a horizontal geodesic. Since $n \le d_{(T, \B)}(\omega_m, \tau_m) \le N$ and the length of $(\omega_{m - n}, \ldots, \tau_{m - n})$ is at most $N$, it follows that $d(x, y) \le N\max_{w \in (T)_{m - N}} \diam{K_w, d} \to 0$ as $m \to \infty$.
Therefore $x = y$.\par
Conversely, if $x = y$, then $x = y \in K_{[\omega]_m} \cap K_{[\omega]_m}$ for any $m \ge 0$. Therefore $d_{(T, \B)}([\omega_m], [\tau]_m) \le 1$ for any $m \ge 0$.
\enddemo

\thm\label{HYP.thm05}
Define an equivalence relation $\sim$ on the collection $\SS$ of the geodesic rays as $\omega \sim \tau$ if and only if $\sup_{n \ge 1} d_{(T, \B)}([\omega]_n, [\tau]_n) < +\infty$. Let $\O_*$ be the natural quotient topology of $\SS/\!\!\sim$ induced by the metric $\rho_*$ on $\SS$. Then
\[
(\SS/\!\!\sim\,, \O_*) = (X, \O),
\]
where we identify $\SS/\!\!\sim$ with $X$ through the map $\s: \SS \to X$.
\endthm

For a Gromov hyperbolic graph, the quotient of the collection of geodesic rays by the equivalence relation $\sim$ is called the hyperbolic boundary of the graph. In our framework, however, the above theorem shows that $\SS/\!\!\sim$ can be always identified with $X$ even if $(T, \B)$ is not hyperbolic.\par
Next we introduce the notion of (Gromov) hyperbolicity of $(T, \B)$. Here we give only basic accounts needed in our work. See \cite{BuyaloSchr}, \cite{MacTyson} and \cite{Vaisala} for details of the general framework of Gromov hyperbolic metric spaces.

\definition\label{HYP.def30}
Define the Gromov product of $w, v \in T$ in $(T, \B)$ with respect to $\phi$ as
\[
(w|v)_{(T, \B), \phi} = \frac{d_{(T, \B)}(\phi, w) + d_{(T, \B)}(\phi, v) - d_{(T, \B)}(w, v)}2.
\]
The graph $(T, \B)$ is called $\eta$-hyperbolic (in the sense of Gromov) if and only if
\[
(w|v)_{(T, \B), \phi} \ge \min\{(w|u)_{(T, \B), \phi}, (u|v)_{(T, \B), \phi}\} - \eta
\]
for any $w, v, u \in T$. $(T, \B)$ is called hyperbolic if and only if it is $\eta$-hyperbolic for some $\eta \in \BbR$.
\enddefinition

It is known that the hyperbolicity can be defined by the thinness of geodesic triangles.
\def\bb{{\bf b}}

\definition\label{HYP.def100}
We say that all the geodesic triangles in $(T, \B)$ are $\d$-thin if and only if for any $w, v, u \in T$, if $\bb(a , b)$ is geodesic between $a$ and $b$ for each $(a, b) \in \{(w, v), (v, u), (u, w)\}$, then $\bb(u, w)$ is contained in $\d$-neighborhood of $\bb(w, v) \cup \bb(v, u)$ with respect to $d_{(T, \B)}$.
\enddefinition

The following theorem is one of the basic facts in the theory of Gromov hyperbolic spaces. A proof can be seen in \cite{Vaisala} for example.

\thm\label{HYP.thm00}
$(T, \B)$ is $\eta$-hyperbolic for some $\eta > 0$ if and only if all the geodesic triangles in $(T, \B)$ are $\d$-thin for some $\d > 0$.
\endthm

The next theorem gives a criterion of the hyperbolicity of the resolution $(T, \B)$. It has explicitly stated and proven by Lau and Wang in \cite{LauWang}. However, Elek had already used essentially the same idea in \cite{Elek1} to construct a hyperbolic graph which is quai-isometric to the hyperbolic cone of a compact metric space. In fact, we are going to recover their works as a part of our general framework later in this section. 

\thm\label{HYP.thm100}
The resolution $(T, \B)$ of $X$ is hyperbolic if and only if there exists $L \ge 1$ such that
\begin{equation}\label{HYP.eq10}
d_{(T, \B)}(w, v) \le L
\end{equation}
for any horizontally minimal pair $(w, v) \in \cup_{m \ge 1} ((T)_m \times (T)_m)$.
\endthm

\remark
As is shown in the proof, if all the geodesic triangles in $(T, \B)$ are $\d$-thin, then $L$ can be chosen as $4\d + 1$. Conversely, if \eqref{HYP.eq10} is satisfied, then $(T, \B)$ is $\frac 32L$-hyperbolic.
\endremark

Since our terminologies and notations differ much from those in \cite{Elek1} and \cite{LauWang}, we are going to present a proof of Theorem~\ref{HYP.thm100} for reader's sake.

\demo
Assume that all the geodesic triangles in $(T, \B)$ are $\d$-thin. Let $(w, v) \in (T)_m$ be horizontally minimal. Consider the geodesic triangle consists of $p_1 = (w, \pi(w), \ldots, \pi^m(w))$, which is the vertical geodesic between $w$ and $\phi$, $p_2 = (v, \pi(v), \ldots, \pi^m(v))$, which is the vertical geodesic between $v$ and $\phi$, and $p_3 = (u(1), \ldots, u(k + 1))$, which is the horizontal geodesic between $w$ and $v$. Since all the geodesic triangles in $(T, \B)$ are $\d$-thin, for any $i$, either there exists $w' \in p_1$ such that $d_{(T, \B)}(w', u(i)) \le \d$ or there exists $v' \in p_2$ such that $d_{(T, \B)}(v', u(i)) \le \d$. Suppose that the former is the case. Since $d_{(T, \B)}(w, w')  = |w| - |w'|$ is the smallest steps from the level $|w| = |u(i)|$ to $|w'|$, it follows that $d_{(T, \B)}(w, w') \le d_{(T, \B)}(w', u(i))$. Hence
\[
d_{(T, \B)}(w, u(i)) \le d_{(T, \B)}(w, w') + d_{(T, \B)}(w', u(i)) \le 2d_{(T, \B)}(w', u(i)) \le 2\d.
\]
Considering the latter case as well, we conclude  that either $d_{(T, \B)}(w, u(i))) \le 2\d$ or $d_{(T, \B)}(v, u(i)) \le 2\d$ for any $i$. This shows that $d_{(T, \B)}(w, v) \le 4\d + 1$.\par
Conversely, assume \eqref{HYP.eq10}. We are going to show that $(T, \B)$ is $\frac 32L$-hyperbolic, namely, 
\begin{equation}\label{HYP.eq20}
(w(1)|w(2)) \ge \min\{(w(2)| w(3)), (w(3)| w(1))\} - \frac 32L
\end{equation}
for any $w(1), w(2), w(3) \in T$. For $(i, j) \in \{(1, 2), (2, 3), (3, 1)\}$, let $\bb_{ij}$ be a bridge between $w(i)$ and $w(j)$ and let $m_{ij}$, $l_{ij}$ and $m_{ji}$ be the lengths of the ascending part, the horizontal part and the descending part respectively. Also set $h_{ij}$ be the height of the bridge $\bb_{ij}$. Then
\[
(w(i)| w(j)) = \frac{h_{ij} + m_{ij} + h_{ij} + m_{ji} - (m_{ij} + m_{ji} + l_{ij})}{2} = h_{ij} - \frac{l_{ij}}2.
\]
Without loss of generality, we may assume that $h_{23} \ge h_{31}$. Then we have three cases;\\
Case 1: $h_{12} \ge h_{23} \ge h_{31}$,\\
Case 2: $h_{23} \ge h_{12} \ge h_{31}$,\\
Case 3: $h_{23} \ge h_{31} \ge h_{12}$.\\
In Case 1 and Case 2, since $h_{31} - h_{12} \le 0$ and $l_{12} \le L$, it follows that 
\[
(w(3)| w(1)) - (w(1)| w(2)) = h_{31} - h_{12} + \frac{l_{12}}2 - \frac{l_{31}}2 \le \frac L2
\]
Thus \eqref{HYP.eq20} holds. In Case 3, let $v(1) \in (T)_{h_{31}}$ belong to the ascending part of $\bb_{12}$ and let $v(2) \in (T)_{h_{31}}$ belong to the descending part of $\bb_{12}$. Moreover, let $\bb_{31}^h$ and $\bb_{23}^h$ be the horizontal parts of $\bb_{31}$ and $\bb_{23}$ respectively. Then the combination of $\bb_{31}^h$ and $\pi^{h_{23} - h_{31}}(\bb_{23}^h)$ gives a chain between $v(1)$ and $v(2)$ whose length is no greater than $l_{31} + l_{23}$. Since the segment of $\bb_{12}$ connecting $v(1)$ and $v(2)$ is a geodesic, we have
\[
2(h_{31} - h_{12}) + l_{12} \le l_{31} + l_{23} \le 2L.
\]
Therefore, it follows that $h_{31} - h_{12} \le L$. This implies
\[
(w(3)| w(1)) - (w(1)| w(2)) = h_{31} - h_{12} + \frac{l_{12}}2 - \frac{l_{31}}2 \le L + \frac L2 = \frac 32L.
\]
Thus we have obtained \eqref{HYP.eq20} in this case as well.
\enddemo

Note that so far weight functions play no role in the statements of results in this section. In order to take weight functions into account, we are going to introduce an rearranged resolution $(\tT^{g, r}, \B_{\tT^{g, r}})$ associated with a weight function $g$ and give the definition of hyperbolicity of the weight function $g$ in terms of the rearranged resolution. 

\definition\label{HYP.def20}
Let $g \in \G(T)$ and let $r \in (0, 1)$. For $m \ge 0$, define $(\tT^{g, r})_m = \LL_{r^m}^g$ and
\[
\tT^{g, r} = \bigcup_{m \ge 0}(\tT^{g, r})_m.
\]
$\tT^{g, r}$ is naturally equipped with a tree structure inherited from $T$. Define $K_{\tT^{g, r}}: \tT^{g, r} \to \C(X, \O)$ by $K_{\tT^{g, r}} = K|_{\tT^{g, r}}$. The collection of geodesic rays of the tree $\tT^{g, r}$ starting from $\phi$ is denoted by $\SS_{\tT^{g, r}}$. Define $\s_{\tT^{g, r}}: \SS_{\tT^{g, r}} \to X$ by $\s_{\tT^{g, r}}(\omega) = \cap_{m \ge 0} K_{\omega(m)}$ for any $\omega = (\phi, \omega(1), \ldots) \in \SS_{\tT^{g, r}}$. For any $w \in \LL_{r^{m + 1}}^g$, the unique $v \in \LL_{r^m}^g$ satisfying $w \in T_v$ is denoted by $\pi^{g, r}(w)$. Also we set $S^{g, r}(w) = \{v | v \in \LL_{r^{m + 1}}^d, v \in T_w\}$ for $w \in \LL_{r^m}^g$. Define the horizontal edges of $\tT^{g, r}$ as
\[
E^h_{g, r} = \bigcup_{n \ge 1}\{(w, v)| w, v \in (\tT^{g, r})_n, K_w \cap K_v \neq \emptyset\}.
\]
Moreover, we define the totality of the horizontal and vertical edges $\B_{\tT^{g, r}}$ by
\[
\B_{\tT^{g, r}} = \{(w, v)| (w, v) \in E^h_{g, r}\,\,\text{or}\,\,w = \pi^{g, r}(v)\,\,\text{or}\,\,v = \pi^{g, r}(w)\}.
\]
The graph $(\tT^{g, r}, \B_{\tT^{g, r}})$ is called the rearranged resolution of $X$ associated with the weight function $g$.
\enddefinition

\remark
Even if $m \neq n$, it may happen that $\LL_{r^m}^g \cap \LL_{r^n}^g \neq \emptyset$. In such a case, for $w \in \LL_{r^m}^g \cap \LL_{r^n}^g$, we regard $w \in (\tT^{g, r})_m$ and $w \in (\tT^{g, r})_n$ as different elements in $\tT^{g, r}$. More precisely, the exact definition of $\tT^{g, r}$ should be $\tT^{g, r} = \cup_{m \ge 0} (\{m\} \times \LL_{r^m}^g)$ and the associated partition $K_{\tT^{g, r}}: \tT^{g, r} \to \C(X,\O)$ is defined as $K_{\tT^{g, r}}((m, w)) = K_w$.
\endremark

\remark
$\SS_{\tT^{g, r}}$ and $\s_{\tT^{g, r}}$ can be naturally identified with $\SS$ and $\s$ respectively. 
\endremark

\definition\label{HYP.def25}
A weight function $g$ is said to be hyperbolic if and only if the rearranged resolution $(\tT^{g, r}, \B_{\tT^{g, r}})$ is hyperbolic for some $r \in (0, 1)$.
\enddefinition

The next theorem shows that the hyperbolicity of a weigh function $g$ is equivalent to the existence of a ``visual metric'' associated with $g$. It also implies that the quantifier ``for some $r \in (0, 1)$'' in Definition~\ref{HYP.def25} can be replaced  by ``for any $r \in (0, 1)$''.

\thm\label{HYP.thm20}
Let $g$ be a weight function. Then the following three conditions are equivalent:\\
{\rm (1)}\,\,There exists $M \ge 1$ such that $({\rm EV})_M$ is satisfied, i.e. there exist $d \in \D(X, \O)$ and $\a > 0$ such that $d$ is $M$-adapted to $g^{\a}$.\\
{\rm (2)}\,\,The weight function $g$ is hyperbolic.\\
{\rm (3)}\,\,$(\tT^{g, r}, \B_{\tT^{g, r}})$ is hyperbolic for any $r \in (0, 1)$.\par
Moreover, if any of the above conditions is satisfied, then there exist $c_1, c_2 > 0$ such that
\begin{equation}\label{HYP.eq30}
c_1\d^g_M(x, y) \le r^{(x|y)_{\tT^{g, r}}} \le c_2\d^g_M(x, y).
\end{equation}
for any $x, y \in X$, where 
\begin{multline*}
(x|y)_{\tT^{g, r}} = \sup\Big\{\lim_{n, m \to \infty} (\omega(n)|\tau(m))_{(\tT^{g, r}, \B_{\tT^{g, r}}), \phi}\Big|\\ \omega = (\phi, \omega(1), \ldots), \tau = (\phi, \tau(1), \ldots) \in \SS_{\tT^{g, r}}, \s_{\tT^{g, r}}(\omega) = x, \s_{\tT^{g, r}}(\tau) = y\Big\}.
\end{multline*}
\endthm

\remark
The proof of Theorem~\ref{HYP.thm20} shows that if every geodesics triangle of $(\tT^{g, r}, \B_{\tT^{g, r}})$ is $\eta$-thin, then $\rm (EV)_M$ is satisfied for $M = \min\{m| m \in \BbN, 4\eta + 1 \le m\}$.
\endremark

By \eqref{HYP.eq30}, if $d$ is $M$-adapted to $g^{\a}$, then there exist $c_1, c_2 > 0$ such that
\[
c_1d(x, y) \le (r^{\a})^{(x|y)_{\tT}} \le c_2d(x, y)
\]
for any $x, y \in X$. Then the metric $d$ is called a visual metric on the hyperbolic boundary $X$ of $(\tT^{g, r}, \B_{\tT^{g, r}})$ in the framework of Gromov hyperbolic metric spaces. See \cite{BuyaloSchr} and \cite{MacTyson} for example.\par
About the original resolution $(T, \B)$, we have the following corollary.

\cor\label{HYP.cor10}
For $r \in (0, 1)$, define a weight function $h_r$ by $h_r(w) = r^{-|w|}$ for $w \in T$. Then $(T, \B)$ is hyperbolic if and only if there exist a metric $d \in \D(X, \O)$, $M \ge 1$ and $r \in (0, 1)$ such that $d$ is $M$-adapted to $h_r$.
\endcor

To show the hyperbolicity of a weigh function $g$, the existence of adapted metric is (an equivalent condition as we have seen in Theorem~\ref{HYP.thm20} but) too restrictive in some cases. In fact, the notion of  ``weakly adapted'' metric is often more useful as we will see in Example~\ref{HYP.ex10} and \ref{HYP.ex20}. 

\definition\label{HYP.def50}
Let $d \in \D(X, \O)$. For $r \in (0, 1]$, $s > 0$ and $x \in X$, define
\begin{multline*}
\widetilde{B}_d^r(x, s) = \{y| y \in B_d(x, s), \text{there exists a horizontal chain $(w(1), \ldots, w(k))$}\\ \text{in $\LL_r^g$ between $x$ and $y$ such that $K_{w(i)} \cap B_d(x, s) \neq \emptyset$ for any $i = 1, \ldots, k$}\}.
\end{multline*}
A metric $d \in \D(X, \O)$ is said to be weakly $M$-adapted to a weight function $g$ if and only if there exist $c_1, c_2 > 0$ such that
\[
\widetilde{B}_d^r(x, c_1r) \subseteq U^g_M(x, r) \subseteq B_d(x, c_2r)
\]
for any $x \in X$ and $r \in (0, 1]$.
\enddefinition

Since $\widetilde{B}_d^r(x, cr) \subseteq B_d(x, cr)$, we immediately have the following fact.

\prop\label{HYP.prop25}
If $d \in \D(X, \O)$ is $M$-adapted to a weight function $g$, then it is weakly $M$-adapted to $g$.
\endprop

The next proposition gives a sufficient condition for a metric being weakly adapted, which will be applied in Examples~\ref{HYP.ex10} and \ref{HYP.ex20}.

\prop\label{HYP.prop30}
If there exist $c_1, c_2 > 0$ and $M \in \BbN$ such that
\begin{equation}\label{HYP.eq80}
\diam{K_w, d} \le c_1r
\end{equation}
for any $w \in \LL_r^g$ and 
\begin{equation}\label{HYP.eq90}
\#(\{w| w \in \LL_r^g, B_d(x, c_2r) \cap K_w \neq \emptyset\}) \le M + 1
\end{equation}
for any $x \in X$ and $r \in (0, 1]$, then $d$ is weakly $M$-adapted to $g$.
\endprop

\demo
Assume that $y \in U_M^g(x, r)$. Then there exists a chain $(w(1), \ldots, w(M + 1))$ between $x$ and $y$ in $\LL^g_r$. Choose $x_i \in K_{w(i)} \cap K_{w(i + 1)}$ for any $i$. Set $x_0 = x$ and $x_{M + 1} = y$ Then by \eqref{HYP.eq80}, it follows that
\[
d(x, y) \le \sum_{i = 0}^M d(x_i, x_{i + 1}) \le (M + 1)c_1r.
\]
Hence $y \in B_d(x, (M + 2)c_1r)$. This implies $U^g_M(x, r) \subseteq B_d(x, (M + 2)c_1r)$. \par
Next, let $y \in \widetilde{B}^r_d(x, c_2r)$. Then there exists  a chain $(w(1), \ldots, w(k))$ between $x$ and $y$ in $\LL^g_r$ such that $K_{w(i)} \cap B_d(x, c_2r) \neq \emptyset$ for any $i = 1, \ldots, k$. We may assume that $w(i) \neq w(j)$ if $i \neq j$. Then by \eqref{HYP.eq90}, we see that $k \le M + 1$. This yields that $y \in U_M^g(x, r)$. Thus we have shown that $d$ is weakly $M$-adapted to $g$.
\enddemo

\prop\label{HYP.prop20}
Let $g$ be a weight function. If there exists a metric $d \in \D(X, \O)$ that is weakly $M$-adapted to $g^{\a}$ for some $M \ge 1$ and $\a > 0$, then $(\tT^{g, r}, \B_{\tT^{g, r}})$ is hyperbolic for any $r \in (0, 1]$.
\endprop

 There have been several works on the construction of a hyperbolic graph whose hyperbolic boundary coincides with a given compact metric space. For example, Elek\cite{Elek1} has studied the case for arbitrary compact subset of $\BbR^n$ and Lau-Wang\cite{LauWang} has considered self-similar sets satisfying the open set condition. Due to above proposition, we may integrate these works into our framework. See Example~\ref{HYP.ex10} and \ref{HYP.ex20} for details.\par

For ease of notations, we use $\tT$, $\tpi$, $\SS_{\tT}$, $\s_{\tT}$ and $\B_{\tT}$ to denote $\tT^{g, r}$, $\pi^{g, r}$, $\SS_{\tT^{g, r}}$, $\s_{\tT^{g, r}}$ and $\B_{\tT^{g, r}}$ respectively. Moreover, we write $d_{\tT} = d_{(\tT, \B_{\tT})}$, which is the geodesic metric of $(\tT, \B_{\tT})$.

\demo
Assume that there exist a metric $d \in \D(X, \O)$, $M \ge 1$ and $\a > 0$ such that $d$ is weakly $M$-adapted to $g^{\a}$. Then there exist $c_1, c_2 > 0$ such that
\[
\widetilde{B}_d^r(x, c_1r) \subseteq U^{g^{\a}}_M(x, r) \subseteq B_d(x, c_2r)
\]
for any $x\in X$ and $r \in (0, 1]$. In particular, $d(x, y) \le c_2\d^{g^{\a}}_M(x, y)$ for any $x, y \in X$. Suppose that  $w, v \in (\tT)_n$, $x_1 \in K_w$, $x_2 \in K_v$ and $(w, v) \in E^h_{g, r^n}$. Since $\d_1^{g}(x_1, x_2) \le r^n$, it follows that 
\begin{equation}\label{HYP.eq70}
d(x_1, x_2) \le c_2r^{{\a}n}.
\end{equation}
Let $m \ge 1$ and fix $r \in (0, 1)$. Suppose that there exists $w_*, v_* \in (\tT)_n$ such that $(w_*, v_*)$ is horizontally minimal  and $d_{\tT}(w_*, v_*) = 3m + 1$. Let $(w(1), \ldots, w(3m + 1))$ be the horizontal geodesic of $(\tT, \B_{\tT})$ between $w_*$ and $v_*$. Let $x \in K_{w_*}$ and let $y \in K_{v_*}$. By \eqref{HYP.eq70} for any $z \in \cup_{i = 1}^{3m + 1}K_{w_*}$, 
\[
d(x, z) \le c_2(3m + 1)r^{\a{n}}.
\]
If $m$ is sufficiently large, then $c_2r^{{\a}m}(3m + 1) < c_1$. Hence
\[
d(x, z) < c_1r^{\a(n - m)}.
\]
Set $v(i) = \tpi^m(w(i))$ for $i = 1, \ldots, 3m + 1$, then $(v(1), \ldots, v(3m + 1))$ is a horizontal chain in $(T)_{n - m}$ between $x$ and $y$ and $K_{v(i)} \cap B_d(x, c_1r^{\a(n - m)}) \supseteq K_{w(i)} \cap B_d(x, c_1r^{\a(n - m)}) \neq \emptyset$. Therefore 
\[
y \in \widetilde{B}^r_d(x, c_1r^{\a(n - m)}) \subseteq U^{g^{\a}}_M(x, r^{\a(n - m)}) = U^g_M(x, r^{n - m}).
\]
So there exists a horizontal chain $(u(1), \ldots, u(M + 1))$ in $(\tT)_{n - m}$ between $x$ and $y$. Combining this horizontal chain with vertical geodesics $(w_*, \ldots, \tpi^{m}(w_*))$ and $(\tpi^m(v_*), \ldots, v_*)$, we have a chain of $(\tT, \B_{\tT})$ between $w_*$ and $v_*$ whose length is $M + 2 + 2m$. Therefore, 
\[
M + 2 + 2m \ge d_{\tT}(w, v) = 3m.
\]
Hence $M + 2 \ge m$.  Applying Theorem~\ref{HYP.thm100} to $(\tT, \B_{\tT})$, we verify that $(\tT, \B_{\tT})$ is hyperbolic.
\enddemo

\demo[Proof of Theorem~\ref{HYP.thm20}]
(1) $\Rightarrow$ (3)\,\, Proposition~\ref{HYP.prop20} suffices.\\
(3) $\Rightarrow $ (2)\,\,This is immediate.\\
(2) $\Rightarrow$ (1)\,\,
Assume that all the geodesic triangles in $(\tT, \B_{\tT})$ are $\d$-thin. Set $L = \min\{m| m \in \BbN, 4\d + 1 \le m\}$. For ease of notation, we use $(w|v)_{\tT}$ to denote the Gromov product of $w$ and $v$ in $(\tT, \B_{\tT})$ with respect to $\phi$. Let $x \neq y \in X$ and let $\omega = (\phi, \omega(1), \ldots), \tau = (\phi, \tau(1), \ldots) \in \SS_{\tT}$ satisfy $\s_{\tT}(\omega) = x$ and $\s_{\tT}(\tau) = y$. Applying Proposition~\ref{HYP.prop10} to $\tT$, we see that there exists $m_* \in \BbN$ such that $d_{\tT}(\omega(m), \tau(m)) > L$ for any $m \ge m_*$. Let $\bb$ be a bridge between $\omega(m_*)$ and $\tau(m_*)$. If $k_*$ is the height of $\bb$ and $(\omega(k_*), w(1). \ldots, w(l - 1), \tau(k_*))$ is the horizontal part of $\bb$, then $\bb$ is the concatenation of  $(\omega(m_*), \ldots, \omega(k_*))$, $(\omega(k_*), w(1). \ldots, w(l - 1), \tau(k_*))$ and $(\tau(k_*), \ldots, \tau(m_*))$. If $m, n \ge m_*$, then $(\omega(m), \ldots, \omega(k_*), w(1), \ldots, w(l - 1), \tau(k_*), \ldots, \tau(n))$ is a bridge between $\omega(m)$ and $\tau(n)$. Therefore,
\[
(\omega(m)|\tau(n))_{\tT} = k_* - \frac l2.
\]
Hence, if we define
\[
(\omega, \tau)_{\tT} = \lim_{m, n \to \infty} (\omega(m)|\tau(n))_{\tT},
\]
then $(\omega|\tau)_{\tT} = k_* - \frac l2$. \par
 Applying Theorem~\ref{HYP.thm100} to $(\tT, \B_{\tT})$, we see that the length of any horizontal geodesic is no greater than $L$. Since the length of the horizontal part of $\bb$ is $l$, it follows that $l \le L$. Therefore letting $s_* = \d_L^g(x, y)$, then we see that
\begin{equation}\label{HYP.eq40}
s_* \le r^{k_*} \le r^{k_* - \frac l2} = r^{(\omega|\tau)_{\tT}}.
\end{equation}
 Choose $n_*$ so that $r^{k_* + n_*} \ge s_* > r^{k_* + n_* + 1}$. Then there exists a horizontal chain $(v(1), \ldots, v(L + 1))$ in $(\tT)_{k_* + n_*}$ such that $x \in K_{v(1)}$ and $y \in K_{v(L + 1)}$. Hence $(\omega(k_* + n_*), v(1), \ldots, v(L + 1), \tau(k_* + n_*))$ is a chain between $\omega(k_* + n_*)$ and $\tau(k_* + n_*)$. Comparing this chain with $(\omega(k_* + n_*), \ldots, \omega(k_*), w(1), \ldots, w(l - 1), \tau(k_*), \ldots, \tau(k_* + n_*))$, we obtain
\[
2n_* + l \le L + 2.
\]
This implies $k_* + n_* + 1 - \frac L2 - 2 \le k_* - \frac l2$. Therefore
\begin{equation}\label{HYP.eq50}
r^{(\omega|\tau)_{\tT}} = r^{k_* - \frac l2} \le r^{k_* + n_* + 1}r^{-\frac L2 - 2} \le r^{-\frac L2 - 2}s_*.
\end{equation}
Set $c_1 = 1$ and $c_2 = r^{-\frac L2 - 2}$. Then we have
\[
c_1\d_L^g(x, y) \le r^{(\omega|\tau)_{\tT}} \le c_2\d_L^g(x, y).
\]
Define $(x|y)_{\tT} = \sup\{(\omega|\tau)_{\tT}| \omega, \tau \in \SS_{\tT}, \s_{\tT}(\omega) = x, \s_{\tT}(\tau) = y\}$. Then 
\begin{equation}\label{HYP.eq60}
c_1\d_L^g(x, y) \le r^{(x|y)_{\tT}} \le c_2\d_L^g(x, y).
\end{equation}
It is known that if $(\tT, \B_{\tT})$ is hyperbolic, then $r^{(x|y)_{\tT}}$ is a quasimetric. Hence by \eqref{HYP.eq60}, $\d_L^g(x, y)$ is a quasimetric as well. Thus we have obtained ${\rm(EV2)}_M$.
\enddemo

In short, in the above reasonings, we have two steps:\\
(1)\,\,The existence of weakly adapted metric $d$ implies the hyperbolicity of $g$.\\
(2)\,\,The hyperbolicity of $g$ implies the existence of an adapted metric $\rho$.\\
It is notable that the original weakly adapted metric $d$ may essentially differ from the adapted metric $d$. In fact, in Example~\ref{HYP.ex10}, we are going to present an explicit example where no power of the original weakly adapted metric is bi-Lipschitz equivalent to any adapted metric. 

\demo[Proof of Corollary~\ref{HYP.cor10}]
Note that $(T, \B) = (\tT^{h_r, r}, \B_{\tT^{h_r, r}})$ for any $r \in (0, 1)$. Assume that $(T, \B)$ is hyperbolic. By Theorem~\ref{HYP.thm20}, there exist $d \in \D(X, \O)$, $M \ge 1$ and $\a > 0$ such that $d$ is adapted to $(h_{1/2})^{\a}$. Since $(h_{1/2})^{\a} = h_{2^{-\a}}$, we have the desired statement with $r = 2^{-\a}$. Conversely, with the existence of $d$, $M$ and $r$, Theorem~\ref{HYP.thm20} implies that $(\tT^{h_r, s}, \B_{\tT^{h_r, s}})$ is hyperbolic for any $s \in (0, 1)$. Letting $s = r$, we see that $(T, \B)$ is hyperbolic.
\enddemo

To end this section, we are going to integrate the works by Elek\cite{Elek1} and Lau-Wang\cite{LauWang} into our framework.

\example\label{HYP.ex10}
In Example~\ref{PAS.ex20}, we have obtained a partition of a compact metric space in $\BbR^N$ corresponding to the hyperbolic graph constructed by Elek in \cite{Elek1}. In fact, we have obtained two graphs $(T, \widetilde{\B})$ and $(T, \B)$ satisfying $\widetilde{\B} \supseteq \B$. The former coincides with Elek's graph and the latter is the resolution associated with the partition constructed from the Dyadic cubes. In this example, using Propositions~\ref{HYP.prop30} and ~\ref{HYP.prop20}, we are going to show the hyperbolicity of the graph $(T, \B)$. The hyperbolicity of the original graph $(T, \widetilde{\B})$ may be shown in a similar fashion. \par
Let $X$ be a compact subset of $[0, 1]^n$ and let $(T, \A, \phi)$ be the tree associated with $X$ constructed from the dyadic cubes in Example~\ref{PAS.ex20}. Also let $K: T \to \C(X, \O)$ be the partition of $X$ parametrized by $(T, \A, \phi)$ given in Example~\ref{PAS.ex20}. Set $g(w) = 2^{-m}$ if $w \in (T)_m$. Then $\LL^g_{r} = (T)_m$ if and only if $2^{-m} \le r < 2^{-m + 1}$. Let $d_*$ be the Euclidean metric. Then for any $w \in \LL^g_r$, 
\[
\diam{K_w, d_*} \le 2\sqrt{n}r.
\]
This shows \eqref{HYP.eq80}. Moreover, if $w \in \LL^g_r$ and $K_w \cap B_{d_*}(x, cr) \neq \emptyset$, then $C(w) \subseteq B_{d_*}(x, (c + 2\sqrt{N})r)$. Note that $|C(w)|_n = 2^{-mn}$ for any $w \in (T)_m$, where $|\cdot|_n$ is the $n$-dimensional Lebesgue measure. Therefore, if $2^{-m} \le r < 2^{-m + 1}$, then
\begin{multline*}
\#(\{w| w \in \LL_r^g, B_{d_*}(x, cr) \cap K_w \neq \emptyset\}) \le \frac{|B_{d_*}(x, (c + 2\sqrt{n})r)|_n}{2^{-mn}}\\ = |B_{d_*}(0, 1)|_N(c + 2\sqrt{n})^n(2^mr)^n \le  |B_{d_*}(0, 1)|_n(c + \sqrt{n})^n2^n.
\end{multline*}
Therefore choosing $M \in \BbN$ so that $|B_{d_*}(0, 1)|_N(c + \sqrt{n})^n2^n \le M + 1$, we have \eqref{HYP.eq90}. Hence by Proposition~\ref{HYP.prop30}, (the restriction of) $d_*$ is weakly $M$-adapted to $g$. Since $(\tT^{g, \frac12}, \B_{\tT^{g, \frac12}}) = (T, \B)$,  Proposition~\ref{HYP.prop20} yields that $(T, \B)$ is hyperbolic and its hyperbolic boundary coincides with $X$. \par
As we have mentioned above, in this example, the weakly adapted metric $d_*$ is not necessarily adapted to any power of $g$. For example, let
\begin{multline*}
X = [0, 1] \bigcup \{(t, t)| t \in [0, 1]\} \bigcup \\
\bigg(\bigcup_{m \ge  1}\Big\{\Big(\frac 1{2^m}, s\Big)\Big| s \in \Big[0, \frac 1{2^m}\Big]\backslash\Big(\frac{1 - \e_m}{2^m}, \frac{1 + \e_m}{2^m}\Big)\Big\}\bigg),
\end{multline*}
where $\displaystyle \e_m = \frac 1{2^{m^2}}$. Set $\displaystyle x_m = \Big(\frac 1{2^m}, \frac{1- \e_m}{2^m}\Big)$ and $\displaystyle y_m = \Big(\frac 1{2^m}, \frac{1+ \e_m}{2^m}\Big)$. Then $\displaystyle d_*(x_m, y_m) = \frac{\e_m}{2^{m - 1}}$ and $\displaystyle \d^g_1(x_m, y_m) = \frac 1{2^{m - 1}}$. Then 
\[
\frac{d_*(x_m, y_m)}{\d^{g^{\a}}_1(x_m, y_m)} = 2^{\a(m - 1) - m^2 - m + 1} \to 0
\]
as $m \to \infty$ for any $\a > 0$. Thus for any $\a > 0$, the Euclidean metric is not bi-Lipschitz equivalent to any metric adapted to $g^{\a}$.
\endexample

\example\label{HYP.ex20}
Let $X$ be the self-similar set associated with the collection of contractions $\{F_1, \ldots, F_N\}$ and let $K: T^{(N)} \to X$ be the partition of $K$ parametrized by $(T^{(N)}, \A^{(N)}, \phi)$ introduced in Example~\ref{PAS.ex10}. We write $T = T^{(N)}$ for simplicity. In this example, we further assume that for any $i = 1, \ldots, N$, $F_i: \BbR^n \to \BbR^n$ is a similitude, i.e. there exist an orthogonal matrix $A_i$, $r_i \in (0, 1)$ and $a_i \in \BbR^n$ such that $F_i(x) = r_iA_ix + a_i$. Furthermore, we assume that the open set condition holds, i.e. there exists a nonempty open subset $O$ of $\BbR^n$ such that $F_w(O) \subseteq O$ for any $w \in T$ and $F_w(O) \cap F_v(O) = \emptyset$ if $w, v \in T$ and $T_w \cap T_v = \emptyset$. Define $g(w) = r_{w_1}{\cdots}r_{w_m}$ for any $w = w_1\ldots{w_m} \in T$. In this case, the conditions \eqref{HYP.eq80} and \eqref{HYP.eq90} have been known to hold for the Euclidean metric $d_*$. See \cite[Proposition~1.5.8]{AOF} for example. Hence Proposition~\ref{HYP.prop30} implies that $d_*$ is weakly $M$-adapted to $g$ for some $M \in \BbN$. Using Proposition~\ref{HYP.prop20}, we see that $(\tT^{g, r}, \B_{\tT^{g, r}})$ is hyperbolic for any $r \in (0, 1)$ and hence the self-similar set $X$ is the hyperbolic boundary of $(\tT^{g, r}, \B_{\tT^{g, r}})$. This fact has been shown by Lau and Wang in \cite{LauWang}. As in the previous example, the Euclidean metric is not necessarily a visual metric in this case.
\endexample

\part{Relations of weight functions}
\setcounter{equation}{0}
\section{Bi-Lipschitz equivalence}\label{BLE}
In this section, we define the notion of bi-Lipschitz equivalence of weight functions. Originally the definition, Definition~\ref{BLE.def10}, only concerns the tree structure $(T, \A, \phi)$ and has nothing to do with a partition of a space. Under proper conditions, however, we will show that the bi-Lipschitz equivalence of weight functions is identified with
\begin{itemize}
\item
absolutely continuity with uniformly bounded Radon-Nikodym derivative from below and above between measures in \ref{BLE1}.
\item
usual bi-Lipschitz equivalence between metrics in \ref{BLE2}.
\item
Ahlfors regularity of a measure with respect to a metric in \ref{BLE3}.
\end{itemize}

As in the previous sections, $(T, \A, \phi)$ is a locally finite tree with a reference point $\phi$, $(X, \O)$ is a compact metrizable topological space with no isolated point and $K: T \to \C(X, \O)$ is a partition of $X$ parametrized by $(T, \A, \phi)$.

\definition\label{BLE.def10}
Two weight functions $g, h \in \G(T)$ are said to be bi-Lipschitz equivalent if and only if there exist $c_1, c_2 > 0$ such that
\[
c_1g(w) \le h(w) \le c_2g(w)
\]
for any $w \in T$. We write $g \bl h$ if and only if $g$ and $h$ are bi-Lipschitz equivalent.
\enddefinition

By the definition, we immediately have the next fact.

\prop\label{BLE.prop10}
The relation $\bl$ is an equivalent relation on $\G(T)$.
\endprop

\subsection{bi-Lipschitz equivalence of measures}\label{BLE1}
As we mentioned above, the bi-Lipschitz equivalence between weight functions can be identified with other properties according to classes of weight functions. First we consider the case of weight functions associated with measures.

\definition\label{BLE.def20}
Let $\mu, \nu \in \M_P(X, \O)$. We write $\mu \ac \nu$ if and only if there exist $c_1, c_2 > 0$ such that
\begin{equation}\label{BLE.eq170}
c_1\mu(A) \le \nu(A) \le c_2\mu(A)
\end{equation}
for any Borel set $A \subseteq X$.
\enddefinition

It is easy to see that $\ac$ is an equivalence relation and $\mu \ac \nu$ if and only if $\mu$ and $\nu$ are mutually absolutely continuous and the Radon-Nikodym derivative $\frac{d\nu}{d\mu}$ is uniformly bounded from below and above.

\thm\label{BLE.thm20}
Assume that the partition $K: T \to \C(X, \O)$ is strongly finite. Let $\mu, \nu \in \M_P(X, \O)$. Then $g_{\mu} \bl g_{\nu}$ if and only if $\mu \ac \nu$. Moreover, the natural map  $\M_P(X, P)/\!\!\!\ac\,\,\to\,\,\G(X)/\!\!\!\bl$ given by $[g_{\mu}]_{\bl}$ is injective, where $[\,\cdot\,]_{\bl}$ is the equivalence class under $\bl$.
\endthm

\demo
By \eqref{BLE.eq170}, we see that $\a_1\nu(K_w) \le \mu(K_w) \le \a_2\nu(K_w)$ and hence $g_{\mu} \bl g_{\nu}$. Conversely, if
\[
c_1\mu(K_w) \le \nu(K_w) \le c_2\mu(K_w)
\]\
for any $w \in T$. Let $U \subset X$ be an open set. Assume that $U \neq X$. For any $x \in X$, there exists $w \in T$ such that $x \in K_w \subseteq U$. Moreover, if $K_w \subseteq U$, then there exists $m \in \{1, \ldots, |w|\}$ such that $K_{[w]_m} \subseteq U$ but $\sd{K_{[w]_{m - 1}}}{U} \neq \emptyset$. Therefore, if 
\[
T(U) = \{w| w \in T, K_w \subseteq U, \sd{K_{\pi(w)}}{U} \neq \emptyset\},
\]
then $T(U) \neq \emptyset$ and $U = \cup_{w \in T(U)} K_w$. Now, since $K$ is strongly finite, there exists $N \in \BbN$ such that $\#(\s^{-1}(x)) \le N$ for any $x \in X$. Let $y \in U$. If $w(1), \ldots, w(m) \in T(U)$ are mutually different and $y \in K_{w(m)}$ for any $i = 1, \ldots, m$, then there exists $\omega(i) \in \SS_{w(i)}$ such that $\s(\omega(i)) = y$ for any $i = 1, \ldots, m$. Hence $\#(\s^{-1}(y)) \ge m$ and therefore $m \le N$.  By Proposition~\ref{MEM.prop10}, we see that
\begin{align*}
\nu(U) \le \sum_{w \in T(U)} \nu(K_w) \le\sum_{w \in T(U)} c_2\mu(K_w) \le c_2N\mu(U)\\
\mu(U) \le \sum_{w \in T(U)} \mu(K_w) \le \sum_{w \in T(U)} \frac 1{c_1}\nu(K_w) \le \frac{N}{c_1}\nu(U).
\end{align*}
Hence letting $\a_1 = c_1/N$ and $\a_2 = c_2N$, we have
\[
\a_1\mu(U) \le \nu(U) \le \a_2\mu(U)
\]
for any open set $U \subseteq X$. Since $\mu$ and $\nu$ are Radon measures, this yields \eqref{BLE.eq170}.
\enddemo

\subsection{bi-Lipschitz equivalence of metrics}\label{BLE2}

Under the tightness of weight functions defined below, we will translate bi-Lipschitz equivalence of weight functions to the relations between ``balls'' and ``distances'' associated with weight functions in Theorem~\ref{BLE.thm10}. The tightness of a weight function ensures that $\d_M^g$ is comparable with $g$, i.e the diameter with respect to $\d_M^g$ of $K_w$ is bi-Lipschitz equivalent to $g$.

\definition\label{BLES.def10}
A weight function $g$ is called tight if and only if for any $M \ge 0$, there exists $c > 0$ such that
\[
\sup_{x,  y \in K_w} \d_M^g(x, y) \ge cg(w)
\]
for any $w \in T$.
\enddefinition

\prop\label{BLE.prop300}
Let $g$ and $h$ be weight functions. Assume that $g \bl h$.  If $g$ is tight and $g \bl h$, then $h$ is tight.
\endprop

\demo
Since $g \bl h$, there exist $\c_1, \c_2 > 0$ such that $\c_1g(w) \le h(w) \le \c_2g(w)$ for any $w \in T$. Therefore,
\[
\c_1D_M^g(x, y) \le D_M^h(x, y) \le \c_2D_M^g(x, y)
\]
for any $x, y \in X$ and $M \ge 0$. By Proposition~\ref{MAG.prop10}, for any $M \ge 0$, there exist $c_1, c_2 > 0$ such that
\[
c_1\d_M^g(x, y) \le \d_M^h(x, y) \le c_2\d_M^g(x, y)
\]
for any $x, y \in X$. Hence
\[
\sup_{x, y \in K_w} \d_M^h(x, y) \ge c_1\sup_{x, y \in K_w}\d_M^g(x, y) \ge c_1cg(w) \ge c_1c(\c_2)^{-1}h(w)
\]
for any $w \in T$. Thus $h$ is tight.
\enddemo

Any weight function induced from a metric is tight.

\prop\label{BLE.prop310}
Let $d \in \D(X, \O)$. Then $g_d$ is tight.
\endprop

\demo
Let $x, y \in X$ and let $(w(1), \ldots, w(M + 1)) \in \CH_K(x, y)$. Set $x_0 = x$ and $x_{M + 1} = y$. For each $i = 1, \ldots, M$, choose $x_i \in K_{w(i)} \cap K_{w(i + 1)}$. Then
\[
\sum_{i = 1}^{M + 1} g_d(x) \ge \sum_{i = 1}^{M + 1} d(x_{i - 1}, x_i) \ge d(x, y).
\]
Using this inequality and Proposition~\ref{MAG.prop10}, we obtain
\[
(M + 1)\d_M^g(x, y) \ge D_M^{g_d}(x, y) \ge d(x, y)
\]
and therefore $(M + 1)\sup_{x, y \in K_w}\d_M^{g_d}(x, y) \ge g_d(w)$ for any $w \in T$. Thus $g_d$ is tight.
\enddemo
Now we give geometric conditions which are equivalent to bi-Lipschitz equivalence of tight weight functions. The essential point is that bi-Lipschitz condition between weight function $g$ and $h$ are equivalent to that between $\d^g_M(\cdot, \cdot)$ and $\d^h_M(\cdot, \cdot)$ in the usual sense as is seen in (BL2) and (BL3). 

\thm\label{BLE.thm10}
Let $g$ and $h$ be weight functions. Assume that both $g$ and $h$ are tight. Then the following conditions are equivalent:\\
{\rm (BL)}\,\,
$g \bl h$.\\
{\rm (BL1)}\,\,
There exist $M_1, M_2$ and $c > 0$ such that
\[
\d^g_{M_1}(x, y) \le c\d^h_0(x, y)\quad\text{and}\quad \d^h_{M_2}(x, y) \le c\d^g_0(x, y)
\]
for any $x, y \in X$.\\
{\rm (BL2)}\,\,
There exist $c_1, c_2 > 0$ and $M \ge 0$ such that
\[
c_1\d^g_M(x, y) \le \d^h_M(x, y) \le c_2\d^g_M(x, y)
\]
for any $x, y \in X$.\\
{\rm (BL3)}\,\,
For any $M \ge 0$, there exist $c_1, c_2 > 0$ such that
\[
c_1\d^g_M(x, y) \le \d^h_M(x, y) \le c_2\d^g_M(x, y)
\]
for any $x, y \in X$.
\endthm

Before a proof of this theorem, we state two notable corollaries of it. The first one is the case when weight functions are induced from adapted metrics. In such a case bi-Lipschitz equivalence of weight functions exactly corresponds to the usual bi-Lipschitz equivalence of metrics.

\definition\label{BLE.def30}
(1)\,\,
Let $d, \rho \in \D(X, \O)$. $d$ and $\rho$ are said to be bi-Lipschitz equivalent, $d \bl \rho$ for short, if and only if there exist $c_1, c_2 > 0$ such that 
\[
c_1d(x, y) \le \rho(x, y) \le c_2d(x, y)
\]
for any $x, y \in X$.\\
(2)\,\,Define
\[
\D_{A}(X, \O) = \{d | d \in \D(X, \O), \text{$d$ is adapted.}\}
\]
\enddefinition

\cor\label{BLE.cor10}
Let $d, \rho \in \D_A(X, \O)$.  Then $g_d \bl g_{\rho}$ if and only if $d \bl \rho$. In particular, the correspondence of $[d]_{\bl}$ with $[g_d]_{\bl}$ gives an well-defined injective map $\D_A(X, \O)/\!\!\!\bl\,\,\to\,\,\G(X)/\!\!\!\bl$.
\endcor

The next corollary shows that an adapted metric is adapted to a weight function if and only if they are bi-Lipschitz equivalent in the sense of weight functions.

\cor\label{BLE.cor20}
Let $d \in \D(X, \O)$ and let $g$ be a weight function. Then $d$ is adapted to $g$ and $g$ is tight if and only if $g_d \bl g$ and $d \in \D_A(X, \O)$.
\endcor

Now we start to prove Theorem~\ref{BLE.thm10} and its corollaries.

\lemma\label{ADD.lemma100}
Let $h$ be a weight function. If $x \in K_w$ and $\sd{K_w}{U^h_0(x, s)} \neq \emptyset$, then $s \le h(w)$.
\endlemma

\demo
If $\pi^n(w) \in \LL_{s, 0}^h(x)$ for some $n \ge 0$, then $U^h_0(x, s) \supseteq K_{\pi^n(w)} \supseteq K_w$. This contradicts the assumption and hence $\pi^n(w) \notin \LL_{s, 0}^h(x)$ for any $n \ge 0$. Therefore there exists $v \in T_w \cap \LL_{s, 0}^h(x)$ such that $|v| > |w|$. Then we have $h(w) \ge h(\pi(v)) > s$.
\enddemo

\prop\label{ADD.prop310}
Let $g$ and $h$ be weight functions. Assume that $g$ is tight. Let $M \ge 0$. If there exists $\a > 0$ such that 
\begin{equation}\label{ADD.eq300}
\a\d_M^g(x, y) \le \d_0^h(x, y)
\end{equation}
for any $x, y  \in X$. Then there exists $c > 0$ such that
\[
cg(w) \le h(w)
\]
for any $w \in T$.
\endprop

\demo
Since $g$ is tight, there exists $\b > 0$ such that,  for any $w \in T$, 
\[
\sd{K_w}{U^g_M(x, {\b}g(w))} \neq \emptyset
\]
for some $x \in K_w$. On the other hand, by \eqref{ADD.eq300}, there exists $\c > 0$ such that $U^g_M(x, s) \supseteq U^h_0(x, {\c}s)$ for any $x \in X$ and $s \ge 0$. Therefore, 
\[
\sd{K_w}{U^h_0(x, {\b\c}g(w))} \neq \emptyset.
\]
By Lemma~\ref{ADD.lemma100}, we have ${\b\c}g(w) \le h(w)$.
\enddemo

\lemma\label{ADD.lemma400}
Let $g$ and $h$ be weight functions. Assume that $g$ is tight. Then the following conditions are equivalent:\\
{\rm (A)}\,\,
There exists $c > 0$ such that $g(w) \le ch(w)$ for any $w \in T$.\\
{\rm (B)}\,\,
For any $M, N \ge 0$ with $N \ge M$, there exists $c > 0$ such that
\[
\d^g_N(x, y) \le c\d^h_M(x, y)
\]
for any $x, y \in X$.\\
{\rm (C)}\,\,
There exist $M, N \ge 0$ and $c > 0$ such that $N \ge M$ and
\[
\d^g_N(x, y) \le c\d^h_M(x, y)
\]
for any $x, y \in X$.
\endlemma

\demo
(A) implies
\begin{equation}\label{ADD.eq100}
D^g_M(x, y) \le cD^h_M(x, y)
\end{equation}
for any $x, y \in X$ and $M \ge 0$. By Proposition~\ref{MAG.prop10}, we see 
\[
\d_M^g(x, y) \le c(M + 1)\d_M^h(x, y)
\]
for any $x, y \in X$.
Since $\d^g_N(x, y) \le \d^g_M(x, y)$, if $N \ge M$, then we have (B). Obviously (B) implies (C). Now assume (C). Then we have $\d_N^g(x, y) \le c\d_0^h(x, y)$. Hence Proposition~\ref{ADD.prop310} yields (A).
\enddemo

\demo[Proof of Theorem~\ref{BLE.thm10}]
Lemma~\ref{ADD.lemma400} immediately implies the desired statement.
\enddemo

\demo[Proof of Corollary~\ref{BLE.cor10}]
Since $d$ and $\rho$ are adapted, by \eqref{MAG.eq10}, there exist $M \ge 1$ and $c > 0$ such that
\begin{align}
c\d^d_M(x, y) &\le d(x, y) \le \d^d_M(x, y), \label{BLE.eq150}\\
c\d^{\rho}_M(x, y) &\le \rho(x, y) \le \d^{\rho}_M(x, y)\label{BLE.eq160}
\end{align}
for any $x, y \in X$. Assume $g_d \bl g_{\rho}$. Since $g_d$ and $g_{\rho}$ are tight, we have (BL3) by Theorem~\ref{BLE.thm10}. Hence by \eqref{BLE.eq150} and \eqref{BLE.eq160},  $d(\cdot, \cdot)$ and $\rho(\cdot, \cdot)$ are bi-Lipschitz equivalent as metrics. The converse direction is immediate.
\enddemo

\demo[Proof of Corollary~\ref{BLE.cor20}]
If $d$ is $M$-adapted to $g$ for some $M \ge 1$, then by (ADa), there exists $c > 0$ such that $d_w \le cg(w)$ for any $w \in K_w$. Moreover, \eqref{MAG.eq10} implies
\[
d(x, y) \ge c_2\d_M^g(x, y)
\]
for any $x, y \in X$, where $c_2$ is independent of $x$ and $y$. Hence the tightness of $g$ shows that there exists $c' > 0$ such that
\[
d_w \ge c_2\sup_{x, y} \d_M^g(x, y) \ge c'g(w).
\]
Thus $g_d \bl g$.  Moreover, by Proposition~\ref{PAS.prop40}, $d$ is adapted. Conversely, assume that $d$ is $M$-adapted and $g_d \bl g$. Then Theorem~\ref{BLE.thm10} implies (BL3) with $h = g_d$. At the same time, since $d$ is $M$-adapted, we have \eqref{MAG.eq20} with $g = g_d$.  Combining these two, we deduce \eqref{MAG.eq20}. Hence $d$ is $M$-adapted to $g$.
\enddemo

\subsection{bi-Lipschitz equivalence between measures and metrics}\label{BLE3}
Finally in this section, we consider what happens if the weight function associated with a measure is bi-Lipschitz equivalent to the weight function associated with a metric.\par
To state our theorem, we need the following notions.

\definition\label{BLE.def40}
(1)\,\,
A weight function $g: T \to (0, 1]$ is said to be uniformly finite if $\sup\{\#(\LL^{g}_{s, 1}(w))| s \in (0, 1],  w \in \LL^g_s\} < +\infty$.\\
(2)\,\,A function $f: T \to (0, \infty)$ is called sub-exponential if and only if there exist $m \ge 0$ and $c_1 \in (0, 1)$ such that $f(v) \le c_1f(w)$ for any $w \in T$ and $v \in T_w$ with $|v| \ge |w| + m$. $f$ is called super-exponential if and only if there exists $c_2 \in (0, 1)$ such that $f(v) \ge c_2f(w)$ for any $w \in T$ and $v \in S(w)$. $f$ is called exponential if it is both sub-exponential and super-exponential.
\enddefinition

The following proposition and the lemma are immediate consequences of the above definitions.

\prop\label{VDP.prop10}
Let $h$ be a weight function. Then $h$ is super-exponential if and only if there exists $c \ge 1$ such that $ch(w) \ge s \ge h(w)$ whenever $w \in \LL^h_s$.
\endprop

\demo
Assume that $h$ is super-exponential. Then there exists $c_2 < 1$ such that $h(w) \ge c_2h(\pi(w))$ for any $w \in T$. If $w \in \LL^h_s$, then $h(\pi(w)) > s \ge h(w)$. This implies $(c_2)^{-1}h(w) \ge s \ge h(w)$. \par
Conversely, assume that $ch(w) \ge s \ge h(w)$ for any $w$ and $s$ with $w \in \LL^h_s$. If $h(\pi(w)) > ch(w)$, then $w \in \LL^h_{t}$ for any $t \in (ch(w), h(\pi(w))$. This contradicts the assumption that $ch(w) \ge t \ge h(w)$. Hence $h(\pi(w)) \le ch(w)$ for any $w \in T$. Thus $h$ is super-exponential.
\enddemo

\lemma\label{VDP.lemma20}
If a weight function $g: T \to (0, 1]$ is uniformly finite, then
\[
 \sup\{\#(\LL_{s, M}(x)) | x \in X, s \in (0, 1]\} < +\infty
\]
for any $M \ge 0$.
\endlemma

\demo
Let $C = \sup\{\#(\LL_{s, 1}(w))| s \in (0, 1], w \in \LL_s\}$. Then $\#(\LL_{s, M}(x)) \le C + C^2 + \ldots + C^{M + 1}$.
\enddemo

\definition\label{BLE.def50}
Let $\a > 0$. A radon measure $\mu$ on $X$ is said to be $\a$-Ahlfors regular with respect to $d \in \D(X, \O)$ if and only if there exist $C_1, C_2 > 0$ such that
\begin{equation}\label{MEM.eq100}
C_1r^{\a} \le \mu(B_d(x, r)) \le C_2r^{\a}
\end{equation}
for any $r \in [0, \diam{X, d}]$.
\enddefinition

\definition\label{VDP.def30}
Let $g:T \to (0, 1]$ be a weight function. We say that $K$ has thick interior with respect to $g$, or $g$ is thick for short, if and only if there exist $M \ge 1$ and $\a > 0$ such that $K_w \supseteq U^{g}_{M}(x, {\a}s)$ for some $x \in K_w$ if $s \in (0, 1]$ and $w \in \LL^g_s$.
\enddefinition

The value of the integer $M \ge 1$ is not crucial in the above definition. In Proposition~\ref{ADD.prop20}, we will show if the condition of the above definition holds for a particular $M \ge 1$, then it holds for all $M \ge 1$. \par
The thickness is invariant under the bi-Lipschitz equivalence of weight functions as follows.

\prop\label{BLE.prop30}
Let $g$ and $h$ be weight functions. If $g$ is thick and $g \bl h$, then $h$ is thick.
\endprop

Since we need further results on thickness of weight functions, we postpone a proof of this proposition until the next section.\par
Now we give the main theorem of this sub-section.

\thm\label{MEM.thm100}
Let $d \in \D_A(X, \O)$ and let $\mu \in \M_P(X, \O)$. Assume that $K$ is minimal and $g_d$ is super-exponential and thick. Then $(g_d)^{\a} \bl g_{\mu}$ and $g_d$ is uniformly finite if and only if $\mu$ is $\a$-Ahlfors regular with respect to $d$. Moreover, if either/both of the these two conditions is/are satisfied, then $g_{\mu}$ and $g_d$ are exponential.
\endthm

By the same reason as Proposition~\ref{BLE.prop30}, a proof of this theorem will be given at the end of Section~\ref{VDP}.

\setcounter{equation}{0}
\section{Thickness of weight functions}\label{TGF}

In this section, we study conditions for a weight function being thick and relation between the notions ``thick'' and ``tight''.  For instance in Theorem~\ref{AAA.thm100} we present topological condition (TH1) ensuring that all super-exponential weight functions are thick.  In particular, this is the case for partitions of $S^2$ discussed in Section~\ref{SMR} because partitions satisfying \eqref{SMR.eq20} are minimal and the condition (TH) in  Section~\ref{SMR} yields the condition (TH1).  Moreover in this case, all super-exponential weight functions are tight as well by Corollary~\ref{BLES.cor10}.

\prop\label{ADD.prop20}
$g$ is thick if and only if for any $M \ge 0$, there exists $\b > 0$ such that, for any $w \in T$, $K_w \supseteq U^g_M(x, {\b}g(\pi(w)))$ for some $x \in K_w$.
\endprop

\demo
Assume that $g$ is thick. By induction, we are going to show  the following claim $(C)_M$ holds for any $M \ge 1$:\\
$(C)_M$\,\,\,There exists $\a_M > 0$ such that, for any $s \in (0, 1]$ and $w \in \LL^g_s$, one find $x \in K_w$ satisfying $K_w \supseteq U_M(x, \a_Ms)$. \\
Proof of $(C)_M$.\,\,Since $g$ is thick,  $(C)_M$ holds for some $M \ge 1$. Since $U^g_1(x, s) \subseteq U^g_M(x, s)$ if $M \ge 1$, $(C)_1$ holds as well. Now, suppose that $(C)_M$ holds. Let $w \in \LL^g_s$ and choose $x$ as in $(C)_M$. Then there exists $v \in \LL^g_{\a_Ms}$ such that $v \in T_w$ and $x \in K_v$. Applying $(C)_M$ again, we find $y \in K_v$ such that $K_v \supseteq U^g_M(y, (\a_M)^2s)$. Since $M \ge 1$, it follows that $U^g_{M + 1}(y, (\a_M)^2s) \subseteq U^g_M(x, \a_Ms) \subseteq K_w$. Therefore, letting $\a_{M + 1} = (\a_M)^2$, we have obtained $(C)_{M + 1}$. Thus we have shown $(C)_M$ for any $M \ge 1$. \qed\par
Next, fix $M \ge 1$ and write $\a = \a_M$. Note that $w \in \LL^g_s$ if and only if $g(w) \le s < g(\pi(w))$. Fix $\e \in (0, 1)$. Assume that $g(\pi(w)) > g(w)$. There exists $s_*$ such that $g(w) \le s_* < g(\pi(w))$ and $s_* > \e{g(\pi(w))}$. Hence we obtain
\[
K_w \supseteq U^g_M(x, {\a}s_*) \supseteq U^g_M(x, {\a}\e{g(\pi(w))}).
\]
If $g(w) = g(\pi(w))$, then there exists $v \in T_w$ such that $g(v) < g(\pi(v)) = g(w) = g(\pi(w))$. Choosing $s_*$ so that $g(v) \le s_* < g(\pi(v)) = g(\pi(w))$ and ${\e}g(\pi(w)) < s_*$, we obtain
\[
K_w \supseteq K_v \supseteq U^g_M(x, \a{s_*}) \supseteq U^g_M(x, \a\e{g(\pi(w))}).
\]
Letting $\b = \a\e$, we obtain the desired statement.\\
Conversely, assume that for any $M \ge 0$, there exists $\b > 0$ such that, for any $w \in T$, $K_w \supseteq U^g_M(x, {\b}g(\pi(w)))$ for some $x \in K_w$. If $w \in \LL_s$, then $g(w) \le s < g(\pi(w))$. Therefore $K_w \supseteq U_M^g(x, {\b}s)$. This implies that $g$ is thick.
\enddemo

\prop\label{ADD.prop21}
Assume that $K$ is minimal. Let $g: T \to (0, 1]$ be a weight function. Then $g$ is thick if and only if, for any $M \ge 0$, there exists $\c > 0$ such that, for any $w \in T$, $O_w \supseteq U_M^g(x, {\c}g(\pi(w)))$ for some $x \in O_w$. 
\endprop

\demo
Assume that $g$ is thick. By Proposition~\ref{ADD.prop20}, for any $M \ge 0$, we may choose $\a > 0$ so that for any $w \in T$, there exists $x \in K_w$ such that  $K_w \supseteq U_{M + 1}^g(x, \a{g(\pi(w))})$. Set $s_w = g(\pi(w))$. Let $y \in \sd{U^g_M(x, \a{s_w})}{O_w}$. There exists $v \in (T)_{|w|}$ such that $y \in K_v$ and $w \neq v$. Then we find $v_* \in T_v \cap \LL^g_{\a{s_w}}$ satisfying $y \in K_{v_*}$. Since $K_{v_*} \cap U^g_M(x, \a{s_w}) \neq \emptyset$, we have
\[
K_{v_*} \subseteq U^g_{M + 1}(x, \a{s_w}) \subseteq K_w.
\]
Therefore, $K_{v_*} \subseteq \cup_{w' \in T_w, |w'| = |v_*|} K_{w'}$. This implies that $O_{v_*} = \emptyset$, which contradicts the fact that $K$ is minimal. So, $\sd{U^g_M(x, \a{s_w})}{O_w} = \emptyset$ and hence $U^g_M(x, {\a}s_w) \subseteq O_w$. \\
The converse direction is immediate.
\enddemo

Using the above proposition, we give a proof of Proposition~\ref{BLE.prop30}.

\demo[Proof of Proposition~\ref{BLE.prop30}]
By Proposition~\ref{ADD.prop20}, there exists $\b > 0$ such that for any $w \in T$, $K_w \supseteq U^g_M(x, {\b}g(\pi(w)))$ for some $x \in K_w$. On the other hand, since there exist $c_1, c_2 > 0$ such that $c_1h(w) \le g(w) \le c_2h(w)$ for any $w \in T$. It follows that $D^g_M(x, y) \le c_2D^h_M(x, y)$ for any $x, y \in X$. Proposition~\ref{MAG.prop10} implies that there exists $\a > 0$ such that $\a\d_M^g(x, y) \le \d_M^h(x, y)$ for any $x, y \in X$. Hence $U^h_M(x, \a{s}) \subseteq U^g_M(x, s)$ for any $x \in X$ and $s \in (0, 1]$. Combining them, we see that
\[
K_w \supseteq U^g_M(x, {\b}g(\pi(w))) \supseteq U^h_M(x, \a\b{g(\pi(w))}) \supseteq U^h(x, \a\b{c_2}h(\pi(w))).
\]
Thus by Proposition~\ref{ADD.prop20}, $h$ is thick.
\enddemo

\thm\label{AAA.thm100}
Assume that $K$ is minimal. Define $h_*: T \to (0, 1]$ by $h_*(w) = 2^{-|w|}$ for any $w \in T$. Then the following conditions are equivalent:\\
{\rm (TH1)}\,\,
\[
\sup_{w \in T} \min\big\{|v_*| - |w| \big| v_* \in T_w,K_{v_*} \subseteq O_w\big\} < \infty.
\]
{\rm (TH2)}\,\,
Every super-exponential weight function is thick.\\
{\rm (TH3)}\,\,
There exists a sub-exponential weight function which is thick.\\
{\rm (TH4)}\,\,
The weight function $h_*$ is thick.
\endthm

\demo
(TH1) $\Rightarrow$ (TH2):\,\,Assume (TH1). Let $m$ be the supremum in (TH1). Let $g$ be a super-exponential weight function. Then there exists $\lambda \in (0, 1)$ such that $g(w) \ge {\lambda}g(\pi(w))$ for any $w \in T$. Let $w \in \LL_s^g$. By (TH1), there exists $v_* \in T_w \cap (T)_{|w| + m}$ such that $K_{v_*} \subseteq O_w$. For any $v \in T_w \cap (T)_{|w| + m}$,
\[
g(v) \ge \lambda^mg(w) \ge \lambda^{m + 1}g(\pi(w)) > \lambda^{m + 1}s
\]
Choose $x \in O_{v_*}$. Let $u \in \LL_{\lambda^{m + 1}s, 1}(x)$. Then there exists $v' \in \LL_{\lambda^{m + 1}s, 0}(x)$ such that $K_{v'} \cap K_u \neq \emptyset$. Since $g(v_*) > \lambda^{m + 1}s$ and $x \in O_{v_*}$, it follows that $v' \in T_{v_*}$. Therefore $K_{u} \cap O_w \supseteq K_{u} \cap K_{v_*} \neq \emptyset$. This implies that either $u \in T_w$ or $w \in T_{u}$. Since $g(w) > \lambda^{m + 1}s$, it follows that $u \in T_w$. Thus we have shown that $\LL_{\lambda^{m + 1}s, 1}(x) \subseteq T_w$. Hence
\[
U_1^g(x, \lambda^{m + 1}s) \subseteq K_w.
\]
This shows that $g$ is thick.\\
(TH2) $\Rightarrow$ (TH4):\,\,Apparently $h_*$ is an exponential weight function. Hence by (TH2), it is thick.\\
(TH4) $\Rightarrow$ (TH3):\,\,Since $h_*$ is exponential and thick, we have (TH3).\\
(TH3) $\Rightarrow$ (TH1):\,\,Assume that $g$ is a sub-exponential weight function which is thick. Proposition~\ref{ADD.prop21} shows that there exist $\c \in (0, 1)$ and $M \ge 1$ such that for any $w \in T$, $O_w \supseteq U_M^g(x, \c{g(\pi(w))})$. Choose $v_* \in \LL_{\c{g(\pi(w))}, 0}^g(x)$. Then $K_{v_*} \subseteq U_M^g(x, \c{g(\pi(w))}) \subseteq O_w$ and $g(\pi(v_*)) > {\c}g(\pi(w)) \ge {\c}g(w)$. Since $g$ is sub-exponential, there exists $k \ge 1$ and $\eta \in (0, 1)$ such that $g(u) \le \eta{g(v)}$ if $v \in T_v$ and $|u| \ge |w| + k$.  Choose $l$ so that $\eta^l < \c$ and set $m = kl + 1$. Since $g(\pi(v_*)) > \eta^{l}g(w)$, we see that $|\pi(v_*)| \le |w| + m - 1$. Therefore, $|v_*| \le |w| + m$ and hence we have (TH1).
\enddemo

\thm\label{BLES.thm10}
Assume that $K: T \to \C(X, \O)$ is minimal, that there exists $\lambda \in (0, 1)$ such that if $B_w = \emptyset$, then $\#(T_w \cap \LL_{\lambda{g(w)}}^g) \ge 2$ and that $g$ is thick. Then $g$ is tight.
\endthm
\demo
By Proposition~\ref{ADD.prop21}, there exists $\c$ such that, for any $v \in T$, $O_v \supseteq U_M^g(x, \c{g(\pi(v))})$ for some $x \in K_v$. First suppose that $B_w \neq \emptyset$. Then there exists $x \in K_w$ such that $O_w \supseteq U_M^g(x, \c{g(\pi(w))})$. For any $y \in B_w$, it follows that $\d_M^g(x, y) > \c{g(\pi(w))}$. Thus 
\[
\sup_{x, y \in K_w} \d_M^g(x, y) \ge \c{g(\pi(w))}.
\]
Next if $B_w = \emptyset$, then there exists $u \neq v \in T_w \cap \LL_{\lambda{g(w)}}^g$. If $B_u \neq \emptyset$, then the above discussion implies
\[
\sup_{x, y \in K_w} \d_M^g(x, y) \ge \sup_{x, y \in K_v} \d_M^g(x, y) \ge \c{g(\pi(v))} \ge \c\lambda{g(w)}.
\]
If $B_u = \emptyset$, then $\d_M^g(x, y) \ge \lambda{g(w)}$ for any $(x, y) \in K_u \times K_v$. Thus for any $w \in T$, we conclude that
\[
\sup_{x, y \in K_w} \d_M^g(x, y) \ge \c\lambda{g(w)}.
\]
\enddemo

The above theorem immediately implies the following corollary.

\cor\label{BLES.cor10}
Assume that $(X, \O)$ is connected and $K$ is minimal. If $g$ is thick, then $g$ is tight.
\endcor

\setcounter{equation}{0}
\section{Volume doubling property}\label{VDP}

In this section, we introduce the notion of a relation called ``gentle'' written as $\gen$ between weight functions. This relation is not an equivalence relation in general. In Section~\ref{GAE}, however, it will be shown to be an equivalence relation among exponential weight functions. As was the case of the bi-Lipschitz equivalence, the gentleness will be identified with other properties in classes of weight functions. In particular, we are going to show that the volume doubling property of a measure with respect to a metric is equivalent to the gentleness of the associated weight functions.\par
As in the previous sections, $(T, \A, \phi)$ is a locally finite tree with a reference point $\phi$, $(X, \O)$ is a compact metrizable topological space with no isolated point and $K: T \to \C(X, \O)$ is a partition of $X$ parametrized by $(T, \A, \phi)$. \par
 The notion of gentleness of a weight function to another weight function is defined as follows.
 \remark
In the case of the natural partition of a self-similar set in Example~\ref{PAS.ex10}, the main results of this section, Theorems~\ref{VDP.thm10} and \ref{VDP.thm20} have been obtained in \cite{Ki13}.
\endremark

\definition\label{VDP.def10}
Let $g: T \to (0, 1]$ be a weight function. A function $f: T \to (0, \infty)$ is said to be gentle with respect to $g$ if and only if there exists $c_G > 0$ such that $f(v) \le c_Gf(w)$ whenever $w, v \in \LL^{g}_{s}$ and $K_w \cap K_v \neq \emptyset$ for some $s \in (0, 1]$. We write $f \gen g$ if and only if $f$ is gentle with respect to $g$.
\enddefinition

Alternatively, we have a simpler version of the definition of gentleness under a mild restriction.

\prop\label{VDP.prop20}
Let $g: T \to (0, 1]$ be an exponential weight function. Let $f : T \to (0, \infty)$. Assume that $f(w) \le f(\pi(w))$ for any $w \in T$ and $f$ is super-exponential. Then $f$ is gentle with respect to $g$ if and only if there exists $c > 0$ such that $f(v) \le cf(w)$ whenever $g(v) \le g(w)$ and $K_v \cap K_w \neq \emptyset$.
\endprop

\demo
By the assumption, there exist $c_1, c_2 > 0$ and $m \ge 1$ such that $f(v) \ge c_2f(w)$, $g(v) \ge c_2g(w)$ and $g(u) \le c_1g(w)$ for any $w \in T$, $v \in S(w)$ and $u \in T_w$ with $|u| \ge |w| + m$. \par
First suppose that $f$ is gentle with respect to $g$. Then there exists $c > 0$ such that $f(v') \le cf(w')$ whenever $w', v' \in \LL^g_s$ and $K_{w'} \cap K_{v'} \neq \emptyset$ for some $s \in (0, 1]$. Assume that $g(v) \le g(w)$ and $K_v \cap K_w \neq \emptyset$. There exists $u \in T_w$ such that $K_{u} \cap K_v \neq \emptyset$ and $g(\pi(u)) > g(v) \ge g(u)$. Moreover, $g(\pi([v]_m)) > g([v]_m) = g(v)$ for some $m \in [0, |v|]$. Then $[v]_m, u \in \LL_{g(v)}^g$ and hence $f(v) \le f([v]_m) \le cf(u) \le cf(w)$. \par
Conversely, assume that $f(v') \le cf(w')$ whenever $g(v') \le g(w')$ and $K_{v'} \cap K_{w'} \neq \emptyset$. Let $w, v \in \LL^g_s$ with $K_w \cap K_v \neq \emptyset$. If $g(v) \le g(w)$, then $f(v) \le cf(w)$. Suppose that $g(v) > g(w)$. Since $g$ is super-exponential, we see that
\[
s \ge g(w) \ge c_2g(\pi(w)) \ge c_2s \ge c_2g(v).
\]
Set $N = \min\{n| c_2 \ge c_1^n\}$. Choose $u \in T_v$ so that $K_{u} \cap K_w \neq \emptyset$ and $|u| = |v| + Nm$. Then $g(w) \ge c_2g(v) \ge (c_1)^Ng(v) \ge g(u)$. This implies $f(u) \le cf(w)$. Since $f(u) \ge (c_2)^{Nm}f(v)$, we have $f(v) \le c(c_2)^{-Nm}f(w)$. Therefore, $f$ is gentle with respect to $g$.
\enddemo

The following is the standard version of the definition of the volume doubling property.

\definition\label{VDP.def20}
Let $\mu$ be a radon measure on $(X, \O)$ and let $d \in \D(X, O)$. $\mu$ is said to have the volume doubling property with respect to the metric $d$ if and only if there exists $C > 0$ such that
\[
\mu(B_d(x, 2r)) \le C\mu(B_d(x, r))
\]
for any $x \in X$ and any $r > 0$. 
\enddefinition

Since $(X, \O)$ has no isolated point, if a Radon measure $\mu$ has the volume doubling property with respect to some $d \in \D(X, \O)$, then the normalized version of $\mu$, $\mu/\mu(X)$, belongs to $\M_P(X, O)$. Taking this fact into account, we are mainly interested in (normalized version of) a Radon measure in $\M_P(X, \O)$.\par
 The main theorem of this section is as follows.
 
 \thm\label{VDP.thm00}
 Let $d \in \D(X, \O)$ and let $\mu \in \M_P(X, \O)$. Assume that $d$ is adapted, that $g_d$ is thick, exponential and uniformly finite and that $\mu$ is exponential. Then $\mu$ has the volume doubling property with respect to $d$ if and only if $g_d \gen g_{\mu}$.
 \endthm
 
So, this says that the volume doubling property equals the gentleness in the world of weight functions having certain regularities. This theorem is an immediate corollary of Theorem~\ref{VDP.thm10} and \ref{VDP.thm20}.\par
To describe a refined version of Theorem~\ref{VDP.thm00}, we define the notion of volume doubling property of a measure with respect to a weight function $g$ as well by means of balls ``$U^g_M(x, s)$''.

\definition\label{VDP.def25}
Let $\mu \in \M_P(X, \O)$ and let $g$ be a weight function. For $M \ge 1$, we say $\mu$ has $M$-volume doubling property with respect to $g$ if and only if there exist $\c \in (0, 1)$ and $\b \ge 1$ such that $\mu(U^g_M(x, s)) \le \b\mu(U^g_M(x, {\c}s))$ for any $x \in X$ and any $s \in (0, 1]$.
\enddefinition

It is rather annoying that the notion of ``volume doubling property'' of a measures with respect to a weight function depends on the value $M \ge 1$. Under certain conditions including exponentiality and the thickness, however, we will show that if $\mu$ has $M$-volume doubling property for some $M \ge 1$, then it has $M$-volume doubling property for all $M \ge 1$ in Theorem~\ref{VDP.thm20}. \par
Naturally, if a metric is adapted to a weight function, the volume doubling with respect to the metric and that with respect to the weight function are virtually the same as is seen in the next proposition.

\prop\label{VDP.prop00}
Let $d \in \D(X, \O)$, let $\mu \in M_P(X, \O)$ and let $g$ be a weight function. Assume that $d$ is adapted to $g$. Then $\mu$ has the volume doubling property with respect to $d$ if and only if there exists $M_* \ge 1$ such that $\mu$ has $M$-volume doubling property with respect to $g$ for any $M \ge M_*$.
\endprop

\demo
Since $d$ is adapted to $g$, for sufficiently large $M$, there exist $\a_1, \a_2 > 0$ such that 
\[
U^g_M(x, \a_1{s}) \subseteq B_d(x, s) \subseteq U^g_M(x, \a_2{s})
\]
for any $x \in X$ and $s \in (0, 1]$. Suppose that $\mu$ has the volume doubling property with respect to $d$. Then there exists $\lambda  > 1$ such that
\[
\mu(B_d(x, 2^mr)) \le  \lambda^m\mu(B_d(x, r))
\]
for any $x \in X$ and $r \ge 0$. Hence
\[
\mu(U^g_M(x, \a_12^mr)) \subseteq \lambda^m\mu(U^g_M(x, \a_2{r})).
\]
Choosing $m$ so that $\a_12^m > \a_2$, we see that $\mu$ has $M$-volume doubling property with respect to $g$ if $M$ is sufficiently large. Converse direction is more or less similar.
\enddemo 

By the above proposition, as far as we confine ourselves to adapted metrics, it is enough to consider the volume doubling property of a measure with respect to a weight function.
Thus we are going to investigate relations between ``the volume doubling property with respect to a weight function'' and other conditions like a weight function being exponential, a weight function being uniformly finite, a measure being super-exponential, and a measure being gentle with respect to a weight function. To begin with, we show that the last four conditions imply the volume doubling property of $\mu$ with respect to $g$.

\thm\label{VDP.thm10}
Let $g: T \to (0, 1]$ be a weight function and let $\mu \in \M_P(X, \O)$. Assume that $g$ is exponential, that $g$ is uniformly finite, that $\mu$ is gentle with respect to $g$ and that $\mu$ is super-exponential. Then $\mu$ has $M$-volume doubling property with respect to $g$ for any $M \ge 1$.
\endthm

Hereafter in this section, we are going to omit $g$ in notations if no confusion may occur. For example, we write $\LL_s, \LL_{s, M}(w), \LL_{s, M}(w)$ and $U_M(x, s)$ in place of $\LL^g_s, \LL^{g}_{s, M}(w), \LL^{g}_{s, M}(x)$ and $U_{M}^g(x, s)$ respectively.\par
The following lemma is a step to prove the above theorem.

\lemma\label{VDP.lemma30}
Let $g:T \to (0, 1]$ be a weight function and let $\mu \in \M_P(X, \O)$. For $s \in (0, 1]$, $\lambda > 1$ and $c > 0$, define
\[
\Theta(s, \lambda, k, c) = \{v| v \in \LL_s, \mu(K_u) \le c\mu(K_v)\,\,\, \text{for any $u \in \LL_{\lambda{s}, k}((v)_{\lambda{s}})$}\},
\]
where $(v)_{\lambda{s}}$ is the unique element of $\{[v]_n| 0 \le n \le |v|\} \cap \LL_{\lambda{s}}$. Assume that $g$ is uniformly finite and that there exist $N \ge 1, \lambda > 1$ and $c > 0$ such that $\LL_{s, N}(w) \cap \Theta(s, \lambda, 2N + 1, c) \neq \emptyset$ for any $s \in (0, 1]$ and $w \in \LL_s$. Then $\mu$ has the $N$-volume doubling property with respect to $g$.
\endlemma

\demo
Let $w \in \LL_{s, 0}(x)$ and let $v \in \LL_{s, N}(w) \cap \Theta(s, \lambda, 2N + 1, c)$. If $u \in \LL_{\lambda{s}, N}(x)$, then $u \in \LL_{\lambda{s}, 2N + 1}((v)_{\lambda{s}})$. Moreover,  since $v \in \LL_{s, N}(x)$, we see that
\[
\mu(K_u) \le c\mu(K_v) \le c\mu(U_N(x, s)).
\]
Therefore,
\[
\mu(U_N(x, \lambda{s})) \le \sum_{u \in \LL_{\lambda{s}, N}(x)} \mu(K_u) \le \#(\LL_{\lambda{s}, N}(x))c\mu(U_N(x, s)).
\]
Since $g$ is uniformly finite, Lemma~\ref{VDP.lemma20} shows that $\#(\LL_{\lambda{s}, N}(x))$ is uniformly bounded with respect to $x \in X$ and $s \in (0, 1]$. 
\enddemo

\demo[Proof of Theorem~\ref{VDP.thm10}]
Fix $\lambda > 1$. By Proposition~\ref{VDP.prop10}, there exists $c \ge 1$ such that $cg(w) \ge s \ge g(w)$ if $w \in \LL_{s}$. Since $g$ is sub-exponential, there exist $c_1 \in (0, 1)$ and $m \ge 1$ such that $c_1g(w) \ge g(v)$ whenever $v \in T_w$ and $|v| \ge |w| + m$. Assume that $w \in \LL_s$. Set $w_* = (w)_{\lambda{s}}$. Then $\lambda{s} \ge g(w_*)$. If $|w| \ge |w_*| + nm$, then $(c_1)^ng(w_*) \ge g(w)$ and hence $(c_1)^n\lambda{s} \ge g(w) \ge g(w)/c$. This shows that $(c_1)^n\lambda{c} \ge 1$. Set $l = \min\{n|n \ge 0, (c_1)^n\lambda{c} < 1\}$. Then we see that $|w| < |w_*| + lm$.\par
On the other hand, since $\mu$ is super-exponential, there exists $c_2 > 0$ such that $\mu(K_u) \ge c_2\mu(K_{\pi(u)})$ for any $u \in T$. This implies that $\mu(K_{w_*}) \le (c_2)^{-ml}\mu(K_w)$. Since $\mu$ is gentle, there exists $c_* > 0$ such that $\mu(K_{w(1)}) \le c_*\mu(K_{w(2)})$ whenever $w(1), w(2) \in \LL_s$ and $K_{w(1)} \cap K_{w(2)} \neq \emptyset$ for some $s \in (0, 1]$. Therefore for any $u \in \LL_{{\lambda}s, M}(w_*)$,
\[
\mu(K_u) \le (c_*)^M\mu(K_{w_*}) \le (c_*)^M(c_2)^{-ml}\mu(K_w).
\]
Thus we have shown that
\[
\LL_s = \Theta(s, \lambda, M, (c_*)^M(c_2)^{-ml})
\]
for any $s \in (0, 1]$. Now by Lemma~\ref{VDP.lemma30}, $\mu$ has $M$-volume doubling property with respect to $g$ for any $M \ge 1$.
\enddemo

In order to study the converse direction of Theorem~\ref{VDP.thm10}, we need the thickness  of $K$ with respect to the weight function in question.

\thm\label{VDP.thm20}
Let $g: T \to (0, 1]$ be a weight function and let $\mu \in \M_P(X, \O)$. Assume that $g$ is thick.\\
{\rm (1)}\,\,
Suppose that $g$ is exponential and uniformly finite. Then the following conditions are equivalent:\\
{\rm (VD1)}\,\,
$\mu$ has $M$-volume doubling property with respect to $g$ for some $M \ge 1$. \\
{\rm (VD2)}\,\,
$\mu$ has $M$-volume doubling property with respect to $g$ for any $M \ge 1$.\\
{\rm (VD3)}\,\,
$\mu$ is gentle with respect to $g$ and $\mu$ is super-exponential.\\
{\rm (2)}\,\,
Suppose that $K$ is minimal and $g$ is super-exponential. Then {\rm (VD1)}, {\rm (VD2)} and the following condition {\rm (VD4)} are equivalent:\\
{\rm (VD4)}\,\,
$g$ is sub-exponential and uniformly finite, $\mu$ is gentle with respect to $g$ and $\mu$ is super-exponential. \\
Moreover, if any of the above conditions {\rm (VD1), (VD2)} and {\rm (VD4)} hold, then $\mu$ is exponential and
\[
\sup_{w \in T}\#(S(w))  < +\infty.
\]
\endthm

In general, the statement of Theorem~\ref{VDP.thm20} is false if $g$ is not thick. In fact, in Example~\ref{ESS.ex10}, we will present an example without thickness where $d$ is adapted to $g$, $g$ is exponential and uniformly finite, $\mu$ has the volume doubling property with respect to $g$ but $\mu$ is neither gentle to $g$ nor super-exponential.\par
As for a proof of Theorem~\ref{VDP.thm20}, it is enough to show the following theorem.

\thm\label{VDP.thm30}
Let $g: T \to (0, 1]$ be a weight function and let $\mu \in \M_P(X, \O)$. Assume that $\mu$ has $M$-volume doubling property with respect to $g$ for some $M \ge 1$.\\
{\rm (1)}\,\,
If $g$ is thick, then $\mu$ is gentle with respect to $g$.\\
{\rm (2)}\,\,
If $g$ is thick and $g$ is super-exponential, then $\mu$ is super-exponential.\\
{\rm (3)}\,\,
If $g$ is thick and $K$ is minimal, then $g$ is uniformly finite.\\
{\rm (4)}\,\,
If $g$ is thick, $K$ is minimal, and $\mu$ is super-exponential, then 
\[
\sup_{w \in T} \#(S(w)) < +\infty.
\]
and $\mu$ is sub-exponential.\\
{\rm (5)}\,\,
If $g$ is uniformly finite, $\mu$ is gentle with respect to $g$, $\mu$ is sub-exponential, then $g$ is sub-exponential.
\endthm

To prove Theorem~\ref{VDP.thm30}, we need several lemmas.

\lemma\label{VDP.lemma50}
Let $g: T \to (0, 1]$ be a weight function. Assume that $K$ is minimal and $g$ is thick. Let $\mu \in \M_P(X, \O)$. If $\mu$ has $M$-volume doubling property with respect to $g$ for some $M \ge 1$, then there exists $c > 0$ such that $\mu(O_w) \ge c\mu(K_w)$ for any $w \in T$.
\endlemma

\demo
By Proposition~\ref{ADD.prop21}, there exists $\c > 0$ such that $O_v \supseteq U^g_M(x, {\c}s)$ for some $x \in K_v$ if $v \in \LL_s$. Let $w \in T$.  Choose $u \in T_w$ such that $u \in \LL_{g(w)/2}$. Then
\[
\mu(O_w) \ge \mu(O_u) \ge \mu(U^g_M(x, {\c}g(w)/2)).
\]
Since $\mu$ has $M$-volume doubling property with respect to $g$, there exists $c > 0$ such that
\[
\mu(U^g_M(y, {\c}r/2)) \ge  c\mu(U^g_M(y, r))
\]
for any $y \in X$ and $r > 0$. Since $U_M(x, g(w)) \supseteq K_w$, it follows that
\[
\mu(O_w) \ge \mu(U^g_M(x, {\c}g(w)/2)) \ge c\mu(U_M(x, g(w))) \ge c\mu(K_w).
\]
\enddemo

\lemma\label{VDP.lemma60}
Let $g: T \to (0, 1]$ be a weight function. Assume that $\mu \in \M_P(X, \O)$ is gentle with respect to $g$ and that $g$ is uniformly finite. Then there exists $c > 0$ such that
\[
c\mu(K_w) \ge \mu(U_M(x, s))
\]
if $w \in \LL_{s, 0}(x)$. 
\endlemma

\demo
Since $\mu$ is gentle with respect to $g$, there exists $c_1 > 0$ such that $\mu(K_v) \le c_1\mu(K_w)$ if $w \in \LL_s$ and $v \in \LL_{s, 1}(w)$. Hence if $v \in \LL_{s, {M + 1}}(w)$, then it follows that $\mu(K_v) \le (c_1)^{M + 1}\mu(K_w)$. Since $\LL_{s, M}(x) \subseteq \LL_{s, M + 1}(w)$, 
\begin{multline*}
\mu(U_M(x, s)) \le \sum_{v \in \LL_{s, M}(x)} \mu(K_v)\\ \le \sum_{v \in \LL_{s, M}(x)}(c_1)^{M + 1}\mu(K_w) = (c_1)^{M + 1}\#(\LL_{s, M}(x))\mu(K_w).
\end{multline*}
By Lemma~\ref{VDP.lemma20}, we obtain the desired statement.
\enddemo

\demo[Proof of Theorem~\ref{VDP.thm30}]
{\rm (1)}\,\,
Since $g$ is thick, there exists $\b \in (0, 1)$ such that, for any $s \in (0, 1]$ and $w \in \LL_s$, $K_w \supseteq U_M(x, {\b}s)$ for some $x \in K_w$. By $M$-volume doubling property of $\mu$, there exists $c > 0$ such that $\mu(U_M(x, {\b}s)) \ge c\mu(U_M(x, s))$ for any $s \in (0, 1]$ and $x \in X$. Hence
 \begin{equation}\label{VDP.eq10}
\mu(K_w) \ge \mu(U_M(x, {\b}s)) \ge c\mu(U_M(x, s)).
\end{equation}
 If $v \in \LL_s$ and $K_v \cap K_w \neq \emptyset$, then $U_M(x, s) \supseteq K_v$. \eqref{VDP.eq10} shows that $\mu(K_w) \ge c\mu(K_v)$. Hence $\mu$ is gentle with respect to $g$.\\
{\rm (2)}\,\,
Let $v \in \sd{T}{\{\phi\}}$. Choose $u \in T_v$ so that $u \in \LL_{g(v)/2}$. Applying \eqref{VDP.eq10} to $u$ and using the volume doubling property repeatedly, we see that there exists $x \in K_u$ such that
\begin{equation}\label{VDP.eq20}
\mu(K_v) \ge \mu(K_u) \ge \mu(U_M(x, {\b}g(v)/2)) \ge c^n\mu(U_M(x, \b^{1 - n}g(v)/2))
\end{equation}
for any $n \ge 0$. Since $g$ is super-exponential, there exists $n \ge 0$, which is independent of $v$, such that $\b^{1 - n}g(v) /2 > g(\pi(v))$. By \eqref{VDP.eq20}, we obtain $\mu(K_v) \ge c^n\mu(K_{\pi(v)})$. Thus $\mu$ is super-exponential. \\
{\rm (3)}\,\,
Let $w \in \LL_s$. Then $\{O_v\}_{v \in \LL_{s, 1}(w)}$ is mutually disjoint by Lemma~\ref{PAS.lemma20}-(2). By \eqref{VDP.eq10} and Lemma~\ref{VDP.lemma50},
\[
\mu(K_w) \ge c\mu(U_M(x, s)) \ge c\sum_{v \in \LL_{s, 1}(w)} \mu(O_v) \ge c^2\sum_{v \in \LL_{s, 1}(w)} \mu(K_v)
\]
As $\mu$ is gentle with respect to $g$ by (1), there exists $c_* > 0$, which is independent of $w$ and $s$, such that $\mu(K_v) \ge c_*\mu(K_w)$ for any $v \in \LL_{s, 1}(w)$. Therefore,
\[
\mu(K_w) \ge c^2\sum_{v \in \LL_{s, 1}(w)} \mu(K_v) \ge c^2c_*\#(\LL_{s, 1}(w))\mu(K_w)
\]
Hence $\#(\LL_{s, 1}(w)) \le c^{-2}(c_*)^{-1}$ and $g$ is uniformly finite.\\
{\rm (4)}\,\,
By Lemma~\ref{VDP.lemma50}, for any $w \in T$, we have
\[
\mu(K_w) \ge \mu(\cup_{v \in S(w)} O_v) = \sum_{v \in S(w)} \mu(O_v) \ge c\sum_{v \in S(w)} \mu(K_v).
\]
Since $\mu$ is super-exponential, there exists $c' > 0$ such that $\mu(K_v) \ge c'\mu(K_w)$ if $w \in T$ and $v \in S(w)$. Hence
\[
\mu(K_w) \ge c\sum_{v \in S(w)}\mu(K_v) \ge cc'\#(S(w))\mu(K_w).
\]
Thus $\#(S(w)) \le (cc')^{-1}$, which is independent of $w$. By the above arguments, \begin{equation}\label{VDP.eq30}
\mu(O_v) \ge c\mu(K_v) \ge c_*\mu(K_w) \ge c_*\mu(O_w)
\end{equation}
for any $w \in T$ and $v \in S(w)$, where $c_* = cc'$. Let $v_* \in S(w)$. If $\mu(O_{v_*}) = (1 - a)\mu(O_w)$, then
\[
\mu(O_w) \ge \sum_{v \in S(w)} \mu(O_v) = (1 - a)\mu(O_{w}) + \sum_{v \in S(w), v \neq v_*} \mu(O_v).
\]
This implies $a\mu(O_w) \ge \mu(O_v)$ for any $v \in \sd{S(w)}{\{v_*\}}$. By \eqref{VDP.eq30}, $a \ge c_*$. Therefore, $\mu(O_v) \le (1 - c_*)\mu(O_w)$ for any $v \in S(w)$. This implies
\[
c\mu(K_v) \le \mu(O_v) \le (1 - c_*)^m\mu(O_w) \le (1 - c_*)^m\mu(K_w)
\]
if $v \in T_w$ and $|v| = |w| + m$. Choosing $m$ so that $(1 - c_*)^m < c$, we see that $\mu$ is sub-exponential.\\
{\rm (5)}\,\,
As $\mu$ is sub-exponential, there exist $\a \in (0, 1)$ and $m \ge 0$ such that $\mu(K_v) \le \a\mu(K_w)$ if $u \in T_w$ and $|u| \ge |w| + m$. Since $\mu$ has $M$-volume doubling property with respect to $g$, there exist $\lambda, c \in (0, 1)$ such that $\mu(U_M(x, \lambda{s})) \ge c\mu(U_M(x, s))$ for any $x \in X$ and $s > 0$. Let $\b \in (\lambda, 1)$. Assume that $g$ is not sub-exponential. Then for any $n \ge 0$, there exist $w \in T$ and $u \in T_w$ such that $|u| \ge |w| + nm$ and $g(u) \ge {\b}g(w)$. In case $g(w) = g(\pi(w))$, we may replace $w$ by $v = [w]_m$ for some $m \in \{0, 1, \ldots, |w|\}$ satisfying $g(\pi(v)) > g(v) = g(w)$ or $g(v) = g(w) = 1$.  Consequently we may assume $w \in \LL_{g(w)}$. Set $s = g(w)$. Since $\b > \lambda$, there exists $u_* \in T_u \cap \LL_{s\lambda}$. Let $x \in K_{u_*}$. Then by the volume doubling property,
\[
\mu(U_M(x, \lambda{s})) \ge c\mu(U_M(x, s)) \ge c\mu(K_w).
\]
By Lemma~\ref{VDP.lemma60}, there exists $c_* > 0$ which is independent of $n, w$ and $u$ such that
\[
c_*\mu(K_{u_*}) \ge \mu(U_M(x, \lambda{s})).
\]
Since $\mu$ is sub-exponential, 
\[
\a^nc_*\mu(K_w) \ge c_*\mu(K_{u_*}) \ge \mu(U_M(x, \lambda{s})) \ge c\mu(K_w).
\]
This implies $\a^nc_* \ge c$ for any $n \ge 0$ which is a contradiction.
 \enddemo
 
 At the end of this section, we give a proof of Theorem~\ref{MEM.thm100} by using Theorem~\ref{VDP.thm30}.
 
 \demo[Proof of Theorem~\ref{MEM.thm100}]
If is enough to show the case where $\a = 1$. Assume that $g_d \bl g_{\mu}$ and $d$ is uniformly finite. Since $d$ is adapted, there exist $M \ge 1$ and $\a_1, \a_2 > 0$ such that
\[
U_M^d(x, \a_1{r}) \subseteq B_d(x, r) \subseteq U_M^d(x, \a_2{r})
\]
for any $x \in X$ and $r > 0$. \par
Write $d_w = g_d(w)$ and $\mu_w = g_{\mu}(w)$. Assume that $g_d \bl g_{\mu}$.  Then there exist $c_1, c_2 > 0$ such that
\[
c_1d_w \le \mu_w \le c_2d_w
\]
for any $w \in T$. For any $x \in X$, choose $w \in T$ so that $x \in K_w$ and $w \in \LL^d_{\a_1{r}}$. Then since $d$ is super-exponential, there exists $\lambda$ which is independent of $x, r$ and $w$ such that
\[
\mu(B_d(x, r)) \ge \mu(U^d_M(x, {\a_1}r)) \ge \mu(K_w) \ge c_1d_w \ge c_1{\lambda}d_{\pi(w)} \ge c_1\lambda\a_1r.
\]
On the other hand, since $d$ is uniformly finite, Lemma~\ref{VDP.lemma20} implies
\begin{multline*}
\mu(B_d(x, r)) \le \mu(U^d_M(x, \a_2{r})) \le C\sum_{w \in \LL^d_{\a_2r, M}(x)} \mu(K_w) \\
\le Cc_2\sum_{w \in \LL^d_{\a_2r, M}(x)}d_w \le Cc_2\#(\LL^d_{\a_2r, M})\a_2r \le C_2r
\end{multline*}
Conversely, assume \eqref{MEM.eq100}. For any $w \in T$ and $x \in K_w$, 
\[
K_w \subseteq U^d_M(x, d_w) \subseteq B_d(x, d_w/\a_1).
\] 
Hence
\[
\mu(K_w) \le \mu(B_d(x, d_w/c_1)) \le C_2d_w/\a_1.
\]
By Proposition~\ref{ADD.prop20}, there exists $z \in K_w$ such that
\[
K_w \supseteq U^d_M(z, \b{d_{\pi(w)}}) \supseteq B_d(z, \b{d_{\pi(w)}}/\a_2).
\]
By \eqref{MEM.eq100},
\[
\mu(K_w) \ge \mu(B_d(z, \b{d_{\pi(w)}}/\a_2)) \ge C_1\b{d_{\pi(w)}}/\a_2 \ge C_1\b{d_w}/\a_2.
\]
Thus we have shown that $g_d \bl g_{\mu}$. Furthermore, since $d$ is $M$-adapted for some $M \ge 1$, $\mu$ has $M$-volume doubling property with respect to the weight function $g_d$. Applying Theorem~\ref{VDP.thm30}-(3), we see that $g_d$ is uniformly finite.  In the same way, by Theorem~\ref{VDP.thm30}, both $g_d$ and $g_{\mu}$ are exponential.
\enddemo

\def\wI{\widehat{I}}
\setcounter{equation}{0}
\section{Example: subsets of the square}\label{ESS}
\begin{figure}
\centering
\setlength{\unitlength}{20mm}

\begin{picture}(4.3,2.5)(0.7,-0.5)
\linethickness{1pt}
\thinlines
\multiput(0,0)(0.81,0){4}{\drawline(0, 0)(0, 2.43)}
\multiput(0, 0)(0,0.81){4}{\drawline(0,0)(2.43,0)}
\multiput(0,0)(0.27,0){4}{\drawline(0,0)(0,0.81)}
\multiput(0,0)(0,0.27){4}{\drawline(0,0)(0.81,0)}
\multiput(0, 0)(1.215,0){3}{\put(0, 0){\circle*{0.07}}}
\multiput(0, 1.215)(1.215,0){3}{\put(0, 0){\circle*{0.07}}}
\multiput(0, 2.43)(1.215,0){3}{\put(0, 0){\circle*{0.07}}}
\put(0.405, 1.215){\makebox(0,0){\large$Q_8$}}
\put(0.405, 2.025){\makebox(0,0){\large$Q_7$}}
\put(1.215, 2.025){\makebox(0,0){\large$Q_6$}}
\put(2.025, 2.025){\makebox(0,0){\large$Q_5$}}
\put(1.215, 1.4){\makebox(0,0){\large$Q_9$}}
\put(2.025, 1.215){\makebox(0,0){\large$Q_4$}}
\put(2.025, 0.405){\makebox(0,0){\large$Q_3$}}
\put(1.215, 0.405){\makebox(0,0){\large$Q_2$}}
\put(0.135, 0.405){\makebox(0,0){\small$Q_{18}$}}
\put(0.405, 0.405){\makebox(0,0){\small$Q_{19}$}}
\put(0.675, 0.405){\makebox(0,0){\small$Q_{14}$}}
\put(0.405, 0.135){\makebox(0,0){\small$Q_{12}$}}
\put(0.405, 0.675){\makebox(0,0){\small$Q_{16}$}}
\put(0.135, 0.675){\makebox(0,0){\small$Q_{17}$}}
\put(0.675, 0.675){\makebox(0,0){\small$Q_{15}$}}
\put(0.675, 0.135){\makebox(0,0){\small$Q_{13}$}}
\put(0.135, 0.135){\makebox(0,0){\small$Q_{11}$}}

\put(-0.12, -0.12){\makebox(0,0){\large$p_1$}}
\put(-0.1,1.215){\makebox(0,0)[r]{\large$p_8$}}
\put(1.215, -0.1){\makebox(0,0)[t]{\large$p_2$}}
\put(2.63, -0.1){\makebox(0,0)[t]{\large$p_3$}}
\put(2.53, 1.215){\makebox(0,0)[l]{\large$p_4$}}
\put(1.215, 1.1){\makebox(0,0)[t]{\large$p_9$}}
\put(2.58, 2.58){\makebox(0,0){\large$p_5$}}
\put(1.215, 2.53){\makebox(0,0)[b]{\large$p_6$}}
\put(-0.1, 2.5){\makebox(0,0)[b]{\large$p_7$}}

\put(1.25,-0.4){\makebox(0, 0){The square $Q$ as a self-similar set}}

\multiput(3, 0)(0,2.43){2}{\drawline(0,0)(2.43,0)}
\multiput(3, 0)(2.43,0){2}{\drawline(0,0)(0,2.43)}
\shade\path(3.5,1.2)(4.3,1.2)(4.3,1.5)(3.5,1.5)(3.5,1.2)
\put(3.9, 1.35){\makebox(0,0){\large$R_3$}}
\shade\path(4.5,0.3)(5.2,0.3)(5.2,1.9)(4.5,1.9)(4.5,0.3)
\put(4.85, 1.1){\makebox(0,0){\large$R_1$}}
\shade\path(3.3,0.2)(3.8,0.2)(3.8,1.0)(3.3,1.0)(3.3,0.2)
\put(3.55,0.6){\makebox(0,0){\large$R_2$}}
\shade\path(3.9, 0.3)(4.3,0.3)(4.3, 0.6)(3.9,0.6)(3.9,0.3)
\put(4.1, 0.45){\makebox(0,0){\large$R_4$}}
\shade\path(3.7,1.8)(4.0,1.8)(4.0,2.2)(3.7,2.2)(3.7,1.8)
\put(3.85,2.0){\makebox(0,0){\large$R_5$}}
\shade\path(3.2,1.7)(3.5,1.7)(3.5,2.0)(3.2,2.0)(3.2,1.7)
\put(3.35, 1.85){\makebox(0,0){\large$R_6$}}
\shade\path(4.8,2.0)(5.1,2.0)(5.1,2.2)(4.8,2.2)(4.8,2.0)
\put(4.95,2.1){\makebox(0,0){\large$R_7$}}

\put(4.25,-0.4){\makebox(0, 0){$\displaystyle X = Q \backslash \Bigg(\bigcup_{i \ge 1}\inte{R_j}\Bigg)$}}

\end{picture}
\caption{The square $Q$ and its subset $X$}\label{SQ1}
\end{figure}
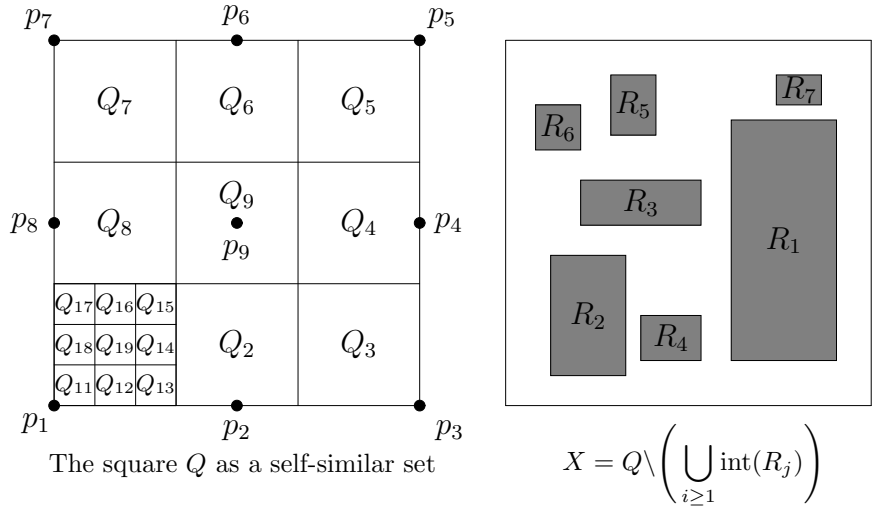

In this section, we give illustrative examples of the results in the previous sections. For simplicity, our examples are subsets of the square $[0, 1]^2$ denoted by $Q$ and trees parametrizing partitions are sub-trees of $(T^{(9)}, \A^{(9)}, \phi)$ defined in Example~\ref{TWR.ex10}. Note that $[0, 1]^2$ is divided into $9$-squares with the length of the sides $\frac 13$. As in Example~\ref{PAS.ex10}, the tree $(T^{(9)}, \A^{(9)}, \phi)$ is naturally appears as the tree parametrizing the natural partition associated with this self-similar division. Namely, let $p_1 = (0, 0), p_2 = (1/2, 0), p_3 = (1, 0), p_4 = (1, 1/2), p_5 = (1, 1), p_6 = (1/2, 1), p_7 = (0, 1), p_8 = (0, 1/2)$ and $p_9 = (1/2, 1/2)$. Set $W = \{1, \ldots, 9\}$. Define $F_i: Q \to Q$ by
\[
F_i(x) =\frac 13 (x - p_i) + p_i
\]
for any $i \in W$. Then $F_i$ is a similitude for any $i \in W$ and
\[
Q = \bigcup_{i \in W} F_i(Q).
\]
See Figure~\ref{SQ1}. In this section, we write $(W_*, \A_*, \phi) = (T^{(9)}, \A^{(9)}, \phi)$, which is a locally finite tree with a reference point $\phi$. Set $W_m = \{1, \ldots, 9\}^m$.  Then $(W_*)_m = W_m$ and $\pi^{(W_*, \A_*, \phi)}(w) = \word w{m - 1}$ for any $w = \word wm \in W_m$. For simplicity, we use $|w|$ and $\pi$ in place of $|w|_{(W_*, \A_*, \phi)}$ and $\pi^{(W_*, \A_*, \phi)}$ respectively hereafter. Define $g: W_* \to (0, 1]$ by $g(w) = 3^{-|w|}$ for any $w \in W_*$. Then $g$ is an exponential weight function.\par
As for the natural associated partition of $Q$, define $F_w = F_{w_1}\circ\ldots\circ{F_{w_m}}$ and $Q_w = F_w(Q)$ for any $w = \word wm \in W_m$. Set $Q_*(w) = Q_w$ for any $w \in W_*$. (If $w = \phi$, then $F_{\phi}$ is the identity map and $Q_{\phi} = Q$.) Then $Q_*: W_* \to \C(Q, \O)$ is a partition of $Q$ parametrized by $(W_*, \A_*, \phi)$, where $\O$ is the natural topology induced by the Euclidean metric. In fact, $\cap_{m \ge 0} Q_{[\omega]_m}$ for any $\omega \in \SS$, where $\SS = W^{\BbN}$, is a single point. Define $\s:\SS \to Q$ by $\{\s(\omega)\} = \cap_{m \ge 0} Q_{[\omega]_m}$.\par
It is easy to see that the partition $Q_*$ is minimal, $g$ is uniformly finite, $g$ is thick with respect to the partition $Q_*$, and the (restriction of) Euclidean metric $d_E$ on $Q$ is $1$-adapted to $g$.\par
In order to have more interesting examples, we consider certain class of subsets of $Q$ whose partition is parametrized by a subtree $(T, \A_*|_{T \times T}, \phi)$ of $(W_*, \A_*, \phi)$. Let $\{I_m\}_{m \ge 0}$ be a sequence of subsets of $W_*$ satisfying the following conditions (SQ1), (SQ2) and (SQ3):\newline
(SQ1)\,\,\,
For any $m \ge 0$, $I_m \subseteq W_m$ and if $\wI_{m + 1} = \{wi|w \in I_m, i \in W\}$, then $I_{m + 1} \supseteq \wI_{m + 1}$.\newline
(SQ2)\,\,\,
$Q_w \cap Q_v = \emptyset$ if $w \in \wI_{m + 1}$ and $v \in \sd{I_{m + 1}}{\wI_{m + 1}}$.\newline
(SQ3)\,\,\,
For any $m \ge 0$, the set $\cup_{w \in I_m} Q_w$ is a disjoint union of rectangles $R^m_j = [a^m_j, b^m_j] \times [c^m_j, d^m_j]$ for $j = 1, \ldots, k_m$. \par
See Figure~\ref{SQ1}. By (SQ2), we may assume that $k_m \le k_{m + 1}$ and $R^m_j = R^{m + 1}_j$ for any $m$ and $j = 1, \ldots, k_m$ without loss of generality. Under this assumption, we may omit $m$ of $R^m_j, a^m_j, b^m_j, c^m_j$ and $d^m_j$ and simply write $R_j, a_j, b_j, c_j$ and $d_j$ respectively.\par

\notation
As a topology of $Q = [0, 1] \times [0, 1]$, we consider the relative topology induced by the Euclidean metric. We use $\inte{A}$ and $\partial{A}$ to denote the interior and the boundary, respectively, of a subset $A$ of $Q$ with respect to this topology.
\endnotation

Note that $\inte{\cup_{w \in I_{m}} Q_w} = \cup_{j = 1, \ldots, k_m} \inte{R_j}$. 

\prop\label{ESS.prop10}
{\rm (1)}\,\,\,
Define 
\[
X^{(m)} = \sd{Q}{\Bigg(\bigcup_{j = 1, \ldots, k_m} \inte{R_j}\Bigg     )}.
\]
then $X^{(m)} \supseteq X^{(m + 1)}$ for any $m \ge 0$ and $X = \cap_{m \ge 0} X^{(m)}$ is a non-empty compact set. Moreover, $\partial{R_j} \subseteq X$ for any $j \ge 1$. \newline
{\rm (2)}\,\,\,
 Define $(T)_m = \{w|w \in W_m, \inte{Q_w} \cap X \neq \emptyset\}$ for any $m \ge 0$. If $T = \cup_{m \ge 0} (T)_m$ and $\A = \A_*|_{T \times T}$, then $(T, \A, \phi)$ is a locally finite tree with the reference point $\phi$ and $\#(S(w)) \ge 3$ for any $w \in T$. Moreover, let 
\[
\SS_T = \{\omega|\omega \in \SS, [\omega]_m \in (T)_m\,\,\,\text{for any $m \ge 0$}\}
\]
Then $X = \s(\SS_T)$. \newline
{\rm (3)}\,\,\,
Define $K_w = Q_w \cap X$ for any $w \in T$. Then $K_w \neq \emptyset$ and $K: T \to \C(X)$ defined by $K(w) = K_w$ is a minimal partition of $X$ parametrized by $(T, \A, \phi)$. Moreover, $g|_T$ is exponential and uniformly finite.
\endprop

To prove the above proposition, we need the following lemma.

\lemma\label{ESS.lemma10}
If $w \in T$, then $\cup_{i \in W, wi \notin T} Q_{wi}$ is a disjoint union of rectangles and $\#(\{i | i \in W, wi \in T\}) \in \{3, 5, 7, 8, 9\}$.
\endlemma

\demo
Set $I = \{i | i \in W, wi \notin T\}$. For each $i \in I$, there exists $k_i \ge 1$ such that $Q_{wi} \subseteq R_{k_i}$. Hence $\cup_{i \in I} Q_{wi} = \cup_{i \in I} (Q_w \cap R_{k_i})$. Since $\{R_j\}_{j \ge 1}$ are mutually disjoint, we have the desired conclusion. Assume that $I = W$. Suppose $|i - j| = 1$. Since $Q_{wi} \cap Q_{wj} \neq \emptyset$, we see that $R_{k_i} = R_{k_j}$. Hence $R_{k_1} = \ldots = R_{k_9}$ and $Q_w \subseteq R_{k_1}$. This contradicts the fact that $\inte{Q_w} \cap X \neq \emptyset$. Thus $I \neq W$. Considering all the possible shapes of $\cup_{i \in W, wi \notin T} Q_{wi}$, we conclude $\#(\{i | i \in W, wi \in T\}) \in \{3, 5, 7, 8, 9\}$.
\enddemo

\demo[Proof of Proposition~\ref{ESS.prop10}]
(1)\,\,
Since $\{X^{(m)}\}_{m \ge 0}$ is a decreasing sequence of compact sets, $X$ is a nonempty compact set. By (SQ2), $R_j \cap R_i = \emptyset$ for any $i \neq j$. Therefore, $\partial{R_j} \subseteq X^{(m)}$ for any $m \ge 0$. Hence $\partial{R_j} \subseteq X$.\\
(2)\,\,
If $w \in (T)_m$, then $\inte{Q_{\pi(w)}} \cap X \supseteq \inte{Q_w} \cap X \neq \emptyset$. Hence $\pi(w) \in (T)_{m - 1}$. Using this inductively, we see that $[w]_k \in (T)_k$ for any $k \in \{0, 1, \ldots, m\}$. This implies that $(T, \A, \phi)$ is a locally finite tree with a reference point $\phi$. By Lemma~\ref{ESS.lemma10}, we see that $\#(\{i | i \in W, wi \in (T)_{m + 1}\}) \ge 3$. Next if $\omega \in \SS_T$, then for any $m \ge 0$, there exists $x_m \in \inte{Q_{[\omega]_m}} \cap X$. Therefore, $x_m \to \s(\omega)$ as $m \to \infty$. Since $X$ is compact, it follows that $\s(\omega) \in X$. \par
Conversely, assume that $x \in X$. Set $W_{m, x} = \{w| w \in W_m, x \in Q_w\}$. Note that $\#(\s^{-1}(x)) \le 4$ and $\cup_{w \in W_{m, x}} Q_w$ is a neighborhood of $x$. Suppose that $(T)_m \cap W_{m, x} \neq \emptyset$ for any $m \ge 0$. Then there exists $w_m \in (T)_m \cap W_{m, x}$ such that  $x \in Q_{w_m}$. Since $W_{m, x} = \{[\omega]_m| \omega \in \s^{-1}(x)\}$, there exists $\omega \in \s^{-1}(x)$ such that $[\omega]_m = w_m$ for infinitely many $m$. As $\inte{Q_{[\omega]_m}}$ is monotonically decreasing, it follows that $[\omega]_m \in (T)_m$ for any $m \ge 0$. This implies $x \in \s(\SS_T)$. Suppose that there exists $m \ge 0$ such that $W_{m, x} \cap (T)_m = \emptyset$. By this assumption, $\inte{Q_w} \cap X = \emptyset$ for any $w \in W_{m, x}$ and hence there exists $j_w \ge 1$ such that $Q_w \subseteq R_{j_{w}}$. Note that $Q_w \cap Q_{w'} \neq \emptyset$ for any $w, w' \in W_{m, x}$ and hence $R_{j_w} = R_{j_{w'}}$. Therefore, $ \cup_{w \in W_{m, x}} Q_w \subseteq R_j$ for some $j \ge 1$. Since $\cup_{w \in W_{m, x}}Q_w$ is a neighborhood of $x$, it follows that $x \notin X$. This contradiction concludes the proof.\\
(3)\,\,
The fact that $K$ is a partition of $X$ parametrized by $(T, \A|_{T \times T}, \phi)$ is straightforward from (1) and (2). As $\sd{K_w}{(\cup_{v \in (T)_m, v \neq w} K_v)}$ is contained in the sides of the square $Q_w$, the partition $K$ is minimal. Since $\LL^{g|_T}_s = (T)_m$ if and only if $\frac 1{3^m} \le s < \frac 1{3^{m - 1}}$, it follows that $g|_T$ is exponential. Furthermore, $\LL^{g|_T}_{s, 1}(w) \subseteq \{v| v \in W_m, Q_v \cap Q_w \neq \emptyset\}$ for any $w \in (T)_m$. Hence $\#(\LL^{g|_T}_{s, 1}(w)) \le 8$. This shows that $g|_T$ is uniformly finite.
\enddemo

Now, we consider when the restriction of the Euclidean metric is adapted. 

\definition\label{ESS.def10}
Let $R = [a, b] \times [c, d]$ be a rectangle. The degree of distortion of $R$, $\kappa(R)$, is defined by
\[
\kappa(R) = \max\bigg\{1, (1 - \delta_{c0})(1 - \delta_{d1})\frac {|b - a|}{|d - c|}, (1 - \delta_{a0})(1 - \delta_{b1})\frac {|d - c|}{|b - a|}\bigg\},
\]
where $\delta_{xy}$ is the Kronecker delta defined by $\delta_{xy} = 1$ if $x = y$ and $\delta_{xy} = 0$ if $x \neq y$. Moreover, for $\kappa \ge 1$, we define
\[
\R_{\kappa}^0 = \{R| \text{$R$ is a rectangle, $R \subseteq Q$ and $\kappa(R) \le \kappa$}\}
\]
and
\begin{multline*}
\R_{\kappa}^1 = \{R| R \subseteq Q, \text{$R$ is a rectangle, there exists $w \in T$ such that $\sd{Q_w}{\inte{R}}$}\\
\text{ has two connected components and $\kappa(Q_w \cap R) \le \kappa$}\}
\end{multline*}
\enddefinition

The extra factors $(1 - \delta_{c0}), (1 - \delta_{d1}), (1 - \delta_{a0})$ and $(1 - \delta_{b1})$ become effective if the rectangle $R$ has an intersection with the boundary of the square $Q$.

\thm\label{ESS.thm10}
Let $d$ be the restriction of the Euclidean metric on $X$. Then $d$ is adapted to $g|_{T}$ if and only if the following condition {\rm (SQ4)} holds:\\
{\rm (SQ4)}\,\,
There exists $\kappa \ge 1$ such that $R_j \in \R^0_{\kappa} \cup \R^1_{\kappa}$ for any $j \ge 1$.
\endthm

Several lemmas are needed to prove the above theorem.

\lemma\label{ESS.lemma20}
Define $N(x, y) = \min\{[-\frac{\log{|x_1 - y_1|}}{\log 3}],[-\frac{\log{|x_2 - y_2|}}{\log 3}]\}$ for any $x = (x_1, x_2), y = (y_1, y_2) \in Q$, where $[a]$ is the integer part of a real number $a$.\\
{\rm (1)}\,\,
\[
\frac 1{3^{N(x, y) + 1}} < d(x, y) \le \frac{\sqrt{2}}{3^{N(x, y)}}
\]
{\rm (2)}\,\,If $x, y \in X$, then there exist $w, v, u \in W_{N(x, y)}$ such that $w, v \in T$, $x \in Q_w$, $y \in Q_u$, $Q_w \cap Q_v \neq \emptyset$ and $Q_v \cap Q_v \neq \emptyset$.
\endlemma

\demo
Set $N = N(x, y)$. Let $n_i = [-\frac{\log{|x_i - y_i|}}{\log 3}]$ for $i = 1, 2$. Then $N = \min\{n_1, n_2\}$ and
\[
\frac 1{3^{N + 1}}< |x_j - y_j| \le \frac 1{3^N}
\]
if $n_j = N$. This yields (1). Since $x, y \in X$, then there exist $w, u \in (T)_m$ such that $x \in K_w$ and $y \in K_u$. Since $|x_1 - y_1| \le 1/3^N$ and $|x_2 - y_2| \le 1/3^N$, we find $v \in W_m$ satisfying $Q_w \cap Q_v \neq \emptyset$ and $Q_v \cap Q_u \neq \emptyset$.\enddemo

\notation
For integers $n, k, l \ge 0$, we set
\[
Q(n, k, l) = \bigg[\frac k{3^n}, \frac{(k + 1)}{3^n}\bigg] \times \bigg[\frac l{3^n}, \frac{(l + 1)}{3^n}\bigg]
\]
\endnotation

\lemma\label{ESS.lemma30}
Assume {\rm (SQ4)}. Let $M = [\log{(2\kappa)}/\log 3] + 1$ and $L = 2[2\kappa] + 9$. If $w, v \in (T)_m$ and $Q_w \cap Q_v \neq \emptyset$, then there exists a chain $(w(1), \ldots, w(L))$ of $K$ such that $w \in T_{w(1)}$, $v \in T_{w(L)}$ and $|w(k)| \ge m -M$.
\endlemma
\demo
Case 1:\,\,Assume that $Q_w \cap Q_v$ is a line segment. Without loss of generality, we may assume that $Q_w = Q(m, k - 1, l)$ and $Q_v = Q(m, k, l)$.\\
Case 1a: $K_w \cap K_v \neq \emptyset$, then $(w, v)$ is a desired chain of $K$. \\
Case 1b: In case $K_w \cap K_v = \emptyset$, $Q_w \cap Q_v \cap K_w$ and $Q_w \cap Q_v \cap K_v$ are disjoint closed subsets of $Q_w \cap Q_v$. Since $Q_w \cap Q_v$ is connected, there exists $a \in Q_w \cap Q_v$ such that $a \notin K_w \cap K_v$. Since $K_w \cup K_v$ is closed, there exists an open neighborhood of $a$ which has no intersection with $K_w \cap K_v$. This open neighborhood must be contained in $R_j$ for some $j$. So, we see that $R_j \cap \inte{Q_w \cap Q_v} \neq \emptyset$ and $(k - 1)/3^m \le a_j \le k/3^m \le b_j \le (k + 1)/3^m$. Assume $c_i > l/3^m$. Then since the line segment $[a_j, b_j] \times \{c_j\}$ is contained in $X$, we see that $K_w \cap K_v \neq \emptyset$. Therefore $c_j \le l/3^m$. By the same argument we have $d_j \ge (l + 1)/3^m$. Now if $R_j \in \R^0_{\kappa}$, it follows that $|d_j - c_j| \le 2\kappa/3^m$. Hence the line segment $[a_j, b_j] \times \{c_j\}$ and $[a_j, b_j] \times \{d_j\}$ are covered by at most $4$ pieces of $K_u$'s for $u \in (T)_m$ and the line segment $\{a_j\} \times [c_j, d_j]$ and $\{b_j\} \times [c_j, d_j]$ is covered by at most $2\kappa + 2$ pieces of $K_u$'s for $u \in (T)_m$. Since $K_w$ and $K_v$ are pieces of these coverings, we obtain a chain $(w(1), \ldots, w(k))$ of $K$ from these coverings where $w(1) = w, w(k) = v$ and $l \le 2\kappa + 5$. Next assume $R_j \in \R^1_{\kappa}$. Note that $2\kappa/3^m \le 1/3^{m - M}$. By the definition of $\R_{\kappa}^1$, there exists $u \in (T)_{m - M}$ such that $\sd{Q_u}{R_j}$ has two connected component. Sifting $Q_u$ up and down, we may find $u' \in (T)_{m - M}$ such that $Q_w \cup Q_v \subseteq Q_{u'}$. Then $(u')$ is a desired chain of $K$.\\
Case 2:\,\,Assume that $Q_w \cap Q_v$ is a single point. Without loss of generality, we may assume that $Q_w = Q(m, k - 1, l - 1)$ and $Q_v = Q(m, k, l)$. Choose $u(1), u(2) \in W_m$ so that $Q_{u(1)} = Q(m, k - 1, l)$ and $Q_{w(2)} = Q(m, k + 1, l - 1)$. If neither $u(1)$ nor $u(2)$ belongs to $T$, then there exist $i, j \ge 1$ such that $Q_{u(1)} \subseteq R_i$ and $Q_{u(2)} \subseteq R_j$. Since $Q_{u(1)} \cap Q_{u(2)} \neq \emptyset$, it follows that $R_i = R_j$ and hence $Q_w \cup Q_v \subseteq R_i$. This contradicts the fact that $w, v \in T$. Hence $u(1) \in T$ or $u(2) \in T$. Let $u(1) \in T$. Then $Q_w \cap Q_{u(1)}$ and $Q_{u(1)} \cap Q_v$ are line segments. By using the method in (1), we find  a chain between $w$ and $u(1)$ and a chain between $u(1)$ and $v$. Connecting these two chains, we obtain the desired chain $(w(1), \ldots, w(L))$. 
\enddemo

\demo[Proof of Theorem~\ref{ESS.thm10}]
Assume (SQ4). Let $x, y \in X$. Define $N = N(x, y)$ and choose $w, v, u \in W_N$ as in Lemma~\ref{ESS.lemma20}. We fix the constants $M$ and $L$ as in Lemma~\ref{ESS.lemma30}. There are two cases.\newline
{\bf Case 1}: Suppose $v \in T$. Applying Lemma~\ref{ESS.lemma30} to two pairs $\{w, v\}$ and $\{v, u\}$ and connecting the two resultant chains, we obtain a chain $(w(1), \ldots, w(2L - 1)) \in \CH_K(x, y)$ satisfying $w \in T_{w(1)}, u \in \in T_{w(2L - 1)}$ and $|w(i)| \ge N - M$ for any $i$. This concludes Case 1.\\
{\bf Case 2}: Suppose $v \notin T$. If $Q_w \cap Q_u \neq \emptyset$, then we have a chain $(w(1), \ldots, w(L))$ between $x$ and $y$ satisfying $w \in T_{w(1)}$, $u \in T_{w(L)}$ and $|w(i)| \ge N - M$ for any $i$ by Lemma~\ref{ESS.lemma30}. Assume $Q_w \cap Q_u = \emptyset$. Without loss of generality, we may assume one of the following tree situations (a), (b) and (c):\\
(a)\,\,
$Q_w = Q(N, k - 1, l - 1)$ and $Q_u = Q(N, k + 1, l - 1)$.\\
(b)\,\,
$Q_w = Q(N, k - 1, l - 1)$ and $Q_u = Q(N, k + 1, l)$.\\
(c)\,\,
$Q_w = Q(N, k - 1, l - 1)$ and $Q_u = Q(N, k + 1, l + 1)$.\par
Set $Q_{v(1)} = Q(N, k, l - 1)$ and $Q_{v(2)} = Q(N, k, l)$. In each case,  $x_1 = k/3^N$ and $y_1 = (k + 1)/3^N$.\par
First consider cases (a) and (b). If either $v(1)$ or $v(2)$ belongs to $T$, then replacing $v$ by either $v(1)$ or $v(2)$, we end up with Case 1. So we assume that neither $v(1)$ nor $v(2)$ belongs to $T$. Then there exists $j \ge 1$ such that $Q_{v(1)} \cup Q_{v(2)} \subseteq R_j$. Since $x_1 = k/3^N$ and $y_1 = (k + 1)/3^N$, we have $a_j = k/3^N$ and $b_j = (k + 1)/3^N$. Then by the same argument as in the proof of Lemma~\ref{ESS.lemma30}, there exists a chain $(w(1), \ldots, w(L)) \in \CH_K(x, y)$ such that $w \in T_{w(1)}$, $u \in T_{w(L)}$ and $|w(i)| \ge N - M$ for any $i$.\par
Next in the situation of (c), $x = (k/3^N, l/3^N)$, $y = ((k + 1)/3^N, (l + 1)/3^N)$ and $v = v(2)$. Since $v = v(1) \notin T$, there exists $j \ge 1$ such that $Q_v \subseteq R_j$. Note that $x, y \in X  \cap Q_v$. Hence $Q_v = R_j$. Choose $v(3), v(4) \in W_N$ so that $Q_{v(3)} = Q(N, k + 1, l - 1)$ and $Q_{v(4)} = Q(N, k + 1, l)$. Then $v(3), v(4) \in T$ and therefore $(w, v(1), v(3), v(4), u)$ is a chain of $K$ between $x$ and $y$. This concludes Case 2.\par
As a consequence, we may always find a chain $(w(1), \ldots, w(2L - 1)) \in \CH_K(x, y)$ satisfying $|w(i)| \ge N(x, y) - M$ for any $i$. By Lemma~\ref{ESS.lemma20}-(1),
\[
3^{M + 1}d(x, y) \ge 3^M\frac 1{3^N} \ge \frac 1{3^{w(i)}} = g(w(i)).
\]
Thus we have verified the conditions (ADa) and $\rm (ADb)_{2L - 2}$ in Theorem~\ref{PAS.thm10}. Hence $d$ is $(2L - 2)$-adapted to $g|_T$ by Theorem~\ref{PAS.thm10}.\par
Conversely, assume that $d$ is $J$-adapted to $g|_T$. By $\rm (ADb)_J$, there exists $C \ge 0$ such that for any $x, y \in X$, there exists a chain $(w(1), \ldots, w(J + 1)) \in \CH_K(x, y)$ satisfying
\begin{equation}\label{ESS.eq100}
Cd(x, y) \ge \frac 1{3^{|w(i)|}}
\end{equation}
for any $i = 1, \ldots, J + 1$. Set $M = [\log{(\sqrt{2}C)}/\log 3] + 1$. Suppose that (SQ4) does not hold; for any $\kappa \ge 1$, there exists $R_j \notin \R^0_{\kappa} \cup \R^1_{\kappa}$. In particular, we choose $\kappa \ge 3^{M + 2}$. Write $R = R_j$ and set $R = [a, b] \times [c, d]$. Define $\partial{R}_L = \{a\} \times [c, d]$ and $\partial{R}_R = \{b\} \times [c, d]$. (The symbols ``L'' and ``R'' correspond to the words ``Left'' and ``Right'' respectively.) Without loss of generality, we may assume that $|a - b| \le |c - d|$. Since $R \notin \R^0_{\kappa}$, we have $\kappa|b - a| \le |d - c|$. Let $x = (a, (c + d)/2)$ and let $y = (b, (c + d)/2)$. Set $N = N(x, y)$. There exists $(w(1), \ldots, w(J + 1)) \in \CH_K(x, y)$ such that \eqref{ESS.eq100} holds for any $i = 1, \ldots, J + 1$. By Lemma~\ref{ESS.lemma20}-(1),  
\begin{equation}\label{ESS.eq110}
|w(i)| \ge N - M
\end{equation}
for any $i = 1, \ldots, J + 1$. Define $A = [0, 1] \times (c, d)$. If $Q_{w(i)} \subseteq A$, $Q_{w(i)} \cap \partial{R}_L \neq \emptyset$ and $Q_{w(i)} \cap \partial{R}_R \neq \emptyset$, then the fact that $R \notin \R^1_{\kappa}$ along with Lemma~\ref{ESS.lemma20}-(1) shows 
\begin{equation}\label{ESS.eq120}
\frac1{3^{|w(i)|}} \ge \kappa|b - a| = \kappa{d(x, y)} \ge \frac {\kappa}{3^{N + 1}} \ge \frac 1{3^{N + M - 1}}.
\end{equation}
This contradicts \eqref{ESS.eq110} and hence we verify the following claim (I):\\
(I)\,\,
If $Q_{w(i)} \subseteq A$, then $Q_{w(i)} \cap \partial{R}_L = \emptyset$ or $Q_{w(i)} \cap \partial{R}_R = \emptyset$. \par
Next we prove that there exists $j \ge 1$ such that $\sd{Q_{w(j)}}{A} \neq \emptyset$. Otherwise, $Q_{w(i)} \subseteq A$ for any $i = 1, \ldots, J + 1$. Let $A_L = [0, a] \times (c, d)$ and let $A_R = [b, 1] \times (c, d)$. Define $I_L = \{i | i = 1, \ldots, J + 1, Q_{w(i)} \cap A_L \neq \emptyset\}$ and $I_R = \{i | i = 1, \ldots, J + 1, Q_{w(i)} \cap A_R \neq \emptyset\}$. Since $K_{(w(i))} \subseteq X \cap A \subseteq A_L \cup A_R$, it follows that $\{1, \ldots, J + 1\} = I_L \cup I_R$. Moreover, the claim (I) implies $I_L \cap I_R = \emptyset$. Hence $I_L = \{i | i = 1, \ldots, J + 1, K_{w(i)} \subseteq A_L\}$ and $I_R = \{i | i = 1, \ldots, J + 1, K_{w(i)} \subseteq A_R\}$. This contradicts the fact that $(w(1), \ldots, w(J + 1))$ is a chain of $K$ between $x$ and $y$. Thus there exists $j \ge 1$ such that $\sd{Q_{w(j)}}A \neq \emptyset$. Define $i_* = \min\{i| i = 1, \ldots, J + 1, \sd{Q_{w(i)}}A \neq \emptyset\}$. Without loss of generality, we may assume that $Q_{w(i_*)} \cap [0, 1] \times \{d\} \neq \emptyset$. Set 
\[
\partial{R}_L^T = \{a\} \times \bigg[\frac{c + d}2 , d - \frac 1{3^{|w(i_*)|}}\bigg].
\]
Shifting $Q_{w(i)}$'s for $i = 1, \ldots, i_* - 1$ horizontally towards $\partial{R}_L$, we obtain a covering of $\partial{R}_L^T$. Note that the length of $\partial{R}_L^T$ is $|d - c|/2 - 1/3^{|w(i_*)|}$ and
\[
\frac {|d - c|}2 - \frac 1{3^{|w(i_*)|}} \ge \frac{\kappa|b - a|}2 - \frac 1{3^{N - M}} = \frac {\kappa}2d(x, y) - \frac 1{3^{N - M}} \ge \frac {\kappa}2\frac 1{3^{N + 1}} - \frac 1{3^{N - M}}.
\]
On the other hand, the lengths of the sides of $Q_{w(i)}$'s are no less that $1/3^{N - M}$ by \eqref{ESS.eq110}. Hence
\[
i_* - 1 \ge 3^{N - M}\bigg(\frac {\kappa}2\frac 1{3^{N + 1}} - \frac 1{3^{N - M}}\bigg) \ge \frac{\kappa}2\frac 1{3^{M + 1}} - 1.
\]
Since $J + 1 \ge i_*$, it follows that
\[
2(J + 1)3^{M + 1} \ge \kappa.
\]
This contradicts the fact that $\kappa$ can be arbitrarily large. Hence we conclude that (SQ4) holds.
\enddemo

In the followings, we give four examples. The first one has infinite connected components but still the restriction of the Euclidean metric is adapted.

 \example[Figure~\ref{SQ2}]\label{ESS.ex00}
Let $X$ be the self-similar set associated with the contractions $\{F_1, F_3, F_4, F_5, F_7, F_8\}$, i.e. $X$ is the unique nonempty compact set which satisfies
\[
X = \bigcup_{i \in S} F_j(X),
\]
where $S = \{1, 3, 4, 5, 7, 8\}$. Then $X = C_3 \times [0, 1]$, where $C_3$ is the ternary Cantor set. Define $(T)_m = S^m$ and $T = \cup_{m \ge 1} (T)_m$. If $K_w = F_w(X)$ for any $w \in T$, then $K$ is a partition of  $X$ parametrized by $(T, \A|_T, \phi)$. Define
\[
I_{\phi} = \Big[\frac 13, \frac 23\Big] \times [0, 1]\quad\text{and}\quad I_{i_1, \ldots, i_n} = \bigg[\sum_{k = 1}^n\frac{i_k}{3^k} + \frac 1{3^{n + 1}}, \sum_{k = 1}^n\frac{i_k}{3^k} + \frac 2{3^{n + 1}}\bigg] \times [0, 1]
\]
for any $1 \ge 0$ and $i_1, \ldots, i_n \in \{0, 2\}$. Then
\[
\{R_j\}_{j \ge 1} = \{I_{\phi}, I_{i_1, \ldots, i_n}| n \ge 1, i_1, \ldots, i_n \in \{0, 2\}\}.
\]
Set $J_{i_1, \ldots, i_n} = [\sum_{k = 1}^n\frac{i_k}{3^k}, \sum_{k = 1}^n\frac{i_k}{3^k} + \frac 1{3^n}] \times [0, \frac 1{3^n}]$. Then there exists $w \in (T)_n$ such that $J_{i_1, \ldots, i_n} = Q_{w}$, $\sd{Q_w}{\inte{I_{i_1, \ldots, i_n}}}$ has two connected component and $\kappa({Q_w} \cap {I_{i_1, \ldots, i_n}}) = 3$. Therefore, $\{R_j\}_{j \ge 1} \subseteq \R^1_{3}$ and hence $d$ is adapted to $g|_T$.
\endexample

\begin{figure}
\centering
\setlength{\unitlength}{20mm}

\begin{picture}(4.3,2.5)(0.7,-0.5)
\linethickness{1pt}
\thinlines
\multiput(0,0)(0.81,0){4}{\drawline(0,0)(0, 2.43)}
\multiput(0,0)(0,2.43){2}{\drawline(0,0)(2.43,0)}
\shade\path(0.81,0)(0.81,2.43)(1.62,2.43)(1.62,0)(0.81,0)
\put(1.215,2.6){\makebox(0,0){\large$I_{\phi}$}}
\multiput(0,0)(1.62,0){2}{\shade\path(0.27,0)(0.27,2.43)(0.54,2.43)(0.54,0)(0.27,0)}
\put(0.405,2.6){\makebox(0,0){$I_{0}$}}
\put(2.025,2.6){\makebox(0,0){$I_{2}$}}
\multiput(0,0)(1.62,0){2}{\multiput(0,0)(0.54,0){2}{\shade\path(0.09,0)(0.18,0)(0.18,2.43)(0.09,2.43)(0.09,0)}}
\put(0.135,2.6){\makebox(0,0){\small$I_{00}$}}
\put(0.675,2.6){\makebox(0,0){\small$I_{02}$}}
\put(1.775,2.6){\makebox(0,0){\small$I_{20}$}}
\put(2.295,2.6){\makebox(0,0){\small$I_{22}$}}
\path(0, 0)(0,0.81)(0.81,0.81)(0.81,0)(0,0)
\multiput(0,0)(0,0.03){27}{\drawline(0,0)(0.81,0)}
\multiput(0,0)(0.03,0){27}{\drawline(0,0)(0,0.81)}
\put(0.405,-0.15){\makebox(0,0){\large$J_0$}}
\path(1.62,0)(1.89,0)(1.89,0.27)(1.62,0,27)(1.62,0)
\multiput(1.62,0)(0,0.03){9}{\drawline(0,0)(0.27,0)}
\multiput(1.62,0)(0.03,0){9}{\drawline(0,0)(0,0.27)}
\put(1.775,-0.15){\makebox(0,0){$J_{20}$}}
\path(2.16,0)(2.25,0)(2.25,0.09)(2.16,0.09)(2.16,0)
\multiput(2.16,0)(0,0.03){3}{\drawline(0,0)(0.09,0)}
\multiput(2.16,0)(0.03,0){3}{\drawline(0,0)(0,0.09)}
\put(2.225,-0.15){\makebox(0,0){\small$J_{220}$}}
\put(1.25,-0.4){\makebox(0, 0){Example~\ref{ESS.ex00}}}

\multiput(3,0)(0,2.43){2}{\drawline(0,0)(2.43,0)}
\multiput(3,0.27)(0,0.27){8}{\dottedline[\circle*{0.005}]{0.03}(0,0)(2.43,0)}
\multiput(3,0.81)(0,0.81){2}{\dottedline[\circle*{0.02}]{0.06}(0,0)(2.43,0)}
\multiput(3.27,0)(0.27,0){8}{\dottedline[\circle*{0.005}]{0.03}(0,0)(0,2.43)}
\multiput(3.81,0)(0.81,0){2}{\dottedline[\circle*{0.02}]{0.06}(0,0)(0,2.43)}
\multiput(3, 0)(2.43,0){2}{\drawline(0,0)(0,2.43)}
\shade\path(3.54,0)(4.08,0)(4.08,2.43)(3.54,2.43)(3.54,0)
\put(3.81,-0.15){\makebox(0,0){$R_1$}}
\shade\path(3.24,0)(3.3,0)(3.3,2.43)(3.24,2.43)(3.24,0)
\put(3.32,-0.15){\makebox(0,0){\small$R_2$}}
\shade\path(3.02889,0)(3.03111,0)(3.03111,2.43)(3.02889,2.43)(3.02889,0)
\put(3.08,-0.15){\makebox(0,0){\small$R_3$}}
\put(4.25,-0.4){\makebox(0, 0){Example~\ref{ESS.ex20}}}
\end{picture}
\caption{Examples~\ref{ESS.ex00} and \ref{ESS.ex20}}\label{SQ2}
\end{figure}

The second example is the case where the restriction of the Euclidean metric is not adapted.

\example[Figure~\ref{SQ2}]\label{ESS.ex20}
Set  $x_j = \frac 1{3^j} - \frac 1{3^{2j}}$, $y_j = \frac 1{3^j} + \frac 1{3^{2j}}$ and $R_j = [x_j, y_j] \times [0, 1]$ for any $j \ge 1$. Define $X = \sd Q{(\cup_{j \ge 1} \inte{R_j})}$. Let $T = \{w| w \in W_*, \inte{Q_w} \cap X \neq \emptyset\}$ and let $K_w = X \cap Q_w$ for any $w \in T$. Then $K: T \to \C(X)$ is a partition of $X$ parametrized by $(T, \A|_{T \times T}, \phi)$ by Proposition~\ref{ESS.prop10}. In this case, we easily see the following facts:
\begin{itemize}
\item
$\kappa(R_j) = 3^{2j}/2$ for any $j \ge 1$,
\item
If $w \in \cup_{m \ge j} (T)_m$, then $\sd{Q_w}{\inte{R_j}}$ is a rectangle,
\item
Set $(1)^n = \underset{\text{$n$ times}}{1\cdots{1}} \in (T)_n$. Then $\sd{Q_{(1)^{j - 1}}}{\inte{R_j}}$ has two connected components and $\kappa(Q_{(1)^{j - 1}} \cap R_j) = 2\cdot{3^{j + 1}}$.
\end{itemize}
These facts yield that $R_j \notin \R_{2\cdot{3^j}}^0 \cup \R_{2\cdot{3^j}}^1$ for sufficiently large $j$. By Theorem~\ref{ESS.thm10}, $d$ is not adapted to $g|_T$. In fact, $D^g_M((x_j, 0), (y_j, 0)) = 3^{-(j - 1)}$ for any $j \ge 1$ while $d((x_j, 0), (y_i, 0)) = 2{\cdot}3^{-2j}$. Hence the ratio between $D^g_M(\cdot, \cdot)$ and $d(\cdot, \cdot)$ is not bounded for any $M \ge 0$.\par
Furthermore, let $d_*(x, y) = \max\{|x_1 - y_1|, |x_2 - y_2|\}$ for any $x = (x_1, x_2), y = (y_1, y_2) \in X$. Then $g|_T = g_{d_*}$. Note that $d$ and $d_*$ are bi-Lipschitz equivalent. Since $d$ is not adapted to $g|_T$, it follows that $d_*$ is not adapted to $g|_{T} = g_{d_*}$ as well. Thus $d$ and  $d_*$ are not adapted.
\endexample

The third one is the case when the restriction of the Euclidean metric is not $1$-adapted but $2$-adapted. 

\example[Figure~\ref{SQ3}]\label{COM.ex10} 
Define
\[
w_*(j) = (1)^{j - 1}9(1)^j,\quad
R_j = Q_{w_*(j)}
\quad\text{and}\quad
k_m = \Big[\frac m2\Big]
\]
for $j \in \BbN$ and $m \in \BbN$. Note that $(1)^n = \underset{\text{$n$-times}}{1\ldots{1}}$ as is defined in Example~\ref{ESS.ex20}. Then it follows that $T = \sd{T^{(9)}}{\cup_{j \in \BbN} T^{(9)}_{w_*(j)}}$, where $T^{(9)}_w = \{wi_1i_2\ldots| i_1, i_2, \ldots \in \{1, \ldots, 9\}\}$. Let $g(w) = 3^{-|w|}$ for any $w \in T$. Define $w(m) = (1)^{m - 1}9$ and $v(m) = (1)^m$. Then $(w(m), (1)^{m - 1}8(3)^{k}, v(m))$ is a chain for $k = 0, 1, \ldots, m$. See Figure~\ref{SQ3}. Therefore, $w(m)$ and $v(m)$ are $1$-separated in $\LL_{3^{-m}}^g$ but not $2$-separated in $\LL_{3^{-2m}}^g$. This means that the condition ${\rm (EV5)}_M$ for $M = 1$ does not hold. Therefore, there exists no metric which is $1$-adapted to $g^{\a}$ for any $\a > 0$. On the other hand, since $\kappa(R_j) = 1$ for any $j \in \BbN$, the restriction of the Euclidean metric to $X$, which is denoted by $d$, is adapted to $g$. In fact, it is easy to see that $d$ is $2$-adapted to $g$. As a consequence, $d$ is not $1$-adapted but $2$-adapted to $g$.
\endexample

\begin{figure}
\centering
\setlength{\unitlength}{20mm}

\begin{picture}(4.3,2.5)(0.7,-0.5)
\linethickness{1pt}
\thinlines
\multiput(0,0)(1.26,0){3}{\drawline(0,0)(0,2.52)}
\multiput(0,0)(0,1.26){3}{\drawline(0,0)(2.52,0)}
\shade\path(1.26,1.26)(1.26,1.4)(1.4,1.4)(1.4,1.26)(1.26,1.26)
\multiput(0,0)(1.26,1.26){2}{\path(0,0)(1.26,0)(1.26,1.26)(0,1.26)(0,0)}
\put(0.63,0.63){\makebox(0,0){$v(m) = (1)^m$}}
\put(1.75,1.99){\makebox(0,0){$w(m)$}}
\put(1.95,1.79){\makebox(0,0){$= (1)^{m - 1}9$}}
\path(1.26,1.26)(1.12,1.26)(1.12,1.4)(1.26,1.4)(1.26,1.26)
\put(0.92,1.6){\vector(1,-1){0.2}}
\put(0.75,1.75){\makebox(0,0){$(1)^{m-1}8(3)^m$}}
\put(1.89,0.63){\makebox(0,0){$(1)^{m - 1}2$}}
\put(1.53,1.06){\vector(-1,1){0.2}}
\put(1.64,0.97){\makebox(0,0){$w_*(m)$}}
\put(1.25,-0.2){\makebox(0, 0){The chain $(v(m), (1)^{m-1}8(3)^m, w(m))$}}
\put(1.25,-0.5){\makebox(0, 0){Example~\ref{COM.ex10}}}

\multiput(3, 0)(2.43,0){2}{\drawline(0,0)(0,2.43)}
\multiput(3,0)(0,2.43){2}{\drawline(0,0)(2.7,0)}
\put(5.6,1.4){\vector(0,1){1.03}}
\put(5.6,1.03){\vector(0,-1){1.03}}
\shade\path(3.27,0.27)(5.16,0.27)(5.16,2.16)(3.27,2.16)(3.27,0.27)
\dottedline[\circle*{0.005}]{0.03}(5.16, .27)(5.43,0.27)
\dottedline[\circle*{0.005}]{0.03}(5.16, .27)(5.16,0)
\dottedline[\circle*{0.005}]{0.03}(3.27, 0.27)(3.27,0)
\dottedline[\circle*{0.005}]{0.03}(3, .27)(3.27,0.27)
\dottedline[\circle*{0.005}]{0.03}(3.27,2.16)(3.27,2.43)
\dottedline[\circle*{0.005}]{0.03}(3,2.16)(3.27,2.16)
\dottedline[\circle*{0.005}]{0.03}(5.16,2.16)(5.16,2.43)
\dottedline[\circle*{0.005}]{0.03}(5.16,2.16)(5.43,2.16)
\put(5.295,0.135){\vector(1,0){0.135}}
\put(5.295,0.135){\vector(-1,0){0.135}}
\put(4.215,1.215){\makebox(0,0){\Large$R(v)$}}
\put(5.7,1.215){\makebox(0,0){\Large$\big(\frac1{3}\big)^m$}}
\put(5.4,-0.2){\makebox(0,0){$\big(\frac{1}{3}\big)^{2m}$}}
\put(4.25,2.6){\makebox(0,0){{\Large$v9$}}}
\put(4.25,-0.2){\makebox(0,0){$v \in \{1, 3, 5, 7\}^{m - 1}$}}
\put(4.25,-0.5){\makebox(0,0){Example~\ref{ESS.ex10}}}
\end{picture}
\caption{Example~\ref{COM.ex10} and \ref{ESS.ex10}}\label{SQ3}
\end{figure}

In the fourth example, we do not have thickness while the restriction of the Euclidean metric is adapted.

\example[Figure~\ref{SQ3}]\label{ESS.ex10}
Define $\Delta{Q} = (\sd{\BbR^2}{\inte{Q}}) \cap Q$, which is the topological boundary of $Q$ as a subset of $\BbR^2$. Let $I_0 = \emptyset$ and let $E = \{1, 3, 5, 7\}$. Define $\{I_n\}_{n \ge 0}$ inductively by $I_{2m - 1} = \wI_{2m - 1}$ and $I_{2m} = J_m \cup \wI_{2m}$ for $m \ge 1$, where
\[
J_{m} = \{v9w | v \in E^{m - 1}, w \in W_m, Q_w \cap \Delta{Q} = \emptyset\}.
\]
$\{I_m\}_{m \ge 0}$ satisfies (SQ1), (SQ2) and (SQ3). In fact, if $J_{m, v} = \{v9w| w \in W_m, Q_w \cap \Delta{Q} = \emptyset\}$ for any $v \in E^{m - 1}$, $J_{m, v}$ is a collection of $(3^m - 2)^2$-words in $W_{2m}$. Set $R(v) = \cup_{u \in J_{m, v}}Q_u$ for any $m \ge 1$ and $v \in E^{m - 1}$.  See Figure~\ref{SQ3}. Then $\{R_j\}_{j \ge 1} = \{R(v)| m \ge 1, v \in E^{m - 1}\}$. More precisely $R(v) \subseteq Q_{v9}$ and $R(v)$ is a square which has the same center, i.e. the intersection of two diagonals, as $Q_{v9}$ and the length of the sides is $\frac 1{3^m}(1 - \frac 2{3^{m}})$. Note that the length of the sides of $Q_{v9}$ is $\frac 1{3^m}$. Hence the relative size of $R(v)$ in comparison with $Q_{v9}$ is monotonically increasing and convergent to $1$ as $m \to \infty$. The corresponding tree $(T, \A|_T, \phi)$ and the partition $K: T \to \C(X)$ of $X = \sd{Q}{\cup_{j \ge 1} \inte{R_j}}$ have the following properties:\\
Let $d$ be the restriction of the Euclidean metric to $X$. Then 
\begin{itemize}
\item[({\bf a})] 
$d$ is adapted to $g|_T$.
\item[({\bf b})]
$g|_T$ is exponential and uniformly finite.
\item[({\bf c})]
Let $\mu_*$ be the restriction of the Lebesgue measure on $X$. Then $\mu_*$ has the volume doubling property with respect to $d$.
\item[({\bf d})]
$\mu_*$ is not gentle with respect to $g|_T$.
\item[({\bf e})]
$\mu_*$ is not super-exponential.
\item[({\bf f})]
$g|_T$ is not thick.
\end{itemize}

In the rest, we present proofs of the above claims.\\
\noindent({\bf a})\,\,
Since $\kappa(R_m) = 1$ for any $m \ge 1$, we see that $\{R_m\}_{m \ge 1} \subseteq \R^0_1$. Hence Theorem~\ref{ESS.thm10} shows that $d$ is adapted to $g|_T$. In fact, $d$ is $1$-adapted to $g|_T$ in this case.\\
({\bf b})\,\,
This is included in the statement of Proposition~\ref{ESS.prop10}-(3).\\
({\bf c})\,\,\,
If $v \in \LL^{g|_T}_s$ and $Q_v = K_v$, then $\mu_*(K_v) = 9^{-|v|}$ and hence $\mu_*(K_u) \le 9^{-|u|} = 9^{-|v| + 1} \le 9\mu_*(K_v)$ for any $u \in \LL^{g|_T}_{3s}$. Therefore, $v \in \Theta(s, 3, k, 9)$ for any $k \ge 1$. On the other hand, for any $w \in T$, there exists $v \in \LL^{g|_T}_{s, 1}(w)$ such that $K_v = Q_v$. Therefore, we see that $\LL_{s, 1}^{g|_T}(w) \cap \Theta(s, 3, 3, 9) \neq \emptyset$. By Lemma~\ref{VDP.lemma30}, we have ({\bf c}).\\
({\bf d}) and ({\bf e})\,\,\,
Set $w(m) = (1)^{m - 1}9$. Then $K_{w(m)} = \sd{Q_{w(m)}}{\inte{R_m}}$, where $R_m = \cup_{w \in J_m} Q_w$. Then $\mu_*(K_{w(m)}) = 4(3^m - 1)3^{-4m}$. On the other hand, if $v(m) = (1)^{m - 1}8$, then $\mu_*(K_{v(m)}) = 3^{-2m}$. Since $K_{w(m)} \cap K_{v(m)} \neq \emptyset$, $\mu_*$ is not gentle with respect to $g|_T$. Moreover, since $K_{\pi(w(m))}$ contains $Q_{v(m)}$, we have $\mu_*(K_{\pi(w(m)}) \ge 3^{-2m}$. This implies that $\mu_*$ is not super-exponential.\newline
({\bf f})\,\,\,
To clarify the notation, we use $B(x, r) = \{y|y \in Q, |x - y| < r\}$ and $B_*(x, r) = B(x, r) \cap X$. This means that $B_*(x, r)$ is the ball of radius $r$ with respect to the metric $d$ on $X$. Assume that $g|_T$ is thick. Since $K$ is minimal, Proposition~\ref{ADD.prop21} implies that $K_{w(m)} \supseteq B_*(x, c3^{-m})$ for some $x \in K_{w(m)}$, where $c$ is independent of $m$ and $x$. However, for any $x \in K_{w(m)}$, there exists $y \in \sd{X}{K_{w(m)}}$ such that $|x - y| \le 2\cdot3^{-2m}$. This contradiction shows that $g|_T$ is not thick.
\endexample

\section{Gentleness and exponentiality}\label{GAE}
In this section, we show that the gentleness ``$\gen$'' is an equivalence relation among exponential weight functions. Moreover, the thickness of the interior, tightness, the uniformly finiteness and the existence of visual metric will be proven to be invariant under the gentle equivalence.\par
As in the section~\ref{VDP}, $(T, \A, \phi)$ is a locally finite tree with a reference point $\phi$, $(X, \O)$ is a compact metrizable topological space with no isolated point and $K: T \to \C(X, \O)$ is a partition of $X$ parametrized by $(T, \A, \phi)$. 

\definition\label{GAE.def05}
Define $\G_e(T)$ as the collection of exponential weight functions. 
\enddefinition

\thm\label{GAE.thm10}
The relation $\gen$ is an equivalence relation on $\G_e(T)$.
\endthm

Several steps of preparation are required to prove the above theorem.

\definition\label{GAE.def10}
(1)\,\,
Let $A \subseteq T$. For $m \ge 0$, we define $S^m(A) \subseteq T$ as
\[
S^m(A) = \bigcup_{w \in A} \{v| v \in (T)_{m + |w|}, [v]_{|w|} = w\}.
\]
(2)\,\,\,Let $g: T \to (0, 1]$ be a weight function. For any $w \in T$, define 
\[
N_g(w) = \min\{n| n \ge 0, \pi^n(w) \in \LL^g_{g(w)}\}
\]
and $\pi^*_g(w) = \pi^{N_g(w)}(w)$.\newline
(3)\,\,\,
$(u, v) \in T \times T$ is called an ordered pair if and only if $u \in T_v$ or $v \in T_u$. Define $|u, v| = ||u| - |v||$ for an ordered pair $(u, v)$.
\enddefinition

Note that if $g(w) < 1$, then we have
\[
N_g(w) = \min\{n| n \ge 0, g(\pi^{n + 1}(w)) > g(w)\}.
\]
Therefore, if $g(\pi(w)) > g(w)$ for any $w \in T$, then $N_g(w) = 0$ and $\pi^*_g(w) = w$ for any $w \in T$.\par
The following lemma is immediate from the definitions.

\lemma\label{GAE.lemma03}
Let $g: T \to (0, 1]$ be a super-exponential weight function, i.e. there exists $\c \in (0, 1)$ such that $g(w) \ge {\c}g(\pi(w))$ for any $w \in T$. If $(u, v)$ is an ordered pair, then $g(u) \le \c^{-|u, v|}g(v)$.
\endlemma

\lemma\label{GAE.lemma05}
Let $g: T \to (0, 1]$ be a weight function. If $g$ is sub-exponential, then $\sup_{w \in T}N_g(w) < +\infty$.
\endlemma

\demo
Since $g$ is sub-exponential, there exist $c \in (0, 1)$ and $m \ge 0$ such that $cg(w) \ge g(u)$ if $w \in T$, $u \in T_w$ and $|u, v| \ge m$. This immediately implies that $N_g(w) \le m$.
\enddemo

\lemma\label{GAE.lemma10}
Assume that $g, h \in \G_e(T)$ and $h$ is gentle with respect to $g$. Then there exist $M$ and $N$ such that  if $s \in (0, 1]$, $w \in \LL^h_s$ and $u \in S^M(\LL^g_{g(w), 1}(\pi^*_g(w)))$, then one can choose $n(u) \in [0, N]$ so that $\pi^{n(u)}(u) \in \LL^h_s$. Moreover, define $\eta^{g, h}_{s, w}: S^M(\LL^g_{g(w), 1}(\pi^*_g(w))) \to \LL^h_s$ by $\eta^{g, h}_{s, w}(u) = \pi^{n(u)}(u)$. Then $\LL^h_{s, 1}(w) \subseteq \eta^{g, h}_{s, w}(S^M(\LL^g_{g(w), 1}(\pi^*_g(w))))$. In particular, for any $s \in (0, 1]$, $w \in \LL^h_s$ and $v \in \LL^h_{s, 1}(w)$, there exists $u \in \LL^g_{g(w), 1}(\pi^*_g(w))$ such that $(u, v)$ is an ordered pair and $|u, v| \le \max\{M, N\}$.
\endlemma

\demo
Since $h$ is sub-exponential, there exist $c_1 \in (0, 1)$ and $m \ge 0$ such that $c_1h(w) \ge h(u)$ for any $w \in T$ and $u \in S^m(w)$. Let $w \in \LL_s^h$ and let $w' = \pi^*_g(w)$. Set $t = g(w)$.  Let $v \in \LL^g_{t, 1}(w')$. As $h$ is gentle with respect to $g$, there exists $c \ge 1$ such that
\[
h(w')/c \le h(v) \le ch(w'),
\]
where $c$ is independent of $s, w$ and $v$. By Lemma~\ref{GAE.lemma05}  and the fact that $h$ is super-exponential, there exists $c' \ge 1$ such that 
\[
h(w)/c \le h(v) \le c'h(w)
\]
for any $s$, $w$ and $v$. Using this, $h$ being sub-exponential and Proposition~\ref{VDP.prop10}, we see that there exist $c'' > 0$ and $M$ which are independent of $s$ and $w$ such that $c''s \le h(u) \le s$ for any $u \in S^M(\LL_{t, 1}^g(w'))$. Choose $k$ so that $c''(c_1)^{-k} > 1$. Then $h(\pi^{km}(u)) \ge (c_1)^{-k}h(u) \ge c''(c_1)^{-k}s > s$. Set $N = km - 1$. Then, for any $u \in S^M(\LL_{t, 1}^g(w'))$, there exists $n(u)$ such that $n(u) \le N$ and $\pi^{n(u)}(u) \in \LL^h_s$. Now for any $\rho \in \LL^h_{s, 1}(w)$, there exists $v \in \LL^g_{t, 1}(w')$ such that $(\rho, v)$ is an ordered pair. Since $\pi^{n(u)}(u) = \rho$ for any $u \in S^M(v)$, it follows that $\eta^{g, h}_{s, w}(S^M(\LL^g_{g(w), 1}(\pi^*_g(w)))) \supseteq \LL^h_{s, 1}(w)$. The rest is straightforward.
\enddemo

Finally we are ready to give a proof of Theorem~\ref{GAE.thm10}.

\demo[Proof of Theorem~\ref{GAE.thm10}]
Let $g, h, \xi \in \G_e(T)$. Then there exists $\c \in (0, 1)$ such that $g(w) \ge {\c}g(\pi(w))$, $h(w) \ge {\c}h(\pi(w))$ and $\xi(w) \ge {\c}\xi(\pi(w))$ for any $w \in T$.\par
First we show $g \gen g$. By Proposition~\ref{VDP.prop10}, there exists $c \in (0, 1)$ such that if $w \in \LL^g_s$, then $cg(w) \le s \le g(w)$. As a consequence, if $w, v \in \LL^g_s$, then $g(w) \le s/c \le g(v)/c$. Thus $g \gen g$.\par
Next assume $g \gen h$. Suppose that $w, v \in \LL^h_s$ and $K_w \cap K_v \neq \emptyset$. Since $v \in \LL^h_{s, 1}(w)$, Lemma~\ref{GAE.lemma10} implies that there exists $u \in \LL^{g}_{g(w), 1}(\pi^*_g(w))$ such that $(u, v)$ is an ordered pair and $|u, v| \le L$, where $L = \max\{M, N\}$. By Lemma~\ref{GAE.lemma03}, $g(v) \ge {\c}^{L}g(u) \ge \c^{L}g(w)$. Hence $h \gen g$.\par
Finally assume that $g \gen h$ and $h \gen \xi$. Suppose that $w, v \in \LL^{\xi}_s$ and $K_w \cap K_v \neq \emptyset$. Since $v \in \LL^{\xi}_{s, 1}(w)$, Lemma~\ref{GAE.lemma10} implies that there exists $u \in \LL^h_{h(w), 1}(\pi^*_h(w))$ such that $(u, v)$ is an ordered pair and $|u, v| \le L$. By Lemma~\ref{GAE.lemma03}, it follows that $g(v) \ge \c^{L}g(u)$. Set $s' = h(w)$ and $w' = \pi^*_h(w)$. Note that $w' \in \LL^h_s$ and $u \in \LL^h_{s', 1}(w')$. Again by Lemma~\ref{GAE.lemma10}, there exists $a \in \LL^g_{g(w'), 1}(\pi^*_g(w'))$ such that $(u, a)$ is an ordered pair and $|a, u| \le L$. Lemma~\ref{GAE.lemma03} shows that $g(u) \ge \c^{L}g(a) \ge \c^Lg(\pi^*_h(w))$. By Lemma~\ref{GAE.lemma05}, $N_h(w)$ is uniformly bounded and hence there exists $c_* > 0$ which is independent of $s$, $w$ and $v$ such that $g(\pi^*_h(w)) \ge c_*g(w)$. Combining these, we obtain $g(v) \ge \c^{2L}g(\pi^*_h(w)) \ge \c^{2L}c_3g(w)$. Hence $\xi \gen g$. Consequently we verify $g \gen \xi$ by the above arguments.
\enddemo

Next, we show the invariance of thickness, tightness and uniform finiteness under the equivalence relation $\gen$.

\thm\label{GAE.thm20}
Let $g, h \in \G_e(T)$. Suppose $g \gen h$. \newline
{\rm (1)}\,\,\,
Suppose that $\sup_{w \in T}\#(S(w)) < +\infty$. If $g$ is uniformly finite then so is $h$.\newline 
{\rm (2)}\,\,\,
If $g$ is thick, then so is $h$.\\
{\rm (3)}\,\,\,
If $g$ is tight, then so is $h$.
\endthm

We need the next lemma to prove Theorem~\ref{GAE.thm20}.

\lemma\label{COM.lemma300}
Let $g, h \in \G_e(T)$. Assume that $g$ is gentle with respect to $h$. Then for any $\a \in (0, 1]$ and $M \ge 0$, there exists $\c \in (0, 1)$ such that 
\[
U_M^g(x, {\a}g(w)) \supseteq U_M^h(x, {\c}h(w))
\]
for any $w \in T$ and $x \in K_w$.
\endlemma
\demo
Since $g$ and $h$ are exponential, there exist $c_1, c_2 \in (0, 1)$ and $m \ge 1$ such that $h(w) \ge c_2h(\pi(w)), g(w) \ge c_2g(\pi(w))$, $h(v) \le c_1h(w)$ and $g(v) \le c_1g(w)$ for any $w \in T$ and $v \in S^m(w)$. Moreover, since $g$ is gentle with respect to $g$, there exists $c > 1$ such that $g(w) \le cg(u)$ whenever $w, u \in \LL_s^h$ and $K_w \cap K_v \neq \emptyset$. Note that  $N_g(w) \le m$ and $N_h(w) \le m$ for any $w \in T$.\par
Let $w \in T$ and let $x \in K_w$.  Assume that $\c < (c_2)^{lm}$. Let $v \in \LL_{{\c}h(w), 0}^h(x)$. Then $h(\pi(v)) > {\c}h(w) \ge h(v)$. There exists $k \ge 0$ such that $\pi^k(v) \in \LL_{h(w)}^h$. Then $h(\pi^{k + 1}(v)) > h(w) \ge h(\pi^k(v))$. Thus we have
\[
{\c}h(\pi^{k + 1}(v)) \ge h(v)
\]
Therefore, it follows that $k + 1 \ge lm$. Let $w_* = \pi^{N_h(w)}(w)$. Then we see that $x \in K_{\pi^{k + 1}(v)} \cap K_{w_*}$. Therefore, $c^{-1}g(w_*) \le g(\pi^{k + 1}(v)) \le cg(w_*)$. Since $k + 1 \ge lm$ and $N_h(w) \le m$, it follows that
\[
g(v) \le (c_1)^{l}g(\pi^{k + 1}(v)) \le c(c_1)^lg(w_*) \le c(c_1)^l(c_2)^{-m}g(w).
\]
Now suppose that $(w(1), \ldots, w(M + 1))$ is a chain in $\LL_{{\c}h(w)}^h$ with $w(1) \in \LL_{{\c}h(w), 0}^h(x)$. Using the above arguments, we obtain
\[
g(w(i)) \le c^{i - 1}g(w(1)) \le c^i(c_1)^l(c_2)^{-m}g(w) \le c^{M + 1}(c_1)^l(c_2)^{-m}g(w)
\]
for any $i = 1, \ldots, M + 1$. Choosing $l$ so that $c^{M + 1}(c_1)^l(c_2)^{-m} < \a$, we see that $U_M^h(x, \c{h(w)}) \subseteq U_M^g(x, \a{g(w)})$.
\enddemo

\demo[Proof of Theorem~\ref{GAE.thm20}]
(1)\,\,\,
Set $L = \sup_{w \in T} \#(S(w))$. By Lemma~\ref{GAE.lemma10}, it follows that $\#(\LL^h_{s, 1}(w)) \le L^M\#(\LL^g_{g(w), 1}(\pi^*_g(w)))$. This suffices to the desired conclusion.\newline
(2)\,\,\,
By the thickness of $g$ and Proposition~\ref{ADD.prop20}, for any $M \ge 0$, there exists $\b > 0$ such that, for any $w \in T$,
\[
K_w \supseteq U_M^g(x, {\b}g(\pi(w)))
\]
for some $x \in K_w$. By Lemma~\ref{COM.lemma300}, there exists $\c \in (0, 1)$ such that 
\[
U_M^g(x, {\b}g(\pi(w))) \supseteq U_M^h(x, {\c}h(\pi(w)))
\]
 for any $w \in T$. Thus making use of Proposition~\ref{ADD.prop20} again, we see that $h$ is thick.\\
 (3)\,\,
 Since $g$ is tight, for any $M \ge 0$, there exists $\a > 0$ such that, for any $w \in T$, $\sd{K_w}{U_M^g(x, {\a}g(w))} \neq \emptyset$ for some $x \in K_w$. By Lemma~\ref{COM.lemma300}, there exists $\c \in (0, 1)$ such that $U_M^g(x, {\a}g(w)) \supseteq U_M^h(x, {\c}h(w))$ for any $w \in T$ and $x \in K_w$. Hence 
 \[
 \sup_{x, y \in K_w} \d_M^h(x, y) \ge \c{h(w)}
 \]
 for any $w \in T$. Thus we have shown that $h$ is tight.
\enddemo

Finally, the existence of visual metric is also invariant under $\gen$ as follows.

\thm\label{COM.thm30}
Assume that the partition $K: T \to \C(X, \O)$ is minimal. Let $g, h \in \G_e(T)$ and let $M \in \BbN$. Assume that $g \gen h$. Then $g$ is hyperbolic if and only if $h$ is hyperbolic.
\endthm
\demo
Since $g$ and $h$ are exponential, there exists $\lambda \in (0, 1)$ and $m \ge 1$ such that
\begin{align*}
g(w') &\le \lambda{g(w)} \le g(w'')\\
h(w') &\le \lambda{h(w)} \le h(w'')
\end{align*}
if $w \in T$, $w', w'' \in T_w$, $|w'| - |w| \ge m$ and $|w''| - |w| = 1$. Moreover, since $g \gen h$, there exists $\eta > 1$ such that if $w, v \in \LL_s^g$ and $K_w \cap K_v \neq \emptyset$, then $h(w) \le \eta{h(v)}$ and if $w, v \in \LL_s^h$ and $K_w \cap K_v \neq \emptyset$, then $g(w) \le \eta{g(v)}$. Fix $k \in \BbN$ satisfying $\eta^M\lambda^k < 1$.\par
Now assume that $g$ is hyperbolic. Then by Theorem~\ref{COM.thm10} and \ref{HYP.thm20}, $g$ satisfies ${\rm(EV5)}_M$. Let $w, v \in \LL_s^h$ and assume that $(w, v)$ is $M$-separated in $\LL_s^h$. Set $t = g(v)$. Suppose that $(w, v)$ is not $M$-separated in $\LL_{\lambda^{km}t}^g$. Then there exists a chain $(w_*(1), \ldots, w_*(M - 1))$ in $\LL_{\lambda^{km}t}^g$ such that $(w, w_*(1), \ldots, w_*(M - 1), v)$ is a chain. Choose $v_* \in \LL_{\lambda^{km}t}^g \cap T_v$ so that $K_{w_*(M - 1)} \cap K_{v_*} \neq \emptyset$.  Since $g(v_*) \le \lambda^{km}t = \lambda^{km}g(v)$, it follows that $|v_*| - |v| \ge km$. Then we have
\[
h(w_*(i)) \le \eta^Mh(v_*) \le \eta^M\lambda^kh(v) < h(v).
\]
Hence there exists a chain $(w(1), \ldots, w(M - 1))$ in $\LL_s^h$ such that $w_*(i) \in T_{w(i)}$ for any $i = 1, \ldots, M - 1$. This implies that $(w, v)$ is not $M$-separated in $\LL_s^h$. This contradiction implies that $(w, v)$ is $M$-separated in $\LL_{\lambda^{km}t}^g$. \par
Since ${\rm(EV5)}_M$ holds for $g$, we see that $(w, v)$ is $(M + 1)$-separated in $\LL_{\tau\lambda^{km}t}^g$. Set $t_* = \tau\lambda^{km}t$. Choose $v' \in \LL_{t_*}^g \cap T_v$. Then exchanging $g$ and $h$ and using the same argument as above, we see that $(w, v)$ is $(M + 1)$-separated in $\LL_{\lambda^{km}h(v')}^h$.\par
Since $h$ is exponential, Proposition~\ref{VDP.prop10} shows that there exists $c > 0$ such that $cr \le g(u) \le r$ for any $r \in (0, 1]$ and $u \in \LL_r^g$. Choose $n_*$ so that $\lambda^{n_*} < c\tau$. Suppose $|v'| - |v| \ge (km + n_*)m$. Then
\[
\lambda^{km + n_*}g(v) < c\tau\lambda^{km}g(v) \le ct_* \le g(v') \le \lambda^{km + n_*}g(v).
\]
This contradiction yields that $|v'| - |v| < (km + n_*)m$. Therefore, $h(v') \ge \lambda^{(km + n_*)m}h(v) \ge \lambda^{(km + n_*)m}s$. Thus $\lambda^{km}h(v') \ge \lambda^{(km + n_* + k)m}s$. Set $\tau_* = \lambda^{(km + n_* + k)m}$. Then $(w, v)$ is $(M + 1)$-separated in $\LL_{\tau_*s}^h$. Thus we have shown that ${\rm(EV5)}_M$ is satisfied for $h$. Using Theorem~\ref{COM.thm10} and \ref{HYP.thm20}, we see that $h$ is hyperbolic.
\enddemo

\setcounter{equation}{0}
\section{Quasisymmetry}\label{QSY}

In this section, we are going to identify the equivalence relation, gentleness ``$\gen$'' with the quasisymmetry ``$\qs$'' among the metrics under certain conditions.
As in the last section, $(T, \A, \phi)$ is a locally finite tree with a reference point $\phi$, $(X, \O)$ is a compact metrizable topological space with no isolated point and $K: T \to \C(X, \O)$ is a partition of $X$ parametrized by $(T, \A, \phi)$ throughout this section.

\definition[Quasisymmetry]\label{INT.def10}
A metric $\rho \in \D(X, \O)$ is said to be quasisymmetric to a metric $d \in \D(X, \O)$ if and only if there exists a homeomorphism $h$ from $[0, +\infty)$ to itself such that $h(0) = 0$ and, for any $t > 0$, $\rho(x, z) < h(t)\rho(x, y)$ whenever $d(x, z) < td(x, y)$. We write $\rho \qs d$ if $\rho$ is quasisymmetric to $d$.
\enddefinition

It is known that $\qs$ is an equivalence relation on $\D(X, \O)$.

\definition\label{QSY.def100}
Let $d \in \D(X, \O)$. We say that $d$ is (super-, sub-)exponential if and only if $g_d$ is (super-, sub-)exponential.
\enddefinition

Under the uniformly perfectness of a metric space defined below, we can utilize a useful equivalent condition for quasisymmetry obtained in \cite{Ki16}. See the details in the proof of Theorem~\ref{QSY.lemma10}.

\definition\label{QSY.def20}
A metric space $(X, d)$ is called uniformly perfect if and only if there exists $\e > 0$ such that $\sd{B_d(x, (1 + \e)r)}{B_d(x, r)} \neq \emptyset$ unless $B_d(x, r) = X$.
\enddefinition

\lemma\label{QSY.lemma10}
Let $d \in \D(X, \O)$. If $d$ is super-exponential, then $(X, d)$ is uniformly perfect.
\endlemma
 
\demo
Write $d_w = g_d(w)$ for any $w \in T$. Since $d$ is super-exponential, there exists $c_2 \in (0, 1)$ such that $d_w \ge c_2d_{\pi(w)}$ for any $w \in T$. Therefore, $s \ge d_w > c_2s$ if $w \in \LL^d_s$. For any $x \in X$ and $r \in (0, 1]$, choose $w \in \LL^d_{r/2, 0}(x)$. Then $d(x, y) \le d_w \le r/2$ for any $y \in K_w$. This shows $K_w \subseteq B_d(x, r)$. Since $\diam{B_d(x, c_2r/4), d} \le c_2r/2 < d_w$, it follows that $\sd{K_w}{B_d(x, c_2r/2)} \neq \emptyset$. Therefore $\sd{B_d(x, r)}{B_d(x, c_2r/2)} \neq \emptyset$. This shows that $(X, d)$ is uniformly perfect.
\enddemo

\definition\label{QSY.def120}
Define
\[
\D_{A, e}(X, \O) = \{d| d \in \D(X, \O), \text{$d$ is adapted and exponential.}\}
\]
\enddefinition

The next theorem is the main result of this section.

\thm\label{QSY.thm10}
Let $d \in \D_{A, e}(X, \O)$ and let $\rho \in \D(X, \O)$. Then $d \qs \rho$ if and only if $\rho \in \D_{A, e}(X, \O)$ and $d \gen \rho$. Moreover, if $d$ is $M$-adapted and $d \qs \rho$, then $\rho$ is $M$-adapted as well.
\endthm

\remark
In the case of natural partitions of self-similar sets introduced in Example~\ref{PAS.ex10}, the above theorem has been obtained in \cite{Ki18}.
\endremark

The following corollary is straightforward from the above theorem.

\cor\label{QSY.cor10}
Let $d, \rho \in \D_{A, e}(X, \O)$. Then $d \qs \rho$ if and only if $d \gen g$.
\endcor

The rest of this section is devoted to a proof of the above theorem.

\demo[Proof of Theorem~\ref{QSY.thm10}: Part 1]
Assume that $d$ and $\rho$ belong to $\D_{A, e}(X, \O)$. We show that if $d \gen \rho$, then $d \qs \rho$. By Lemma~\ref{QSY.lemma10}, both $(X, d)$ and $(X, \rho)$ are uniformly perfect. By \cite[Theorems~11.5 and 12.3]{Ki16}, $d \qs \rho$ is equivalent to the facts that there exists $\delta \in (0, 1)$ such that
\begin{equation}\label{SQD.eq30}\begin{split}
B_d(x, r) &\supseteq B_{\rho}(x, \delta\orho_d(x, r))\\
B_{\rho}(x, r) &\supseteq B_d(x, \delta\od_{\rho}(x, r))
\end{split}\end{equation}
and
\begin{equation}\label{SQD.eq40}\begin{split}
\orho_d(x, r/2) &\ge \delta\orho_d(x, r)\\
\od_{\rho}(x, r/2) &\ge \delta\od_{\rho}(x, r)
\end{split}\end{equation}
for any  $x \in X$ and $r > 0$, where $\orho_d(x, r) = \sup_{y \in B_d(x, r)} \rho(x, y)$ and $\od_d(x, r) = \sup_{y \in B_{\rho}(x, r)}d(x, y)$. We are going to show \eqref{SQD.eq30} and \eqref{SQD.eq40}. Since $d$ and $\rho$ are adapted, there exist $\b \in (0, 1)$, $\c > 1$ and $M \ge 1$ such that
\begin{align*}
U^d_M(x, \b{r}) &\subseteq B_d(x, r) \subseteq U^d_M(x, \c{r})\\
U^{\rho}_M(x, \b{r}) &\subseteq B_{\rho}(x, r) \subseteq U^{\rho}_M(x, \c{r})
\end{align*}
for any $x \in X$ and $r \in (0, 1]$.  By Lemma~\ref{COM.lemma300}, there exists $\a \in (0, 1)$ such that $U_M^{\rho}(x, \rho_w) \supseteq U_M^d(x, {\a}d_w)$ and $U_M^d(x, d_w) \supseteq U_M^{\rho}(x, {\a}\rho_w)$ for any $w \in T$ and $x \in K_w$. If $w \in \LL^d_{\c{r}/\a, 0}(x)$, then
\begin{equation}\label{SQD.eq90}
B_d(x, r) \subseteq U^d_M(x, {\c}r) \subseteq U^d_M(x, {\a}d_w) \subseteq U^{\rho}_M(x, \rho_w),
\end{equation}
where $w \in \LL^d_{x, {\c_1}r/\a}$. Hence for any $y \in B_d(x, r)$, there exists $(w(1), \ldots, w(k)) \in \CH_K(x, y)$ such that $k \le M + 1$ and $w(i) \in \LL^{\rho}_{\rho_w}$. Since $\rho(x, y) \le \sum_{i = 1}^k \rho_{w(i)} \le (M + 1)\rho_w$, we have
\[
\orho_d(x, r) \le (M + 1)\rho_w.
\]
Let $w \in \LL^d_{\c{r}/\a, 0}(x)$ as above. Since $\b/2 < 1 < \c/\a$, there exists $v \in T_w$ such that $v \in \LL^d_{\b{r}/2, 0}(x)$. Note that $\b{r}/2 \ge d_v$. Hence we have
\begin{equation}\label{SQD.eq100}
B_d\Big(x, \frac r2\Big) \supseteq U^d_M\Big(x, \frac{\b{r}}2\Big) \supseteq U^d_M(x, d_v) \supseteq U^{\rho}_M(x, \a\rho_{v}).
\end{equation}
Since $d$ is sub-exponential, the fact that $w \in \LL_{\c{r}/\a, 0}^d(x)$ and $v \in \LL_{\b{r}/2, 0}^d(x) \cap T_w$ implies that $|v| - |w|$ is uniformly bounded with respect to $x, r$ and $w$. This and the fact that $\rho$ is super-exponential imply that there exists $c > 0$ which is independent of $x, r$ and $w$ such that $\rho_v \ge c\rho_w$. Now we see that $\a\rho_v \ge \eta\orho_d(x, r)$, where $\eta = \a{c}/(M + 1)$. Hence
\[
B_d\Big(x, \frac r2\Big) \supseteq U^{\rho}_M(x, \eta\orho_d(x, r)) \supseteq B_{\rho}\Big(x, \frac{\eta}{\c}\orho_d(x, r)\Big).
\]
By the fact that $(X, \rho)$ is uniformly perfect, there exists $c_* \in (0, 1)$ such that $\sd{B_{\rho}(y,  t)}{B_{\rho}(y, c_*t)} \neq \emptyset$ unless $B_{\rho}(y, c_*t) = X$. Set $\delta = c_*\eta/\c$. In case $B_{\rho}(x, \delta\orho_d(x, r)) = X$, we have $\orho_d(x, r/2) = \orho_d(x, r)$. Otherwise, there exists $z \in B_d(x, r/2)$ such that $\rho(x, z) \ge \delta\orho_d(x, r)$. In each case, we have $\orho_d(x, r/2) \ge \delta\orho_d(x, r)$. Furthermore, $B_d(x, r) \supseteq B_{\rho}(x, \eta\orho_d(x, r)/\c) \supseteq B_{\rho}(x, \delta\orho_d(x, r))$. Thus we have obtained halves of \eqref{SQD.eq30} and \eqref{SQD.eq40}. Exchanging $d$ and $\rho$, we have the other halves of \eqref{SQD.eq30} and \eqref{SQD.eq40}.
\enddemo

\lemma\label{QSY.lemma30}
Let $d \in \D_{A, e}(X, \O)$ and let $\rho \in \D(X, \O)$. Assume that $d \qs \rho$. Let $\delta \in (0, 1)$ be the constant appearing in \eqref{SQD.eq30} and \eqref{SQD.eq40}.\\
{\rm (1)}\,\,
For any $w \in T$ and $x, y \in K_w$,
\[
\orho_d(x, d_w) \le \delta^{-1}\orho_d(y, d_w).
\]
{\rm (2)}\,\,
There exists $c > 0$ such that
\[
c\orho_d(x, d_w) \le \rho_w \le \delta^{-1}\orho_d(x, d_w)
\]
for any $w \in T$ and $x \in K_w$.
\endlemma

\demo
Assume $d \qs \rho$. Lemma~\ref{QSY.lemma10} implies that $(X, d)$ is uniformly perfect. Since $d \qs \rho$, $(X, \rho)$ is uniformly perfect as well. Hence \eqref{SQD.eq30} and \eqref{SQD.eq40} hold.\\
(1)\,\,
Since $B_d(x, d_w) \subseteq B_d(y, 2d_w)$, it follows that $\orho_d(x, d_w) \le \orho_d(y, 2d_w)$. Applying \eqref{SQD.eq40}, we obtain the desired inequality.\\
(2)\,\,
For any $x \in K_w$, $K_w \subseteq B_d(x, 2d_w)$. Hence $\rho_w \le \orho_d(x, 2d_w)$. By \eqref{SQD.eq40}, we see that
\[
\rho_w \le \delta^{-1}\orho_d(x, d_w).
\]
Set $s = d_w/2$ and choose $v \in T_w \cap \LL^d_s$. Since $d$ is adapted and tight, there exists ${\c} > 0$ which is independent of $w, v$ and $s$ such that
\[
\sd{K_v}{B_d(z, {\c}d_v)} \neq \emptyset
\]
for some $z \in K_v$. By \eqref{SQD.eq30},
\[
\sd{K_v}{B_{\rho}(z, \delta\orho_d(z, {\c}d_v))} \neq \emptyset.
\]
Hence $\rho_w \ge \delta\orho_d(z, {\c}d_v)$. Since $d$ is super-exponential, there exists $\c' > 0$ which is independent of $w, v$ and $s$ such that $\c{d_v} \ge \c'd_w$.  Choose $n \ge 1$ so that $2^{n - 1}\c' \ge 1$. Using \eqref{SQD.eq40} $n$-times,  we have
\[
\rho_w \ge \delta\orho_d(z, {\c'}d_w) = \delta^{n + 1}\orho_d(z, d_w).
\]
By (1), if $c = \delta^{n + 2}$, then $\rho_w \ge c\orho_d(x, d_w)$.
\enddemo

\demo[Proof of Theorem~\ref{QSY.thm10}: Part 2]
Assume that $d \in \D_{A, e}(X, \O)$. We show that if $d \qs \rho$, then $\rho \in \D_{A, e}(X, \O)$ and $d \gen \rho$. As in the proof of Lemma~\ref{QSY.lemma30}, \eqref{SQD.eq30} and \eqref{SQD.eq40} hold.\\
{\bf Claim 1}\,\,
$\rho$ is super-exponential.\newline
Proof of Claim 1:\,\,
Since $d$ is super-exponential, there exists $c' \in (0, 1)$ such that $d_w \ge c'd_{\pi(w)}$ for any $w \in T$. Choose $l \ge 1$ so that $2^{l}c' \ge 1$. By Lemma~\ref{QSY.lemma30}-(2) and \eqref{SQD.eq40}, if $x \in K_w$, then
\[
\rho_{w} \ge c\orho_d(x, d_{w}) \ge c\delta^{l}\orho_d(z, 2^ld_w) \ge c\delta^l\orho_d(x, d_{\pi(w)}) \ge c\delta^{l + 1}\rho_{\pi(w)}.
\]
{\bf Claim 2}\,\,
$\rho$ is sub-exponential.\newline
Proof of Claim 2:\,\,
Since $d$ is sub-exponential, there exist $c_1 \in (0, 1)$ and $m \ge 1$ such that 
\[
d_{v'} \le c_1d_w
\]
for any $w \in T$ and $v' \in S^m(w)$. Let $w \in T$. If $v \in S^{mj}(w)$ for $j \ge 1$ and $x \in K_v$, then by Lemma~\ref{QSY.lemma30}-(1)
\begin{equation}\label{SQD.eq60}
\rho_{v} \le \delta^{-1}\orho_d(x, d_v) \le \delta^{-1}\orho_d(x, (c_1)^jd_w).
\end{equation}

On the other hand, by \cite[Proposition 11.7]{Ki16}, there exist $\lambda \in (0, 1)$ and $c'' > 0$ such that
\[
\orho_d(x, c_1s) \le c''\lambda\orho_d(x, s)
\]
for any $x \in X$ and $s \in (0, 1]$. This together with \eqref{SQD.eq60} and Lemma~\ref{QSY.lemma30}-(2) yields
\[
\rho_{v} \le \delta^{-1}\orho_d(x, (c_1)^jd_w) \le \delta^{-1}c''\lambda^j\orho_d(x, d_w) \le \delta^{-1}c''\lambda^jc^{-1}\rho_w.
\]
Choosing $j$ so that $\delta^{-1}c''\lambda^jc^{-1} < 1$, we see that $\rho$ is sub-exponential.\\
{\bf Claim 3}\,\,
$d \gen \rho$.\newline
Proof of Claim 3:\,\,
Since $d$ is super-exponential, there exists $c_2 \in (0, 1)$ such that
\begin{equation}\label{QSY.eq70}
s \ge d_w > c_2s
\end{equation}
for any $s \in (0, 1]$ and $w \in \LL^d_s$. Let $w, v \in \LL_s^d$ with $K_w \cap K_v \neq \emptyset$. Then $d_w \le d_v/c_2$. Choose $k \ge 1$ so that $2^kc_2 \ge 1$. If $x \in K_w \cap K_v$, then by Lemma~\ref{QSY.lemma30}-(2) and \eqref{SQD.eq40},
\[
\rho_w \le \delta^{-1}\orho_d(x, d_w) \le \delta^{-1}\orho(x, d_v/c_2) \le \delta^{-(k + 1)}\orho(x, d_v) \le c^{-1}\delta^{-(k + 1)}\rho_v.
\]
Hence $d \gen \rho$.\\
{\bf Claim 4}\,\,
$\rho$ is adapted. More precisely, if $d$ is $M$-adapted, then so is $\rho$.\newline
Proof of Claim 4:\,\,
Assume that $d$ is $M$-adapted. Let $x \in X$ and let $s \in (0, 1]$. Then there exists $\a > 0$ which is independent of $x$ and $s$ such that $U^d_M(x, {\a}s) \supseteq B_d(x, s)$.  Let $w \in \LL^{\rho}_{s, 0}(x)$. Since $\rho$ is super-exponential, there exists $b \in (0, 1)$ which is independent of $w$ and $s$ such that $\rho_w \ge bs$. By Lemma~\ref{COM.lemma300}, there exists $\c > 0$ such that $U_M^{\rho}(x, \rho_w) \supseteq U_M^d(x, {\c}d_w)$ for any $w \in T$ and $x \in K_w$. Choose $p \ge 1$ so that $2^p\c/\a \ge 1$. Then by Lemma~\ref{QSY.lemma30}-(2), \eqref{SQD.eq30} and \eqref{SQD.eq40},
\begin{multline*}
U^{\rho}_M(x, s) \supseteq U^{\rho}_M(x, \rho_w) \supseteq U^d_M(x, {\c}d_w)\\
 \supseteq B_d\Big(x, \frac{\c}{\a}d_w\Big) \supseteq B_{\rho}\Big(x, \delta\orho_d\big(x, \frac{\c}{\a}d_w\big)\Big) \supseteq B_{\rho}(x, \delta^{p + 1}\orho_d(x, d_w))\\
  \supseteq B_{\rho}(x, \delta^{p + 2}\rho_w) \supseteq B_{\rho}(x, \delta^{p + 2}bs).
\end{multline*}
On the other hand, let $x \in K$ and let $r \in (0, 1]$. Then for any $y \in U^{\rho}_M(x, r)$, there exists $(w(1), \ldots, w(M + 1)) \in \CH_K(x, y)$ such that $w(i) \in \LL^{\rho}_r$ for any $i$. It follows that
\[
\rho(x, y) \le \sum_{i = 1}^{M + 1} \rho_{w(i)} \le (M + 1)r.
\]
This shows that $U_M^{\rho}(x, r) \subseteq B_{\rho}(x, (M + 1)r)$.
Thus we have shown that $\rho$ is $M$-adapted.\par
Using Theorem~\ref{GAE.thm20}-(2), we see that $g_{\rho}$ is thick and hence $\rho \in \D_{A, e}(X, \O)$.  Thus we have shown the desired statement.
\enddemo

\part{Characterization of Ahlfors regular conformal dimension}

\section{Construction of adapted metric I}\label{CAM}
In this section, we present a sufficient condition for the existence of an adapted metric to a given weight function. The sufficient condition obtained in this section is useful to construct an Ahlfors regular metric later.\par
Let $(T, \A, \phi)$ be a locally finite tree with a reference point $\phi$ and let $(X, \O)$ be a compact metrizable topological space with no isolated point. Moreover let $K: T \to \C(X, \O)$ be a minimal partition.

\definition\label{CAM.def10}
Let $M \ge 1$. A chain $(w(1), \ldots, w(M + 1))$ of $K$ is called a horizontal $M$-chain of $K$ if and only if $|w(i)| = |w(i + 1)|$ and $K_{w(i)} \cap K_{w(i + 1)} \neq \emptyset$ for any $i = 1, \ldots, M$. Define 
\begin{multline*}
\GG_M(w, K) = \{u | u \in (T)_{|w|},\,\, \text{there exists a horizontal $M$-chain}\\
\text{$(w(1), \ldots, w(M + 1))$ of $K$ such that $w(1) = w$ and $w(M + 1) = u$}\}
\end{multline*}
for $w \in T$,
\[
J_{M, n}^h(K) = \{(w, u)| w, u \in (T)_n, u \in \GG_M(w, T)\}
\]
for $n \ge 0$,
\[
J_M^h(K) = \bigcup_{n \ge 0} J_{M, n}^h(K),
\]
\[
J^v_M(K) = \{(w, u)| w, u \in T, w \in \pi(\GG_M(u, K)) \,\,\text{or}\,\,u \in \pi(\GG_M(w, T))\},
\]
and
\[
J_M(K) = J_M^v(K) \cup J_M^h(K).
\]
A sequence $(w(1), \ldots, w(m)) \in T$ is called an $M$-jumping path, or an $M$-jpath for short, (resp. horizontal $M$-jumping path, or horizontal $M$-jpath) of $K$ if and only if $(w(i), w(i + 1)) \in J_M(K)$ (resp. $(w(i), w(i + 1)) \in J^h_M(K)$) for any $i = 1, \ldots, m - 1$. Furthermore define
\[
U_M(w, K) = \bigcup_{v \in \GG_M(w, K)} K_v
\]
for any $w \in T$.
\enddefinition

\remark
Define a weight function $h_* : T \to (0, 1]$ as $h_*(w) = 2^{-|w|}$ for any $w \in T$. Then $(T)_m = \LL^{h_*}_{2^{-m}}$ and $\GG_M(w, K) = \LL^{h_*}_{2^{-|w|}, M}(w)$.
\endremark

\remark
Note that the horizontal vertices $E^h_m$ defined in Definition~\ref{HF.def10} is equal to $J^h_{1, m}$. On the contrary the collection $J_1^v(K)$ of the vertices in the vertical direction in $J_1(K)$ is strictly larger than that of vertical edges in the resolution $(T, \B)$ in general.
\endremark

If no confusion may occur, we are going to omit $K$ in $\GG_M(w, K)$, $U_M(w, K)$, $J_{M,n}^h(K)$, $J_M^v(w, K)$, $J_M^h(K)$ and $J_M(K)$ and write $\GG_M(w)$, $U_M(w)$, $J_{M, n}^h$, $J_M^v$, $J_M^h$ and $J_M$ respectively. Moreover, in such a case, we simply say a (horizontal) $M$-jpath instead of a (horizontal) $M$-jpath of $K$.

\definition\label{CAM.def20}
(1)\,\,
For $w \in T$ and $M \in \BbN$, define
\begin{multline*}
\C_w^M = \{(w(1), \ldots, w(m))|\pi(w(i)) \in \GG_M(w)\,\,\text{for any $i = 1, \ldots, m$}, \\
\text{there exist $w(0) \in S(w)$ and $w(m + 1) \in \sd{(T)_{|w| + 1}}S(\GG_M(w))$ such that}\\\text{$(w(0), w(1), \ldots, w(m), w(m + 1))$ is a horizontal $M$-jpath.}\}
\end{multline*}
(2) A function $\vp: T \to (0, \infty)$ is called $M$-balanced if and only if 
\[
\sum_{i = 1}^{m} \vp(w(i)) \ge \vp(\pi(w(m)))
\]
for any $w \in T$ and $(w(1), \ldots, w(m)) \in \C_w^M$.
\enddefinition

\remark
If $J^h_M = \emptyset$, then $\C_w^M = \emptyset$ for any $w \in T$ as well. Therefore, in this case, every $\vp : T \to (0, \infty)$ is $M$-balanced. This happens if and only if the original set is (homeomorphic to) the Cantor set.
\endremark

\thm\label{CAM.thm10}
Define $h_*: T \to (0, 1]$ by $h_*(w) = 2^{-|w|}$. Let $g \in \G_e(T)$. Assume that $g \gen h_*$. If there exists $\vp:T \to (0, \infty)$ such that $\vp$ is $M$-balanced and $\vp \bl g$, i.e. there exist $c_1, c_2 > 0$ such that
\[
c_1g(w) \le \vp(w) \le c_2g(w)
\]
for any $w \in T$, then there exists a metric $\rho \in \D(X, \O)$ which is $M$-adapted to $g$.
\endthm

The rest of this section is devoted to a proof of the above theorem. Throughout this section, $g$ and $\vp$ are assumed to satisfy the conditions required in Theorem~\ref{CAM.thm10}.\par
To begin with, we are going to summarize useful facts on $g$ following immediately from the assumptions.

\prop\label{CAM.prop10}
Let $g$ be a weight function. Assume that $g \gen h_*$ and that $g$ is exponential. \\
{\rm (1)}\,\,
There exist $\kappa \in (0, 1)$ and $N_0 \in \BbN$ such that if $|w| = |v|$ and $K_w \cap K_v \neq \emptyset$, then
\begin{equation}\label{CAM.eq10}
g(w) \ge \kappa{g(v)}.
\end{equation}
and if $w, v \in \LL_s^g$ and $K_w \cap K_v \neq \emptyset$, then
\begin{equation}\label{CAM.eq20}
|w| \le |v| + N_0.
\end{equation}
{\rm (2)}\,\,There exist $\eta \in (0, 1)$ and $n_0$ such that
\begin{equation}\label{CAM.eq30}
{\eta}g(\pi(w)) \le g(w)\quad\text{and}\quad g(v) \le {\eta}g(w)
\end{equation}
if $w \in T$, $v \in T_w$ and $|v| \ge |w| + n_0$. Moreover,
\begin{equation}\label{CAM.eq40}
\eta{s} < g(w) \le s
\end{equation}
for any $s \in (0, 1]$ and $w \in \LL_s^g$. In other words, if $g(w) \le {\eta}s$, then $w \in \LL_{t}^g$ for some $t < s$.
\endprop

The constants $\kappa, N_0, n_0$ and $\eta$ are fixed in the rest of this paper.

\lemma\label{CAM.lemma10}
Assume that $\vp: T \to (0, \infty)$ is $M$-balanced. Let $m \ge 2$ and let $\bp  = (w(1), \ldots, w(m))$ be an $M$-jpath satisfying $|w(1)| = |w(m)| = |w(i)| - 1$ for any $i = \{2, \ldots, m - 1\}$. Then $w(m) \in \GG_M(w(1))$ or there exists a horizontal $M$-jpath $\bp' = (v(1), \ldots, v(n))$ such that $w(1) = v(1)$, $w(m) = v(n)$ and 
\[
\sum_{i = 2}^{m - 1} \vp(w(i)) \ge \sum_{j = 2}^{n - 1}\vp( v(j)).
\]
\endlemma

\demo
Assume that $w(m) \notin \GG_M(w(1))$. This implies that $m \ge 3$. We use an induction on $m$. Let $m \ge 3$. Since $w(1) \in \pi(\GG_M(w(2)))$, we have $\pi(w(2)) \in \GG_M(w(1))$. \\
Case 1:\,\,Suppose $\pi(w(2)), \ldots, \pi(w(m - 1)) \in \GG_M(w(1))$. \\
In this case, since $(w(m - 1), w(m)) \in J_M^v$, there exists $u \in \GG_M(w(m - 1))$ such that $\pi(u) = w(m)$. By the fact that $\pi(u) = w(m) \notin \GG_M(w(1))$, we confirm $(w(2), \ldots, w(m - 1)) \in \C_w^M$. Since $\vp$ is $M$-balanced, we see
\[
\vp(w(2)) + \cdots + \vp(w(m - 1)) \ge \vp(\pi(w(m - 1))).
\]
Hence $(w(1), \pi(w(m - 1)), w(m))$ is the desired horizontal $M$-jpath. Note that if $m = 3$, then Case 1 applies.\\
Case 2:\,\,Suppose that $\pi(w(2)), \ldots, \pi(w(i)) \in \GG_M(w(1))$ and $\pi(w(i + 1)) \notin \GG_M(w(1))$ for some $i \in \{2, \ldots, m - 2\}$.\\
In this case, set $\widetilde{p} = (w(1), \pi(w(i)), w(i + 1), \ldots, w(m))$. Since $(w(1), w(2)) \in J^v_M$, there exists $v \in T$ such that $\pi(v) = w(1)$ and $v \in \GG_M(w(2))$. Therefore, $(v, w(2), \ldots, w(i), w(i + 1))$ is a horizontal $M$-jpath and $(w(2), \ldots, w(i)) \in \C^M_{w(1)}$. Since $\vp$ is $M$-balanced, it follows that
\[
\vp(w(2)) + \ldots, \vp(w(i)) \ge \vp(\pi(w(i)).
\]
Therefore 
\[
\sum_{k \ge 2}^{m - 1}\vp(w(k)) \ge \vp(\pi(w(i))) + \sum_{j = i + 1}^{m - 1} \vp(w(j)).
\]
Applying induction hypothesis to $(\pi(w(i)), w(i + 1), \ldots, w(m - 1), w(m))$, we obtain the desired result in this case. \\
\enddemo

Repeated use of the above lemma yields the following fact.
\lemma\label{CAM.lemma20}
Assume that $\vp$ is $M$-balanced. Let $\bp = (w(1), \ldots, w(m))$ is an $M$-jpath satisfying $|w(1)| = |w(m)|$. Then $w(m) \in \GG_M(w(1))$ or there exists an $M$-jpath $\bp' = (v(1), \ldots, v(k))$ such that $v(1) = w(1)$, $v(k) = w(m)$, $|v(i)| \le |v(1)|$ for any $i = 1, \ldots, k$, $k \le m$ and
\[
\sum_{i = 2}^{m -1} \vp(w(i)) \ge \sum_{j = 2}^{k - 1} \vp(v(j)).
\]
\endlemma

\lemma\label{CAM.lemma30}
Assume that $\vp$ is $M$-balanced. Let $\bp = (w(1), \ldots, w(m))$ be an $M$-jpath. Suppose that $w, v \in T$, $|w| = |v|$, $w \notin \GG_M(v)$, $w(1) \in T_w$ and $w(m) \in T_v$.\\
(1)\,\,
There exists an $M$-jpath $(v(1), \ldots, v(k))$ such that $v(1) = w, w(k) = v$, $|v(j)| \le |w|$ for any $j = 1, \ldots, k$ and
\[
\sum_{i = 2}^{m - 1} \vp(w(i)) \ge \sum_{j = 2}^{k - 1} \vp(v(j))
\]
(2) Assume that there exists $\kappa_0 \in (0, 1)$ such that $\vp(w) \ge \kappa_0\vp(v)$ for any $(w, v) \in J_M$. Then
\[
\sum_{i = 2}^{m - 1} \vp(w(i)) \ge \kappa_0\max\{\vp(w), \vp(v)\}.
\]
\endlemma

\demo
(1)\,\,Since $w(1) \in T_w$ and $w(m) \in T_v$, there exist $n_1, n_2 \ge 0$ such that $\pi^{n_1}(w(1)) = w$ and $\pi^{n_2}(w(m)) = v$. We use an induction on $n_1 + n_2$. If $n_1 + n_2 = 0$, then Lemma~\ref{CAM.lemma20} suffices. Assume that $n_1 + n_2 \ge 1$. Suppose that there exists $(w(k), w(k + 1), \ldots, w(k + l))$ such that $|w(k)| = |w(k + l)|$ and $|w(k + i)| = |w(k)| + 1$ for any $i = 1, \ldots, l - 1$. We call such a sequence $(w(k), \ldots, w(k + l))$ as a plateau. If $w(k + l) \in \GG_M(w(k))$, then we replace this part by $(w(k), w(k + l))$. Otherwise using  Lemma~\ref{CAM.lemma10}, we can replace it by a horizontal $M$-jpath $(w(k), w'(1), \ldots, w'(k'), w(k + l))$ satisfying
\[
\sum_{i = 1}^{l - 1} \vp(w(k + i)) \ge \sum_{j = 1}^{k'} \vp(w'(j)).
\]
Making iterated use of this procedure, we may assume that $(w(1), \ldots, w(m))$ contains no plateau without loss of generality.  Suppose $n_1 \ge n_2$. Then $n_1 \ge 1$. If $|w(1)| \ge |w(2)|$, then $(\pi(w(1)),  w(2)) \in J_M$. Hence  $(\pi(w(1)), w(2), \ldots, w(m))$ is an $M$-jpath to which the induction hypothes applies. Hence we assume that $|w(1)| + 1 = |w(2)|$. Then no plateau assumption yields $|w(1)| < |w(2)| \le \ldots \le |w(m)|$. Therefore, we have $n_1 < n_2$ and this contradict the fact that $n_1 \ge n_2$. If $n_1 \le n_2$, then same argument works by replacing $w(1)$ with $w(m)$. Thus we have shown (1).\\
(2)\,\,By (1), there exists an $M$-jpath $(v(1), \ldots, v(k))$ such that $v(1) = w, v(k) = v$, $|v(i)| \le |w|$ for any $i = 1, \ldots, k$ and
\[
\sum_{i = 2}^{m - 1}\vp(w(i)) \ge \sum_{j = 2}^{k - 1} \vp(v(j)).
\]
If $m = 3$, we have $(w, v(2)), (v(2), v) \in J_M$. Hence 
\[
\sum_{i = 2}^{m - 1} \vp(w(i)) \ge \vp(v(2)) \ge \kappa_0\max\{\vp(w), \vp(v)\}.
\]
If $m > 3$, then $(w, v(1)), (v(k - 1), v) \in J_M$ and hence
\[
\sum_{i = 2}^{m - 1}\vp(w(i)) \ge \vp(v(2)) + \vp(v(k - 1)) \ge \kappa_0(\vp(w) + \vp(v)) \ge \kappa_0\max\{\vp(w), \vp(v)\}.
\]
\enddemo

\definition\label{CAM.def30}
Let $\vp: T \to (0, \infty)$.\\
(1)\,\, For an $M$-jpath $\bp = (w(1), \ldots, w(m))$, we define 
\[
\ell_M^{\vp}(\bp) = \sum_{i = 1}^m \vp(w(i)).
\]
(2)\,\,For a chain $\bp = (w(1), \ldots, w(m))$ of $K$, define $L_{\vp}(\bp)$ by
\[
L_{\vp}(\bp) = \sum_{i = 1}^m \vp(w(i)).
\]
\enddefinition

Since there is no jump in a $1$-jpath, we say $1$-path instead of $1$-jpath in the followings.

\lemma\label{CAM.lemma35}
Let $g \in \G_e(T)$. For any chain $\bp = (w(1), \ldots, w(m))$ of $K$, there exists a $1$-path $\widehat{\bp} = (v(1), \ldots, v(k))$ of $K$ such that $w(1) = v(1)$, $w(m) = v(k)$ and 
\[
L_g(\bp) \ge c\ell_1^g(\widehat{\bp}),
\]
where $c > 0$ is independent of $\bp$ and $\widehat{\bp}$.
\endlemma

\demo
By \eqref{CAM.eq30}, there exists $c \ge 1$ and $\lambda \in (0, 1)$ such that $g(w) \le c\lambda^kg(\pi^k(w))$ for any $w \in T$ and any $k \ge 0$. \par
Now we start to construct a $1$-path $\widehat{\bp}$ by inserting a sequence between $w(i)$ and $w(i + 1)$ for each $i$ with $|w(i)| \neq |w(i + 1)|$. If $|w(i)| > |w(i + 1)|$, then there exists $v \in T_{w(i + 1)}$ such that $|v| = |w(i)|$ and $K_{w(i)} \cap K_{v} \neq \emptyset$. Let $k_i = |w(i)| - |w(i + 1)|$. Then $\pi^{k_i}(v) = w(i + 1)$ and $(w(i), v, \pi(v), \ldots, \pi^{k_i}(v))$ is a $1$-path. Moreover,
\begin{equation}\label{CAM.eq100}
\sum_{j = 0}^{k_i - 1} g(\pi^j(v)) \le (c\lambda^{k_i} + c\lambda^{k_i - 1} + \ldots + c\lambda)g(w(i + 1)) \le \frac{c\lambda}{1 - \lambda}g(w(i + 1)).
\end{equation}
Next suppose that $|w(i)| < |w(i + 1)|$. Then using a similar discussion as above, we find a $1$-path $(\pi^{k_i}(v), \pi^{k_i - 1}(v), \ldots, v, w(i + 1))$ satisfying $\pi^{k_i}(v) = w(i)$ and a counterpart of \eqref{CAM.eq100}. Inserting sequences in this manner, we obtain the desired $\widehat{\bp}$. By \eqref{CAM.eq100}, it follows that
\[
\Big(1 + \frac{2c\lambda}{1 - \lambda}\Big)L_g(\bp) \ge \ell_1^g(\widehat{\bp}).
\]
\enddemo



\demo[Proof of Theorem~\ref{CAM.thm10}]
Fix $M \in \BbN$. Write $\d(x, y) = \delta_M^g(x, y)$ for any $x, y \in X$. For $A \subseteq X$, we set
\begin{align*}
\LL_s^g(A) &= \{w| w \in \LL_s^g, K_w \cap A \neq \emptyset\}\\
(T)_n(A) &= \{w| w \in (T)_n, K_w \cap A \neq \emptyset\}
\end{align*}
In particular, we write $\LL_s^g(x, y) = \LL_s^g(\{x, y\})$ and $(T)_n(x, y) = (T)_n(\{x, y\})$.\par
Since $\vp \bl g$, using Proposition~\ref{CAM.prop10}, we see that there exists $\kappa_0 \in (0, 1)$ such that  $\vp(w) \ge \kappa_0\vp(v)$ for any $(w, v) \in J_M$. \\
{\bf Claim 1}: There exists $N_1 \in \BbN$ such that
\[
||w| - |v|| \le N_1
\]
for any $x, y \in X$ and $w, v \in \LL_{\kappa^M\d(x, y)}^g(x, y)$.
 \demo[Proof of Claim 1] 
 By \eqref{CAM.eq30}, there exists $N' \in \BbN$ such that if $w \in \LL_{\kappa^Ms}^g$, $w' \in \LL_s^g$ and $w \in T_{w'}$, then $|w| \le |w'| + N'$. Let $w, v \in \LL_{\kappa^M\d(x, y)}^g(x, y)$. If $w, v \in \LL_{\kappa^M\d(x, y)}^g(x)$ or $w, v \in \LL_{\kappa^M\d(x, y)}^g(y)$, then we have $K_w \cap K_v \neq \emptyset$ and hence $||w| - |v|| \le N_0$ by \eqref{CAM.eq20}. Otherwise, choose $w', v' \in \LL_{\d(x, y)}^g(x, y)$ so that $w \in T_{w'}$ and $v \in T_{v'}$. Since $y \in U_M^g(x, \delta(x, y))$ by Porpositon~\ref{COM.prop10}, there exists a chain $(w(1), \ldots, w(M + 3))$ such that $w(1) = w', w(M + 3) = v'$ and $w(j) \in \LL_{\d(x, y)}^g$ for any $j = 1, \ldots, M + 3$. By \eqref{CAM.eq20},  it follows that
\[
|w| - N' \le |w'| \le |v'| + N_0(M + 2) \le |v| + N_0(M + 2)
\]
Hence letting $N_1 = N_0(M + 2) + N'$, we obtain the desired claim.
\enddemo

\noindent {\bf Claim 2}: For any $N \in \BbN$, if $w \in \LL_s^g$, $v \in \LL_{\eta^{N + 1}s}^g$ and $v \in T_w$, then $|v| \ge |w| + N$.
\demo[Proof of Claim 2]
By \eqref{CAM.eq30} and \eqref{CAM.eq40}, it follows that
\[
\eta^{|v| - |w|} g(w) \le g(v) \le \eta^{N + 1}s \le \eta^Ng(w).
\]
Hence $|v| - |w| \ge N$.
\enddemo

\noindent
{\bf Claim 3}: Set $N_2 = N_1 + 2$ and $r_* = \eta^{N_2}\kappa^{M}$. For any $x, y \in X$, there exists $l_* = l_*(x, y) \in \BbN$ such that  $|w| \ge l_*$ for any $w \in \LL_{r_*\d(x, y)}(x, y)$ and if $v \in (T)_{l_*}(x, y)$, then $v \in \LL_{s'}^g$ for some $s' < \kappa^M\delta(x, y)$.
\demo[Proof of Claim 3]
Set $s = \kappa^M\d(x, y)$. Let
\[
N_3 = \min_{u \in \LL_{s}^g(x, y)} |u|.
\]
Set $l_* = l_*(x, y) = N_3 + N_1 + 1$. Note that $r_*\d(x, y) = \eta^{N_2}s < s$. For any $w \in \LL_{r_*\d(x, y)}^g(x, y)$, choose $w_* \in \LL_s^g(x, y)$ so that $w \in T_{w_*}$. By Claim 1 and 2, 
\[
|w| \ge |w_*| + N_1 + 1 \ge N_3 + N_1 + 1 = l_* > N_3 + N_1 \ge |u|.
\]
for any $u \in \LL_s^g(x, y)$. Let $v \in (T)_{l_*}(x, y)$. There exists $v' \in \LL_{s}^g(x, y)$ such that $v \in T_{v'}$. Since $|v'|  < l_*$, there exists $s' < s$ such that $v \in \LL_{s'}^g$.
\enddemo

Let $\bp = (w(1), \ldots, w(m))$ be a chain of $K$. Assume that $x \in K_{w(1)}$ and $y \in K_{w(m)}$. If $g(w(1)) \ge r_*\d(x, y)$ or $g(w(m)) \ge r_*\d(x, y)$, then 
\begin{equation}\label{CAM.eq530}
L_g(\bp) \ge r_*\d(x, y).
\end{equation}
Assume that $g(w(1)) < r_*\d(x, y)$ and $g(w(m)) < r_*\d(x, y)$. Set $l_* = l_*(x, y)$. Then by Claim 3, there exist $w, v \in (T)_{l_*}$ such that $w(1) \in T_w$, $w(m) \in T_v$, $w \in \LL_{s_1}^g$ for some $s_1 < \kappa^M\d(x, y)$ and $v \in \LL_{s_2}^g$ for some $s_2 < \kappa^M\d(x, y)$. Note that $x \in K_w$ and $y \in K_v$. Suppose that $v \in \GG_M(w)$. Then there exists a horizontal $M$-chain $(u(1), \ldots, u(M + 1))$ such that $u(1) = w$ and $u(M + 1) = v$. Let $s_* = \max\{g(u(i))| i = 1, \ldots, M + 1\}$. Then $(\pi^{m_1}(u(1)), \ldots, \pi^{m_{M + 1}}(u(M + 1)))$ is a $M$-chain between $x$ and $y$ in $\LL_{s_*}^g$ for some $m_1, \ldots, m_{M + 1} \ge 0$. Therefore, $\d(x, y) \le s_*$. On the other hand, since $\max\{s_1, s_2\} < \kappa^M\d(x, y)$, by \eqref{CAM.eq10} we see that
\[
g(u(i)) \le \kappa^{-\max\{i - 1, M - i + 1\}}\max\{s_1, s_2\} < \d(x, y)
\]
for any $i = 1, \ldots, M + 1$. This implies that $s_* < \d(x, y)$. Thus it follows that $v \notin \GG_M(w)$. By Lemma~\ref{CAM.lemma35}, there exists a $1$-path $\bp_1 = (v(1), \ldots, v(k))$ such that $v(1) = w(1)$, $v(k) = w(m)$ and
\begin{equation}\label{CAM.eq540}
 L_g(\bp) \ge c_0\ell_1^g(\bp_1),
\end{equation}
where $c_0$ is independent of $\bp$. Note that a $1$-path of $K$ is an $M$-jpath of $K$ for any $M \ge 1$ and $\ell_M^g(\bp_1) = \ell_1^g(\bp_1)$. Since $\vp \bl g$, we have
\begin{equation}\label{CAM.eq550}
\ell_M^g(\bp_1) \ge c_1\ell_M^{\vp}(\bp_1),
\end{equation}
where $c_1 > 0$ is independent of $\bp$. Applying Lemma~\ref{CAM.lemma30}, we obtain
\begin{equation}\label{CAM.eq560}
\ell_M^{\vp}(\bp_1) \ge \kappa_0\max\{\vp(w), \vp(v)\} \ge c_2\kappa_0\max\{g(w), g(v)\},
\end{equation}
where $c_2$ is independent of $\bp$. By Claim 3, there exist $w' \in T_w$ and $v' \in T_v$ such that $w', v' \in \LL_{r_*\d(x, y)}^g$. Hence by \eqref{CAM.eq40}, we have $\eta{r_*}\d(x, y) < g(w') \le g(w)$ and $\eta{r_*}\d(x, y) < g(v') \le g(v)$. So by \eqref{CAM.eq560},
\begin{equation}\label{CAM.eq570}
\ell_M^{\vp}(\bp_1) \ge c_3\d(x, y),
\end{equation}
where $c_3$ is independent of $\bp$. Finally combining \eqref{CAM.eq530}, \eqref{CAM.eq540}, \eqref{CAM.eq550} and \eqref{CAM.eq570}, we conclude that there exists $c_4 > 0$ such that if $\bp = (w(1), \ldots, w(m))$ is a chain of $K$, $x \in K_{w(1)}$ and $y \in K_{w(m)}$, then
\[
L_g(\bp) \ge c_4\d(x, y).
\]
This immediately implies
\[
c_4\d_M^g(x, y) \le D^g(x, y) \le D_M^g(x, y) \le (M + 1)\d_M^g(x, y)
\]
for any $x, y \in X$. Thus, the metric $D^g$ is $M$-adapted to $g$.
\enddemo

\setcounter{equation}{0}
\section{Construction of Ahlfors regular metric I}\label{EAR}

In this section, we discuss a condition for a weight function to induce an Ahlfors regular metric, whose definition is given in Definition~\ref{EAR.def10}. \par
As in the last section, $(X, \O)$ is a compact metrizable topological space with no isolated point, $(T, \A, \phi)$ is a locally finite tree and $K: T \to \C(X, \O)$ is a minimal partition. Furthermore, we assume that the partition $K: T \to \C(X, \O)$ is strongly finite, i.e. $\sup_{w \in T}\#(S(w)) < +\infty$ throughout this section. 

\definition\label{EAR.def10}
A metric $d \in \D(X, \O)$ is called Ahlfors regular if there exists a Borel regular probability measure $\mu$ on $X$ which is $\a$-Ahlfors regular with respect to $d$ for some $\a$.
\enddefinition

\thm\label{thm10}
Let $d \in \D_{A, e}(X, \O)$ and assume that $d \gen h_*$, $d$ is thick and uniformly finite. Let $\a > 0$. There exist a metric $\rho \in \D(X, \O)$ and a measure $\mu \in \M_P(X, \O)$ such that $\rho \qs d$ and $\mu$ is $\a$-Ahlfors regular with respect to $\rho$ if and only if there exists $g \in \G_{e}(X)$ such that
\begin{itemize}
\item
$g \gen g_d$,
\item
there exist $c > 0$ and $M \ge 1$ such that
\begin{equation}\label{EAR.eq05}
cD_M^g(x, y) \le D^g(x, y)
\end{equation}
for any $x, y \in X$,
\item
there exists $c > 0$ such that 
\begin{equation}\label{EAR.eq10}
c^{-1}g(w)^{\a} \le \sum_{v \in S^n(w)} g(v)^{\a} \le cg(w)^{\a}
\end{equation}
 for any $w \in T$ and $n \ge 0$, where $S^n(w) = (T)_{|w| + n} \cap T_w$ by definition. (See Definition~\ref{GAE.def10} to recall the definiton of $S^n(w)$.)
\end{itemize}
\endthm
\remark
The condition~\eqref{EAR.eq05} is equivalent to the existence of a metric $\rho'$ which is $M$-adapted to $g$ for some $M \ge 1$.
\endremark
\demo
By Proposition~\ref{BLE.prop310}, $g$ is tight. Hence Theorem~\ref{GAE.thm20} shows that $h_*$ is tight, thick and uniformly finite.\par
Suppose that there exist a metric $\rho$  and a measure $\mu$ such that $\rho \qs d$ and $\mu$ is $\a$-Ahlfors regular with respect to $\rho$. By Theorem~\ref{QSY.thm10}, setting $g = g_{\rho}$, we see that $g \gen g_d$, $g$ is exponential and $\rho$ is adapted. Hence by Theorem~\ref{GAE.thm20}, $g$ is thick and uniformly finite. Using Proposition~\ref{PAS.prop40}, we verify \eqref{EAR.eq05}. Since $\mu$ has the volume doubling property with respect to $\rho$, Lemma~\ref{VDP.lemma50} implies that
\[
\sum_{v \in S^n(w)} c\mu(K_v) \le \sum_{v \in S^n(w)} \mu(O_v) \le \mu(K_w) \le \sum_{v \in S^n(w)} \mu(K_v).
\]
Furthermore, by Theorem~\ref{MEM.thm100}, it follows that $g^{\a} \bl g_{\mu}$.This yields \eqref{EAR.eq10}.\par
Conversely, assume that $g \in \G_e(X)$, $g \gen g_d$, \eqref{EAR.eq05} and \eqref{EAR.eq10}.  Since $g, g_d \in \G_e(T)$, $g \gen g_d$ and $g_d$ is tight and thick, Theorem~\ref{GAE.thm20} shows that $g$ is thick and tight. Define $\rho(x, y) = D^g(x, y)/\sup_{a, b \in X}D^g(a, b)$. By \eqref{EAR.eq05}, it follows that $\rho$ is adapted to $g$. Moreover, Corollary~\ref{BLE.cor20} implies that $g_{\rho} \bl g$. Therefore $g_{\rho} \gen g_d$ and hence by Theorem~\ref{QSY.thm10}, we see that $\rho \qs d$.
Choose $x_w \in O_w$ for each $w \in T$. Define 
\[
\mu_n = \frac{1}{\sum_{w \in (T)_n} g(w)^{\a}}\sum_{w \in (T)_n} g(w)^{\a}\d_{x_w},
\]
where $\d_x$ is Dirac's point mass at $x$. Note that \eqref{EAR.eq10} implies that
\begin{equation}\label{EAR.eq20}
c^{-1} = c^{-1}g(\phi) \le \sum_{w \in (T)_n}g(w)^{\a} \le cg(\phi) = c
\end{equation}
Since $(X, \O)$ is compact, there exist a sub-sequence $\{\mu_{n_i}\}_{i \ge 1}$ and a Borel regular probability measure $\mu$ on $X$ such that $\{\mu_{n_i}\}_{i \ge 1}$ converges weakly to $\mu$ as $i \to \infty$. For $w \in T$, let $A_w$ be the $\e$-neighborhood of $K_w$ with respect to $d$, i.e. $A_{w, \e} = \{y|y \in X, \inf_{x \in K_w} d(x, y) < \e\}$ and let $f_{w, \e}: X \to [0, 1]$ be a continuous function satisfying $f_{w, \e}|_{K_w} = 1$ and $f_{w, \e}|_{\sd X{A_{w, \e}}} = 0$. Set 
\[
U_w = K_w \bigcup \Bigg(\bigcup_{v \in \GG_1(w)} O_v\Bigg).
\]
Then $U_w$ is an open neighborhood of $K_w$. Therefore for sufficiently small $\e > 0$, $A_{w, \e} \subseteq U_w$. By this fact, we have
\[
\sum_{v \in S^n(w)} g(v)^{\a} \le \int_X f_{w, \e}d\mu_{|w| + n} \le \sum_{u \in \GG_1(w)}\sum_{v \in S^n(u)} g(v)^{\a}
\]
By \eqref{EAR.eq10},
\[
c^{-1}g(w)^{\a} \le \int_X f_{w, \e}d\mu_{|w| + n} \le c\sum_{u \in \GG_1(w)} g(u)^{\a}
\]
Since $h_*$ is uniformly finite, there exists $c_1 > 0$ which is independet of $w$ such that $\#(\GG_1(w)) \le c_1$. Moreover since $g \gen d \gen h_*$, using Proposition~\ref{CAM.prop10}, we have
\[
c^{-1}g(w)^{\a} \le \int_X f_{w, \e}d\mu_{|w| + n} \le cc_1\kappa^{-\a}g(w)^{\a}.
\]
Choosing the proper subsequence of $|w| + n$ and taking the limit, we obtain
\[
c^{-1}g(w)^{\a} \le \int_X f_{w, \e}d\mu \le cc_1\kappa^{-\a}g(w)^{\a}.
\]
Letting $\e \downarrow 0$, we see 
\[
c^{-1}g(w)^{\a} \le \mu(K_w) \le cc_1\kappa^{-\a}g(w)^{\a}.
\]
This implies that $\mu \in \M_P(X, \O)$ and $g^{\a} \bl g_{\mu}$. Moreover by Theorem~\ref{GAE.thm20}, $g_{\rho}$ is uniformly finite. Hence by Theorem~\ref{MEM.thm100}, $\mu$ is $\a$-Ahlfors regular with respect to $d$.
\enddemo

\setcounter{equation}{0}

\section{Basic framework}\label{BFM}

From this section, we start proceeding towards the characterization of Ahlfors regular conformal dimension. To begin with, we fix our framework in this section and keep it until the end. \par
As in the previous sections, $(T, \A, \phi)$ is a locally finite tree with the root $\phi$, $(X, \O)$ is a compact metrizable topological space with no isolated point, $K: T \to \C(X, \O)$ is a minimal partition. We also assume that $(T, \A, \phi)$ is strongly finite, i.e. $\sup_{w \in T} \#(S(w)) < +\infty$.\par
Our standing assumptions in the following sections are as follows:
\vspace{5pt}\par
\noindent {\bf Basic Framework}\,\,
Let $d \in \D(X, \O)$. For $r \in (0, 1)$, define $h_r: T \to (0, 1]$ by
\[
h_r(w) = r^{|w|}
\]
for any $w \in T$. We assume the following conditions (BF1) and (BF2) are satisfied: \\
(BF1)\,\,$d$ is $M_*$-adapted for some $M_* \ge 1$, exponential, thick, and uniformly finite.\\
(BF2)\,\,There exists $r \in (0, 1)$ such that $h_r \bl d$. \par
\vspace{5pt}

\remark
$h_r$ is an exponential weight function. 
\endremark

Our goal of the rest of this paper is to obtain characterizations of the Ahlfors regular conformal dimension of $(X, d)$ satisfying the above conditions (BF1) and (BF2).\par
The condition (BF2) may be too restrictive. Replacing the original tree $T$ by its subtree $\tT^{g_d, r}$ associated with $d$ defined in Definition~\ref{HYP.def20}, however, we can realize (BF2) providing (BF1) is satisfied.\par
To simplify the notation, we write $\tT^{d, r}$ in place of $\tT^{g_d, r}$ hereafter.

\prop\label{BFM.prop10}
Assume that $d$ is $M_*$-adapted, exponential, thick and uniformly finite. For any $r > 0$, if we replace $T$ and $K: T \to \C(X, \O)$ by $\wT^{d, r}$ and $K_{\tT^{d,r}}: \tT^{d, r} \to \C(X, \O)$ respectively, then (BF1) and (BF2) are satisfied. 
\endprop

\demo
Since we should handle two different structures $(T, K)$ and $(\tT^{d, r}, K_{\tT^{d, r}})$ here, we denote $\LL_s^g$ and $U_M^g(x, s)$ by $\LL_s^g(T, K)$ and $U_M^g(x, s; T, K)$ respectively to emphasize the dependency on a tree structure $T$ and a partition $K$. Then it follows that
\[
\LL_{r^n}^{d}(\tT^{d, r}, K_{\tT^{d, r}}) = \LL_{r^n}^{d}(T, K)
\]
for any $n \ge 0$ and
\[
U_M^d(x, r^n; T, K) = U_M^{\wth_r}(x, r^n;\tT^{d, r}, K_{\tT^{d, r}})
\]
Since $d$ is $M_*$-adapted, there exist $c_1, c_2 > 0$ such that
\begin{equation}\label{BF.eq10}
B_d(x, c_1s)  \subseteq U_{M_*}^d(x, s; T, K) \subseteq B_d(x, c_2s)
\end{equation}
for any $x \in X$ and $s \in (0, 1]$. Let $r^m  > s \ge r^{m + 1}$.  Then
\begin{multline*}
U_{M_*}^{\wth_r}(x, s; \tT^{d, r}, K_{\tT^{d, r}}) \subseteq U_{M_*}^{\wth_r}(x, r^m; \tT^{d, r}, K_{\tT^{d, r}})\\
 = U_{M_*}^d(x, r^m;T, K) \subseteq B_d(x, c_2r^m) \subseteq B_d(x, c_2r^{-1}s).
\end{multline*}
In the same manner, we also obtain
\[
B_d(x, c_1rs) \subseteq U^{\wth_r}_{M_*}(x, s;\tT^{d, r}, K_{\tT^{d, r}}).
\]
Thus $d$ is $M_*$ adapted with respect to $(\tT^{d, r}, K_{\tT^{d, r}})$. It is straightforward to see that $d$ is exponential, tight and uniformly finite with respect to $(\tT^{d, r}, K_{\tT^{d, r}})$. Since $d$ is exponential, there exists $c > 0$ such that
\[
cd(w) \ge s  \ge d(w)
\]
if $w \in \LL_s^d$. This implies that $\wth_r \bl d$ as a weight function of $\tT^{d, r}$. Thus (BF1) and (BF2) with respect to $(\tT^{d, r}, K_{\tT^{d, r}})$ is satisfied. 
\enddemo

Due to this proposition, if $d$ is $M_*$-adapted, exponential, thick and uniformly finite, then  we replace $(T, K), \pi$ and $S$ by $(\tT^{d, r}, K_{\tT^{d, r}}), \pi^{d, r}$ and $S^{d, r}$ respectively and assume that (BF1) and (BF2) are satisfied hereafter. For ease of notations, we use $(T, K), \pi$ and $S$ to denote $(\tT^{d, r}, K_{\tT^{d, r}}), \pi^{d, r}$ and $S^{d, r}$. \par
Note that even after the modification the condition $\sup_{w \in T} \#(S(w)) < +\infty$ still holds as $d$ is exponential. Moreover, since $d$ is uniformly finite, the following notion is well-defined.

\definition\label{BFM.def10}
Define
\[
L_* = \sup_{w \in T} \#(\GG_1(w))
\]
and
\[
 N_* = \sup_{w \in T} \#(S(w))
\]
\enddefinition

\notation
We write $U_M(x, s) = U^{h_r}_M(x, s)$ for any $M \ge 0$, $x \in X$ and $s \in (0, 1]$.
\endnotation

By the above definition, it is straightforward to deduce the next lemma.

\lemma\label{BF.lemma20}
{\rm (1)}\,\,For any $w \in T$ and $k \ge 1$,
\begin{equation}\label{BF.eq40}
\#(S^k(w)) \le (N_*)^k
\end{equation}
{\rm (2)}\,\,For any $w \in T$ and $M \ge 1$,
\begin{equation}\label{BF.eq50}
\#(\GG_M(w)) \le (L_*)^M
\end{equation}
\endlemma

We present two useful propositions. The first one is an observation on the geometry of $\GG_M(w)$'s which holds without (BF1) and (BF2).

\prop\label{BF.prop20}
 Let $M_1, M_2 \in \BbN$. Suppose that 
\begin{equation}\label{BF.eq60}
\pi(\GG_{M_1 + M_2}(v)) \subseteq \GG_{M_2}(\pi(v))
\end{equation}
for any $v \in T$. Then
\[
\pi(\GG_{M_1 + M_2}(v)) \subseteq \GG_{M_1 + M_2}(w)
\]
if $\pi(v) \in \GG_{M_1}(w)$.
\endprop

\demo
Assume that $\pi(v) \in \GG_{M_1}(w)$, i.e. there exists a horizontal $1$-jpath $(w(1), \ldots, w(M_1 + 1))$ such that $w(1) = w$ and $w(M_1 + 1) = \pi(v)$. By \eqref{BF.eq60}, we see $\pi(u) \in \GG_{M_2}(\pi(v))$ for any $u \in \GG_{M_1 + M_2}(v)$. Hence there exists  a horizontal $1$-path $(v(1), \ldots, v(M_2 + 1))$ such that $v(1) = \pi(v)$ and $v(M_2 + 1) = \pi(u)$. As $(w(1), \ldots, w(M_1), v(1), \ldots, v(M_2 + 1))$ is a horizontal $1$-path, we see that $\pi(u) \in \GG_{M_1 + M_2}(w)$.
\enddemo

The second observation requires (BF1) and (BF2).

\prop\label{BFM.prop40}
Assume {\rm (BF1)} and {\rm (BF2)}. For any $M \ge 1$, there exists $m_0 \in \BbN$ such that, for any $m \ge m_0$ and $w \in T$, $\GG_M(v) \subseteq S^m(w)$ for some $v \in S^m(w)$.
\endprop

To prove this proposition, we need the following fact.

\lemma\label{BFM.lemma10}
Suppose that a partition $K$ is minimal. Let $A \subseteq T$ and let $v \in T$. If $|v| \ge |w|$ for any $w \in A$ and $K_v \subseteq \cup_{w \in A}K_w$, then $v \in \cup_{w \in A} T_w$.
\endlemma

\demo
Set $A' = \cup_{w \in A} S^{|v| - |w|}(w)$. Then $A' \subseteq (T)_{|v|}$ and
\[
K_v \subseteq \bigcup_{w \in A} K_w = \bigcup_{u \in A'} K_u.
\]
Since $O_v \neq \emptyset$, we see that $v \in A' \subseteq  \cup_{w \in A} T_w$.
\enddemo 

\demo[Proof of Proposition~\ref{BFM.prop40}]
Since $d$ is thick, so does $h_r$. By Proposition~\ref{ADD.prop20}, there exists $\a \in (0, 1)$ such that for any $w \in T$, 
\[
U_{M}(x, {\a}r^{|w|}) \subseteq K_w
\]
for some $x \in K_w$. Set $m_0 = \min\{m| r^m < \a\}$ and let $m \ge m_0$. Then $U_{M}(x, r^{|w| + m}) \subseteq K_w$. Therefore, if $v \in S^m(w)$ and $x \in K_v$, then $U_{M}(v) \subseteq U_{M}(x, r^{|w| + m}) \subseteq K_w$. Since the partition is minimal, Lemma~\ref{BFM.lemma10} implies that $\GG_{M}(v) \subseteq S^m(w)$.
\enddemo

\setcounter{equation}{0}
\section{Construction of adapted metric II}\label{WFG}

In this section, we study a sufficient condition for the existence of an adapted metric to a given weigh function $d$ under the basic framework presented in Section~\ref{BFM}. This is the continuation of the study of Section~\ref{CAM}.\par
As in the previous sections, $(T, \A, \phi)$ is a locally finite tree with the root $\phi$, $(X, \O)$ is a compact metrizable topological space with no isolated point, $K: T \to \C(X, \O)$ is a minimal partition. We also assume that $\sup_{w \in T} \#(S(w)) < +\infty$. Moreover, we assume that $d \in \D_{A, e}(X, \O)$ and the basic framework given in Section~\ref{BFM}, i.e. the conditions (BF1) and (BF2) are satisfied.\par
In this section, we fix $g \in \G_e(T)$ satisfying $g \gen d$.

\remark
The modification of $T$ in the previous section does not affect the relation $\gen$ and the exponentiality of $g$.
\endremark

\definition~\label{WFG.def10}
Let $\vp : T \to (0, \infty)$. For $M \ge 1$ and $w \in T$, define 
\[
\inn{\vp}_M(w) = \min_{v \in \GG_M(w)} \vp(v).
\]
and
\[
\Pi_M^\vp(w) = \min_{v \in \GG_M(w)} \frac{\vp(v)}{\vp(\pi(v))}
\]
\enddefinition

The next lemma, which holds without (BF1) and (BF2), gives a sufficient condition for a non-negative function on $T$ to be balanced. 

\lemma\label{WFG.lemma10}
Let $\vp: T \to (0. \infty)$ and let $M_1, M_2 \in \BbN$. Suppose that \eqref{BF.eq60} holds for any $v \in T$. \\
(1)\,\,For any $v \in T$,
\[
\inn{\vp}_{M_1 + M_2}(v) \ge \Pi_{M_1 + M_2}^{\vp}(v)\max\{\inn{\vp}_{M_1 + M_2}(u)| u \in \GG_{M_1}(\pi(v))\}.
\]
(2)\,\,If
\begin{equation}\label{WFG.eq200}
 \sum_{i = 1}^m \Pi^{\vp}_{M_1 + M_2}(w(i)) \ge 1
\end{equation}
for any $w \in T$ and $(w(1), \ldots, w(m)) \in \C_w^{M_1}$, then $\inn{\vp}_{M_1 + M_2}$ is $M_1$-balanced. 
\endlemma

\demo
For simplicity, set $\vp_* = \inn{\vp}_{M_1 + M_2}$. \\
(1)\,\,There exists $v' \in \GG_{M_1 + M_2}(v)$ such that $\vp_*(v) = \vp(v')$. Let $u \in \GG_{M_1}(\pi(v))$. Since $\pi(v) \in \GG_{M_1}(u)$, Proposition~\ref{BF.prop20} shows that $\pi(\GG_{M_1 + M_2}(v)) \subseteq \Gamma_{M_1 + M_2}(u)$. Therefore, $\pi(v') \in \Gamma_{M_1 + M_*}(u)$. Hence
\[
\vp_*(v) = \vp(v') = \frac{\vp(v')}{\vp(\pi(v'))}\vp(\pi(v')) \ge \Pi_{M_1 + M_2}^{\vp}(v)\max\{\vp_*(u)| u \in \GG_{M_1}(\pi(v))\}
\]
(2)\,\,Let $(w(1), \ldots, w(m)) \in C_w^{M_1}$. Then by (1), 
\begin{equation}\label{WFG.eq100}
\sum_{i = 1}^{m} \vp_*(w(i)) \ge \sum_{i = 1}^{m}\Pi_{M_1 + M_2}^g(w(i))\max\{\vp_*(u)| u \in \GG_{M_1}(\pi(w(i))\}.
\end{equation}
If $\pi(w(m)) \in \GG_{M_1}(\pi(w(i))$ for any $i = 1, \ldots, m$, then \eqref{WFG.eq100} and \eqref{WFG.eq200} imply
\begin{equation}\label{WFG.eq110}
\sum_{i = 1}^{m} \vp_*(w(i)) \ge \sum_{i = 1}^{m}\Pi_{M_1 + M_2}^{\vp}(w(i))\vp_*(\pi(w(m))) \ge \vp_*(\pi(w(m))).
\end{equation}
Suppose that there exists $j \in \{1, \ldots, m - 1\}$ such that $\pi(w(m)) \notin \GG_{M_1}(\pi(w(j)))$ and $\pi(w(m)) \in \GG_{M_1}(\pi(w(i)))$ for any $i \in \{j + 1, \ldots, m - 1\}$. Note that $j \le m - 2$ because $w(m) \in \GG_{M_1}(w(m - 1))$. Then $(w(m - 1), w(m -2), \ldots, w(j + 1)) \in \C_{\pi(w(m))}^{M_1}$. Hence by \eqref{WFG.eq200}
\[
\sum_{k = j + 1}^{m - 1} \Pi_{M_1 + M_1}^{\vp}(w(k)) \ge 1.
\]
So using \eqref{WFG.eq100}, we see that
\[
\sum_{i = 1}^{m} \vp_*(w(i)) \ge \sum_{i = j + 1}^{m - 1}\Pi_{M_1 + M_*}^{\vp}(w(i))\max_{u \in \GG_{M_1}(\pi(w(i)))}\vp_*(u) \ge \vp_*(\pi(w(m))).
\]
Hence we obtain \eqref{WFG.eq110} in this case as well. Thus we have shown that $\vp_*$ is $M_1$-balanced.
\enddemo

Using the last lemma, we have a sufficient condition for the existence of a metric which is adapted to a given weight function.

\thm\label{WFG.thm10}
Let $M_1, M_2 \in \BbN$. Suppose that \eqref{BF.eq60} holds for any $v \in T$. Assume that $g \in \G_e(T)$ and $g \gen d$. If 
\begin{equation}\label{WFG.eq201}
 \sum_{i = 1}^m \Pi^g_{M_1 + M_2}(w(i)) \ge 1
\end{equation}
for any $w \in T$ and $(w(1), \ldots, w(m)) \in \C_w^{M_1}$, then there exists a metric $\rho \in \D_{\A, e}(X)$ which is $M_1$-adapted to $g$ and quasisymmetric to $d$.
\endthm

\demo
By Lemma~\ref{WFG.lemma10}, $\inn{g}_{M_1 + M_2}$ is $M_1$-balanced. Since $g \gen d$, we have $g \bl \inn{g}_{M_1 + M_2}$. Therefore Theorem~\ref{CAM.thm10} shows that there exists a metric $\rho \in \D(X)$ such that $\rho$ is $M_1$-adapted to $g$. Since $g \bl \rho$, we see that $\rho$ is exponential. Moreover, the fact that $\rho \gen d$ implies $\rho \qs d$ by Theorem~\ref{QSY.thm10}.
\enddemo

To utilize Theorem~\ref{WFG.thm10}, we need to find $M_1$ and $M_2$ satisfying \eqref{BF.eq60};
\[
\pi(\GG_{M_1 + M_2}(v)) \subseteq \GG_{M_2}(\pi(v))
\]
for any $v \in T$. Since $d$ is $M_*$-adapted and the metric $\rho$ obtained in the above theorem is quaisymmetric to $d$, Theorem~\ref{QSY.thm10} implies that $\rho$ is $M_*$-adapted, and conversely, $d$ is $M_1$-adapted as well. If $M_* = \min\{M| \text{$d$ is $M$-adapted.}\}$, then it follows $M_1 \ge M_*$. This requirement on $M_1$ can make it hard to find $M_1$ and $M_2$. Replacing $\pi$ by $\pi^{k}$, however, we have the following fact.

\prop\label{WFG.prop10}
Let $M_1 \in \BbN$. There exists $k_{M_1} \ge 1$ such that if $v \in T$ and $k \ge k_{M_1}$, then $\pi^k(\GG_{M_1 + M_*}(v)) \subseteq \GG_{M_*}(\pi^k(v))$.
\endprop

\demo
If $|k| \le k$, then it is obvious. So, let $|v| \ge k$. Since $d$ is $M_*$-adapted, it is $(M_* + M_1)$-adapted as well. Therefore, there exist $c_1, c_2 > 0$ such that
\[
U_{M_* + M_1}^{h_r}(x, r^n) \subseteq B_d(x, c_1r^n)
\quad\text{and}\quad
B_d(x, c_2r^n) \subseteq U_{M_*}^{h_r}(x, r^n)
\]
for any $x \in X$ and $n \ge 0$. Choose $k \in \BbN$ so that $c_1r^k < c_2$. Then 
\[
U_{M_* + M_1}^{h_r}(x, r^{|v|}) \subseteq B_d(x, c_1r^{|v|}) \subseteq B_d(x, c_2r^{|v| - k}) \subseteq U_{M_*}^{h_r}(x, r^{|v| - k}).
\]
If $x \in O_v$, then this implies that $K_u \subseteq \cup_{w \in \GG_{M_*}(\pi^k(v))} K_w$ for any $u \in \GG_{M_* + M_1}(v)$.  By Lemma~\ref{BFM.lemma10}, it follows that $\pi^k(u) \in \GG_{M_*}(\pi^k(v))$.
\enddemo

\definition\label{CAM.def40}
Let $q \in \BbN$. Define $T^{(q)} = \cup_{m \ge 0} (T)_{mq}$. Define $\pi_q: T^{(q)} \to T^{(q)}$ by $\pi_q = \pi^q$, which is the $q$-th iteration of $\pi$. We consider $T^{(q)}$ as a tree with the reference point $\phi$ under the natural tree structure inherited from $T$. Then $(T^{(q)})_m = (T)_{mq}$. Moreover, set $K^{(q)} = K|_{T^{(q)}}$. 
\enddefinition

Note that a horizontal $M$-jpath of $K^{(q)}$ is a horizontal $M$-jpath of $K$. Similarly, $\GG_M(w, K^{(q)}) = \GG_M(w, K)$ and $U_M(w, K^{(q)}) = U_M(w, K)$ for any $w \in T^{(q)}$ and $J_{M, mq}^h(K) = J_{M, m}^h(K^{(q)})$. \par

\definition\label{WFG.def100}
(1)\,\,For $w \in T$, define
\begin{multline*}
\C_{w, k}(N_1, N_2, N) = \{(w(1), \ldots, w(m))|\text{$(w(1), \ldots, w(m))$ is}\\\text{ a horizontal $N$-jpath,} \,\,w(j) \in S^k(\GG_{N_2}(w))\,\,\text{for any $j = 1, \ldots, m$}, \\
\GG_N(w(1)) \cap S^k(\GG_{N_1}(w)) \neq \emptyset\,\,\text{and}\,\,\sd{\GG_N(w(m))}{S^k(\GG_{N_2}(w))} \neq \emptyset.\}
\end{multline*}
(2)\,\,
For a weight function $g \in \G(T^{(k)})$ on $T^{(k)}$, define
\[
 \Pi_M^{g, k}(w) = \min_{v \in \GG_M(w)} \frac{g(v)}{g(\pi^k(v))}
\]
for any $w \in T^{(k)}$.

\enddefinition

Note that $\C_w^M = \C_{w, 1}(0, M, M)$.\\
Replacing $T$ by $T^{(k)}$ and applying Theorem~\ref{WFG.thm10}, we obtain the following corollary.

\cor\label{WFG.cor10}
Let $M_1 \in \BbN$ and let $k \ge k_{M_1}$, where $k_{M_1}$ is the constant appearing in Proposition~\ref{WFG.prop10}. Assume that $g \in \G_e(T^{(k)})$ and $g \gen d$ as weight functions on $T^{(k)}$. If
\[
\sum_{i = 1}^m \Pi^{g, k}_{M_1 + M_*}(w(i)) \ge 1
\]
for any $w \in T^{(k)}$ and $(w(1), \ldots, w(m)) \in \C_{w, k}(0, M_1, M_1)$, then there exists a metric $\rho \in \D_{\A, e}(X)$ such that $\rho$ is $M_*$-adapted to $g$ and $\rho$ is quasisymmetric to $d$.
\endcor

\remark
It is easy to see that any weight function $g \in \G_e(T^{(k)})$ can be extended to a weight function $\widetilde{g} \in \G_e(T)$, i.e. there exists a weight function $\widetilde{g} \in \G_e(T)$ such that $\widetilde{g}|_{T^{(k)}} = g$. Since $\widetilde{g}$ is eponential, we can see that for any $M \ge 1$, there exist $c_1, c_2 > 0$ such that
\[
c_1D^g_M(x, y) \le D^{\widetilde{g}}_M(x, y) \le c_2D^g_M(x, y)
\]
for any $x, y \in X$. Therefore, the metric $\rho$ obtained in the above corollary is $M_*$-adapted to $\widetilde{g}$ as well.
\endremark

\section{Construction of Ahlfors regular metric II}\label{CAR}

In this section, making use of Theorem~\ref{thm10} and Corollary~\ref{WFG.cor10}, we are going to establish a sufficient condition for the existence of an adapted metric $\rho$ and a measure $\mu$ where $\mu$ is Ahlfors regular with respect to the metric $\rho$.\par
As in the previous sections, $(T, \A, \phi)$ is a locally finite tree with the root $\phi$, $(X, \O)$ is a compact metrizable topological space with no isolated point, $K: T \to \C(X, \O)$ is a minimal partition. We also assume that $\sup_{w \in T} \#(S(w)) < +\infty$.  Furthermore, we continue to employ the basic framework (BF1) and (BF2) in Section~\ref{BFM}.\par
Our main theorem of this section is as follows:

\thm\label{CAR.thm10}
Let $M_1 \in \BbN$. Assume that $k \ge \max\{m_0, k_{M_1}, k_{M_*}\}$, where $m_0$ is the constant appearing in Proposition~\ref{BFM.prop40} with $M = 1$ and $k_{M_1}$ and $k_{M_*}$ are the constants appearing in Proposition~\ref{WFG.prop10}. If there exists $\vp: T^{(k)} \to (0, 1]$ such that 
\[
\sum_{i = 1}^m \vp(w(i)) \ge 1
\]
for any $w \in T^{(k)}$ and $(w(1), \ldots, w(m)) \in \C_{w, k}(0, M_1, M_1)$ and
\[
\sum_{v \in S^k(w)} \vp(v)^p < \frac 12(L_*)^{-2(M_1 + 2M_*)},
\]
for any $w \in T^{(k)}$, then there exist a metric $\rho \in \D_{\A, e}(X)$ and a Borel regular probability measure $\mu$ on $X$ such that $\rho \qs d$ and $\mu$ is Ahlfors $p$-regular with respect to the metric $\rho$.
\endthm

\remark
By the choice of $k$ in the above theorem, it follows that for any $v \in T$,
\[
\pi^k(\GG_{M_1 + M_*}(v)) \cup \pi^k(\GG_{2M_*}(v)) \subseteq \GG_{M_*}(\pi^k(v)).
\]
\endremark

The rest of this section is devoted to a proof of this theorem. First we present two key lemmas.

\lemma\label{MCA.lemma20}
Let $(V, E)$ be a non-directed graph. Assume that $(v, v) \in E$ for any $v \in V$. For $m \ge 1$ and $A \subseteq V$, define 
\begin{multline*}
V_m(A) = \{u|\text{there exists $(x(1), \ldots, x(m + 1))$ such that}\\
\text{ $x(1) \in A$, $x(k) = u$ and $(x(i), x(i + 1)) \in E$ for any $i = 1, \ldots, m$}\}
\end{multline*}
Write $V_m(x) = V_m(\{x\})$ for $x \in V$. For any $f: V \to [0, \infty)$, there exists $\s: V \to [0, \infty)$ such that 
\[
f(v) \le \min\{\s(u)| u \in V_m(v)\} \le \s(v) \le \max_{u \in V_m(v)} f(u)
\]
 for any $v \in V$ and
\begin{equation}\label{MCA.eq150}
\sum_{v \in U} \s(v)^p \le \bigg(\max_{v \in V} \#(V_m(v))\bigg)\sum_{v \in V_m(U)} f(v)^p.
\end{equation}
for any $U \subseteq V$.
\endlemma

\demo
Define $\s(v) = \max\{f(u)| u \in V_m(v)\}$. Since $v \in V_m(u)$ if and only if $u \in V_m(v)$, it follows that $f(v) \le \s(u)$ for any $u \in V_m(v)$. Hence $f(v) \le \min\{\s(u)| u \in V_m(v)\}$. Moreover
\begin{multline*}
\sum_{v \in U} \s(v)^p \le \sum_{v \in U} \sum_{u \in V_m(v)} f(u)^p\\ = \sum_{u \in V_m(U)}\sum_{v \in V_m(u)}f(u)^p = \sum_{u \in V_m(U)} \#(V_m(u))f(u)^p.
\end{multline*}
Hence \eqref{MCA.eq150} holds.
\enddemo

\lemma\label{CON.lemma10}
Let $k \ge k_{M_*}$. Let $\kappa_0 \in (0, 1)$ and let $f: \sd{T^{(k)}}{\{\phi\}} \to [\kappa_0, 1)$. 
Then there exists $g: T^{(k)} \to (0, 1]$ such that
\begin{equation}\label{CON.eq10}
g(u) \ge \kappa_0g(v)
\end{equation}
if $(u, v) \in J^h_{M_*}$,
\begin{equation}\label{CON.eq20}
f(u) \le \frac{g(u)}{g(\pi^k(u))} \le \max_{v \in \GG_{M_*}(u)} f(v)
\end{equation}
for any $u \in \sd{T^{(k)}}{\{\phi\}}$ and
\begin{equation}\label{CON.eq30}
\sum_{v \in S^k(w)} \bigg(\frac{g(v)}{g(\pi^k(v))}\bigg)^p \le (L_*)^{2M_*}\sup\bigg\{\sum_{u \in S^k(w')} f(u)^p\bigg| w' \in \GG_{M_*}(w)\bigg\}.
\end{equation}
for any $p > 0$ and $w \in T^{(k)}$.
\endlemma

\demo
First we are going to construct $g: \cup_{n \ge 0} (T)_{kn} \to (0, 1]$ satisfying \eqref{CON.eq10} and \eqref{CON.eq20} inductively. Set $g(\phi) = 1$ and $g(w) = f(w)$ for any $w \in (T)_k$. Then \eqref{CON.eq10} and \eqref{CON.eq20} are satisfied for any $w \in (T)_k$. Assume that there exists $g: \cup_{n = 0}^m (T)_{kn} \to (0, 1]$ satisfying \eqref{CON.eq10} and \eqref{CON.eq20} up to the $m$-th level. Define 
\begin{align*}
g_1(v) &= f(v)g(\pi^k(v)),\\
g_2(v) &= \kappa_0\max_{u \in \GG_{M_*}(v)} g_1(u),\\
g(v) &= \max\{g_1(v), g_2(v)\}.
\end{align*}
for any $v \in (T)_{k(m + 1)}$. First we are going to show \eqref{CON.eq20} for $u \in (T)_{k(m + 1)}$. If $g(v) = g_1(v)$, then $f(v) = \frac{g(v)}{g(\pi^k(v))} \le \max_{u \in \GG_{M_*}(v)}f(u)$. Hence we have \eqref{CON.eq20}. If $g(v) = g_2(v)$, then $g_1(v) \le g(v)$ and hence $f(v) \le \frac{g(v)}{g(\pi^k(v))}$. There exists $u \in \GG_{M_*}(v)$ such hat $g(v) = \kappa_0f(u)g(\pi^k(u))$. Since $(\pi^k(u), \pi^k(v)) \in J^h_{M_*}$, we have $g(\pi^k(u)) \le (\kappa_0)^{-1}g(\pi^k(v))$ by \eqref{CON.eq10} for $(T)_{km}$. Therefore,
\[
\frac{g(v)}{g(\pi^k(v))} \le f(u) \le \max_{w \in \GG_{M_*}(v)}f(w).
\]
Thus we have shown \eqref{CON.eq20} for $(T)_{k(m + 1)}$.\par
Next to show \eqref{CON.eq10} for $(T)_{k(m + 1)}$, we need the following fact.\\
{\bf Claim}: If $(u, v), (v, v') \in J^h_{M_*, k(m + 1)}$, then $g_1(u) \ge (\kappa_0)^2g_1(v')$. \\
Proof of Claim:\,\,Note that $\pi^k(\GG_{2M_*}(u)) \subseteq \GG_{M_*}(\pi^k(u))$ because $k \ge k_{M_*}$. Hence $\pi^k(v') \in \GG_{M_*}(\pi^k(u))$ and so we have $(\pi^k(v'), \pi^k(u)) \in J^h_{M_*, km}$. Therefore, by \eqref{CON.eq10}, if $g_1(u) < (\kappa_0)^2g_1(v')$, then
\[
(\kappa_0)^2g_1(v') > g_1(u) = f(u)g(\pi^k(u)) \ge \kappa_0g(\pi^k(u)) \ge (\kappa_0)^2g(\pi^k(v')).
\]
Therefore $g_1(v') = f(v')g(\pi^k(v')) > g(\pi^k(v'))$. This contradicts the fact that $f(v') \le 1$ and hence we have confirmed the claim.\qed\\
Now let $(u, v) \in J^h_{M_*, k(m + 1)}$. If $g(v) = g_2(v)$, then there exists $v' \in \GG_{M_*}(v)$ such that $g(v) = \kappa_0g_1(v')$. By Claim, it follows that
\[
g(u) \ge g_1(u) \ge (\kappa_0)^2g_1(v') = \kappa_0g(v).
\]
If $g(v) = g_1(v)$, then $g(u) \ge g_2(u) \ge \kappa_0g_1(v) = \kappa_0g(v)$. Thus \eqref{CON.eq10} holds for $(T)_{k(m + 1)}$.\par
Using this construction of $g$ inductively, we obtain $g: T \to (0, 1]$ satisfying \eqref{CON.eq10} and \eqref{CON.eq20} at every level. Next we are going to proof \eqref{CON.eq30}. Note that $\underset{v \in S^k(w)}{\bigcup} \GG_{M_*}(v) \subseteq \underset{w' \in \GG_{M_*}(w)}{\bigcup} S^k(w')$. By \eqref{CON.eq20},
\begin{multline*}
\sum_{v \in S^k(w)} \bigg(\frac{g(v)}{g(\pi^k(v))}\bigg)^p \le \sum_{v \in S^k(w)}\sum_{u \in \GG_{M_*}(v)} f(u)^p\\
= \sum_{u \in \underset{v \in S^k(w)}{\bigcup} \GG_{M_*}(v)} \#(\GG_{M_*}(u) \cap S^k(w))f(u)^p \le (L_*)^{M_*}\sum_{u \in \underset{v \in S^k(w)}{\bigcup} \GG_{M_*}(v)} f(u)^p\\
 \le (L_*)^{M_*}\sum_{w' \in \GG_{M_*}(w)}\sum_{u \in S^k(w')} f(u)^p\\
 \le (L_*)^{2M_*}\sup\Big\{\sum_{u \in S^k(w')} f(u)^p\Big| w' \in \GG_{M_*}(w)\Big\}
\end{multline*}
\enddemo

\demo[Proof of Theorem~\ref{CAR.thm10}]
Set $\eta = \frac 12(L_*)^{-2(2M_* + M_1)}$. Assume that $\vp: T^{(k)} \to [0, 1]$ satisfies
\begin{equation}\label{PEN.eq100}
\sum_{i = 1}^m \vp(w(i)) \ge 1
\end{equation}
for any $w \in T$ and $(w(1), \ldots, w(m)) \in \C_{w, k}(0, M_1, M_1)$ and
\[
\sum_{v \in S^k(w)} \vp(v)^p < \eta
\]
for any $v \in T^{(k)}$.
Define
\[
\tvp(v) = \bigg(\vp(v)^p + \frac{\eta}{(N_*)^k}\bigg)^{1/p}
\]
for any $v \in T$. Then \eqref{PEN.eq100} still holds if we replace $\vp$ by $\tvp$. Moreover,
\[
\sum_{v \in S^k(w)} \tvp(v)^p < 2\eta
\]
and
\[
\bigg(\frac{\eta}{(N_*)^k}\bigg)^{1/p} \le \tvp(v) \le \bigg(\frac{1 + (N_*)^k}{(N_*)^k}\eta\bigg)^{1/p}.
\]
Set $M_3 = M_1 + M_*$. Letting $V = (T)_{km}$ and $f = \tvp$ and applying Lemma~\ref{MCA.lemma20}, we obtain $\p: T \to [0, \infty)$ satisfying
\[
\tvp(v) \le \inn{\p}_{M_3}(v) \le \p(v) \le \bigg(\frac{1 + (N_*)^k}{(N_*)^k}\eta\bigg)^{1/p}
\]
for any $v \in T$ and
\begin{multline*}
\sum_{v \in S^k(w)} \p(v)^p \le (L_*)^{M_3}\sum_{u \in \underset{v \in S^k(w)}{\bigcup} \GG_{M_3}(v)} \tvp(u)^p\\ \le (L_*)^{M_3} \sum_{w' \in \GG_{M_3}(w)}\sum_{u \in S^k(w')} \tvp(u)^p < 2(L_*)^{2M_3}\eta.
\end{multline*}
Next step is to use Lemma~\ref{CON.lemma10}. Set $\kappa_0 = \big(\frac{\eta}{(N_*)^k}\big)^{1/p}$ and $\kappa_1 = \big(\frac{1 + (N_*)^k}{(N_*)^k}\eta\big)^{1/p}$. Note that since $\kappa_1 < 1$ because $\eta < \frac 12$. Applying Lemma~\ref{CON.lemma10} with $f = \p$, we obtain $g: T^{(k)} \to (0, 1]$ satisfying \eqref{CON.eq10}, \eqref{CON.eq20} and \eqref{CON.eq30}. Define $\tau(w) = g(w)/g(\pi^k(w))$ for any $w \in \sd{T^{(k)}}{\{\phi\}}$. Then by \eqref{CON.eq20}, for any $w \in \sd{T^{(k)}}{\{\phi\}}$,
\[
\kappa_0 \le \p(w) \le \tau(w) \le \kappa_1
\]
and by \eqref{CON.eq30}
\[
\sum_{v \in S^k(w)} \tau(v)^p < 1
\]
for any $w \in T^{(k)}$. To construct desired weight function, we need to modify $\tau$ once more. Since $k \ge m_0$, Proposition~\ref{BFM.prop40} shows that, for any $w \in T^{(k)}$, there exists $v_w \in S^k(w)$ such that $\GG_1(v_w) \subseteq S^k(w)$. Define
\[
\s(v) = \begin{cases}
\tau(v)\quad&\text{if $v \neq v_{\pi^k(v)}$},\\
\displaystyle\bigg(1 - \sum_{u \in \sd{S^k(w)}{\{v\}}} \tau(u)^p\bigg)^{1/p}\quad&\text{if $v = v_{\pi^k(v)}$.}
\end{cases}
\]
Then
\begin{equation}\label{CAR.eq10}
\kappa_0 \le \tau(v) \le \s(v) \le \max\{\kappa_1, \big(1 - (\kappa_0)^p\big)^{1/p}\} < 1
\end{equation}
and 
\begin{equation}\label{CAR.eq20}
\sum_{v \in S^k(w)} \s(v)^p = 1
\end{equation}
for any $w \in T^{(k)}$. Since
\[
\tvp(v) \le \inn{\p}_{M_3}(v)\quad\text{and}\quad\p(v) \le \tau(v) \le \s(v),
\]
 it follows that $\tvp(v) \le \inn{\s}_{M_3}(v)$. Hence
\begin{equation}\label{CAR.eq30}
\sum_{i = 1}^m \inn{\s}_{M_3}(w(i)) \ge 1
\end{equation}
for any $w \in T^{(k)}$ and $(w(1), \ldots, w(m)) \in \C_{w, k}(0, M_1, M_1)$. Define $\wtg(w)$ inductively by $\wtg(\phi) = 1$ and
\[
\wtg(w) = \s(w)\wtg(\pi^k(w)).
\]
Suppose that $u, v \in (T)_{kn}$, $u \neq v$ and $u \in \GG_1(v)$. Then $\pi^{kl}(u) \neq \pi^{kl}(v)$ and $\pi^{k(l + 1)}(u) = \pi^{k(l + 1)}$ for some $l \ge 0$. Note that $\pi^{kj}(u) \in T_{\pi^{kl}(u)}$, $\pi^{kj}(v) \in T_{\pi^{kl}(v)}$ and $\pi^{kj}(u) \in \LL_1(\pi^{kl}(v))$ for any $j = 0, \ldots, l - 1$. Hence we see that
\[
\frac{\wtg(u)}{\wtg(v)} = \frac{\s(u)\s(\pi^k(u))\cdots\s(\pi^{kl}(u))}{\s(v)\s(\pi^k(v))\cdots\s(\pi^{kl}(v))} = \frac{\tau(u)\cdots\tau(\pi^{k(l - 1)}(u))\s(\pi^{kl}(u))}{\tau(v)\cdots\tau(\pi^{k(l - 1)}(v))\s(\pi^{kl}(v))}
\]
On the other hand,
\[
\frac{g(u)}{g(v)} = \frac{\tau(u)\cdots\tau(\pi^{k(l - 1)}(u))\tau(\pi^{kl}(u))}{\tau(v)\cdots\tau(\pi^{k(l - 1)}(v))\tau(\pi^{kl}(v))}.
\]
Thus if $\kappa_2 = \max\{\kappa_1, \big(1 - (\kappa_0)^p\big)^{1/p}\}$, then
\begin{equation}\label{CAR.eq40}
\frac{\wtg(u)}{\wtg(v)} = \frac{g(u)}{g(v)}\frac{\s(\pi^{kl}(u))}{\s(\pi^{kl}(v))}\frac{\tau(\pi^{kl}(v))}{\tau(\pi^{kl}(u))} \ge \kappa_0\frac{\kappa_0}{\kappa_2}\frac{\kappa_0}{\kappa_1}.
\end{equation}
By \eqref{CAR.eq10}, the weight function $\wtg$ on $T^{(k)}$ is exponential. By \eqref{CAR.eq40}, $\wtg \gen d$ as weight functions on $T^{(k)}$. By \eqref{CAR.eq30}, using Corollary~\ref{WFG.cor10}, we deduce that there exists a metric $\rho \in \D_{\A, e}(X)$ such that $\rho$ is $M_1$-adapted to $\wtg$, $\rho \bl \wtg$ and $\rho \qs d$. Moreover, by \eqref{CAR.eq20}, $\wtg$ satisfies the $T^{(k)}$-version of \eqref{EAR.eq10}. Hence applying Theorem~\ref{thm10} to $\wtg$ on $T^{(k)}$,  we verify the existence of a Borel regular probability measure $\mu$ on $X$ satisfying $\mu(K_w) = \wtg(w)^p$ for any $w \in T^{(k)}$. This implies $\mu \bl \wtg^p \bl \rho^p$. Hence $\mu$ is Ahlfors $p$-regular with respect to $\rho$.
\enddemo

\section{Critical index of $p$-energies and the Ahlfors regular conformal dimension}\label{PEN}
Finally in this section, we establish the characterization of the Ahlfors regular conformal dimension as the critical index of $p$-energies. \par
As in the previous sections, $(T, \A, \phi)$ is a locally finite tree with the root $\phi$, $(X, \O)$ is a compact metrizable topological space with no isolated point, $K: T \to \C(X, \O)$ is a minimal partition. We also assume that $\sup_{w \in T} \#(S(w)) < +\infty$. \par
Throughout this and the following sections, we fix $d \in \D_{A, \e}(X, \O)$ satisfying the basic framework, i.e. (BF1) and (BF2) in Section~\ref{BFM}.\par
First we recall the definition of Ahlfors regular conformal dimension.

\definition\label{PEN.def00}
Let $(X, d)$ be a metric space. The Ahlfors regular conformal dimension, AR conformal dimension for short, of a metric space $(X, d)$ is defined as 
\begin{multline*}
\dim_{AR}(X, d) =\\
 \inf\{\a| \text{there exist a metric $\rho$ on $X$ and a Borel regular measure $\mu$ on $X$}\\
\text{ such that $\rho \qs d$ and $\mu$ is $\a$-Ahlfors regular with respect to $\rho$}\}.
\end{multline*}
\enddefinition

 The definition of $p$-energy $\E_p(f|V, E)$ of a function $f$ on a graph $(V, E)$ is as follows.\par

\definition\label{PEN.def10}
Let $G = (V, E)$ be a (non-directed) graph. For $f: V \to \BbR$, define
\[
\E_p(f|V, E) = \frac 12\sum_{(x, y) \in E} |f(x) - f(y)|^p.
\]
If $E = \emptyset$, then we define $\E_p(f, V, E) = 0$ for any $f: V \to \BbR$.
Let $V_1, V_2 \subseteq V$. Assume that $V_1 \cap V_2 = \emptyset$. Define
\[
\F_F(V, E, V_1, V_2) = \{f| f: V \to [0, \infty), f|_{V_1} \ge 1, f|_{V_2} \equiv 0\}
\]
and
\[
\E_p(V, E, V_1, V_2) = \inf\{\E_p(f|V, E)|f \in \F_F(V, E, V_1, V_2)\}
\]
\enddefinition

For $p = 2$, on the analogy of electric circuits, the quantity $\E_2(V, E, V_1, V_2)$ is considered as the conductance (and its reciprocal is considered as the resistance) between $V_1$ and $V_2$. In the same way, we may regard $\E_p(V, E, V_1, V_2)$ as the $p$-conductance between $V_1$ and $V_2$. \par
Applying the above definition to the horizontal graphs $((T)_m, J_{N, m}^h)$, we define the critical index $I_{\E}(N_1, N_2, N)$ of $p$-energies.

\definition\label{PEN.def20}
Let $N_1, N_2$ and $N$ be integers satisfying $N_1 \ge 0$, $N_2 > N_1$ and $N \ge 1$. Define
\[
\E_{p, k}(N_1, N_2, N)  = \sup_{w \in T}\E_p((T)_{|w| + k}, J^h_{N, |w| + k}, S^k(\GG_{N_1}(w)), S^k(\GG_{N_2}(w))^c),
\]
\[
\overline{\E}_p(N_1, N_2, N) = \limsup_{k \to \infty} \E_{p, k}(N_1, N_2, N)
\]
and
\[
\underline{\E}_p(N_1, N_2, N) = \liminf_{k \to \infty} \E_{p, k}(N_1, N_2, N).
\]
Furthermore, define
\[
I_{\E}(N_1, N_2, N) = \inf\{p | \underline{\E}_p(N_1, N_2, N) = 0\}.
\]
\enddefinition

The last quantity $I_{\E}(N_1, N_2, N)$ is called the critical index of $p$-energies. Two values $\overline{\E}_p(N_1, N_2, N)$ and $\underline{\E}_p(N_1, N_2, N)$ represents the asymptotic behavior of the $p$-conductance between $\GG_{N_1}(w)$ and the complement of $\GG_{N_2}(w)$ as we refine the graphs between those two sets.

\thm\label{PEN.thm10}
For any $N \ge 1$,
\[
I_{\E}(N_1, N_2, N) = \dim_{AR}(X, d)
\]
if $N_1+ M_* \le N_2$.
\endthm

\remark
As is shown in Theorem~\ref{RCI.thm10}, even if we replace $\underline{\E}_p$ by $\overline{\E}_p$, the value of $I_{\E}$ is the same.
\endremark

Up to now we have considered the critical exponent for $p$-energies associated with simple graphs $\{((T)_m, J^h_{N, m})\}_{m \ge 0}$. In fact, the critical exponent is robust with respect to certain class of modifications of graphs as will be seen in Theorem~\ref{RCI.thm10}. The admissible class of modified graphs is called a proper system of horizontal networks.

\definition\label{RCI.def10}
A sequence of graphs $\{(\Omega_m, E_m)\}_{m \ge 0}$ is called a proper system of horizontal networks with indices $(N, L_0, L_1, L_2)$ if and only if the following conditions (N1), (N2), (N3), (N4) and (N5) are satisfied:\\
(N1)\,\,
For every $m \ge 0$, $\Omega_m = A_m \cup V_m$ where $A_m \subseteq (T)_m$ and $V_m \subseteq X$.\\
(N2)\,\,
For any $m \ge 0$ and $w \in (T)_m$, $\Omega_{m, w} \neq \emptyset$, where $\Omega_{m, w}$ is defined as
\[
\Omega_{m, w} = (\{w\} \cap A_m) \cup (V_m \cap K_w).
\]
(N3)\,\,
If we define 
\[
E_m(u, v) = \{(x, y)| (x, y) \in E_m, x \in \Omega_{m, u}, y \in \Omega_{m, v}\},
\]
then
\[
\#(E_m(u, v)) \le L_0
\]
for any $m \ge 0$ and $u, v \in (T)_m$\\
(N4)\,\,
For any $(x, y) \in E_m$, $x \in \Omega_{m, u}$ and $y \in \Omega_{m, v}$ for some $(u, v) \in J^h_{N, m}$.\\
(N5)\,\,
For any $u, v \in J^h_{L_1}$, $x \in \Omega_{|u|, u}$ and $y \in \Omega_{|u|, v}$, there exist $(x_1, \ldots, x_n)$ and $(w(1), \ldots, w(n))$ such that $w(i) \in \GG_{L_2}(u)$ for any $i \ge 1, \ldots, n$, $(x_i, x_{i + 1}) \in E_{|v|}(w(i), w(i + 1))$ for any $i = 1, \ldots, n - 1$ and $x_1 = x, x_n = y, w(1) = u, w(n) = v$.
\enddefinition

\example\label{RCI.ex10}
Let $\Omega^{(N)}_* = \{((T)_m, J^h_{N, m})\}_{m \ge 0}$. Then $\Omega^{(N)}_*$ is a proper system of horizontal networks with indices $(N, 1, 1, 1)$. 
\endexample

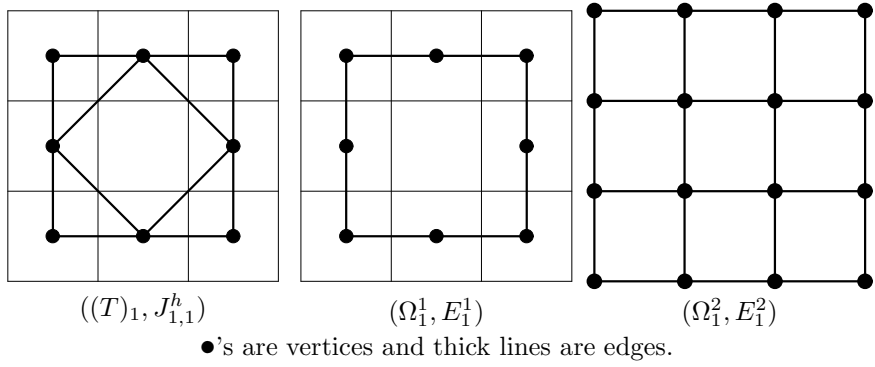
\begin{figure}
\centering
\setlength{\unitlength}{30mm}

\begin{picture}(4, 2)(-1, -0.3)
\drawline (-1, 0)(0.2, 0)(0.2, 1.2)(-1, 1.2)(-1, 0)
\drawline(-0.6, 0)(-0.6, 1.2)
\drawline(-0.2, 0)(-0.2, 1.2)
\drawline(-1, 0.4)(0.2, 0.4)
\drawline(-1, 0.8)(0.2, 0.8)
\put(-0.8, 0.2){\circle*{0.06}}
\put(-0.8, 0.6){\circle*{0.06}}
\put(-0.8, 1){\circle*{0.06}}
\put(-0.4, 0.2){\circle*{0.06}}
\put(-0.4, 1){\circle*{0.06}}
\put(0, 1){\circle*{0.06}}
\put(0, 0.6){\circle*{0.06}}
\put(0, 0.2){\circle*{0.06}}
\thicklines
\drawline(-0.8, 0.2)(-0.8, 1)(0, 1)(0, 0.2)(-0.8, 0.2)
\drawline(-0.4, 0.2)(-0.8, 0.6)(-0.4,1)(0, 0.6)(-0.4, 0.2)
\put(-0.4, -0.2){\makebox(0,0)[b]{$((T)_1, J^h_{1, 1})$}}

\thinlines
\drawline (0.3, 0)(1.5, 0)(1.5, 1.2)(0.3, 1.2)(0.3, 0)
\drawline(0.7, 0)(0.7, 1.2)
\drawline(1.1, 0)(1.1, 1.2)
\drawline(0.3, 0.4)(1.5, 0.4)
\drawline(0.3, 0.8)(1.5, 0.8)
\put(0.5, 0.2){\circle*{0.06}}
\put(0.5, 0.6){\circle*{0.06}}
\put(0.5, 1){\circle*{0.06}}
\put(0.9, 0.2){\circle*{0.06}}
\put(0.9, 1){\circle*{0.06}}
\put(1.3, 1){\circle*{0.06}}
\put(1.3, 0.6){\circle*{0.06}}
\put(1.3, 0.2){\circle*{0.06}}
\thicklines
\drawline(0.5, 0.2)(0.5, 1)(1.3, 1)(1.3, 0.2)(0.5, 0.2)
\put(0.9, -0.2){\makebox(0,0)[b]{$(\Omega^1_1, E^1_1)$}}

\path (1.6, 0)(2.8, 0)(2.8, 1.2)(1.6, 1.2)(1.6, 0)
\drawline(2.0, 0)(2.0, 1.2)
\drawline(2.4, 0)(2.4, 1.2)
\drawline(1.6, 0.4)(2.8, 0.4)
\drawline(1.6, 0.8)(2.8, 0.8)
\put(1.6, 0){\circle*{0.06}}
\put(2, 0){\circle*{0.06}}
\put(2.4, 0){\circle*{0.06}}
\put(2.8, 0){\circle*{0.06}}
\put(1.6, 0.4){\circle*{0.06}}
\put(2, 0.4){\circle*{0.06}}
\put(2.4, 0.4){\circle*{0.06}}
\put(2.8, 0.4){\circle*{0.06}}
\put(1.6, 0.8){\circle*{0.06}}
\put(2, 0.8){\circle*{0.06}}
\put(2.4, 0.8){\circle*{0.06}}
\put(2.8, 0.8){\circle*{0.06}}
\put(1.6, 1.2){\circle*{0.06}}
\put(2, 1.2){\circle*{0.06}}
\put(2.4, 1.2){\circle*{0.06}}
\put(2.8, 1.2){\circle*{0.06}}
\put(2.2, -0.2){\makebox(0,0)[b]{$(\Omega^2_1, E^2_1)$}}

\put(0.9, -0.3){\makebox(0, 0){{\large $\bullet$}'s are vertices and thick lines are edges.}}
\end{picture}
\caption{Proper systems of horizontal networks: the Sierpinski carpet}\label{sc}
\end{figure}

\example[the Sierpinski carpet; Figure~\ref{sc}]\label{RCI.ex20}
Let $d$ be the Euclidean metric (divided by $\sqrt{2}$ so that the diameter of $[0, 1]^2$ is one.) Then $h_{1/3} \bl d$. Obviously, $d$ is $1$-adapted to the weight function $h_{1/3}$, exponential and uniformly finite. In this case, the original edges of the horizontal graph $((T)_m, J^h_{1, m})$ contain slanting edges, which are $(w, v) \in (T)_m \times (T_m)$ with $K_w \cap K_v$ being a single point. Even if all the slanting edges are deleted, we still have a proper system of horizontal networks $\{(\Omega^1_m, E^1_m)\}_{m \ge 0}$ given by
\[
\Omega^1_m = (T)_m
\]
and
\[
E^1_m  = \{(w, v) | w, v \in (T)_m,\dim_H(K_w \cap K_v, d) = 1\}.
\]
$\{(\Omega^1_m, E^1_m)\}_{m \ge 0}$ is a proper system of horizonal networks with indices $(1, 1, 1, 2)$. \par
There is another natural proper system of horizontal networks. Note that the four points $p_1 = (0, 0), p_3 = (1, 0), p_5 = (1, 1)$ and $p_7 = (0, 1)$ are the corners of the square $[0, 1]^2$. Define $\Omega^2_0 = \{p_1, p_3, p_5, p_7\}$ and 
\begin{multline*}
E^2_0 = \{(p_i, p_j)| \text{the line segment $p_ip_j$ is one of the four line segments}\\\text{ of the boundary of $[0, 1]^2$.}\}
\end{multline*}
For $m \ge 1$, we define
\[
\Omega^2_m = \cup_{w \in (T)_m} F_w(\Omega^2_0)
\]
and
\[
E^2_m = \{(F_w(p_j), F_w(p_j))| w \in (T)_m, (p_i, p_j) \in E^2_0\}
\]
Then $\{(\Omega^2_m, E^2_m)\}_{m \ge 0}$ is a proper system of horizontal networks with indices $(1, 5, 1, 1)$. In this case all the vertices are the points in the Sierpinski carpet, and the length between the end points of an edge in $E_m^2$ is $3^{-m}$.
\endexample

\notation
Let $\Omega = \{(\Omega_m, E_m)\}_{m \ge 0}$ be a proper system of horizontal networks. For any $U \subseteq (T)_m$, we define
\begin{equation}\label{RCI.eq20}
\Omega_m(U) = \bigcup_{v \in U} \Omega_{m, v}.
\end{equation}
Furthermore, for $w \in T$, $k \ge 0$ and $n \ge 0$, we define
\begin{align}
\Omega^k(w, n) &= \Omega_{|w| + k}(S^k(\GG_n(w)))\label{RCI.eq30}\\
\Omega^{k, c}(w, n) &= \Omega_{|w| + k}(\sd{(T)_{|w| + k}}{S^k(\GG_n(w))})\label{RCI.eq40}
\end{align}
\endnotation

In the same manner as the original case, we define the $p$-conductances and the critical index of $p$-energies for a proper system of horizontal networks as follows.

\definition\label{RCI.def20}
Let $\Omega = \{(\Omega_m, E_m)\}_{m \ge 0}$ be a proper system of horizontal networks. Define
\begin{align*}
\E_{p, k, w}(N_1, N_2, \Omega) &= \E_p(\Omega_{|w| + k}, E_{|w| + k}, \Omega^k(w, N_1), \Omega^{k, c}(w, N_2))\\
\E_{p, k}(N_1, N_2, \Omega) &= \sup_{w \in T} \E_{p, k, w}(N_1, N_2, \Omega)\\
\overline{\E}_p(N_1, N_2, \Omega) &= \limsup_{k \to \infty} \E_{p, k}(N_1, N_2, \Omega),\\
\underline{\E}_p(N_1, N_2, \Omega) &= \liminf_{k \to \infty} \E_{p, k}(N_1, N_2, \Omega)\\
\oI_{\E}(N_1, N_2, \Omega) &= \inf\{p| \oE_p(N_1, N_2, \Omega) = 0\}\\
\uI_{\E}(N_1, N_2, \Omega) &= \inf\{p| \uE_p(N_1, N_2, \Omega) = 0\}
\end{align*}
\enddefinition

Comparing Definitions~\ref{PEN.def20} and \ref{RCI.def20}, we notice that 
\[
 \E_{p, k}(N_1, N_2, N) = \E_{p, k}(N_1, N_2, \Omega^{(N)}_*).
\]
Thus Theorem~\ref{PEN.thm10} is a corollary of the following theorem.

\thm\label{RCI.thm10}
Let $\Omega$ be a proper system of horizontal networks. If $N_2 \ge N_1 + M_*$, then
\[
{\oI}_{\E}(N_1, N_2, \Omega) = \uI_{\E}(N_1, N_2, \Omega) = \dim_{AR}(X, d).
\]
\endthm

Before a proof of this theorem, we are going to present a corollary which ensures the finiteness of $\dim_{AR}(X, d)$. To begin with, we need to define growth rates of volumes.

\definition\label{RCI.def30}
Define
\begin{align*}
\oN_* &= \limsup_{n \to \infty} \big(\sup_{w \in T} \#(S^n(\GG_{N_2}(w)))\big)^{\frac 1n}\\
\uN_* &= \liminf_{n \to \infty} \big(\sup_{w \in T} \#(S^n(\GG_{N_2}(w)))\big)^{\frac 1n}.
\end{align*}
\enddefinition

It is easy to see that
\[
\uN_* \le \oN_* \le N_*.
\]
The quantities $\oN_*$ and $\uN_*$ appear to depend on the value of $N_2$ but they do not  as is shown in the next lemma.

\lemma\label{APP.lemma10}
\begin{align*}
\oN_* &= \limsup_{n \to \infty} \big(\sup_{w \in T} \#(S^n(w))\big)^{\frac 1n}\\
\uN_* &= \liminf_{n \to \infty} \big(\sup_{w \in T} \#(S^n(w))\big)^{\frac 1n}.
\end{align*}
\endlemma

\demo
Since $S^n(w) \subseteq S^n(\GG_{N_2}(w))$, we have
\[
\sup_{w \in T} \#(S^n(w)) \le \sup_{w \in T} \#(S^n(\GG_{N_2}(w))).
\]
On the other hand, by the fact that $\#(S^n(w)) \le (N_*)^n$, there exists $w(n) \in T$ such that $\#(S^n(w(n)))$ attains the supremum. Note that
\[
S^n(\GG_{N_2}(w)) = \bigcup_{v \in \GG_{N_2}(w)} S^n(v).
\]
Therefore, 
\[
\#(S^n(\GG_{N_2}(w))) \le \#(\GG_{N_2}(w))\#(S^n(w(n)) \le (L_*)^{N_2}\sup_{w \in T} \#(S^n(w)).
\]
\enddemo

\cor\label{RCI.cor10}
Let $\Omega$ be a proper system of horizontal networks. Then
\[
\dim_{AR}(X, d) \le -\frac{\log{\uN_*}}{\log r} \le -\frac{\log{N_*}}{\log r}.
\]
\endcor

Now we start proving the theorem and the corollary.\par
Using the condition (N5), one can easily obtain the first lemma.

\lemma\label{RCI.lemma00}
Let $\Omega = \{(\Omega_m, E_m)\}_{m \ge 0}$ is a proper system of horizontal networks with indices $(N, L_0, L_1, L_2)$. Then $\Omega$ is a proper system of horizontal networks with indices $(N, L_0, nL_1, (n - 1)L_1 + L_2)$ for any $n \ge 1$.
\endlemma

\lemma\label{RCI.lemma05}
Let $\Omega = \{(\Omega_m, E_m)\}_{m \ge 0}$ be a proper system of horizontal networks. Assume $N_2 \ge N_1 + M_*$. If $\rho$ is a metric on $X$ with $\diam{X, \rho} = 1$ and $\rho \qs d$. Then for any $p > 0$, there exists $c > 0$ such that
\[
\E_{p, k, w}(N_1, N_2, \Omega) \le c\sum_{u \in S^k(\GG_{N_2}(w))} \bigg(\frac{\rho(u)}{\rho(w)}\bigg)^p
 \]
 for any $w \in T$ and $k \ge 0$.
\endlemma

\demo
Since $\rho \qs d$, Theorem~\ref{QSY.thm10} implies that $\rho \in \D_{A, e}(X)$ and $d \gen \rho$. Moreover, $\rho$ is $M_*$-adapted. For $w \in T$, define
\[
f_w(x) = \begin{cases} \min\bigg\{\displaystyle\frac{\rho(K_x, U_{N_2}(w)^c)}{\rho(U_{N_1}(w), U_{N_2}(w)^c)}, 1\bigg\}\quad&\text{if $x \in A_m$,}\\
\min\bigg\{\displaystyle\frac{\rho(x, U_{N_2}(w)^c)}{\rho(U_{N_1}(w), U_{N_2}(w)^c)}, 1\bigg\}\quad&\text{if $x \in V_m$.}
\end{cases}
\]
for any $x \in \cup_{k \ge 0} \Omega_{|w| + k}$. If is easy to see that $f_w(x) = 1$ for any $x \in \Omega^k(w, N_1)$ and $f_w(x) = 0$ for any $x \in \Omega^{k, c}(w, N_2)$. Since $N_2 \ge N_1 + M_*$, we have
\[
U_{N_2}(w) \supseteq U_{M_*}^d(x, r^{|w|}) \supseteq U_{M_*}^{\rho}(x, \c\rho(w)) \supseteq B_{\rho}(x, \c'\rho(w))
\]
for any $x \in \inte{U_{N_1}(w)}$. Thus 
\[
\rho(U_{N_1}(w), (U_{N_2}(w))^c) \ge \c'\rho(w).
\]
Let $(N, L_0, L_1, L_2)$ be the indices of $\Omega$. Since $\rho \ge d$, there exists $\kappa \in (0, 1)$ such that
\[
\kappa\rho(u) \le \rho(v)
\]
if $|u| = |v|$ and $K_u \cap K_v \neq \emptyset$. For any $(u, v) \in J^h_{N, |w| + k}$ and $(x, y) \in E_{|w| + k}(u, v)$, if $(w(1), \ldots, w(N + 1)$ is a horizontal $N$-chain between $u$ and $v$, then
\begin{multline*}
|f_w(x) - f_w(y)| \le \frac{\sup_{a \in K_v, b \in K_w} \rho(a, b)}{\c'\rho(w)}\\ \le \frac{1}{\c'\rho(w)}\sum_{i = 1}^{N + 1} \rho(w(i)) \le (N + 1)\kappa^{-N}(\c')^{-1}\frac{\rho(u)}{\rho(w)}
\end{multline*}
Hence
\[
\sum_{(x, y) \in E_{|w| + k}(u, v)} |f_w(x) - f_w(y)|^p \le L_0((N + 1)\kappa^{-N}(\c')^{-1})^p\bigg(\frac{\rho(u)}{\rho(w)}\bigg)^p.
\]
Set $c_1 = L_0((N + 1)\kappa^{-N}(\c')^{-1})^p$.  Then the above inequality and the condition (N4) imply that
\begin{multline*}
\E_{p, k, w}(N_1, N_2, \Omega)\\
 \le \frac 12\sum_{u \in S^k(\GG_{N_2}(w))}\sum_{v \in \GG_N(u)}\sum_{(x, y) \in E_{|w| + k}(u, v)} |f_w(x) - f_w(y)|^p \\
 \le c_1(L_*)^{N}\sum_{u \in [S^k(\GG_{N_2}(w))} \bigg(\frac{\rho(u)}{\rho(w)}\bigg)^p.
 \end{multline*}
\enddemo

The next lemma yields the fact that $I_{\E}(N_1, N_2, \Omega) \le \dim_{AR}(X, d)$ if $N_2 \ge N_1 + M_*$.
 
 \lemma\label{RCI.lemma10}
Let $\Omega = \{(\Omega_m, E_m)\}_{m \ge 0}$ be a proper system of horizontal networks. Assume $N_2 \ge N_1 + M_*$. If $\dim_{AR}(X, d) < p$, then $\overline{\E}_p(N_1, N_2, \Omega) = 0$.
\endlemma

\demo
Since $\dim_{AR}(X, d) < p$, there exist $q \in [\dim_{AR}(X, d), p)$, a metric $\rho$,  a Borel regular measure $\mu$ and constants $c_1, c_2 > 0$ such that $d \qs \rho$ and
\[
c_1r^q \le \mu(B_{\rho}(x, r)) \le c_2r^q
\]
for any $x \in X$ and $r > 0$. By Lemma~\ref{RCI.lemma05},
\begin{multline}\label{PEN.eq50}
\E_{p, k, w}(N_1, N_2, \Omega)
 \le c\sum_{u \in S^k(\GG_{N_2}(w))} \bigg(\frac{\rho(u)}{\rho(w)}\bigg)^p\\
  \le c\max_{u \in S^k(\GG_{N_2}(w))} \bigg(\frac{\rho(u)}{\rho(w)}\bigg)^{p - q} \sum_{u \in S^k(\GG_{N_2}(w))} \bigg(\frac{\rho(u)}{\rho(w)}\bigg)^{q}.
 \end{multline}
 Since $\rho$ is exponential, there exist $\lambda \in (0, 1)$ and $c > 0$ such that $\rho(v) \le c\lambda^k\rho(\pi^k(v))$ for any $v \in T$. Choose $\kappa$ as in the proof of Lemma~\ref{RCI.lemma05}. If $u \in S^k(\GG_{N_2}(w))$, then
 \begin{equation}\label{PEN.eq60}
 \rho(u) \le c\lambda^k\rho(\pi^k(u)) \le c\lambda^k\kappa^{-N_2}\rho(w).
 \end{equation}
 On the other hand, by Theorem~\ref{MEM.thm100}(7.21), there exist $\c_1, \c_2 > 0$ such that
\[
\c_2\mu(K_u) \le \rho(u)^q \le \c_3\mu(K_u)
\]
for any $u \in T$. This implies
\[
\sum_{u \in S^k(v)} \rho(u)^q \le \c_3L_*\mu(K_v) \le (\c_2)^{-1}\c_3L_*\rho(u)^q \le (\c_2)^{-1}\c_3L_*(\kappa^{-N_2}\rho(w))^q
\]
for any $v \in (T)_{|w|}$. Thus
\begin{multline}\label{PEN.eq70}
\sum_{u \in S^k(\GG_{N_2}(w))} \bigg(\frac{\rho(u)}{\rho(w)}\bigg)^{q} = \sum_{v \in \GG_{N_2}(w)}\sum_{u \in S^k(v)}  \bigg(\frac{\rho(u)}{\rho(w)}\bigg)^{q} \\
\le \#(\GG_{N_2}(w)) (\c_2)^{-1}\c_3L_*\kappa^{-N_2q} \le (L_*)^{N_2 + 1}(\c_2)^{-1}\c_3\kappa^{-N_2q}.
\end{multline}
Combining \eqref{PEN.eq50}, \eqref{PEN.eq60} and \eqref{PEN.eq70}, we obtain
\begin{equation}\label{RCI.eq10}
\E_{p, k, w}(N_1, N_2, \Omega) \le c'\lambda^{(p - q)k}
\end{equation}
where $c'$ is independent of $w$. Therefore, we conclude that $\overline{\E}_p(N_1, N_2, \Omega) = 0$.
\enddemo

The following lemma enable us to apply Theorem~\ref{CAR.thm10} and to construct desired pair of a metric and a measure with Ahlfors regularity.

\lemma\label{RCI.lemma20}
Let $\Omega = \{(\Omega_m, E_m)\}_{m \ge 0}$ be a proper system of horizontal networks and let $N \in \BbN$. If $\underline{\E}_p(N_1, N_2, \Omega) = 0$, then for any $\eta > 0$ and $k_0 \in \BbN$, there exists $k_* \ge k_0$ and $\vp : \sd{T^{(k_*)}}\{\phi\} \to [0, 1]$ such that, for any $w \in T^{(k_*)}$,
\begin{equation}\label{RCI.eq100}
\sum_{i = 1}^m \vp(w(i)) \ge 1
\end{equation}
for any $(w(1), \ldots, w(m)) \in \C_{w, k_*}(N_1, N_2, N)$ and
\begin{equation}\label{RCI.eq110}
\sum_{v \in S^{k_*}(w)} \vp(v)^p < \eta.
\end{equation}
\endlemma

\demo
As $\uE_p(N_1, N_2, \Omega) = 0$, for any $\eta_0 > 0$ and $k_0 \in \BbN$, there exists $k_* \ge k_0$ such that $\E_{p, k_*, w}(N_1, N_2, \Omega) < \eta_0$ for any $w \in T$. Hence there exists $f_w : (T)_{|w| + k_*} \to [0, 1]$ such that $f_w(x) = 1$ for any $x \in \Omega^{k_*}(w, N_1)$, $f_w(x) = 0$ for any $x \in \Omega^{k_*, c}(w, N_2)$ and
\[
\sum_{(x, y) \in E_{|w| + k_*}} |f_w(x) - f_w(y)|^p < \eta_0
\]
Let $(N_0, L_0, L_1, L_2)$ be the indices of $\Omega$. Set $n_0 = \min\{n| N \le nL_1\}$ and $\oN = (n_0 - 1)L_1 + L_2$. Note that $N \le \oN$ because $L_1 \le L_2$.
Define $E_m(U) = \cup_{u_1, u_2 \in U} E_m(u_1, u_2)$ for $U \subseteq (T)_m$. Define $\vp_w: (T)_{|w| + k_*} \to [0, 1)$ by
\[
\vp_w(v) = \begin{cases}\displaystyle\bigg(\sum_{(x, y) \in E_{|w| + k_*}(\GG_{\oN}(v))} |f_w(x) - f_w(y)|^p\bigg)^{1/p}&\text{if $v \in S^{k_*}(\GG_{N_2}(w))$,}\\
\quad\quad\quad\quad\quad\quad\quad\quad\quad0\,\,&\text{otherwise.}
\end{cases}
\]
Let $(w(1), \ldots, w(m)) \in \C_{w, k_*}(N_1, N_2, N)$. By definition,  there exist $w(0)$ and $w(m + 1) \in (T)_{|w| + k_*}$ such that $w(0) \in S^{k_*}(\GG_{N_1}(w))$, $w(m + 1) \notin S^{k_*}(\GG_{N_2}(w))$ and $(w(0), w(1), \ldots, w(m), w(m + 1))$ is a horizontal $N$-jpath. Choose $x_i \in \Omega_{|w| + k_*, w(i)}$ for $i = 0, \ldots, m + 1$. By Lemma~\ref{RCI.lemma00}, for any $i = 0, \ldots, m$, there exist $(x^i_1, \ldots, x^i_{l_i})$ and $(w^i(1), \ldots, w^i(l_i))$ such that $x^i_1 = x_i$, $x^i_{l_i} = x_{i + 1}$, $(x^i_j, x^i_{j + 1}) \in E_{|w| + k_*}(w^i(j), w^i(j + 1))$ and $w^i(j) \in \GG_{\oN}(w(i))$ for any $j = 1, \ldots, l_i - 1$. Concatenating the paths $(x^i_1, \ldots, x^i_{l_i})$ for $i = 0, \ldots, m$ and removing all loops from it, we obtain a path $(z_1, \ldots, z_n)$. By the nature of the construction, for any $i = 1, \ldots, n - 1$, $(z_i, z_{i + 1}) \in E_{|w| + k_*}(\GG_{\oN}(w(j))$ for some $j = 0, \ldots, m$. Hence
\begin{equation}\label{RCI.eq120}
\sum_{i = 1}^m \vp_w(w(i)) \ge \sum_{i = 1}^{n - 1} |f_w(z_i) - f_w(z_{i + 1})| \ge f_w(z_1) - f_w(z_n) = 1.
\end{equation}
Since $\#(\{u| u \in (T)_m, x \in K_u) \le L_*$, we see $\#((u_1, u_2)| (x, y) \in E_m(u_1, u_2)\}) \le (L_*)^2$ for any $(x, y) \in L_m$. Making use of this fact, we see that
\begin{multline*}
\#(\{v|(x, y) \in E_{|w| + k_*}(\GG_{\oN}(v))\}) \le \\\sum_{(u_1, u_2): (x, y) \in E_{|w| + k_*}(u_1, u_2)} \#(\GG_{\oN}(u_1) \cap \GG_{\oN}(u_2))\\
\le \sum_{(u_1, u_2): (x, y) \in E_{|w| + k_*}(u_1, u_2)} (L_*)^{\oN} \le (L_*)^{\oN + 2}.
\end{multline*}
Hence
\begin{multline*}
\sum_{v \in (T)_{|w| + k_*}} \vp_w(v)^p = \sum_{v \in (T)_{|w| + k_*}}\sum_{(x, y) \in E_{|w| + k_*}(\GG_{\oN}(v))} |f_w(x) - f_w(y)|^p\\
\le \sum_{(x,y) \in E_{|w| + K_*}}\#(\{v|(x, y) \in \GG_{\oN}(v)\}) |f_w(x) - f_w(y)|^p < (L_*)^{\oN + 2}\eta_0.
\end{multline*}
Define $\vp: \sd{T^{(k_*)}}{\{\phi\}} \to [0, 1]$ by
\[
\vp(v) = \max\{\vp_w(v)| w \in (T)_{k_*(|v| - 1)}\}.
\]
By \eqref{RCI.eq120}, we obtain \eqref{RCI.eq100}. Since if $\vp_w(v) > 0$, then $\pi^{k_*}(v) \in \GG_{N_2}(w)$, it follows that
\[
\vp(v) = \max\{\vp_w(v) | w \in \GG_{N_2}(\pi^{k_*}(v))\}.
\]
Therefore
\begin{multline*}
\sum_{v \in S^{k_*}(w)} \vp(v)^p \le \sum_{v \in S^{k_*}(w)} \sum_{w' \in \GG_{N_2}(w)} \vp_{w'}(v)^p\\
 < \#(\GG_{N_2}(w))(L_*)^{\oN + 2}\eta_0 \le (L_*)^{N_2 + \oN + 2}\eta_0.
\end{multline*}
So, letting $\eta_0 = (L_*)^{-(N_2 + \oN + 2)}\eta$, we have shown \eqref{RCI.eq110}.
\enddemo

Finally, we are going to complete the proof of Theorem~\ref{RCI.thm10}.

\demo[Proof of Theorem~\ref{RCI.thm10}]
Suppose $N_2 \ge N_1 + M_*$. By Lemma~\ref{RCI.lemma10}, it follows that $\oI_{\E}(N_1, N_2, \Omega) \le \dim_{AR}(X, d)$. To prove the opposite inequality, we assume that $\uI_{\E}(0, N_2, \Omega.) < p$. Set $k_0 = \max\{m_0, k_{N_1}, k_{M_*}\}$. Since $\uE_p(0, N_2, \Omega) = 0$, Lemma~\ref{RCI.lemma20} yields a function $\vp$ satisfying the assumptions of Theorem~\ref{CAR.thm10}. Hence by Theorem~\ref{CAR.thm10}, we find a metric $\rho$ which is quasisymmetric to $d$ and a measure $\mu$ which is Ahlfors $p$-regular with respect to $\rho$. This immediately shows that $\dim_{AR}(X, d) \le p$. Hence we obtain $\dim_{AR}(X, d) \le \uI_{\E}(0, N_2, \Omega) \le \uI_{\E}(N_1, N_2, \Omega)$. Thus we have obtained 
\[
\dim_{AR}(X, d) \le \uI_{\E}(N_1, N_2, \Omega) \le \oI_{\E}(N_1, N_2, \Omega) \le \dim_{AR}(X, d).
\]
\enddemo

\demo[Proof of Corollary~\ref{RCI.cor10}]
Applying Lemma~\ref{RCI.lemma05} in the case where $\rho = d$, we obtain
\[
\E_{p, k, w}(N_1, N_2, \Omega) \le c\sum_{u \in S^k(\GG_{M_*}(w)])} \Bigg(\frac{d(u)}{d(w)}\Bigg)^p \le c'\#(S^k(\GG_{M_*}(w)))r^{pk}.
\]
Set $\oN_k(M) = \sup_{w \in T} \#(S^k(\GG_{M}(w)))$. Then 
\begin{equation}\label{RCI.eq500}
\E_{p, k}(N_1, N_2, \Omega) \le c\oN_k(N_2)r^{pk}
\end{equation}
If $q - \e > \liminf_{k \to \infty} -\frac{\log{\oN_k(N_2)}}{n\log r}$, then there exists $\{k_j\}_{j \ge 1}$ such that
\[
r^{{\e}k_j} \ge \oN_{k_j}(N_2)r^{q{k_j}}.
\]
Hence by \eqref{RCI.eq500}, $I_{\E}(N_1, N_2, \Omega) \le q$. This implies
\[
I_{\E}(N_1, N_2, \Omega) \le \liminf_{k \to \infty} -\frac{\log{\oN_k(N_2)}}{n\log r}.
\]
The rest of the statement follows from the fact that
\[
\#(S^k(\GG_{M}(w))) \le \#(\GG_M(w))(N_*)^k \le (L_*)^M(N_*)^k.
\]

\enddemo

\setcounter{equation}{0}
\section{Relation with $p$-spectral dimensions}\label{APP}

By Theorem~\ref{RCI.thm10}, we see that
\[
\lim_{k \to \infty} \E_{p, k}(N_1, N_2, \Omega) = 0
\]
if $p > \dim_{AR}(X, d)$ and
\[
\liminf_{k \to \infty} \E_{p, k}(N_1, N_2, \Omega) > 0
\]
if $p < \dim_{AR}(X, d)$. So, how about the rate of decrease and/or increase of $\E_{p, k}(N_1, N_2, \Omega)$ as $k \to \infty$?  In this section, we define and investigate the rates and present another characterization of the Ahlfors regular conformal dimension in terms of them.
\par
As in the previous sections, $(T, \A, \phi)$ is a locally finite tree with the root $\phi$, $(X, \O)$ is a compact metrizable topological space with no isolated point, $K: T \to \C(X, \O)$ is a minimal partition. We also assume that $\sup_{w \in T} \#(S(w)) < +\infty$. Furthermore we fix $d \in \D_{A, \e}(X, \O)$ satisfying (BF1) and (BF2) in Section~\ref{BFM}.\par
In this framework, the rates are defined as follows.

\definition\label{APP.def10}
Define
\begin{align*}
\oR_p(N_1, N_2, \Omega) &= \limsup_{n \to \infty} \E_{p, n}(N_1, N_2, \Omega)^{\frac 1n}\\
\uR_p(N_1, N_2, \Omega) &= \liminf_{n \to \infty} \E_{p, n}(N_1, N_2, \Omega)^{\frac 1n}.
\end{align*}
\enddefinition

Ideally, we expect that $\oR_p(N_1, N_2, \Omega) = \uR_p(N_1, N_2, \Omega)$ and that
\[
\frac 1{(R_p)^m}\E_p(f|\Omega, E_m) \to \E_p(f)
\]
as $m \to \infty$ for $f$ belonging to some reasonably large class of functions, where we write $R_p = \oR_p(N_1, N_2, \Omega)$. In particular, for a class of (random) self-similar sets including the Sierpinski gasket and the (generalized) Sierpinski carpets, it is known that 
\begin{equation}\label{APP.eq80}
\oR_2(N_1, N_2, \Omega) = \uR_2(N_1, N_2, \Omega).
\end{equation}
 Moreover, the rate $R_2$ is called the resistance scaling ratio and $\E_2(f)$ has known to induce the ``Brownian motion'' $(\{X_t\}_{t > 0}, \{P_x\}_{x \in X})$ and the ``Laplacian'' $\Delta$ through the formula
\begin{align*}
\E_2(f) &= \int_X f\Delta{f}d\mu\\
E_x(f(X_t)) &= (e^{-t\Delta}f)(x),
\end{align*}
where $E_x(\cdot)$ is the expectation with respect to $P_x(\cdot)$. See \cite{BP}, \cite{BB1}, \cite{Bar4}, \cite{KusZ1} and \cite{AOF} for details\par
Now we start to study the relation between $R_p$'s for different values of $p$.
\lemma\label{APP.lemma00}
Let $\Omega = \{(\Omega_m, E_m)\}_{m \ge 0}$ be a proper system of horizontal networks. If $p < q$, then there exists $c > 0$ such that
\[
\E_{p, k}(N_1, N_2, \Omega)^{\frac 1p} \le c\E_{q, k}(N_1, N_2, \Omega)^{\frac 1q}\sup_{w \in T} \#(S^k(\GG_{N_2}(w)))^{\frac 1p - \frac 1q}.
\]
\endlemma

\demo
Let $(Y, \M, \mu)$ is a measurable space. Assume $\mu(Y) < \infty$. Then by the H{\"o}lder inequality,
\[
\int_Y |u|^pd\mu \le \Bigg(\int_Y |u|^qd\mu\Bigg)^{\frac pq}\mu(Y)^{\frac {q - p}q}.
\]
Applying this to $\E_p(f|\Omega_{|w| + k}, E_{|w| + k})$, we obtain
\begin{multline*}
\E_p(f|\Omega_{|w| + k}, E_{|w| + k})^{\frac 1p} \\
\le \E_q(f|\Omega_{|w| + k}, E_{|w| + k})^{\frac 1q}\Big(\#(\{(x, y)| (x, y) \in E_{|w| + k}, x \in S^k(\GG_{N_2}(w))\}))\Big)^{\frac 1p - \frac 1q} \\
\le \E_q(f|\Omega_{|w| + k}, E_{|w| + k})^{\frac 1q}\Bigg(\sum_{u \in S^k(\GG_{N_2}(w))}\sum_{v \in \GG_N(u)} \#(E_m(u, v))\Bigg)^{\frac 1p - \frac 1q}\\
\le \E_q(f|\Omega_{|w| + k}, E_{|w| + k})^{\frac 1q}\big(L_0(L_*)^N\#(S^k(\GG_{N_2}(w)))\big)^{\frac 1p - \frac 1q}
\end{multline*}
for  any $f \in \F(\Omega_{|w| + k}, E_{|w| + k}, \Omega^k(w, N_1), \Omega^{k, c}(w, N_2))$. Set $c = \big(L_0(L_*)^N\big)^{\frac 1p - \frac 1q}$. Then the above inequality implies that
\[
\E_{p, k, w}(N_1, N_2, \Omega)^{\frac 1p} \le c\E_{q, k, w}(N_1, N_2, \Omega)^{\frac 1q}\#(S^k(\GG_{N_2}(w)))^{\frac 1p - \frac 1q}.
\]
This immediately verifies the desired inequality.
\enddemo

Lemma~\ref{APP.lemma00} immediately implies the following fact.

\lemma\label{APP.lemma20}
Let $\Omega = \{(\Omega_m, E_m)\}_{m \ge 0}$ be a proper system of horizontal networks. If $p < q$, then
\begin{align}\label{APP.eq30}
\oR_p(N_1, N_2, \Omega)^{\frac 1p} &\le \oR_q(N_1, N_2, \Omega)^{\frac 1q}(\oN_*)^{\frac 1p - \frac 1q}\\
\uR_p(N_1, N_2, \Omega)^{\frac 1p} &\le \uR_q(N_1, N_2, \Omega)^{\frac 1q}(\oN_*)^{\frac 1p - \frac 1q}
\end{align}
\endlemma

Using this lemma, we can show the continuity and the monotonicity of $\uR_p$ and $\oR_p$.

\prop\label{APP.prop10}
Let $\Omega = \{(\Omega_m, E_m)\}_{m \ge 0}$ be a proper system of horizontal networks.\\
{\rm (1)}\,\, $p > 0$,
\[
\uR_p(N_1, N_2, \Omega) \le \oR_p(N_1, N_2, \Omega)\le r^p\oN_*
\]
{\rm (2)}\,\,$\oR_p(N_1, N_2, \Omega)$ and $\uR_p(N_1, N_2, \Omega)$ are continuous and monotonically non-increasing as a function of $p$.\\
{\rm (3)}\,\,If $N_2 \ge N_1 + M_*$, then $\oR_p(N_1, N_2, \Omega) < 1$ for any $p > \dim_{AR}(X, d)$.
\endprop

\demo
(1)\,\,By \eqref{RCI.eq500}, 
\[
\E_{p, k}(N_1, N_2, \Omega) \le c'\oN_k(N_2)r^{pk}.
\]
This immediately implies the desired inequality.\\
(2)\,\,
Since $|f(x) - f(y)| \le 1$ if $0 \le f(x) \le 1$ for any $x \in \Omega_m$, we see that
\[
\E_{p, n, w}(N_1, N_2, \Omega) \ge \E_{q, n, w}(N_1, N_2, \Omega) 
\]
whenever $p < q$. Hence $\uR_p$ and $\oR_p$ are monotonically decreasing. Set $R_p = \oR_p(N_1, N_2, \Omega)$. By \eqref{APP.eq30}, if $p < q$, then
\[
R_q \le R_p \le (R_q)^{\frac pq}(\oN_*)^{1 - \frac pq}
\]
This shows that $\lim_{p \uparrow q}R_p = R_q$. Exchanging $p$ and $q$, we obtain
\[
(R_q)^{\frac pq}(\oN_*)^{1 - \frac pq} \le R_p \le R_q
\]
if $q < p$. This implies $\lim_{p \downarrow q} R_p = R_q$. Thus $R_p$ is continuous. The same discussion works for $\uR_p(N_1, N_2, \Omega)$ as well. \\
(3)\,\,
This follows from \eqref{RCI.eq10}.
\enddemo

Consequently, we obtain another characterization of the AR conformal dimension.

\thm\label{APP.thm05}
Let $\Omega = \{(\Omega_m, E_m)\}_{m \ge 0}$ be a proper system of horizontal networks. Assume that $N_2 \ge N_1 + M_*$. Then
\begin{align}\label{APP.eq60}
\dim_{AR}(X, d) &= \inf\{p| \oR_p(N_1, N_2, \Omega) < 1\} = \max\{p | \oR_p(N_1, N_2, \Omega) = 1\}\\
&= \inf\{p| \uR_p(N_1, N_2, \Omega) < 1\} = \max\{p | \uR_p(N_1, N_2, \Omega) = 1\}.
\end{align}
In particular, if $p_* = \dim_{AR}(X, d)$, then 
\begin{equation}\label{APP.eq70}
\lim_{n \to \infty} \E_{p_*, n}(N_1, N_2, \Omega)^{\frac 1n} = 1.
\end{equation}
\endthm

\demo
Write $\uR_p = \uR_p(N_1, N_2, \Omega)$ and $p_* = \dim_{AR}(X, d)$. Since $\uR_p$ is continuous, $\lim_{p \downarrow d_*} \uR_p \le 1$. If this limit is less than $1$, the continuity of $\uR_p$ implies that $\uR_{p_* + \e} < 1$ for sufficiently small $\e > 0$. Then $\E_{p_* + \e}(N_1, N_2, \Omega) = 0$. Hence $p_* + \e \le d_*$. This contradiction shows that $\lim_{p \downarrow p_*} \uR_p = 1$. Consequently, $\uR_{p_*} = 1$. Since $\uR_p < 1$ if $p > d_*$, we may verify \eqref{APP.eq60} for $\uR_p(N_1, N_2, \Omega)$. The same discussion works for $\oR_p(N_1, N_2, \Omega)$ as well. Consequently, $\oR_{p_*}(N_1, N_2, \Omega) = \uR_{p_*}(N_1, N_2, \Omega) = 1$. Hence we have \eqref{APP.eq70}
\enddemo

Next we define the (upper and lower) $p$-spectral dimension $\od^S_p(N_1, N_2, \Omega)$ and $\ud^S_p(N_1, N_2, \Omega)$. 
\definition\label{APP.def100}
Define $\od^S_p(N_1, N_2, \Omega)$ and $\ud^S_p(N_1, N_2, \Omega)$ by
\begin{align*}
\od^S_p(N_1, N_2, \Omega) &= \frac{p\log{\oN_*}}{\log{\oN_*} - \log {\oR_p(N_1, N_2, \Omega)}}\\
\ud^S_p(N_1, N_2, \Omega) &= \frac{p\log{\oN_*}}{\log{\oN_*} - \log {\uR_p(N_1, N_2, \Omega)}}.
\end{align*}
The quantities $\od^S_p(N_1, N_2, \Omega)$ and $\ud^S_P(N_1, N_2, \Omega)$ are called the upper $p$-spectral dimension and the lower $p$-spectral dimension respectively.
\enddefinition
Note that $\od^S_p(N_1, N_2, \Omega)$ and $\ud^S_p(N_1, N_2, \Omega)$ coincide with the unique numbers $\od, \ud \in \BbR$ which satisfy
\begin{equation}\label{APP.eq50}
\oN_*\Big(\frac{\oR_p(N_1, N_2, \Omega)}{\oN_*}\Big)^{\od/p} = 1\quad\text{and}\quad\oN_*\Big(\frac{\uR_p(N_1, N_2, \Omega)}{\oN_*}\Big)^{\ud/p} = 1
\end{equation}
respectively.\par
For the Sierpinski gasket and the generalized Sierpinski carpets, the equality \eqref{APP.eq80} implies $\od^S_2(N_1, N_1, \Omega) = \ud^S_2(N_1, N_2, \Omega)$, which is called the spectral dimension and written as $d^S$. The spectral dimension has been known to represent asymptotic behaviors of the Brownian motion and the Laplacian. See \cite{BB6}, \cite{BP} and \cite{KL1} for example. For example, if $p(t, x, y)$ is the transition density of the Brownian motion, then
\[
c_1t^{-d^S/2} \le p(t, x, x) \le c_2t^{-d^S/2}
\]
for any $t \in (0, 1]$ and $x \in X$. Moreover, let $N(\cdot)$ is the eigenvalue counting function of $\Delta$, i.e.
\[
N(\lambda) = \text{the number of eigenvalues $\le \lambda$ taking the multiplicity into account}.
\]
Then
\[
c_1\lambda^{d^S/2} \le N(\lambda) \le c_2\lambda^{d^S/2}
\]
for any $\lambda \ge 1$.
Immediately by the above definition, we obtain the following lemma.
\lemma\label{APP.lemma30}
{\rm (a)}\,\,$\oR_p(N_1, N_2, \Omega) < 1$ if and only if $\od^S_p(N_1, N_2, \Omega) < p$, \\
{\rm (b)}\,\,$\oR_p(N_1, N_2, \Omega) = 1$ if and only if $\od^S_p(N_1, N_2, \Omega) = p$\\ 
{\rm (c)}\,\,$\oR_p(N_1, N_2, \Omega) > 1$ if and only if $\od^S_p(N_1, N_2, \Omega) > p$.
\endlemma

Finally we present the relation between $p$-spectral dimension and the Ahlfors regular conformal dimension.
\thm\label{APP.thm10}
Let $\Omega = \{(\Omega_m, E_m)\}_{m \ge 0}$ be a proper system of horizontal networks.\\
{\rm (1)}\,\,If $\uR_p(N_1, N_2, \Omega) < 1$, then 
\[
\dim_{AR}(X, d) \le \ud^S_p(N_1, N_2, \Omega)  \le \od^S_p(N_1, N_2, \Omega) < p.
\]
{\rm (2)}\,\,If $\oR_p(N_1, N_2, \Omega) \ge 1$, then 
\[
\dim_{AR}(X, d) \ge \od^S_p(N_1, N_2, \Omega) \ge \ud^S_p(N_1, N_2, \Omega) \ge p.
\]
\endthm

For the case of the Sierpinski gasket and the (generalized) Sierpinski carpets, the above theorem shows that either
\[
\dim_{AR}(X, d_*) \le d^S < 2
\]
or
\[
\dim_{AR}(X, d_*) \ge d^S \ge 2,
\]
where $d_*$ is the restriction of the Euclidean metric. For the standard planar Sierpinski carpet in Example~\ref{RCI.ex20}, it has been shown in \cite{BBS} that 
\[
d^S \le 1.805
\]
by rigorous numerical estimate. This gives an upper estimate of the Ahlfors regular conformal dimension of the Sierpinski carpet.

\demo[Proof of Theorem~\ref{APP.thm10}]
(1)\,\,
Write $\uR_p = \uR_p(N_1, N_2, \Omega)$ and $\ud_p = \ud^S_p(N_1, N_2, \Omega)$.
Suppose that $\ud_p < q < p$. Using \eqref{APP.eq50}, we see that there exists $\e > 0$ such that
\[
(1 + \e)\oN_*\Big(\frac{R_p}{\oN_*}\Big)^{q/p} < 1.
\]
Choose $\{n_j\}_{j \ge 1}$ so that $\uR_p = \lim_{j \to \infty} \E_{p, n_j}(N_1, N_2, \Omega)^{\frac 1{n_j}}$. Then for sufficiently large $j$, we have
\[
\E_{p, n_j}(N_1, N_2, \Omega) \le ((1 + \e)\uR_p)^{n_j}\,\,{and}\,\, \sup_{w \in T} \#(S^{n_j}(\GG_{M_1}(w))) \le ((1 + \e)\oN_*)^{n_j}.
\]
Hence by Lemma~\ref{APP.lemma00}, as $j \to 0$,
\[
\E_{q, n_j}(0, M_1, \Omega) \le c\Bigg((1 + \e)\oN_*\bigg(\frac{\uR_p}{\oN_*}\bigg)^{\frac qp}\Bigg)^{n_j} \to 0.
\]
Therefore $\uI_{\E}(0, M_1, \Omega) \le q$. Using Theorem~\ref{RCI.thm10}, we have $\dim_{AR}(X, d) \le \ud_p$.\\
(2)\,\,Set $\oR_p = \oR(N_1, N_2, \Omega)$. Assume that $\oR_p > 1$. Let $q \in (p, \od_p)$. By \eqref{APP.eq50}, 
\[
\frac{\oN_*}{1 - \e}\bigg((1 - \e)^2\frac{\oR_p}{\oN_*}\bigg)^{\frac qp} > 1
\]
for sufficiently small $\e > 0$. Choose $\{n_j\}_{j \ge 1}$ so that  $\E_{p, n_j}(N_1, N_2, \Omega)^{\frac 1{n_j}}  \to \oR_p$ as $j  \to \infty$. Using Lemma~\ref{APP.lemma00}, we have
\[
((1 - \e)\oR_p)^{n_j\frac qp} \le c\E_{q, n_j}(N_1, N_2, \Omega)\bigg(\frac{\oN_*}{1 - \e}\bigg)^{n_j\frac{q - p}p}
\]
for sufficiently large $j$. This implies
\[
1 \le \Bigg(\frac{\oN_*}{1 - \e}\bigg((1 - \e)^2\frac{\oR_p}{\oN_*}\bigg)^{\frac qp}\Bigg)^{n_j} \le c\E_{q, n_j}(N_1, N_2, \Omega).
\]
Thus we have $\oE_q(N_1, N_2, \Omega) > 0$. Hence $q \le \dim_{AR}(X, d)$. Consequently by Theorem~\ref{RCI.thm10}, we have $\od_p \le \dim_{AR}(X, d)$.\\
If $\oR_p = 1$, then $\od_p = p \le \dim_{AR}(X, d)$ by Theorem~\ref{APP.thm05}.
\enddemo


\section{Combinatorial modulus of curves}\label{MIT}

Originally in \cite{CarPiag}, the characterization of the Ahlfors regular conformal dimension has been given in terms of the critical exponent of $p$-combinatorial modulus of curve families. In this section, we are going to show a direct correspondence between $p$-energies and $p$-combinatorial moduli and reproduce Piaggio's result in \cite{CarPiag} within our framework.\par
As in the previous sections, $(T, \A, \phi)$ is a locally finite tree with the root $\phi$, $(X, \O)$ is a compact metrizable topological space with no isolated point, $K: T \to \C(X, \O)$ is a minimal partition. We also assume that $\sup_{w \in T} \#(S(w)) < +\infty$. Furthermore, we fix $d \in \D_{A, \e}(X, \O)$ satisfying the basic framework, i.e. (BF1) and (BF2) in Section~\ref{BFM}.

\definition\label{MIT.def10}
 Let $(V, E)$ be a non-directed graph.  Set
\begin{multline*}
\P(V, E) = \{(x(1), \ldots, x(n))| x(i) \in V\,\,\text{for any $i = 1, \ldots, n$ and}\\\,\,(x(i), x(i + 1)) \in E\,\,\text{for any $i = 1, \ldots, n - 1$}\}.
\end{multline*}
For $U_1, U_2 \subseteq V$ with $U_1 \cap U_2 = \emptyset$, set
\begin{multline*}
\C(V, E, U_1, U_2) = \{(x(1), \ldots, x(m))|\,\,
\text{there exist $x(0), x(m + 1) \in V$}\\
\text{such that}\,\,(x(0), \ldots, x(m + 1)) \in \P(V, E), x(0) \in U_1, x(m + 1) \in U_2\}.
\end{multline*}
Define
\begin{multline*}
\F_M(V, E, U_1, U_2) = \{f| f: V \to [0, \infty), \sum_{i = 1}^m f(x(i)) \ge 1\\\,\,\text{for any $(x(1), \ldots, x(m)) \in \C(V, E, U_1, U_2)$}\}
\end{multline*}
and
\[
\Mod_p(V, E, U_1, U_2) = \inf\{\sum_{x \in V} |f(x)|^p|f \in \F_M(V, E, U_1, U_2)\},
\]
which is called the $p$-modulus of curves connecting $U_1$ and $U_2$.
\enddefinition

\definition\label{MTP.def20}
Let $(V, E)$ be a non-directed graph. Assume that $U_1, U_2 \subseteq V$ and $U_1 \cap U_2 = \emptyset$.\\
(1)\,\,For $f \in \F_M(V, E, U_1, U_2)$, define
\[
F(f)(x) = \min\{\sum_{i = 1}^k f(x(i))| (x(1), \ldots, x(m)) \in \P(V, E), x(1) \in U_2, x(k) = x\}.
\]
(2)\,\,
For $g \in \F_F(V, E, U_1, U_2)$, define
\[
G(g)(x) = \sum_{(x, y) \in E} |g(x) - g(y)|.
\]
\enddefinition

The following version of discrete H{\"o}lder inequality will be used several times. It is obtained by applying the ordinary H{\"o}lder inequality to a sum of Dirac measures.

\lemma\label{MIT.lemma05}
Let $C_h(p, n) = \max\{n^{p - 1}, 1\}$. For any $a_1, \ldots, a_n \in \BbR$,
\[
\Big(\sum_{i = 1}^n |a_i|\Big)^p \le C_h(p, n)\sum_{i = 1}^n |a_i|^p.
\]
\endlemma

 \lemma\label{MIT.lemma10}
Let $(V, E)$ be a non-directed graph. Assume that $U_1, U_2 \subseteq V$ and $U_1 \cap U_2 = \emptyset$. Define $L(V, E) = \max\{\#(\{y| (x, y) \in E\})| x \in V\}$.\\
(1)\,\,
For any $f \in \F_M(V, E, U_1, U_2)$, $F(f) \in \F_F(V, E, U_1, U_2)$ and
\begin{equation}\label{MIT.eq10}
\E_p(F(f)| V, E) \le C_h(p, 2)L(V, E)\sum_{x \in V} f(x)^p.
\end{equation}
(2)\,\,
For any $g \in \F_F(V, E, U_1, U_2)$, $G(g) \in \F_M(V, E, U_1, U_2)$ and
\begin{equation}\label{MIT.eq20}
\sum_{x \in V} G(g)(x)^p \le 2C_h(p, L(E, V))\E_p(g| V, E).
\end{equation}
\endlemma

\demo
(1)\,\,The claim that $F(f) \in \F_F(V, E, U_1, U_2)$ is immediate by the definition. If $(x(1), \ldots, x(m)) \in \P(V, E)$, $x(0) = x$ and $x(m) = y$, then
\begin{align*}
F(f)(x) + \sum_{i = 1}^m f(x(i)) &\ge F(f)(y)\quad\text{and}\\
F(f)(y) + \sum_{i = 0}^{m - 1}f(x(i)) &\ge F(f)(x).
\end{align*}
Therefore,
\begin{multline*}
|F(f)(x) - F(f)(y)| \le \\
\min\{\sum_{i = 1}^m f(x(i))| (x(0), \ldots, x(m)) \in \P(V, E), x(0) = x, x(m) = y\}.
\end{multline*}
 This implies 
\[
|F(f)(x) - F(f)(y)| \le f(x) + f(y)
\]
if $(x, y) \in E$. Thus by Lemma~\ref{MIT.lemma05},
\begin{multline*}
\E_p(F(f)| V, E) = \frac 12\sum_{(x, y) \in E}|F(f)(x) - F(f)(y)|^p \le \frac 12\sum_{(x, y) \in E} (f(x) + f(y))^p\\
\frac{C_h(p, 2)}2\sum_{(x, y) \in E} (f(x)^p + f(y)^p) \le C_h(p, 2)L(V, E)\sum_{x \in V} f(x)^p.
\end{multline*}
(2)\,\,
Let $(x(1), \ldots, x(m)) \in \C(V, E, U_1, U_2)$. Then $(x(0),x(1)) \in E$ for some $x(0) \in V_1$ and $(x(m), x(m + 1)) \in V_2$ for some $x(m + 1) \in V_2$. Set $j = \min\{i| i \in \{0, \ldots, m + 1\}, x(i) \in U_2\} - 1$. Since $G(g)(x(i)) \ge |g(x(i - 1) - g(x(i))|$ for any $i = 1, \ldots, j - 1$ and $G(g)(x(j)) \ge |g(x(j - 1)) - g(x(j))| + |g(x(j)) - g(x(j + 1))|$, we see that
\begin{multline*}
\sum_{i = 1}^m G(g)(x(i)) \ge \sum_{i = 1}^j G(g)(x(i)) \ge \sum_{i = 1}^{j + 1} |g(x(i)) - g(x(i - 1))|\\ \ge \sum_{i = 1}^{j + 1} g(x(i)) - g(x(i - 1)) \ge 1.
\end{multline*}
Thus $G(g) \in \F_M(V, E, U_1, U_2)$. Moreover, by Lemma~\ref{MIT.lemma05}
\begin{multline*}
\sum_{x \in V} G(g)(x)^p \le C_h(p, L(E, V))\sum_{x \in V}\sum_{y: (x, y) \in E} |g(x) - g(y)|^p\\
\le 2C_h(p, L(E, V))\E_p(g| V, E).
\end{multline*}
\enddemo

Taking infimums in \eqref{MIT.eq10} and \eqref{MIT.eq20}, we obtain the following proposition giving a direct connection between $p$-energy and $p$-modulus.

\prop\label{MIT.prop10}
Let $(V, E)$ be a non-directed graph. Assume that $U_1, U_2 \subseteq V$ and $U_1 \cap U_2 = \emptyset$. Then
\[
\E_p(V, E, U_1, U_2) \le C_h(p, 2)L(V, E)\Mod_p(V, E, U_1, U_2)
\]
and
\[
\Mod_p(V, E, U_1, U_2) \le 2C_h(p, L(V, E))\E_p(V, E, U_1, U_2).
\]
\endprop

Next we give a definition of the critical index of $p$-moduli.

\definition\label{MIT.def30}
Let $\Omega = \{(\Omega_m, E_m)\}_{m \ge 0}$ be a proper system of horizontal networks. Define
\begin{align*}
\M_{p, k}(N_1, N_2, \Omega) &= \sup_{w \in T} \Mod_p(\Omega_{|w| + k}, E_{|w| + k}, \Omega^k(w, N_1), \Omega^{k, c}(w, N_2))\\
\overline{\M}_p(N_1, N_2, \Omega) &= \limsup_{k \to \infty} \M_{p, k}(N_1, N_2, \Omega),\\
\underline{\M}_p(N_1, N_2, \Omega) &= \liminf_{k \to \infty} \M_{p, k}(N_1, N_2, \Omega)\\
\oI_{\M}(N_1, N_2, \Omega) &= \inf\{p| \oM_p(N_1, N_2, \Omega) = 0\}\\
\uI_{\M}(N_1, N_2, \Omega) &= \inf\{p| \uM_p(N_1, N_2, \Omega) = 0\}
\end{align*}
\enddefinition

Due to Proposition~\ref{MIT.prop10}, $\E_{p, n}(N_1, N_2, \Omega)$ and $\M_{p, n}(N_1, N_2, \Omega)$ can be compared in the following way.

\lemma\label{MIT.lemma20}
Let $\Omega = \{(\Omega_m, E_m)\}_{m \ge 0}$ be a proper system of horizontal networks with indices $(N, L_0, L_1, L_2)$. Then 
\[
\E_{p, n}(N_1, N_2, \Omega) \le C_h(p, 2)L_0(L_*)^{N + 1}\M_{p, n}(N_1, N_2, \Omega)
\]
and
\[
\M_{p, n}(N_1, N_2, \Omega) \le 2C_h(p, L_0(L_*)^{N + 1})\E_{p, n}(N_1, N_2, \Omega).
\]
\endlemma

\demo
It is enough to show that $L(\Omega_m, E_m) \le L_0(L_*)^{N + 1}$. Let $x \in \Omega_m$. If $x = K_w \cap \Omega_m$, then 
By (N4),
\[
\{y| y \in \Omega_m, (x, y) \in E_m\} \subseteq \bigcup_{w: x \in K_w}\bigcup_{v \in \GG_N(w)}\bigcup_{(x, y) \in E_m(w, v)}\{y\}
\]
Using (N3), we see that
\begin{multline*}
\#(\{y| y \in \Omega_m, (x, y) \in E_m\}) \le \#(\{w| x \in K_w\})\#(\GG_N(w)))L_0\\
\le L_*(L_*)^{N}L_0.
\end{multline*}
If $x = w \in (T)_m \cap \Omega_m$, similar arguments show that $\#(y| y \in \Omega_m, (x, y) \in E_m\}) \le (L_*)^NL_0$. Thus we have $L(\Omega, E_m) \le (L_*)^{N + 1}L_0$.
\enddemo

The above lemma combined with Theorem~\ref{RCI.thm10} immediately yields the following characterization of the Ahlfors regular conformal dimension by the critical exponents of discrete moduli.

\thm\label{MIT.thm10}
Let $\Omega = \{(\Omega_m, E_m)\}_{m \ge 0}$ be a proper system of horizontal networks. If $N_2 \ge N_1 + M_*$, then
\[
\oI_{\M}(N_1, N_2, \Omega) = \uI_{\M}(N_1, N_2, \Omega) = \dim_{AR}(X, d).
\]
\endthm

\setcounter{equation}{0}
\section{Positivity at the critical value}\label{PCC}

One of the advantages of the use of discrete moduli is to show the positivity of $\uM_p(N_1, N_2, \Omega)$ and $\uE_p(N_1, N_2, \Omega)$ at the critical value $p_* = \dim_{AR}(X, d)$.\par
As in the previous sections, $(T, \A, \phi)$ is a locally finite tree with the root $\phi$, $(X, \O)$ is a compact metrizable topological space with no isolated point, $K: T \to \C(X, \O)$ is a minimal partition. We also assume that $\sup_{w \in T} \#(S(w)) < +\infty$. Furthermore, we fix $d \in \D_{A, \e}(X, \O)$ satisfying the basic framework, i.e. (BF1) and (BF2) in Section~\ref{BFM}.

\thm\label{MIT.thm20}
Let $\Omega = \{(\Omega_m, E_m)\}_{m \ge 0}$ be a proper system of horizontal networks. Suppose $N_2 \ge N_1 + M_*$. Let $p_* = \dim_{AR}(X, d)$. Then
\[
\uM_{p_*}(N_1, N_2, \Omega) > 0\quad\text{and}\quad\uE_{p_*}(N_1, N_2, \Omega) > 0.
\]
\endthm

First step of a proof is to modify the original proper system of horizontal networks.

\lemma\label{MIT.lemma50}
Let $\Omega = \{(\Omega_m, E_m)\}_{m \ge 0}$ be a proper system of horizontal networks with indices $(N, L_0, L_1, L_2)$. \\
{\rm (1)}\,\,
\[
\#(\Omega_{m, w}) \le L_0(L_*)^N
\]
for any $m \ge 0$ and $w \in (T)_m$.\\
{\rm (2)}\,\,
Define 
\[
J^h_{M, m}[\Omega] = \{(x, y)| m \ge 0, v, w \in (T)_m, v \in \GG_M(w), x \in \Omega_{m, v}, y \in \Omega_{m, w}\}.
\]
Set $\oOmega^M = \{(\Omega_m, J^h_{M, m}[\Omega])\}_{m \ge 0}$. Then $\oOmega^M$ is a proper system of horizontal networks with indices $(M, L_0^2(L_*)^{2N}, M, M)$. Moreover, there exists $c > 0$ such that
\begin{equation}\label{MIT.eq30}
\E_p(f|\Omega_m, J^h_{M, m}[\Omega]) \le c\E_p(f|\Omega_m, E_m)
\end{equation}
for any $m \ge 0$ and $f: \Omega_m \to \BbR$.
\endlemma

\demo
(1)\,\,
Note that 
\[
\{(x, y)| (x, y) \in E_m, x \in \Omega_{m, w}\} \subseteq \bigcup_{v \in \GG_N(w)} E_m(u, v).
\]
Using (N3), we have
\begin{multline*}
\#(\Omega_{m, w}) \le \#(\{(x, y)| (x, y) \in E_m, x \in \Omega_{m, w}\})\\ \le \sum_{v \in \GG_N(w)} \#(E_m(w, v)) \le (L_*)^NL_0.
\end{multline*}
(2)\,\,
By (1), 
\begin{multline*}
\#(\{(x, y)| (x, y) \in J^h_{M, m}(\Omega), x \in \Omega_{m, v}, y \in \Omega_{m, w}\})\\ = \#(\Omega_{m, v} \times \Omega_{m, w}) \le (L_0)^2(L_*)^{2M}.
\end{multline*}
This shows that $\oOmega^M$ is a proper system of horizontal networks with indices $(M, (L_0)^2(L_*)^{2M}, M, M)$. Assume that $M \le L_1$ for the moment. Let $(x, y) \in \Omega_{m, v} \times \Omega_{m, u}$ for some $u, v \in (T)_m$ with $u \in \GG_{M}(v)$. Since $M \le L_1$, the condition (N5) implies that there exist $(x_1, \ldots, x_n)$ and $(w(1), \ldots, w(n))$ such that $w(i) \in \GG_{L_2}(u)$ for any $i \ge 1, \ldots, n$, $(x_i, x_{i + 1}) \in \Omega_{m}(w(i), w(i + 1))$ for any $i = 1, \ldots, n - 1$ and $x_1 = x, x_n = y, w(1) = u, w(n) = v$. Since $n - 1$ is no greater than the total number of edges in $\GG_{L_2}(u)$, we have
\[
n - 1 \le \#\Big(\bigcup_{v_1, v_2 \in \GG_{L_2}(u)} E_m(v_1, v_2)\Big)  \le (L_*)^{2L_2}L_0.
\]
For any $f: \Omega_m \to \BbR$, by Lemma~\ref{MIT.lemma05}
\[
|f(x) - f(y)|^p \le C_h(p, n - 1)\sum_{i = 1}^{n - 1} |u(x(i)) - u(x(i + 1))|^p.
\]
Let $(z_1, z_2) \in E_m$. Consider how many $(x, y) \in J^h_{M, m}$ there are for which $(z_1, z_2)$ appears as $(x(i), x(i + 1))$ in the above inequality. We start with counting the number of possible $u$'s. Since $z_1$ and $z_2$ must belong to $\GG_{L_2}(u)$, the possible number of $u$'s is no greater than $\#(\GG_{L_2}(z_1) \cap \GG_{L_2}(z_2)) \le (L_*)^{L_2}$. For each $u$, 
\[
\#(\{(x, y)| x \in \Omega_{m, u}, y \in \cup_{v \in \GG_{L_1}(u)} \Omega_{m, v}\}) \le (L_*)^{L_1}(L_0)^2(L_*)^{2M}.
\]
Combining those facts, we see that the possible number of $(x, y)$ for which $(z_1, z_2)$ appears as $(x(i), x(i + 1))$ is at most $(L_*)^{L_2 + L_1 + 2M}(L_0)^2$, which is denoted by $C_1$. Then it follows that
\[
\E_p(f|\Omega_m, J^h_{M, m}[\Omega]) \le C_1C_h(p, (L_*)^{2L_2}L_0)\E_p(f|\Omega_m, E_m).
\]
So, we have finished the proof if $M \le M_1$. For general situation, choosing $n_0$ so that $M \le n_0L_1$, we see that $\Omega$ is a proper system of horizontal networks with indices $(N, L_0, n_0L_1, (n_0 - 1)L_1 + L_2)$. Thus replacing $L_1$ and $L_2$ by $n_0L_1$ and $(n_0 - 1)L_1 + L_2$ respectively, we complete the proof for general cases.
\enddemo

\lemma\label{MIT.lemma60}
Let $\Omega = \{(\Omega_m, E_m)\}_{m \ge 0}$ be a proper system of horizontal networks with indices $(N, L_0, L_1, L_2)$. Then
\[
\M_{p, k + l}(0, M, \Omega^{(J)}_*) \le C\M_{p, k}(0, M, \oOmega^{2M + J})\M_{p, l}(0, M, \Omega^{(J)}_*),
\]
for any $k, l, M, J \in \BbN$ and $p > 0$, where $C = L_*C_h(p, (L_*)^{N + 1}L_0)$.
\endlemma

\notation
\begin{equation}\label{MIT.eq40}
\Q_{w, k}(M_1, M_2, \Omega) = \F_M(\Omega_{|w| + k}, E_{|w| + k}, , \Omega^k(w, M_1), \Omega^{k, c}(w, M_2))
\end{equation}
and
\begin{equation}\label{MIT.eq50}
\C_{w, k}(M_1, M_2, \Omega) = \C(\Omega_{|w| + k}, E_{|w| + k}, \Omega^k(w, M_1), \Omega^{k, c}(w, M_2))
\end{equation}
for $M_1, M_2 \ge 1$.
\endnotation

\demo
Let $f \in \Q_{w, k}(0, M, \oOmega^{2M + N})$ and let $g_v \in \Q_{v, l}(0, M, \Omega^{(J)}_*)$ for any $v \in (T)_{|w| + k}$. Define
\[
h(u) = \max\{f(x)g_v(u) |x \in \Omega_{|w| + k, v}, v \in \GG_M(\pi^l(u))\}\chi_{S^{k + l}(\GG_M(w))}(u)
\]
{\bf Claim 1.}\,\,$h \in \Q_{w, k + l}(0, M, \Omega^{(J)}_*)$.\\
Proof of Claim 1:\,\,
Let $(u(1), \ldots, u(m)) \in \C_{w, k + l}(0, M, J)$. Set $v(i) = \pi^l(u(i))$ for $i = 1, \ldots, m$. Let $v_*(1) = v(1)$ and let $i_1 = 1$. Define $v_*(n)$ and $i_n$ inductively as
\[
i_{n + 1} = \max\{j| v(j) \in \GG_{2M}(v_*(n))\} + 1
\]
and $v_*(n + 1) = v(i_{n + 1})$ while $\max\{j | v(j) \in \GG_{2M}(v_*(n))\} < m$. In this way, we construct $(v_*(1), \ldots, v_*(n_*))$ satisfying $\max\{j| v(j) \in \GG_{2M}(v_*(n_*))\} = m$. Since $v(i_{n + 1} - 1) \in \GG_{2M}(v_*(n))$, it follows that $v_*(n + 1) \in \GG_{2M + J}(v_*(n))$. Hence $(v_*(1), \ldots, v_*(n_*)) \in \C_{w, k}(0, M, 2M + J)$. Moreover, $\GG_{M}(v(i)) \cap \GG_M(v(j)) = \emptyset$ if $i \neq j$ and there exists $(u(j_n), \ldots, u(j_n + k_n)) \in \C_{v_*(n), l}(0, M, J)$. Choose $x_i \in \Omega_{|w| + k, v}$ for each $i = 1, \ldots, n$. Since $g_{v_*(n)} \in \Q_{v_*(n), l}(0, M, \Omega^{(J)}_*)$, we have
\[
\sum_{i = j_n}^{j_n + k_n} h(u(i)) \ge \sum_{i = j_n}^{j_n + k_n} f(x_n)g_{v_*(n)}(u(i)) \ge f(x_n).
\]
This and the fact that $(x_1, \ldots, x_n) \in \C_{k, w}(0, M, \oOmega^{2M + J})$ yield 
\[
\sum_{i = 1}^m h(u(i)) \ge \sum_{j = 1}^{n_*} f(x_j) \ge 1.
\]
Thus Claim 1 has been verified. \qed\\
Set $C_0 = C_h(p, (L_*)^{N + 1}L_0)$. Then by Lemma~\ref{MIT.lemma05} and Lemma~\ref{CAM.lemma10},
\begin{multline*}
h(u)^p \le \bigg(\sum_{v \in \GG_M(\pi^l(u))}\sum_{x \in \Omega_{|w| + k, v}} f(x)g_v(u)\bigg)^p\\ \le C_0\sum_{v \in \GG_M(\pi^l(u))}\sum_{x \in \Omega_{|w| + k, v}} f(x)^pg_v(u)^p.
\end{multline*}
Set $M_{p, j, w'} = \Mod_p((T)_{|w'| + j}, J^h_{J, |w'| + j}, S^j(w'), (S^j(\GG_{M}(w')))^c)$. The above inequality yields
\[
M_{p, k + l, w} \le \sum_{u \in (T)_{|w| + k + l}}h(u)^p
\le C_0\sum_{v \in (T)_{|w| + k}}\sum_{x \in \Omega_{|w| + k, v}}\sum_{u \in (T)_{|w| + k + l}}f(x)^pg_v(u)^p.
\]
Hence following the process for getting $\M_{p, k}(N_1, N_2, \Omega)$, we have
\begin{multline*}
M_{p, k + l, w} \le C_0\sum_{v \in (T)_{|w| + k}}\sum_{x \in \Omega_{|w| + k, v}}f(x)^pM_{p, l, v}\\
 \le C_0\sum_{v \in (T)_{|w| + k}}\sum_{x \in \Omega_{|w| + k, v}}f(x)^p\M_{p, l}(0, M, \Omega^{(J)}_*)\\ \le C_0L_*\M_{p, k}(0, M, \Omega^{2M + J)})\M_{p, l}(0, M, \Omega^{(J)}_*).
\end{multline*}
Continuing the process, we finally obtain
\[
\M_{p, k + l}(0, M, \Omega^{(J)}_*) \le C\M_{p, k}(0, M, \oOmega^{2M + J})\M_{p, l}(0, M, \Omega^{(J)}_*),
\]
where $C = C_0L_*$. 
\enddemo

\demo[Proof of Theorem~\ref{MIT.thm20}]
Write $\M_{p, j} = \M_{p, j}(0, M, \Omega^{(J)}_*)$ and $\M'_{p, j} =$\\$ \M_{p, j}(0, M, \oOmega^{2M + J})$. By Lemma~\ref{MIT.lemma60},
\begin{equation}\label{MIT.eq100}
\M_{p, mk + l} \le (C\M'_{p, k})^m\M_{p, l}.
\end{equation}
Assume that $C\M'_{p, k} < 1 - \e$ for some $\e \in (0, 1)$. Then for any $w \in T$, there exists $f_w \in \Q_{w, k}(0, M, \oOmega^{2M + J})$ such that
\[
C\sum_{x \in \Omega_{|w| + k}(S^k(\GG_M(w)))} f_w(x)^p < 1 - \e.
\]
Since
\[
\lim_{\delta \to 0}\max_{x \in [0, 1]} (x^{p -\delta} - x^p) = 0,
\]
there exists $\d_* > 0$ such that $x^{p - \delta} \le x^p + C_2^{-1}C^{-1}\e/2$ for any $\delta \in (0, \d_*]$ and $x \in [0, 1]$, where $C_2 = (N_*)^k(L_*)^{M + N}L_0$. By Lemma~\ref{MIT.lemma50}-(1), 
\begin{multline*}
\#(\Omega_{|w| + k}(S^k(\GG_M(w)))) \le \#(S^k(\GG_M(w)))\max_{v \in S^k(\GG_M(w))} \#(\Omega_{|w| + k, v})\\
 \le (L_*)^M(N_*)^k(L_*)^NL_0 = C_2.
\end{multline*}
This implies 
\[
C\sum_{v \in S^k(\GG_M(w))} f_w(v)^{p - \d}
\le C\sum_{x \in \Omega_{|w| + k}(S^k(\GG_M(w)))} f_w(x)^{p - \d} + \frac {\e}2 \le 1 - \frac{\e}2.
\]
Therefore $\M_p(0, M, \Omega^{(J)}_*) = 0$. By Theorem~\ref{MIT.thm10}, if $M \ge M_*$, it follows that $\dim_{AR}(X, d) < p$. Consequently, $C\M'_{p, k} \ge 1$ for any $k \ge 1$. Therefore, if $M \ge M_*$, then
\[
C^{-1} \le \uM_{p_*}(0, M, \oOmega^{2M + J}).
\]
Using Lemma~\ref{MIT.lemma20}, we see that $0 < \uE_{p_*}(0, M, \oOmega^{2M + J})$. Then the inequality \eqref{MIT.eq30} shows that $0 < \uE_{p_*}(0, M, \Omega)$. Since $\uE_p(0, M, \Omega) \le \uE_p(M', M, \Omega)$ for any $M' \in \{0, 1, \ldots, M - M_*\}$, we conclude that $0 < \uE_{p_*}(N_1, N_2, \Omega)$ for any $N_1, N_2 \ge 0$ with $N_2 \ge N_1 + M_*$. Again by Lemma~\ref{MIT.lemma20}, it follows that $0 < \uM_{p_*}(N_1, N_2, \Omega)$ as well.
\enddemo
\newpage
\appendix{\noindent \bf \Large Appendix}
\setcounter{equation}{0}
\renewcommand{\theequation}{\Alph{section}.\arabic{equation}}

\section{Fact from measure theory}

\prop\label{MEM.prop10}
Let $(X, \M, \mu)$ be measurable space and let $N \in \BbN$. If $U_i \in \M$ for any $i \in \BbN$ and 
\begin{equation}\label{MEM.eq10}
\#(\{i| i \in \BbN, x \in U_i\}) \le N
\end{equation}
for any $x \in X$, then
\[
\sum_{i = 1}^{\infty}\mu(U_i) \le N\mu\bigg(\bigcup_{i \in \BbN} U_i\bigg).
\]
\endprop

\demo
Set $U = \cup_{i \in \BbN} U_i$. Define $U_{i_1\ldots{i_m}} = \cap_{j = 1, \ldots, m} U_{i_j}$. By \eqref{MEM.eq10}, if $m > N$, then $U_{i_1\ldots{i_m}} = \emptyset$. Fix $m \ge 0$ and let rearrange  $\{U_{i_1\ldots{i_m}}| i_1 < i_2 < \ldots < i_m\}$ so that
\[
\{Y^m_j\}_{j \in \BbN} = \{U_{i_1\ldots{i_m}}| i_1 < i_2 < \ldots < i_m\}.
\]
Define
\[
X^m_j = \sd{Y^m_j}{\big(\bigcup_{i \in \BbN, i \neq j} Y^m_i\big)}.
\]
Then
\[
U = \bigcup_{m = 0}^N \Big(\bigcup_{j \in \BbN} X^m_j\Big)
\]
and $X^m_j \cap X^k_l = \emptyset$ if $(m, j) \neq (k, l)$. This implies
\[
\mu(U) = \sum_{m = 0}^N\sum_{j \in \BbN} \mu(X^m_j).
\]
Set $I_j = \{(k, l)| U_j \supseteq X^k_l \neq \emptyset\}$. Then by \eqref{MEM.eq10}, we have $\#(\{j| (k, l) \in I_j\}) \le N$ for any $(k, l)$. This implies
\[
\sum_{j = 1}^{\infty} \mu(U_i) \le N\sum_{m = 0}^N\sum_{j \in \BbN} \mu(X^m_j) = N\mu(U).
\]
\enddemo

\section{List of definitions, notations and conditions}\label{LOD}
{\bf Definitions}\\
\noindent adapted -- Definition~\ref{PAS.def35}, Definition~\ref{PAS.def60}\\
Ahlfors regular -- Definition~\ref{BLE.def50}\\
Ahlfors regular conformal dimension -- Definition~\ref{PEN.def00}\\
Ahlfors regular metric -- Definition~\ref{EAR.def10}\\
balanced -- Definition~\ref{CAM.def20}\\
bi-Lipschitz (metrics) -- Definition~\ref{BLE.def30}\\
bi-Lipschitz (weight functions) -- Definition~\ref{BLE.def10}\\
bridge -- Definition~\ref{HYP.def10}\\
chain -- Definition~\ref{PAS.def20}\\
resolution -- Definition~\ref{HF.def10}\\
degree of distortion -- Definition~\ref{ESS.def10}\\
end of a tree -- Definition~\ref{TWR.def20}\\
exponential -- Definition~\ref{BLE.def40}\\
(super-, sub-)exponential for metrics -- Definition~\ref{QSY.def100}\\
gentle -- Definition~\ref{VDP.def10}\\
geodesic -- Definition~\ref{TWR.def10}\\
Gromov product -- Definition~\ref{HYP.def30}\\
height (of a bridge) -- Definition~\ref{HYP.def10}\\
horizontal edge -- Definition~\ref{HF.def10}\\
horizontal $M$-chain -- Definition~\ref{CAM.def10}\\
horizontally minimal -- Definition~\ref{HYP.def10}\\
hyperbolic -- Definition~\ref{HYP.def30}\\
hyperbolicity of a weight function --Definition\ref{HYP.def25}\\
infinite binary tree -- Example~\ref{TWR.ex10}\\
infinite geodesic ray -- Definition~\ref{TWR.def20}\\
jpath -- Definition~\ref{CAM.def10}\\
jumping path -- Definition~\ref{CAM.def10}\\
locally finite -- Definition~\ref{TWR.def10}\\
minimal -- Definition~\ref{PAS.def20}\\
modulus -- Definition~\ref{MIT.def10}\\
$m$-separated -- Definition~\ref{COM.def20}\\
open set condition --Example~\ref{HYP.ex20}\\
partition -- Definition~\ref{PAS.def20}\\
path -- Definition~\ref{TWR.def10}\\
proper system of horizontal networks -- Definition~\ref{RCI.def10}\\
$p$-modulus of curves -- Definition~\ref{MIT.def10}\\
$p$-spectral dimension -- Definition~\ref{APP.def100}\\
quasisymmetry -- Definition~\ref{INT.def10}\\
rearranged resolution -- Definition~\ref{HYP.def20}\\
simple path -- Definition~\ref{TWR.def10}\\
strongly finite -- Definition~\ref{PAS.def200}\\
sub-exponential -- Definition~\ref{BLE.def40}\\
super-exponential -- Definition~\ref{BLE.def40}\\
thick -- Definition~\ref{VDP.def30}\\
tight -- Definition~\ref{BLES.def10}\\
tree -- Definition~\ref{TWR.def10}\\
tree with a reference point -- Definition~\ref{TWR.def20}\\
uniformly finite -- Definition~\ref{BLE.def40}\\
uniformly perfect -- Definition~\ref{QSY.def20}\\
vertical edge -- Definition~\ref{HF.def10}\\
volume doubling property with respect to a metric -- Definition~\ref{VDP.def20}\\
volume doubling property with respect to a weight function -- Definition~\ref{VDP.def25}\\
weakly $M$-adapted -- Definition~\ref{HYP.def50}\\
weight function -- Definition~\ref{PAS.def10}
\par\vspace{5pt}

\noindent{\bf Notations}\\
$B_w$ -- Definition~\ref{PAS.def20}\\
$\widetilde{B}_r^d(x, cr)$ -- Definition~\ref{HYP.def50}\\
$\B$ -- Definition~\ref{HF.def10}\\
$\B_{\tT^{g, r}}$ -- Definition~\ref{HYP.def20}\\
$C_h(p, n)$ -- Lemma~\ref{MIT.lemma05}\\
$\C(X, \O)$, $\C(X)$: the collection of nonempty compact subsets, -- Definition~\ref{PAS.def20}\\
$\CH_K(A, B)$ -- Definition~\ref{PAS.def20}\\
$\C(V, E, U_1, U_2)$ -- Definition~\ref{MIT.def10}\\
$\C_w^M$ -- Definition~\ref{CAM.def20}\\
$\C_{w, k}(N_1, N_2, N)$ -- Definition~\ref{WFG.def100}\\
$\C_{w, k}(M_1, M_2, \Omega)$ -- \eqref{MIT.eq50}\\
$\od^S_p(N_1, N_2, \Omega)$, $\ud^S_p(N_1, N_2, \Omega)$ -- Definition~\ref{APP.def100}\\
$d_{(T, \B)}$ -- Definition~\ref{HF.def10}\\
$D^g_M(x, y)$ -- Definition~\ref{MAG.def10}\\
$\D(X, \O)$ -- Definition~\ref{PAS.def50}\\
$\D_A(X, \O)$ -- Definition~\ref{BLE.def30}\\
$\D_{A, e}(X, \O)$ -- Definition~\ref{QSY.def120}\\
$E^h_{g, r}$ -- Definition~\ref{HYP.def20}\\
$E_m^h, E^h$: horizontal vertices -- Definition~\ref{HF.def10}\\
$E_m(u, v)$ -- Definition~\ref{RCI.def10}\\
$\E_p(f|V, E)$, $\E_p(V, E, V_1, V_2)$  -- Definition~\ref{PEN.def10}\\
$\E_{p, k}(N_1, N_2, N)$, $\oE_p(N_1, N_2, N)$, $\uE_p(N_1, N_2, N)$ -- Definition~\ref{PEN.def20}\\
$\E_{p, k, w}(N_1, N_2, \Omega)$, $\E_{p, k}(N_1, N_2, \Omega)$, $\oE_p(N_1, N_2, \Omega)$, $\uE_p(N_1, N_2, \Omega)$\par
\hspace{200pt} -- Definition~\ref{RCI.def20}\\
$F(f)(x)$ -- Definition~\ref{MTP.def20}\\
$\F_F(V, E, V_1, V_2)$ -- Definition~\ref{PEN.def10}\\
$\F_M(V, E, U_1, U_2)$ -- Definition~\ref{MIT.def10}\\
$g_d, g_{\mu}$ -- Definition~\ref{PAS.def50}\\
$G(g)(x)$ -- Definition~\ref{MTP.def20}\\
$\G(T)$ -- Definition~\ref{PAS.def10}\\
$\G_e(T)$ -- Definition~\ref{GAE.def05}\\
$h_*$ -- Definition~\ref{AAA.thm100}\\
$h_r$ -- Corollary~\ref{HYP.cor10}\\
$\I_{\E}(N_1, N_2, N)$ -- Definition~\ref{PEN.def20}\\
$\oI_{\E}(N_1, N_2, \Omega), \uI_{\E}(N_1, N_2, \Omega)$ -- Definition~\ref{RCI.def20}\\
$\oI_{\M}(N_1, N_2, \Omega)$, $\uI_{\M}(N_1, N_2, \Omega)$ -- Definition~\ref{MIT.def30}\\
$J_{M, n}^h(K)$, $J_M^h(K)$, $J_M^v(K)$, $J_M(K)$ -- Definition~\ref{CAM.def10}\\
$J_{M, n}^h[\Omega]$ -- Lemma~\ref{MIT.lemma50}\\
$K^{(q)}$ -- Definition~\ref{CAM.def40}\\
$K_w$ -- Definition~\ref{PAS.def20}\\
$\ell_M^{\vp}(\bp)$ -- Definition~\ref{CAM.def30}\\
$L_*$ -- Definition~\ref{BFM.def10}\\
$L_g(\bp)$ -- Definition~\ref{CAM.def30}\\
$\Mod_p(V, E, U_1, U_2)$ -- Definition~\ref{MIT.def10}\\
$\M_{p, k}(N_1, N_2, \Omega)$, $\uM_p(N_1, N_2, \Omega)$, $\oM_p(N_1, N_2, \Omega)$ -- Definition~\ref{MIT.def30}\\
$\M_P(X, \O)$ -- Definition~\ref{PAS.def50}\\
$N_*$ -- Definition~\ref{BFM.def10}\\
$\oN_*, \uN_*$ -- Definition~\ref{RCI.def30}\\
$N_g(w)$ -- Definition~\ref{GAE.def10}\\
$O_w$ -- Definition~\ref{PAS.def20}\\
$\P(V, E)$ -- Definition~\ref{MIT.def10}\\
$\Q_{w, k}(M_1, M_2, \Omega)$ -- \eqref{MIT.eq40}\\
$\oR_p(N_1, N_2, \Omega)$, $\uR_p(N_1, N_2, \Omega)$ -- Definition~\ref{APP.def10}\\
$\R^0_{\kappa}$, $\R^1_{\kappa}$ -- Definition~\ref{ESS.def10}\\
$S^m(A)$ -- Definition~\ref{GAE.def10}\\
$S(\cdot)$ -- Definition~\ref{TWR.def20}\\
$(T)_m$ -- Definition~\ref{TWR.def20}\\
$T^{(N)}$, $(T^{(N)}, \A^{(N)}, \phi)$, $T_m^{(N)}$ -- Example~\ref{TWR.ex10}\\
$T^{(q)}$ -- Definition~\ref{CAM.def40}\\
$T_w$ -- Definition~\ref{TWR.def30}\\
$(T, \B)$ -- Definition~\ref{HF.def10}\\
$\tT^{g, r}, (\tT^{g, r}, \B_{\tT^{g, r}})$ -- Definition~\ref{HYP.def20}\\
$U_M^g(x, s)$ -- Definition~\ref{PAS.def30}\\
$U_M(w, K)$ -- Definition~\ref{CAM.def10}\\
$\GG_M(w, K)$ -- Definition~\ref{CAM.def10}\\
$\d_M^g(x, y)$ -- Definition~\ref{CNB.def10}\\
$\kappa(\cdot)$ -- Definition~\ref{ESS.def10}\\
$\LL_s^g$ -- Definition~\ref{PAS.def10}\\
$\LL^g_{s, M}(\cdot)$ -- Definition~\ref{PAS.def30}\\
$\Omega^k(w, n)$ -- \eqref{RCI.eq30}\\
$\Omega^{k, c}(w, n)$ -- \eqref{RCI.eq40}\\
$\Omega_m(U)$ -- \eqref{RCI.eq20}\\
$\Omega_{m, w}$ -- Definition~\ref{RCI.def10}\\
$\Omega^{(N)}_*$ -- Example~\ref{RCI.ex10}\\
$\overline{\Omega}^M$ -- Lemma~\ref{MIT.lemma50}\\
$\pi$ -- Definition~\ref{TWR.def20}\\
$\pi^{(T, \A, \phi)}$ -- Remark after Definition~\ref{TWR.def20}\\
$\Pi_M^\vp(w)$ -- Definition~\ref{WFG.def10}\\
$\Pi_M^{g, k}(w)$ -- Definition~\ref{WFG.def100}\\
$\rho_*$ -- Definition~\ref{TWR.def30}\\
$\SS$: the collection of ends -- Definition~\ref{TWR.def20}\\
$\SS^w$, $\SS^w_v$ -- Definition~\ref{TWR.def20}\\
$\SS$ and $\SS_v$; abbreviation of $\SS^{\phi}$ and $\SS_v^{\phi}$ respectively,\\
$\SS^{(N)}$ -- Example~\ref{TWR.ex20}\\
\noindent $|w, v|$ -- Definition~\ref{GAE.def10}\\
$\overline{wv}$: the geodesic between $w$ and $v$ of a tree, -- Definition~\ref{TWR.def10}\\
$(w|v)_{((T, \B), \phi)}$, $(w|v)$: Gromov product --  Definition~\ref{HYP.def30}\\
$|w|$ -- Definition~\ref{TWR.def20}\\
$|w|_{(T, \A, \phi)}$ -- Remark after Definition~\ref{TWR.def20}\\
$w \wedge v$ -- Definition~\ref{TWR.def30}\\
$[\omega]_m$ -- Definition~\ref{TWR.def30}\\
$\inn{\cdot}_M$ -- Definition~\ref{WFG.def10}\\
\par\vspace{10pt}

\noindent{\bf Equivalence relations}\\
$\ac$ -- Definition~\ref{BLE.def20}\\
$\bl$ relation on weight functions -- Definition~\ref{BLE.def10}\\
$\bl$ relation on metrics -- Definition~\ref{BLE.def30}\\
$\gen$ -- Definition~\ref{VDP.def10}\\
$\qs$ -- Definition~\ref{INT.def10}
\par\vspace{10pt}

\noindent{\bf Conditions}\\
\noindent (ADa), $\rm (ADb)_M$ -- Theorem~\ref{PAS.thm10}\\
(BF1), (BF2) -- Section~\ref{BFM}\\
(BL), (BL1), (BL2), (BL3) -- Theorem~\ref{BLE.thm10}\\
${\rm (EV)}_M$, ${\rm (EV2)}_M$, ${\rm (EV3)}_M$, ${\rm (EV4)}_M$, ${\rm (EV5)}_M$ -- Theorem~\ref{COM.thm10}\\
(G1), (G2), (G3) -- Definition~\ref{PAS.def10}\\
(N1), (N2), (N3), (N4), (N5) -- Definition~\ref{RCI.def10}\\
(P1), (P2) -- Definition~\ref{PAS.def20}\\
(SQ1), (SQ2), (SQ3) -- Section~\ref{ESS}\\
(SF) -- Strongly finite, \eqref{SMR.eq100}\\
(TH) -- Few lines before Definition~\ref{SMR.def100}\\
(TH1), (TH2), (TH3), (TH4) -- Theorem~\ref{AAA.thm100}\\
(VD1), (VD2), (VD3), (VD4) -- Theorem~\ref{VDP.thm20}

\providecommand{\bysame}{\leavevmode\hbox to3em{\hrulefill}\thinspace}
\providecommand{\MR}{\relax\ifhmode\unskip\space\fi MR }
\providecommand{\MRhref}[2]{%
  \href{http://www.ams.org/mathscinet-getitem?mr=#1}{#2}
}
\providecommand{\href}[2]{#2}

\end{document}